\documentclass[11pt]{article}
\usepackage{amssymb}
\usepackage{amsmath}
\usepackage{times}

\title{Potential Wadge classes\indent}
\author{Dominique LECOMTE}
\date{January 2010}

\def\ufootnote#1{\let\savedthfn\thefootnote\let\thefootnote\relax
\footnote{#1}\let\thefootnote\savedthfn\addtocounter{footnote}{-1}}

\newcommand{\Ana}{{\it\Sigma}^{1}_{1}}
\newcommand{\Ca}{{\it\Pi}^{1}_{1}}
\newcommand{\ca}{{\bf\Pi}^{1}_{1}}
\newcommand{\Borel}{{\it\Delta}^{1}_{1}}

\newcommand{\borel}{{\bf\Delta}^{1}_{1}}

\newcommand{\boraone}{{\bf\Sigma}^{0}_{1}}
\newcommand{\boratwo}{{\bf\Sigma}^{0}_{2}}

\newcommand{\boraxi}{{\bf\Sigma}^{0}_{\xi}}
\newcommand{\borone}{{\bf\Delta}^{0}_{1}}

\newcommand{\bormz}{{\bf\Pi}^{0}_{0}} 
\newcommand{\bormone}{{\bf\Pi}^{0}_{1}}
\newcommand{\Bormone}{{\it\Pi}^{0}_{1}}

\newcommand{\bormpzm}{{\bf\Pi}^{0}_{1+\theta_m}}

\newcommand{\bormpztm}{{\bf\Pi}^{0}_{1+\theta_{\tau (m)}}}

\newcommand{\bormlxi}{{\bf\Pi}^{0}_{<\xi}}

\newcommand{\bormep}{{\bf\Pi}^{0}_{\eta +1}}

\newcommand{\bormxi}{{\bf\Pi}^{0}_{\xi}}

\newcommand{\bormpe}{{\bf\Pi}^{0}_{1+\eta}}
\newcommand{\bormpz}{{\bf\Pi}^{0}_{1+\theta}}
\newcommand{\bormpxz}{{\bf\Pi}^{0}_{1+\xi +\theta}}

\newcommand{\borxi}{{\bf\Delta}^{0}_{\xi}}

\newtheorem{thm} {Theorem} [section]
\newtheorem{defi} [thm] {Definition}
\newtheorem{cor} [thm] {Corollary}
\newtheorem{lem} [thm] {Lemma}
\newtheorem{prop} [thm] {Proposition}
\newtheorem{them} {Theorem} [subsection]
\newtheorem{defin} [them] {Definition}
\newtheorem{coro} [them] {Corollary}
\newtheorem{lemm} [them] {Lemma}
\newtheorem{propo} [them] {Proposition}

\begin{document}

\maketitle

\centerline{$\bullet$ Universit\' e Paris 6, Institut de Math\'ematiques de Jussieu, Projet Analyse Fonctionnelle,}

\centerline{Tour 46-0, Bo\^\i te 186,}

\centerline{4, place Jussieu, 75 252 Paris Cedex 05, France}

\centerline{dominique.lecomte@upmc.fr}\bigskip

\centerline{$\bullet$ Universit\'e de Picardie, I.U.T. de l'Oise, site de Creil,}

\centerline{13, all\'ee de la fa\"\i encerie, 60 107 Creil, France}
\bigskip\bigskip\bigskip\bigskip\bigskip\bigskip\bigskip

\ufootnote{{\it 2010 Mathematics Subject Classification.}~Primary: 03E15, Secondary: 54H05, 28A05, 26A21}

\ufootnote{{\it Keywords and phrases.}~Products, Potentially, Wadge classes, Borel classes, Reduction}

\ufootnote{{\bf Acknowledgements.}~I would like to thank A. Louveau for some nice remarks made during the talks I gave about the subject in the descriptive set theory seminar of the University Paris 6, especially a simplification of the proof of Lemma 2.4.}

\noindent {\bf Abstract.} Let $\bf\Gamma$ be a Borel class, or a Wadge class of Borel sets, and 
$2\!\leq\! d\!\leq\!\omega$ a cardinal. We study the Borel subsets of ${\mathbb R}^d$ that can be made $\bf\Gamma$ by refining the Polish topology on the real line. These sets are called potentially $\bf\Gamma$. We give a test to recognize potentially $\bf\Gamma$ sets.

\vfill\eject

\section{$\!\!\!\!\!\!$ Introduction}\indent

 The reader should see [K] for the descriptive set theoretic notation used in this paper. The standard way 
of comparing the topological complexity of subsets of $0$-dimensional Polish spaces is the Wadge reducibility quasi-order $\leq_W$. Recall that if $X$ (resp., $Y$) is a $0$-dimensional Polish space and 
$A$ (resp., $B$) a subset of $X$ (resp., $Y$), then
$$(X,A)\leq_W (Y,B)\ \Leftrightarrow\ \exists f\! :\! X\!\rightarrow\! Y\mbox{ continuous such that }
A\! =\! f^{-1}(B).$$
This is a very natural definition since the continuous functions are the morphisms for the topological structure. So the scheme is as follows: 
$$\begin{tabular}{r}
$X$\\ 
\\
\end{tabular}
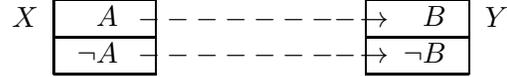
\begin{tabular}{|r|}
\hline
\ \ $A$\ \ \ \ \\
\hline
\ $\neg A$\ \ \ \ \\
\hline
\end{tabular}
\!\!\!\!\!\!\!\begin{tabular}{r}
$--------\!\rightarrow$
\\
$--------\!\rightarrow$
\\
\end{tabular}
\!\!\!\!\!\!\!\!\begin{tabular}{|r|}
\hline
\ \ \ $B$\ \ \\
\hline
\ \ \ $\neg B$\ \ \\
\hline
\end{tabular}
\begin{tabular}{r}
$Y$\\ 
\\
\end{tabular}$$
The ``$0$-dimensional" condition is here to ensure the existence of enough continuous functions (the only continuous functions from $\mathbb{R}$ into $\omega^\omega$ are the constant functions, for example). In the sequel, $\bf\Gamma$ will be a class of Borel subsets of $0$-dimensional Polish spaces. We denote by $\check {\bf\Gamma}\! :=\!\{\neg A\mid A\!\in\! {\bf\Gamma}\}$ the class of complements of elements of 
$\bf\Gamma$. We say that $\bf\Gamma$ is $self$-$dual$ if 
${\bf\Gamma}\! =\!\check {\bf\Gamma}$. We also set 
$\Delta ({\bf\Gamma})\! :=\! {\bf\Gamma}\cap\check {\bf\Gamma}$. Following 4.1 in [Lo-SR2], we give the following definition:

\begin{defi} We say that $\bf\Gamma$ is a $Wadge\ class~of\ Borel\ sets$ if there is a Borel subset $A_0$ of $\omega^\omega$ such that for each $0$-dimensional Polish space $X$, and for each 
$A\!\subseteq\! X$, $A$ is in ${\bf\Gamma}$ if and only if $(X,A)\leq_W (\omega^\omega ,A_0)$. We say that $A_0$ is $\bf\Gamma$-$complete$.\end{defi}
 
 The Wadge hierarchy defined by $\leq_W$, i.e., the inclusion of Wadge classes, is the finest hierarchy of topological complexity in descriptive set theory. The goal of this paper is to study the descriptive complexity of the Borel subsets of products of Polish spaces. More specifically, we are looking for a dichotomy of the following form, quite standard in descriptive set theory: either a set is simple, or it is more complicated than a well-known complicated set. Of course, we have to specify the notions of complexity and comparison we are considering. The two things are actually very much related. The usual notion of comparison between analytic equivalence relations is the Borel reducibility quasi-order $\leq_B$. Recall that if $X$ (resp., $Y$) is a Polish space and $E$ (resp., $F$) an equivalence relation on $X$ (resp., $Y$), then 
$(X,E)\leq_B (Y,F)\ \Leftrightarrow\ \exists f\! :\! X\!\rightarrow\! Y\ \mbox{ Borel such that }\ 
E\! =\! (f\!\times\! f)^{-1}(F)$. 
Note that this makes sense even if $E$ and $F$ are not equivalence relations. The notion of complexity we are considering is a natural invariant for $\leq_B$ in dimension 2. Its definition generalizes Definition 3.3 in [Lo3] to any dimension $d$ making sense in the context of descriptive set theory, and also to any class 
$\bf\Gamma$. So in the sequel $d$ will be a cardinal, and we will have $2\!\leq\! d\!\leq\!\omega$ since 
$2^{\omega_1}$ is not metrizable.

\begin{defi} Let $(X_i)_{i\in d}$ be a sequence of Polish spaces, and $B$ a Borel subset of 
$\Pi_{i\in d}\ X_i$. We say that $B$ is $potentially~in~{\bf\Gamma}$ 
$\big($denoted ${B\!\in\!\hbox{\it pot}({\bf\Gamma})\big)}$ if, for each $i\!\in\! d$, there is a finer 
$0$-dimensional Polish topology $\tau_i$ on $X_i$ such that 
$B\!\in\! {\bf\Gamma}\big(\Pi_{i\in d}\ (X_i,\tau_i)\big)$.\end{defi}

 One should emphasize the fact that the point of this definition is to consider product topologies. Indeed, if $B$ is a Borel subset of a Polish space $X$, then there is a finer Polish topology $\tau$ on $X$ such that $B$ is a clopen subset of $(X,\tau )$ (see 13.1 in [K]). This is not the case in products: if for example 
$\bf\Gamma$ is a non self-dual Wadge class of Borel sets, then there are sets in 
${\bf\Gamma}\big( (\omega^\omega )^2\big)$ that are not $\mbox{pot}(\check {\bf\Gamma})$ (see Theorem 3.3 in [L1]). For example, the diagonal of $\omega^\omega$ is not potentially open.

\vfill\eject

 Note also that since we work up to finer Polish topologies, the ``$0$-dimensional" condition is not a restriction. Indeed, if $X$ is a Polish space, then there is a finer $0$-dimensional Polish topology on $X$ (see 13.5 in [K]). The notion of potential complexity is an invariant for $\leq_B$ in the sense that if $(X,E)\leq_B (Y,F)$ and $F$ is $\mbox{pot}({\bf\Gamma})$, then $E$ is $\mbox{pot}({\bf\Gamma})$ too.\smallskip
 
 The good notion of comparison is not the rectangular version of $\leq_B$. Instead of considering a Borel set $E$ and its complement, we have to consider pairs of disjoint analytic sets. This leads to 
the following notation. Let $(X_i)_{i\in d}$, $(Y_i)_{i\in d}$ be sequences of Polish spaces, and 
$A_0$, $A_1$ (resp., $B_0$, $B_1$) disjoint analytic subsets of $\Pi_{i\in d}\ X_i$ (resp., 
$\Pi_{i\in d}\ Y_i$). Then\smallskip

\leftline{$\big( (X_i)_{i\in d}, A_0, A_1\big)\leq\big( (Y_i)_{i\in d}, B_0, B_1\big)\ \Leftrightarrow\ 
\forall i\!\in\! d\ \ \exists f_i\! :\! X_i\!\rightarrow\! Y_i\mbox{ continuous such that}$}\smallskip

\rightline{$\forall\varepsilon\!\in\! 2\ \ A_\varepsilon\!\subseteq\! 
(\Pi_{i\in d}\ f_i)^{-1}(B_\varepsilon ).$}\smallskip

\noindent So the good scheme of comparison is as follows: 
$$\begin{tabular}{r}
$\Pi_{i\in d}\ X_i$\\ 
\\
\\
\\
\end{tabular}
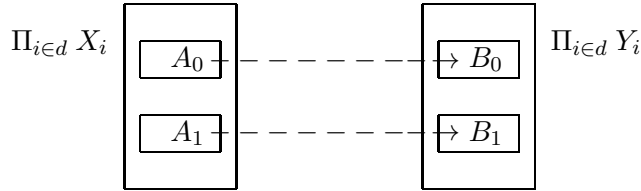
\begin{tabular}{|r|}
\hline
\\ 
\begin{tabular}{|r|}
\hline
\ \ $A_0$\ \\
\hline
\end{tabular}
\\ \\
\begin{tabular}{|r|}
\hline
\ \ $A_1$\ \\
\hline
\end{tabular}
\\ 
\\
\hline
\end{tabular}
\!\!\!\!\!\!\!\!\!\begin{tabular}{r}
$--------\!\rightarrow$
\\
\\
$--------\!\rightarrow$
\\
\end{tabular}
\!\!\!\!\!\!\!\!\!\!\!
\begin{tabular}{|r|}
\hline
\\ 
\begin{tabular}{|r|}
\hline
\ \ $B_0$\ \\
\hline
\end{tabular}
\\
\\
\begin{tabular}{|r|}
\hline
\ \ $B_1$\ \\
\hline
\end{tabular}
\\ 
\\
\hline
\end{tabular}
\begin{tabular}{r}
$\Pi_{i\in d}\ Y_i$\\ 
\\
\\
\\
\end{tabular}$$
The notion of potential complexity was studied in [L1]-[L7] for $d\! =\! 2$ and the non self-dual Borel classes. The main question of this long study was asked by A. Louveau to the author in 1990. A. Louveau wanted to know whether Hurewicz's characterization of $G_\delta$ sets could be extended to 
$\mbox{pot}({\bf\Gamma})$ sets when $\bf\Gamma$ is a Wadge class of Borel sets. The main result of this paper gives a complete and positive answer to this question:
 
\begin{thm} Let $\bf\Gamma$ be a Wadge class of Borel sets, or the class $\borxi$ for some 
$1\!\leq\!\xi\! <\!\omega_1$. Then there are Borel subsets $\mathbb{S}^0$, $\mathbb{S}^1$ of $(d^\omega )^d$ such that for any sequence of Polish spaces $(X_i)_{i\in d}$, and for any disjoint analytic subsets $A_0$, $A_1$ of $\Pi_{i\in d}\ X_i$, exactly one of the following holds:\smallskip

\noindent (a) The set $A_0$ is separable from $A_1$ by a $\mbox{pot}({\bf\Gamma})$ set.\smallskip

\noindent (b) The inequality $\big( (d^\omega )_{i\in d},\mathbb{S}^0,\mathbb{S}^1\big)\leq
\big( (X_i)_{i\in d}, A_0, A_1\big)$ holds.\end{thm}

 It is natural to try to prove Theorem 1.3 since it is a result of continuous reduction, which appears in the very definition of a Wadge class. So it goes beyond a simple generalization. The work in this paper is the continuation of the article [L7], that was announced in [L6]. We generalize the main results of [L7]. The generalization goes in different directions:\smallskip
 
\noindent - It works in any dimension $d$.\smallskip

\noindent - It works for the self-dual Borel classes $\borxi$.\smallskip

\noindent - It works for any Wadge class of Borel sets, which is the hardest part.\smallskip

 We generalize, and also in fact give a new proof of the dimension 1 version of Theorem 1.3 obtained by A. Louveau and J. Saint Raymond (see [Lo-SR1]), which itself was a  generalization of Hurewicz's result. The new proof is without games, and gives a new approach to the study of Wadge classes. Note that A. Louveau and J. Saint Raymond proved that if $\bf\Gamma$ is not self-dual, then the reduction map in (b) can be one-to-one (see Theorem 5.2 in [Lo-SR2]). We will see that there is no injectivity in general in Theorem 1.3. However, G. Debs proved that we can have the $f_i$'s one-to-one when $d\! =\! 2$, 
${\bf\Gamma}\!\in\!\{\bormxi ,\boraxi\}$ and $\xi\!\geq\! 3$. Some injectivity details will be given in the last section.
 
\vfill\eject

 We introduce the following notation and definition in order to specify Theorem 1.3. One can prove that a reduction on the whole product is not possible, for acyclicity reasons (see [L5]-[L7]). We now specify this.  We emphasize the fact that in this paper, there will be a constant identification between $(d^d)^l$ and 
$(d^l)^d$, for $l\!\leq\!\omega$, to avoid as much as possible heavy notation.\bigskip

\noindent\bf Notation.\rm\ If ${\cal X}$ is a set, then $\vec x\! :=\! (x_i)_{i\in d}$ is an arbitrary element of ${\cal X}^d$. If ${\cal T}\!\subseteq\! {\cal X}^d$, then we denote by $G^{\cal T}$ the graph 
with set of vertices ${\cal T}$, and with set of edges 
$\big\{\{\vec x,\vec y\}\!\subseteq\! {\cal T}\mid\vec x\!\not=\!\vec y\ \ 
\hbox{\rm and}\ \ \exists i\!\in\! d\ \ x_i\! =\! y_i\big\}$ 
(see [B] for the basic notions about graphs). So $\vec x\!\not=\!\vec y\!\in\! {\cal T}$ are $G^{\cal T}$-related if they have at least one coordinate in common.

\begin{defi} (a) We say that ${\cal T}$ is $one\mbox{-}sided$ if the following holds:
$$\forall\vec x\!\not=\!\vec y\!\in\! {\cal T}\ \ \forall i\!\not=\! j\!\in\! d\ \ 
(x_i\!\not=\! y_i\vee x_j\!\not=\! y_j).$$
This means that if $\vec x\!\not=\!\vec y\!\in\! {\cal T}$, then they have at most one coordinate in common.\smallskip

\noindent (b) We say that ${\cal T}$ is $almost\ acyclic$ if for every $G^{\cal T}$-cycle 
$(\overrightarrow{x^n})_{n\leq L}$ there are $i\!\in\! d$ and $k\! <\! m\! <\! n\! <\! L$ such that 
$x^k_i\! =\! x^m_i\! =\! x^n_i$. This means that every $G^{\cal T}$-cycle contains a ``flat" subcycle, i.e., a subcycle in a single direction $i\!\in\! d$.\smallskip

\noindent (c) We say that a tree $T$ on $d^d$ is a $tree~with~suitable\ levels$ if the set 
${\cal T}^{l}\! :=\! T\cap (d^d)^l\!\subseteq\! (d^l)^d$ is finite, one-sided and almost acyclic for each integer 
$l$.\end{defi}

 We do not really need the finiteness of the levels, but it makes the proof of Theorem 1.3 much simpler. The following classical property will be crucial in the sequel:
 
\begin{defi} We say that $\bf\Gamma$ has the $separation\ property$ if for each 
$A,B\!\in\! {\bf\Gamma}(\omega^\omega )$ disjoint, there is $C\!\in\!\Delta ({\bf\Gamma})(\omega^\omega )$ separating $A$ from $B$.\end{defi}

 The separation property has been studied in [S] and [vW], where the following is proved:
 
\begin{thm} (Steel-van Wesep) Let ${\bf\Gamma}$ be a non self-dual Wadge class of Borel sets. Then exactly one of the two classes ${\bf\Gamma}$, $\check {\bf\Gamma}$ has the separation property.\end{thm}
 
 We now specify Theorem 1.3.
 
\begin{thm} We can find a tree $T_d$ with suitable levels, together with, for each non self-dual Wadge class of Borel sets ${\bf\Gamma}$,\smallskip

\noindent (1) Some set $\mathbb{S}^d_{\bf\Gamma}\!\in\! {\bf\Gamma}(\lceil T_d\rceil)$ not separable from $\lceil T_d\rceil\!\setminus\!\mathbb{S}^d_{\bf\Gamma}$ by a $\mbox{pot}(\check {\bf\Gamma})$ set.\smallskip

\noindent (2) If moreover ${\bf\Gamma}$ does not have the separation property, and 
${\bf\Gamma}\! =\!\boraxi$ or $\Delta ({\bf\Gamma})$ is a Wadge class, some disjoint sets 
$\mathbb{S}^0_{\bf\Gamma},\mathbb{S}^1_{\bf\Gamma}\!\in\! {\bf\Gamma}(\lceil T_d\rceil)$ not separable by a $\mbox{pot}\big(\Delta ({\bf\Gamma})\big)$ set.\end{thm}
 
\begin{thm} Let $T_d$ be a tree with suitable levels, $\bf\Gamma$ a non self-dual Wadge class of Borel sets, $(X_i)_{i\in d}$ a sequence of Polish spaces, and $A_0$, $A_1$ disjoint analytic subsets of 
$\Pi_{i\in d}\ X_i$.\smallskip

\noindent (1) Let $S\!\in\! {\bf\Gamma}(\lceil T_d\rceil)$ not separable from $\lceil T_d\rceil\!\setminus\! S$ by a $\mbox{pot}(\check {\bf\Gamma})$ set. Then exactly one of the following holds:\smallskip

\noindent (a) The set $A_0$ is separable from $A_1$ by a 
$\mbox{pot}(\check {\bf\Gamma})$ set.\smallskip

\noindent (b) The inequality $\big( (d^\omega )_{i\in d}, S,\lceil T_d\rceil\!\setminus\! S\big)\leq
\big( (X_i)_{i\in d}, A_0, A_1\big)$ holds.

\vfill\eject

\noindent (2) Assume moreover that ${\bf\Gamma}$ does not have the separation property, and that 
${\bf\Gamma}\! =\!\boraxi$ or $\Delta ({\bf\Gamma})$ is a Wadge class. Let 
$S^0\! ,S^1\!\!\in\! {\bf\Gamma}(\lceil T_d\rceil)$ disjoint not separable by a 
$\mbox{pot}\big(\Delta ({\bf\Gamma})\big)$ set. Then exactly one of the following holds:\smallskip

\noindent (a) The set $A_0$ is separable from $A_1$ by a 
$\mbox{pot}\big(\Delta ({\bf\Gamma})\big)$ set.\smallskip

\noindent (b) The inequality $\big( (d^\omega )_{i\in d}, S^0,S^1\big)\leq
\big( (X_i)_{i\in d}, A_0, A_1\big)$ holds.\end{thm}

 We now come back to the new approach to the study of Wadge classes mentioned earlier. There are a lot of dichotomy results in descriptive set theory about equivalence relations, quasi-orders or even arbitrary Borel or analytic sets. So it is natural to ask for common points to these dichotomies. B. Miller's recent work goes in this direction. He proved many known dichotomies without effective descriptive set theory, using variants of the Kechris-Solecki-Todor\v cevi\'c dichotomy about analytic graphs (see [K-S-T]). Here we want to point out another common point, of effective nature. In these dichotomies, the first possibility of the dichotomy is equivalent to the emptyness of some $\Ana$ set. For example, in the Kechris-Solecki-Todor\v cevi\'c dichotomy, the $\Ana$ set is the complement of the union of the $\Borel$ subsets discrete with respect to the $\Ana$ graph considered. We prove a strengthening of Theorem 1.8 in which such a $\Ana$ set appears. We will state in Case (1), unformally. Before that, we need the following notation.\bigskip

\noindent\bf Notation.\rm ~Let $X$ be a recursively presented Polish space. We denote by ${\it\Delta}_{X}$ the topology on $X$ generated by $\Borel (X)$. This topology is Polish (see (iii) $\Rightarrow$ (i) in the proof of Theorem 3.4 in [Lo3]). The topology $\tau_{1}$ on $(\omega^\omega )^d$ will be the product topology ${\it\Delta}_{\omega^\omega}^d$.

\begin{thm} Let $T_d$ be a tree with $\Borel$ suitable levels, $\bf\Gamma$ a non self-dual Wadge class of Borel sets with a $\Borel$ code, $A_0,A_1$ disjoint $\Ana$ subsets of $(\omega^\omega )^d$, and 
$S\!\in\! {\bf\Gamma}(\lceil T_d\rceil)$ not separable from $\lceil T_d\rceil\!\setminus\! S$ by a 
$\mbox{pot}(\check {\bf\Gamma})$ set. Then there is a $\Ana$ subset $R$ of $(\omega^\omega )^d$ such that the following are equivalent:\smallskip

\noindent (a) The set $A_0$ is not separable from $A_1$ by a $\mbox{pot}(\check {\bf\Gamma})$ set.\smallskip

\noindent (b) The set $A_0$ is not separable from $A_1$ by a $\Borel\cap\mbox{pot}(\check {\bf\Gamma})$ set.\smallskip

\noindent (c) The set $A_0$ is not separable from $A_1$ by a $\check {\bf\Gamma}(\tau_1)$ set.\smallskip

\noindent (d) $R\!\not=\!\emptyset$.\smallskip

\noindent (e) The inequality $\big( (d^\omega )_{i\in d}, S,\lceil T_d\rceil\!\setminus\! S\big)\leq
\big( (\omega^\omega )_{i\in d}, A_0, A_1\big)$ holds.\end{thm}

This $\Ana$ set $R$ is build with topologies based on $\tau_1$. The use of the $\Ana$ set $R$ is the new approach to the study of Wadge classes.\bigskip

 We first prove Theorems 1.7 and 1.8 for the Borel classes, self-dual or not. Then we consider the case of the Wadge classes. In Section 2, we start proving Theorem 1.7. We construct a concrete example of a tree with suitable levels, and we give a general condition to get some complicated sets as in the statement of Theorem 1.7. We actually reduce the problem to a problem in dimension one. In Section 3, we prove Theorem 1.7 for the Borel classes. In Section 4, we prove Theorem 1.8 for the Borel classes, using some tools of effective descriptive set theory and the representation theorem of Borel sets proved in [D-SR]. In Section 5, we prove  Theorem 1.7, using the description of Wadge classes in [Lo-SR2]. In Section 6, we prove Theorems 1.3, 1.8 and 1.9. Finally, in Section 7, we give some injectivity complements.

\vfill\eject

\section{$\!\!\!\!\!\!$ A general condition to get some complicated sets}\indent

 We now build an example of a tree with suitable levels. This tree has to be small enough since we cannot have a reduction on the whole product. But as the same time it has to be big enough to ensure the existence of complicated sets, as in Theorem 1.7.\bigskip
 
\noindent\bf Notation.\rm ~Let $b\! :\!\omega\!\rightarrow\!\omega^2$ be the natural bijection. 
More precisely, we set, for $l\!\in\!\omega$,
$$M(l)\! :=\!\hbox{\rm max}\{ m\!\in\!\omega\mid\Sigma_{k\leq m}~k\!\leq\! l\}.$$
Then we define $b(l)\! =\!\big((l)_{0},(l)_{1}\big)\! :=\!\big(M(l)\! -\! l\! +\! 
(\Sigma_{k\leq M(l)}~k),l\! -\! (\Sigma_{k\leq M(l)}~k)\big)$. One can check that 
$<\! n,p\! >:=\! b^{-1}(n,p)\! =\! (\Sigma_{k\leq n+p}~k)\! +\! p$. More concretely, we get 
$$b[\omega ]\! =\!\{(0,0),(1,0),(0,1),(2,0),(1,1),(0,2),\ldots\}.$$

In the introduction, we mentionned the idenfication between $(d^l)^d$ and $(d^d)^l$, for $l\!\leq\!\omega$. More specifically, the bijection we use is given by 
$\vec\alpha\!\mapsto\!\Big(\big(\alpha_i(j)\big)_{i\in d}\Big)_{j\in l}$.

\begin{defi} We say that $E\!\subseteq\!\bigcup_{l\in\omega}~(d^l)^d$ is an $effective~frame$ if\medskip

\noindent (a) $\forall l\!\in\!\omega~\exists !(s^{i}_{l})_{i\in d}\!\in\! E\!\cap\! (d^l)^d$.\medskip

\noindent (b) $\forall p,q,r\!\in\!\omega~\forall t\!\in\! d^{<\omega}~\exists N\!\in\!\omega~
(s^{i}_{q}it0^N)_{i\in d}\!\in\! E$, $(|s^{0}_{q}0t0^N|\! -\! 1)_{0}\! =\! p$ and 
$\big( (|s^{0}_{q}0t0^N|\! -\! 1)_1\big)_0\! =\! r$.\medskip

\noindent (c) $\forall l\! >\! 0~\exists q\! <\! l~\exists t\!\in\! d^{<\omega}\ \forall i\!\in\! d~~
s^{i}_{l}\! =\! s^{i}_{q}it$.\medskip

\noindent (d) The map $l\!\mapsto\! (s^{i}_{l})_{i\in d}$ can be coded by a recursive map from $\omega$ into $\omega^d$.\medskip

 We will call $T_d$ the $tree~on~d^d~associated~with~an~effective~frame$  
$E\! =\!\{(s^{i}_{l})_{i\in d}\mid l\!\in\!\omega\}$: 
$$T_d\! :=\!\big\{\vec s\!\in\! (d^d)^{<\omega}\mid (\forall i\!\in\! d\ \ s_i\! =\!\emptyset )\vee
\big(\exists l\!\in\!\omega ~\exists t\!\in\! d^{<\omega}~\forall i\!\in\! d\ \  s_i\! =\! s^{i}_{l}it\wedge
\forall n\! <\!\vert s_0\vert\ s_0(n)\!\leq\! n\big)\big\}.$$
\end{defi}

 The uniqueness condition in (a) and Condition (c) ensure that $T_d$ is small enough, and also the almost acyclicity. The definition of $T_d$ ensures that $T_d$ has finite levels. Note that 
${\cal T}^l\! =\! T_d\cap (d^d)^l$ can be coded by a $\Bormone$ subset of $(\omega^\omega )^l$ when 
$d\! =\!\omega$. The existence condition in (a) and Condition (b) ensure that $T_d$ is big enough. More specifically, if $(X,\tau )$ is a Polish space and $\sigma$ a finer Polish topology on $X$, then there is a dense $G_{\delta}$ subset of $(X,\tau )$ on which $\tau$ and $\sigma$ coincide. The first part of Condition (b) ensures the possibility to get inside products of dense $G_{\delta}$ sets. The examples in Theorem 1.7 are build using the examples in [Lo-SR1] and [Lo-SR2]. Conditions on verticals are involved, and the second part of Condition (b) gives a control on the choice of verticals. The very last part of Condition (b) is not necessary to get Theorem 1.7 for the Borel classes, but is useful to get Theorem 1.7 for the Wadge classes of Borel sets. Definition 2.1 strengthens Definition 3.1 in [L7], with this very last part of Condition (b), with Condition (d) (ensuring the regularity of the levels of the tree), and also with the last part of the definition of the tree (ensuring the finiteness of the levels of the tree).
 
\begin{prop} The tree $T_d$ associated with an effective~frame is a tree with $\Borel$ suitable levels. In particular, $\lceil T_d\rceil$ is compact.\end{prop}

\noindent\bf Proof.\rm ~Let $l\!\in\!\omega$. Let us prove that ${\cal T}^{l}$ is $\Borel$ and finite. We argue by induction on $l$. The result is clear for $l\!\leq\! 1$ since ${\cal T}^{0}\! =\!\{\vec\emptyset\}$ and 
${\cal T}^{1}\! =\!\{ (i)_{i\in d}\}$. If $l\!\geq\! 1$ and $\vec s\!\in\! (d^d)^{<\omega}$, then 
$$\vec s\!\in\! {\cal T}^{l}\ \Leftrightarrow\ \vert s_0\vert\! =\! l\wedge
\exists q\! <\! l\ \ \exists t\!\in\! d^{<\omega}\ \ \forall i\!\in\! d\ \ s_i\! =\! s^i_qit\wedge
\forall n\! <\! l\ \ s_0(n)\!\leq\! n.$$
But there are only finitely many possibilities for $t$ since $s_0(n)\!\leq\! n$ for each $n\! <\! l$, which implies that $t(m)\! \leq\! q\! +\! 1\! +\! m\! <\! l\! +\! 1\! +\! l$ if $m\! <\!\vert t\vert$. This implies that 
${\cal T}^{l}$ is $\Borel$ and finite.\bigskip

\noindent $\bullet$ Let $\tilde T_d$ be the tree generated by the effective frame:
$$\tilde T_d\! :=\!\big\{\vec s\!\in\! (d^d)^{<\omega}\mid (\forall i\!\in\! d\ \ s_i\! =\!\emptyset )\vee
\big(\exists l\!\in\!\omega ~\exists t\!\in\! d^{<\omega}~\forall i\!\in\! d\ \  s_i\! =\! s^{i}_{l}it\big)\big\}.$$
As $T_d\!\subseteq\!\tilde T_d$, we get, with obvious notation, ${\cal T}^l\!\subseteq\! {\tilde {\cal T}}^l$ for each integer $l$. So it is enough to prove that ${\tilde {\cal T}}^l$ is one-sided and almost acyclic since these properties are hereditary.\bigskip

\noindent $\bullet$ Let us prove that ${\tilde {\cal T}}^l$ is almost acyclic. We argue by induction on $l$. The result is clear for $l\!\leq\! 1$. So assume that $l\!\geq\! 1$. We set, for $j\!\in\! d$, 
$$C_j\! :=\!\big\{ (s^{i}_{q}it)_{i\in d}\!\in\! {\tilde {\cal T}}^{l+1}\mid t\!\not=\!\emptyset\wedge 
t(\vert t\vert\! -\! 1)\! =\! j\big\}.$$ 
We have ${\tilde {\cal T}}^{l+1}\! =\!\{ (s^{i}_li)_{i\in d}\}\cup\bigcup_{j\in d}\ C_j$, and this union is disjoint.\bigskip

 The restriction of $G^{{\tilde {\cal T}}^{l+1}}$ to each $C_j$ is isomorphic to $G^{{\tilde {\cal T}}^{l}}$. The other possible $G^{{\tilde {\cal T}}^{l+1}}$-edges are between $(s^{i}_li)_{i\in d}$ and some vertices in some $C_j$'s. If a $G^{{\tilde {\cal T}}^{l+1}}$-cycle exists, we may assume that it involves only 
$(s^{i}_li)_{i\in d}$ and members of some fixed $C_{j_0}$. But if 
$\vec s\!\in\! C_{j_0}$ is $G^{{\tilde {\cal T}}^{l+1}}$-related to $(s^{i}_li)_{i\in d}$, then we must have 
$s^{j_0}_lj_0\! =\! s_{j_0}$. This implies the existence of $k\! <\! m\! <\! n$ showing that 
${\tilde {\cal T}}^{l+1}$ is almost acyclic.\bigskip

\noindent $\bullet$ Now assume that $\vec x\!\not=\!\vec y\!\in\!{\tilde {\cal T}}^l$, $i,j\!\in\! d$, 
$x_i\! =\! y_i$ and $x_j\! =\! y_j$. Then we can write 
$\vec x\! =\! (s^{i}_{q}it)_{i\in d}$ and $\vec y\! =\! (s^{i}_{q'}it')_{i\in d}$ since 
$\vec x\!\not=\!\vec y$. As $x_i\! =\! y_i$, the reverses $t^{-1}$ and $(t')^{-1}$ of $t$ and $t'$ are compatible. If $t\! =\! t'$, then 
$q\! =\!\vert s^{i}_q\vert\! =\! l\! -\! 1\! -\!\vert t\vert\! =\! l\! -\! 1\! -\!\vert t'\vert\! =\!\vert s^{i}_{q'}\vert\! =\! q'$ and 
$\vec x\! =\!\vec y$, which is absurd. Thus $t\!\not=\! t'$, for example 
$\vert t'\vert\! <\!\vert t\vert$, and $t^{-1}(\vert t'\vert )\! =\! i$. This proves that $i\! =\! j$ and ${\tilde {\cal T}}^l$ is one-sided.\bigskip

\noindent $\bullet$ Let $\pi_l\! :\! {\cal T}^{l+1}\!\rightarrow\! d^d$ defined by 
$\pi_l(\vec s)\! :=\!\big( s_i(l)\big)_{i\in d}$. As ${\cal T}^{l+1}$ is finite, the range $c_l$ of $\pi_l$ is also finite. Thus $\lceil T_d\rceil$ is compact since $\lceil T_d\rceil\!\subseteq\!\Pi_{l\in\omega}\ c_l$.
\hfill{$\square$}\bigskip

 We now give an example of an effective frame.\bigskip

\noindent\bf Notation.\rm ~Let $b_d\! :\!\omega\!\rightarrow\! d^{<\omega}$ be the natural bijection. More specifically,\medskip

\noindent $\bullet$ If $d\! <\!\omega$, then $b_d (0)\! :=\!\emptyset$ is the sequence of length $0$, 
$b_d (1)\! :=\! 0$, ..., $b_d (d)\! :=\! d\! -\! 1$ are the sequences of length $1$ in the lexicographical ordering, and so on.\medskip 

\noindent $\bullet$ If $d\! =\!\omega$, then let $(p_n)_{n\in\omega}$ be the sequence of prime numbers, and ${\cal I}\! :\!\omega^{<\omega}\!\rightarrow\!\omega$ defined by ${\cal I}(\emptyset )\! :=\! 1$, and 
${\cal I}(s)\! :=\! p_0^{s(0)+1}...p_{\vert s\vert -1}^{s(\vert s\vert -1)+1}$ if $s\!\not=\!\emptyset$. Note that 
${\cal I}$ is one-to-one, so that there is an increasing bijection 
$\imath\! :\!\hbox{\rm Seq}\! :=\! {\cal I}[\omega^{<\omega}]\!\rightarrow\!\omega$. We set 
$b_\omega\! :=\! (\imath\circ {\cal I})^{-1}\! :\!\omega\!\rightarrow\!\omega^{<\omega}$.\bigskip
 
 Note that $|b_d (n)|\!\leq\! n$ if $n\!\in\!\omega$. Indeed, this is clear if $d\! <\!\omega$. If $d\! =\!\omega$, then
$${\cal I}\big( b_\omega (n)\vert 0\big)\! <\! {\cal I}\big( b_\omega (n)\vert 1\big)\! <\! ...\! <\! 
{\cal I}\big( b_\omega (n)\big)\hbox{\rm ,}$$
so that $(\imath\circ {\cal I})\big( b_\omega (n)\vert 0\big)\! <\! (\imath\circ {\cal I})
\big( b_\omega (n)\vert 1\big)\! <\! ...\! <\! (\imath\circ {\cal I})\big( b_\omega (n)\big)\! =\! n$. This implies that $\vert b_\omega (n)\vert\!\leq\! n$.

\begin{lem} There is a concrete example of an effective~frame.\end{lem}

\noindent\bf Proof.\rm ~Fix $i\!\in\! d$. We set $s^{i}_{0}\! =\!\emptyset$, and 
$s^{i}_{l+1}\! :=\! s^{i}_{(((l)_{1})_1)_0}~i~b_d\bigg(\Big(\big( (l)_{1}\big)_{1}\Big)_1\bigg)~
0^{l-(((l)_{1})_1)_0-|b_d((((l)_{1})_{1})_1)|}$. Note that 
$(l)_{0}\! +\! (l)_{1}\! =\! M(l)\!\leq\!\Sigma_{k\leq M(l)}~k\!\leq\! l\hbox{\rm ,}$ 
so that $s^{i}_{l}$ is well defined and $|s^{i}_{l}|\! =\! l$, by 
induction on $l$. It remains to check that Condition (b) in the 
definition of an effective~frame is fullfilled. Set $n\! :=\! b_d^{-1}(t)$, 
$s\! := \big<\ r,<q,n>\! \big>$ and $l\! :=<\! p,s\! >$. It remains to put $N\!:=\! l\! -\! q\! -\! |t|$: 
${(s^{i}_{q}it0^N)_{i\in d}\! =\! (s^{i}_{l+1})_{i\in d}}$.\hfill{$\square$}\bigskip

 The previous lemma is essentially identical to Lemma 3.3 in [L7]. Now we come to the lemma crucial for the proof of Theorem 1.7. It strengthens Lemma 3.4 in [L7], even if the proof is essentially the same.\bigskip

\noindent\bf Notation.\rm ~If $s\!\in\!\omega^{<\omega}$ and $q\!\leq\!\vert s\vert$, then $s\! -\! s\vert q$ is defined by $s\! =\! (s\vert q)(s\! -\! s\vert q)$. We extend this definition to the case $s\!\in\!\omega^\omega$ when $q\! <\!\omega$. In particular, we denote $s^*\! :=\! s\! -\! s\vert 1$ when 
$\emptyset\!\not=\! s\!\in\!\omega^{\leq\omega}$. If $\emptyset\!\not=\! s\!\in\!\omega^{<\omega}$, then we define $s^-\! :=\! s\vert (\vert s\vert\! -\! 1)$.\bigskip

\noindent $\bullet$ We now define 
$p\! :\!\omega^{<\omega}\!\setminus\!\{\emptyset\}\!\rightarrow\!\omega$. The definition of $p(s)$ is by induction on $|s|$:
$$p(s)\! :=\!\left\{\!\!\!\!\!\!\!
\begin{array}{ll} 
& s(0)\mbox{ if }|s|\! =\! 1,\cr & \cr
& <\! p(s^-),s(|s|\! -\! 1)\! >~\mbox{otherwise.}
\end{array}
\right.$$
Note that $p_{\vert\omega^n}\! :\!\omega^n\!\rightarrow\!\omega$ is a bijection, 
for each $n\!\geq\! 1$.\bigskip

\noindent $\bullet$ Let $l\!\leq\!\omega$ be an ordinal. The map 
$\Delta\! :\! d^{l}\!\times\! d^{l}\!\rightarrow\! 2^{l}$ is the symmetric difference. So, for $m\!\in\! l$, 
$$(s\Delta t)(m)\! :=\! \Delta (s,t)(m)\! =\! 1~\Leftrightarrow ~s(m)\!\not=\! t(m).$$
$\bullet$ By convention, $\omega\! -\! 1\! :=\!\omega$.

\begin{lem} Let $T_d$ be the tree associated with an effective~frame and, for each $i\!\in\! d$, $G_i$ a dense $G_{\delta}$ subset of $\Pi_i''\lceil T_d\rceil$. Then there are $\alpha_{0}\!\in\! G_0$ and 
$F\! :\! 2^\omega\!\rightarrow\!\Pi_{0<i<d}\ G_i$ continuous such that, for 
$\alpha\!\in\! 2^\omega$,\smallskip

\noindent (a) $\big(\alpha_{0},F(\alpha )\big)\!\in\!\lceil T_d\rceil$.\smallskip

\noindent (b) For each $s\!\in\!\omega^{<\omega}$, and each $m\!\in\!\omega$,\smallskip

\noindent (i) $\alpha\big( p(sm)\big)\! =\! 1~\Rightarrow ~\exists m'\!\in\!\omega ~~
\big(\alpha_{0}\Delta F_0(\alpha )\big)\big( p(sm')\! +\! 1\big)\! =\! 1$.\smallskip

\noindent (ii) $\big(\alpha_{0}\Delta F_0(\alpha )\big)\big( p(sm)\! +\! 1\big)\! =\! 1~\Rightarrow ~
\exists m'\!\in\!\omega ~~\alpha\big( p(sm')\big)\! =\! 1$.\smallskip

 Moreover, there is an increasing bijection 
$$B_\alpha\! :\!\{ m\!\in\!\omega\mid\alpha (m)\! =\! 1\}\!\rightarrow\!
\{ m'\!\in\!\omega\mid \big(\alpha_{0}\Delta F_0(\alpha )\big)(m'\! +\! 1)\! =\! 1\}$$
such that $(m)_0\! =\!\big( B_\alpha (m)\big)_0$ and 
$\big( (m)_1\big)_0\! =\!\Big(\big( B_\alpha (m)\big)_1\Big)_0$ if $\alpha (m)\! =\! 1$.\end{lem}

\noindent\bf Proof.\rm ~Let $(O^i_{q})_{q\in\omega}$ be a non-increasing sequence of dense open subsets of $\Pi_i''\lceil T_d\rceil$ whose intersection is $G_i$. We construct finite approximations of 
$\alpha_{0}$ and $F$. The idea is to linearize the binary tree $2^{<\omega}$. So we will use the bijection $b_2$ defined before Lemma 2.3. To construct $F(\alpha )$ we have to imagine, for each length $l$, the different possibilities for $\alpha\vert l$. More precisely, we construct a map 
$l\! :\! 2^{<\omega}\!\rightarrow\!\omega\!\setminus\!\{ 0\}$ satisfying the following conditions:
$$\begin{array}{ll}
& (1)~\forall t\!\in\! 2^{<\omega}\ \ \forall i\!\in\! d\ \ (i\!\leq\!\vert t\vert\ \Rightarrow\ 
\emptyset\!\not=\! N_{s^i_{l(t)}}\cap\Pi_i''\lceil T_d\rceil\!\subseteq\! O^i_{\vert t\vert})\mbox{,}\cr\cr
& (2)~\exists v_\emptyset\!\in\! d^{<\omega}\ \ \forall i\!\in\! d\ \ s^i_{l(\emptyset )}\! =\! iv_\emptyset
\mbox{,}\cr\cr
& (3)~\forall t\!\in\! 2^{<\omega}\ \ \forall\varepsilon\!\in\! 2\ \ 
\exists v_{t\varepsilon}\!\in\! d^{<\omega}\ \ \forall i\!\in\! d\ \ 
s^i_{l(t\varepsilon )}\! =\! s^i_{l(t)}(i\!\cdot\!\varepsilon )v_{t\varepsilon}\mbox{,}\cr\cr
& (4)~\forall r\!\in\!\omega\ \ s^0_{l(b_2(r))}0\!\subseteq\!s^0_{l(b_2(r+1))}\wedge 
\forall t\!\in\! 2^{<\omega}\ \ \forall n\! <\! l(t)\ \ s^0_{l(t)}(n)\!\leq\! n\mbox{,}\cr\cr
& (5)~\forall t\!\in\! 2^{<\omega}\ \ \big( l(t)\! -\! 1\big)_0\! =\! (\vert t\vert)_0\wedge 
\Big(\big( l(t)\! -\! 1\big)_1\Big)_0\! =\!\big( (\vert t\vert)_1\big)_0.
\end{array}$$ 
$\bullet$ Assume that this construction is done. As $s^0_{l(0^q)}\!\subset_{\not=}\! s^0_{l(0^{q+1})}$ for each integer $q$, we can define $\alpha_{0}\! :=\!\mbox{sup}_{q\in\omega}~s^0_{l(0^q)}$. Similarly, 
as $s^{i+1}_{l(\alpha\vert q)}\!\subset_{\not=}\! s^{i+1}_{l(\alpha\vert (q+1))}$, we can define, for 
$\alpha\!\in\! 2^\omega$ and $i\! <\! d\! -\! 1$,
$$F_i(\alpha )\! :=\!\mbox{sup}_{q\in\omega}~s^{i+1}_{l(\alpha\vert q)}\mbox{,}$$
and $F$ is continuous.\bigskip

\noindent (a) Fix $q\!\in\!\omega$. We have to see that $\big(\alpha_{0}, F(\alpha )\big)\vert q\!\in\! T_d$. Note first that $l(t)\!\geq\!\vert t\vert$ since $l(t\varepsilon )\! >\! l(t)$. Then note that 
$s^0_{l(t)}\!\subseteq\!\alpha_0$ since 
$s^0_{l(0^{\vert t\vert})}\!\subseteq\!s^0_{l(t)}\!\subseteq\!s^0_{l(0^{\vert t\vert +1})}$. Thus 
$\big(\alpha_{0},F(\alpha )\big)\vert l(\alpha\vert q)\! =\! (s^i_{l(\alpha\vert q)})_{i\in d}\!\in\! E$. This implies that $\big(\alpha_{0},F(\alpha )\big)\vert l(\alpha\vert q)\!\in\! T_d$ since 
$s^0_{l(\alpha\vert q)}(n)\!\leq\! n$ if $n\! <\! l(\alpha\vert q)$. We are done since $l(\alpha\vert q)\!\geq\! q$.\bigskip

 Moreover, 
$\alpha_0\!\in\!\bigcap_{q\in\omega}\ N_{s^0_{l(0^q)}}\cap\Pi_0''\lceil T_d\rceil\!\subseteq\!\bigcap_{q\in\omega}\ O^0_q\! =\! G_0$. Similarly, 
$$F_i(\alpha )\!\in\!\bigcap_{q\in\omega}\ N_{s^{i+1}_{l(\alpha\vert q)}}\cap\Pi_{i+1}''\lceil T_d\rceil
\!\subseteq\!\bigcap_{q\geq i+1}\ O^{i+1}_q\! =\! G_{i+1}.$$
(b).(i) We set $t\!:=\!\alpha\vert p(sm)$, so that 
$s^1_{l(t)}\ 1\!\subseteq\! s^1_{l(t1)}\! =\! s^1_{l(\alpha\vert (p(sm)+1))}\!\subseteq\! F_0(\alpha )$. As 
$\big( l(t)\! -\! 1\big)_{0}\! =\! p(s)$ (or $(m)_0$ if $s\! =\!\emptyset$), there is $m'$ with 
$l(t)\! =\! p(sm')\! +\! 1$ (or $l(t)\! =\! m'\! +\! 1$ and $(m')_0\! =\! (m)_0$ if $s\! =\!\emptyset$). But 
$s^0_{l(t)}\ 0\!\subseteq\! s^0_{l(\alpha\vert (p(sm)+1))}\!\subseteq\!\alpha_0$, so that 
$\alpha_{0}\big( l(t)\big)\!\not=\! F_0(\alpha )\big( l(t)\big)$.\bigskip

\noindent (ii) First notice that the only coordinates where $\alpha_{0}$ and $F_0(\alpha )$ can differ are 
$0$ and the $l(\alpha\vert q)$'s. Therefore there is an integer $q$ with 
$p(sm)\! +\! 1\! =\! l(\alpha\vert q)$. In particular, $(q)_{0}\! =\!\big( l(\alpha\vert q)\! -\! 1\big)_{0}\! =\! p(s)$ (or $(m)_0$ if $s\! =\!\emptyset$). Thus there is $m'$ with $q\! =\! p(sm')$ (or $q\! =\! m'$ and 
$(m')_0\! =\! (m)_0$ if $s\! =\!\emptyset$). We have 
$\alpha_{0}\big( l(\alpha\vert q)\big)\! =\! s^0_{l(\alpha\vert (q+1))}\big( l(\alpha\vert q)\big)\! =\! 0\!\not=\! F_0(\alpha )\big( l(\alpha\vert q)\big)\! =\! s^1_{l(\alpha\vert (q+1))}\big( l(\alpha\vert q)\big)\! =\!
\alpha (q)$. So $\alpha(q)\! =\! 1$ and $\alpha\big( p(sm')\big)\! =\! 1$.\bigskip

 Now it is clear that the formula $B_\alpha (m)\! :=\! l(\alpha\vert m)\! -\! 1$ defines the bijection we are looking for.\bigskip

\noindent $\bullet$ So let us prove that the construction is possible. We construct $l(t)$ by induction on 
$b_2^{-1}(t)$.\bigskip

 As $(i0^\infty )_{i\in d}\!\in\!\lceil T_d\rceil$, $0^\infty\!\in\!\Pi_0''\lceil T_d\rceil$ and 
$O^0_0$ is not empty. Thus there is $u^0_0\!\in\! d^{<\omega}\!\setminus\!\{\emptyset\}$ such that 
$\emptyset\!\not=\! N_{u^0_0}\cap\Pi_0''\lceil T_d\rceil\!\subseteq\! O^0_0$. Choose 
$\beta_0\!\in\!N_{u^0_0}\cap\Pi_0''\lceil T_d\rceil$, and $\vec\alpha\!\in\!\lceil T_d\rceil$ such that 
$\alpha_0\! =\!\beta_0$. Then $\vec\alpha\vert\vert u^0_0\vert\!\in\! T_d$ and 
$u_0^0(n)\!\leq\! n$ for each $n\! <\!\vert u^0_0\vert$. Note that $u_0^0(0)\! =\! 0$ and 
$(u_0^0\! -\! u_0^0\vert 1)(n)\! =\! u_0^0(n\! +\! 1)\!\leq\! 1\! +\! n$ for each $n\! <\!\vert u_0^0\vert\! -\! 1$. We choose $N_{\emptyset}\!\in\!\omega$ with 
$\big( i~(u_0^0\! -\! u_0^0\vert 1)~0^{N_{\emptyset}}\big)_{i\in d}\!\in\! E$, 
$(|0~(u_0^0\! -\! u_0^0\vert 1)~0^{N_{\emptyset}}|\! -\! 1)_{0}\! =\! (0)_0$ and 
$\big( (|0~(u_0^0\! -\! u_0^0\vert 1)~0^{N_{\emptyset}}|\! -\! 1)_1\big)_0\! =\!\big( (0)_1\big)_0$. We put  
$v_{\emptyset}\! :=\! (u_0^0\! -\! u_0^0\vert 1)~0^{N_{\emptyset}}$ and 
$l(\emptyset )\! :=\! |0~(u_0^0\! -\! u_0^0\vert 1)~0^{N_{\emptyset}}|$.\bigskip 

 As $(iv_\emptyset 0^\infty )_{i\in d}\!\in\!\lceil T_d\rceil$, $N_{0v_\emptyset 0}\cap\Pi_0''\lceil T_d\rceil$ is a nonempty open subset of $\Pi_0''\lceil T_d\rceil$. Thus there is $u^0_1\!\in\! d^{<\omega}$ such that 
$\emptyset\!\not=\! N_{0v_\emptyset 0u^0_1}\cap\Pi_0''\lceil T_d\rceil\!\subseteq\! O^0_1$. As before we see that $u_1^0(n)\!\leq\! 1\! +\!\vert v_\emptyset\vert\! +\! 1\! +\! n$ for each $n\! <\!\vert u_1^0\vert$. This implies that $(iv_\emptyset 0u^0_10^\infty )_{i\in d}\!\in\!\lceil T_d\rceil$. Thus 
$N_{1v_\emptyset 0u^0_1}\cap\Pi_1''\lceil T_d\rceil$ is a nonempty open subset of 
$\Pi_1''\lceil T_d\rceil$. So there is $u^1_1\!\in\! d^{<\omega}$ such that 
$\emptyset\!\not=\! N_{1v_\emptyset 0u^0_1u^1_1}\cap\Pi_1''\lceil T_d\rceil\!\subseteq\! O^1_1$.  Choose 
$\beta_1\!\in\!N_{1v_\emptyset 0u^0_1u^1_1}\cap\Pi_1''\lceil T_d\rceil$, and 
$\vec\gamma\!\in\!\lceil T_d\rceil$ such that $\gamma_1\! =\!\beta_1$. Then 
$\vec\gamma\vert\vert 1v_\emptyset 0u^0_1u^1_1\vert\!\in\! T_d$ and 
$\gamma_0(n)\!\leq\! n$ for each $n\! <\!\vert 1v_\emptyset 0u^0_1u^1_1\vert$. This implies that 
$\gamma_0(\vert 1v_\emptyset 0u^0_1\vert\! +\! n)\!\leq\!\vert 1v_\emptyset 0u^0_1\vert\! +\! n$ for each 
$n\! <\!\vert u^1_1\vert$. But $u^1_1(n)$ is either $1$, or 
$\gamma_0(\vert 1v_\emptyset 0u^0_1\vert\! +\! n)$. Thus 
$u^1_1(n)\!\leq\!\vert 1v_\emptyset 0u^0_1\vert\! +\! n$ if $n\! <\!\vert u^1_1\vert$. We choose 
$N_{0}\!\in\!\omega$ such that $(s^i_{l(\emptyset )}~0u^0_1u^1_1~0^{N_{0}})_{i\in d}\!\in\! E$, 
$\big( l(\emptyset )\! +\!\vert u^0_1u^1_1\vert\! +\! N_0\big)_{0}\! =\! (1)_0$ and 
$\Big(\big( l(\emptyset )\! +\!\vert u^0_1u^1_1\vert\! +\! N_0\big)_1\Big)_0\! =\!\big( (1)_1\big)_0$. We put 
$v_{0}\! :=\! u^0_1u^1_1~0^{N_{0}}$ and $l(0)\! :=\! l(\emptyset )\! +\! 1\! +\!\vert v_0\vert$.\bigskip

 Assume that $\big(l(t)\big)_{b_2^{-1}(t)\leq r}$ satisfying (1)-(5) have been constructed, which is the case for ${r\! =\! 1}$. Fix $t\!\in\! 2^{<\omega}$ and $\varepsilon\!\in\! 2$ such that  
$b_2(r\! +\! 1)\! =\! t\varepsilon$, with $r\!\geq\! 1$. Note that $b_2^{-1}(t)\! <\! r$, so that 
$l(t)\! <\! l\big( b_2(r)\big)$, by induction assumption.\bigskip

 As $N_{s^0_{l(b_2(r))}}\cap\Pi_0''\lceil T_d\rceil$ is nonempty, 
$N_{s^0_{l(b_2(r))}0}\cap\Pi_0''\lceil T_d\rceil$ is nonempty too. Thus there is 
$u^0_{\vert t\vert +1}$ in $d^{<\omega}$ such that $\emptyset\!\not=\! 
N_{s^0_{l(b_2(r))}0u^0_{\vert t\vert +1}}\cap\Pi_0''\lceil T_d\rceil\!\subseteq\! O^0_{\vert t\vert +1}$. As before we see that $u_{\vert t\vert +1}^0(n)\!\leq\! l\big( b_2(r)\big)\! +\! 1\! +\! n$ for each 
$n\! <\!\vert u_{\vert t\vert +1}^0\vert$. Arguing as in the case $r\! =\! 1$, we prove, for each 
$1\!\leq\! i\!\leq\!\vert t\vert\! +\! 1$, the existence of $u^i_{\vert t\vert +1}\!\in\! d^{<\omega}$ such that 
$\emptyset\!\not=\! N_{s^i_{l(t)}(i\cdot\varepsilon )(s^0_{l(b_2(r))}\! -\! 
s^0_{l(b_2(r))}\vert ( l(t)\! +\! 1))0u^0_{\vert t\vert +1}...u^i_{\vert t\vert +1}}\cap\Pi_i''\lceil T_d\rceil
\!\subseteq\! O^i_{\vert t\vert +1}$ and $u_{\vert t\vert +1}^i(n)\!\leq\! l\big( b_2(r)\big)\! +\! 1\! +\!
\vert u^0_{\vert t\vert +1}...u^{i-1}_{\vert t\vert +1}\vert\! +\! n$ for each 
$n\! <\!\vert u_{\vert t\vert +1}^i\vert$ ($u^i_{\vert t\vert +1}(n)$ can be $i$, in which case we use the fact that $l(t)\!\geq\!\vert t\vert$). We choose $N_{t\varepsilon}\!\in\!\omega$ such that 
$$\bigg( s^i_{l(t)}~(i\!\cdot\!\varepsilon )~
\Big( s^0_{l(b_2(r))}\! -\! s^0_{l(b_2(r))}\vert \big( l(t)\! +\! 1\big)\Big)~0~
u^0_{\vert t\vert +1}...u^{\vert t\vert +1}_{\vert t\vert +1}~0^{N_{t\varepsilon}}\bigg)_{i\in d}\!\in\! E
\mbox{,}$$
$\Big( l\big( b_2(r)\big)\! +\!\vert u^0_{\vert t\vert +1}...u^{\vert t\vert +1}_{\vert t\vert +1}\vert\! +\! N_{t\varepsilon}\Big)_{0}\! =\! (\vert t\vert\! +\! 1)_0$ and 
$$\bigg(\Big( l\big( b_2(r)\big)\! +\!\vert u^0_{\vert t\vert +1}...u^{\vert t\vert +1}_{\vert t\vert +1}\vert\! +\! N_{t\varepsilon}\Big)_1\bigg)_0\! =\!\big( (\vert t\vert\! +\! 1)_1\big)_0.$$ 
We put $l(t\varepsilon )\! :=\! l(t)\! +\! 1\! +\!\vert v_{t\varepsilon}\vert$, where by definition 
$$v_{t\varepsilon}\! :=\!\Big( s^0_{l(b_2(r))}\! -\! s^0_{l(b_2(r))}\vert \big( l(t)\! +\! 1\big)\Big)~0~
u^0_{\vert t\vert +1}...u^{\vert t\vert +1}_{\vert t\vert +1}~0^{N_{t\varepsilon}}.$$ 
This finishes the proof.\hfill{$\square$}\bigskip

 Now we come to the general condition to get some complicated sets as in the statement of Theorem 1.7 announced in the introduction.\bigskip
 
 \noindent\bf Notation.\rm\ The map ${\cal S}\! :\! 2^\omega\!\rightarrow\! 2^\omega$ is the shift map: 
${\cal S}(\alpha )(m)\! :=\!\alpha (m\! +\! 1)$.

\begin{defi} We say that $C\!\subseteq\! 2^\omega$ is $compatible\ with\ comeager\ sets$ (ccs for short) if 
$$\alpha\!\in\! C\ \Leftrightarrow\ {\cal S}\big(\alpha_0\Delta F_0(\alpha )\big)\!\in\! C\mbox{,}$$
for each $\alpha_0\!\in\! d^\omega$ and $F\! :\! 2^\omega\!\rightarrow\! (d^\omega )^{d-1}$ satisfying the conclusion of Lemma 2.4.(b).\end{defi}

\noindent\bf Notation.\rm\ Let $T_d$ be the tree associated with an effective~frame, and 
$C\!\subseteq\! 2^\omega$. We put
$$S^d_C\! :=\!\big\{ \vec\alpha \!\in\!\lceil T_d\rceil\mid 
{\cal S}(\alpha_0\Delta\alpha_1)\!\in\! C\big\}.$$

\begin{lem} Let $T_d$ be the tree associated with an effective~frame, and ${\bf\Gamma}$ a non self-dual Wadge class of Borel sets.\smallskip

\noindent (1) Assume that $C$ is a $\bf\Gamma$-complete ccs set. Then 
$S^d_C\!\in\! {\bf\Gamma}(\lceil T_d\rceil)$ is a Borel subset of $(d^\omega )^d$, and is not separable from $\lceil T_d\rceil\!\setminus\! S^d_C$ by a $\mbox{pot}(\check {\bf\Gamma})$ set.\smallskip

\noindent (2) Assume that $C^0$, $C^1\!\in\! {\bf\Gamma}$ are disjoint, ccs, and not separable by a $\Delta ({\bf\Gamma})$ set. Then 
$S^d_{C^0},S^d_{C^1}$ are in ${\bf\Gamma}(\lceil T_d\rceil )$, disjoint Borel subsets of 
$(d^\omega )^d$, and not separable by a $\mbox{pot}\big(\Delta ({\bf\Gamma})\big)$ set.\end{lem}

\noindent\bf Proof.\rm ~(1) It is clear that $S^d_C\!\in\! {\bf\Gamma}(\lceil T_d\rceil )$ since $\cal S$ and 
$\Delta$ are continuous. So $S^d_C$ is a Borel subset of $(d^\omega )^d$ since $\lceil T_d\rceil$ is a closed subset of $(d^\omega )^d$. Indeed, $\lceil T_\omega\rceil$ is closed:
$$\vec\alpha \!\in\!\lceil T_\omega\rceil\ ~\Leftrightarrow\ ~\forall n\!\in\!\omega\!\setminus\!
\{ 0\}~\ \exists l\! <\! n~\ \forall i\!\in\!\omega\ ~s^i_li\!\subseteq\!\alpha_i\wedge 
(\alpha_i\vert n\! -\! s^i_li)\! =\! (\alpha_0\vert n\! -\! s^0_l0)\wedge\alpha_0(n)\!\leq\! n.$$
We argue by contradiction to see that $S^d_C$ is not separable from $\lceil T_d\rceil\!\setminus\! S^d_C$ by a $\mbox{pot}(\check {\bf\Gamma})$ set: this gives $P\!\in\!\mbox{pot}(\check {\bf\Gamma})$. For each 
$i\!\in\! d$ there is a dense $G_{\delta}$ subset $G_i$ of the compact space 
$\Pi_i''\lceil T_d\rceil$ such that $P\cap (\Pi_{i\in d}\ G_i)\!\in\!\check {\bf\Gamma}(\Pi_{i\in d}\ G_i)$, and 
$S^d_C\cap (\Pi_{i\in d}\ G_i)\!\subseteq\! P\cap (\Pi_{i\in d}\ G_i)\!\subseteq\! (\Pi_{i\in d}\ G_i)\!\setminus\! (\lceil T_d\rceil\!\setminus\! S^d_C)$.\bigskip

 Lemma 2.4 provides $\alpha_{0}\!\in\! G_0$ and $F\! :\! 2^\omega\!\rightarrow\!\Pi_{0<i<d}\ G_i$ continuous. Let 
$$D\! :=\!\big\{\alpha\!\in\! 2^\omega\mid \big(\alpha_0,F(\alpha )\big)\!\in\! P\cap (\Pi_{i\in d}\ G_i)\big\}.$$ 
Then $D\!\in\!\check {\bf\Gamma}$. Let us prove that $C\! =\! D$, which will contradict the fact that 
$C\!\notin\!\check {\bf\Gamma}$. As $C$ is ccs, $\alpha\!\in\! C$ is equivalent to 
${\cal S}\big(\alpha_0\Delta F_0(\alpha )\big)\!\in\! C$. Thus  
$$\alpha\!\in\! C\Rightarrow {\cal S}\big(\alpha_0\Delta F_0(\alpha )\big)\!\in\! C\Rightarrow
\big(\alpha_0, F(\alpha )\big)\!\in\! S^d_C\cap (\Pi_{i\in d}\ G_i)\!\subseteq\! P\cap (\Pi_{i\in d}\ G_i)
\Rightarrow\alpha\!\in\! D.$$
Similarly, $\alpha\!\notin\! C\Rightarrow\alpha\!\notin\! D$, and $C\! =\! D$.\bigskip

\noindent (2) We argue as in (1).\hfill{$\square$}\bigskip

 This lemma reduces the problem of finding some complicated sets as in the statement of Theorem 1.7 to a problem in dimension 1.

\section{$\!\!\!\!\!\!$ The proof of Theorem 1.7 for the Borel classes}\indent 

 The full version of Theorem 1.7 for the Borel classes is as follows:
 
\begin{thm} We can find concrete examples of a tree $T_d$ with $\Borel$ suitable levels, together with, for each $1\!\leq\!\xi\! <\!\omega_{1}$,\smallskip

\noindent (1) Some set $\mathbb{S}^d_\xi\!\in\!\boraxi (\lceil T_d\rceil)$ not separable from 
$\lceil T_d\rceil\!\setminus\!\mathbb{S}^d_\xi$ by a $\mbox{pot}(\bormxi )$ set.\smallskip

\noindent (2) Some disjoint sets $\mathbb{S}^0_\xi,\mathbb{S}^1_\xi\!\in\!\boraxi (\lceil T_d\rceil)$ not separable by a $\mbox{pot}(\borxi )$ set.\end{thm}

 This is an application of Lemma 2.6. We now introduce the objects useful to define the suitable sets $C$'s of this lemma. These objects will also be useful in the general case. The following definition can be found in [Lo-SR2] (see Definition 2.2).
 
\begin{defi} A set $H\!\subseteq\! 2^\omega$ is ${\bf\Gamma}\mbox{-}strategically\ complete$ if\smallskip

\noindent (a) $H\!\in\! {\bf\Gamma}(2^\omega )$.\smallskip

\noindent (b) If $A\!\in\! {\bf\Gamma}(\omega^\omega )$, then Player 2 wins the Wadge game $G(A,H)$ (where Player 1 plays $\alpha\!\in\!\omega^\omega$, Player 2 plays $\beta\!\in\! 2^\omega$ and Player 2 wins if $\alpha\!\in\! A\Leftrightarrow\beta\!\in\! H$).\end{defi}

 The following definition can essentially be found in [Lo-SR1] (see Section 3) and [Lo-SR2] (see Definition 2.3). 
  
\begin{defi} Let $\eta\! <\!\omega_1$. A function 
$\rho\! :\! 2^\omega\!\rightarrow\! 2^\omega$ is an $independent$ $\eta\mbox{-}function$ if\smallskip

\noindent (a) For some function $\pi\! :\!\omega\!\rightarrow\!\omega$, the value $\rho(\alpha )(m)$, for each $\alpha\!\in\! 2^\omega$ and each integer $m$, depends only on the values of $\alpha$ on 
$\pi^{-1}(\{ m\})$.\smallskip

\noindent (b) For each integer $m$, we set  
${\cal C}_m\! :=\!\{\alpha\!\in\! 2^\omega\mid\rho (\alpha )(m)\! =\! 1\}$.\smallskip

\noindent (1) If $\eta\! =\! 0$, then for each integer $m$ the set ${\cal C}_m$ is a $\borone$-complete set.\smallskip

\noindent (2) If $\eta\! =\!\theta\! +\! 1$ is successor, then for each integer $m$ the set ${\cal C}_m$ is a 
$\bormpz$-strategically complete set.\smallskip

\noindent (3) If $\eta$ limit, then for some sequence $(\theta_m)_{m\in\omega}$ with $\theta_m\! <\!\eta$ and $\mbox{sup}_{p\geq 1}\ \theta_{m_p}\! =\!\eta$ for each one-to-one sequence $(m_p)_{p\geq 1}$ of integers, and for each integer $m$ the set ${\cal C}_m$ is a $\bormpzm$-strategically complete set.
\end{defi}

 Note that we added a condition when $\eta\! =\! 0$. Moreover, we do not ask the sequence 
$(\theta_m)_{m\in\omega}$ to be increasing, unlike in [Lo-SR2], Definition 2.3. Note also that an independent $\eta$-function has to be ${\bf\Sigma}^0_{1+\eta}$-measurable. Moreover, if $\rho$ is an independent $\eta$-function, then $\pi$ has to be onto.\bigskip

\noindent\bf Examples.\rm ~In [Lo-SR1], Lemma 3.3, the map 
$\rho_{0}\! :\! 2^\omega\!\rightarrow\! 2^\omega$ defined as follows is introduced:
$$\rho_{0}(\alpha )(m)\! :=\!\left\{\!\!\!\!\!\!\!
\begin{array}{ll}
& 1\mbox{ if }\alpha (<\! m,n\! >)\! =\! 0\mbox{, for each }n\!\in\!\omega\mbox{,}\cr &Ê\cr
& 0\mbox{ otherwise.}
\end{array}\right.$$
Then $\rho_{0}$ is clearly an independent $1$-function, with $\pi (k)\! =\! (k)_0$. In this paper, 
$\rho_{0}^{\eta}\! :\! 2^\omega\!\rightarrow\! 2^\omega$ is also defined for $\eta\! <\!\omega_{1}$ as follows, by induction on $\eta$ (see the proof of Theorem 3.2).

\vfill\eject

 We put\smallskip

\noindent - ${\rho^0_{0}\! :=\!\hbox{\rm Id}_{2^\omega}}$.\smallskip

\noindent - $\rho^{\theta +1}_{0}\!:=\!\rho^{}_{0}\circ\rho^{\theta}_{0}$.\smallskip

\noindent - If $\eta\! >\! 0$ is limit, then fix a sequence $(\theta^\eta_{m})_{m\in\omega}\!\subseteq\!\eta$ of successor ordinals with $\Sigma_{m\in\omega}~\theta^\eta_m\! =\!\eta$. We define 
$\rho^{(m,m+1)}_{0}\! :\! 2^\omega\!\rightarrow\! 2^\omega$ by
$$\rho^{(m,m+1)}_{0}(\alpha )(i)\! :=\!\left\{\!\!\!\!\!\!\!
\begin{array}{ll} 
& \alpha (i)\mbox{ if }i\! <\! m\mbox{,}\cr &Ê\cr
& \rho_{0}^{\theta^\eta_{m}}\big( {\cal S}^m(\alpha )\big)(i\! -\! m)\mbox{ if }i\!\geq\! m.
\end{array}\right.$$ 
We set $\rho^{(0,m+1)}_{0}\! :=\!\rho^{(m,m+1)}_{0}\circ\rho^{(m-1,m)}_{0}\circ
\ldots\circ\rho^{(0,1)}_{0}$ and 
$\rho^\eta_{0}(\alpha )(m)\! :=\!\rho^{(0,m+1)}_{0}(\alpha )(m)$. The authors prove that 
$\rho_0^\eta$ is an independent $\eta$-function (see the proof of Theorem 3.2). In this paper, the set 
$H_{1+\eta}\! :=\! (\rho_0^\eta )^{-1}(\{ 0^\infty\})$ is also introduced, and the authors prove that 
$H_{1+\eta}$ is $\bormpe$-complete (see Theorem 3.2).\bigskip

\noindent\bf Notation.\rm\ Let $1\!\leq\!\xi\! :=\! 1\! +\!\eta\! <\!\omega_1$. We set $C_\xi\! :=\!\neg H_\xi$. If moreover $\varepsilon\!\in\! 2$, then we set
$$C^\varepsilon_\xi\! :=\!\big\{\alpha\!\in\! 2^\omega\mid\exists m\!\in\!\omega\ \ 
\rho_0^{\eta}(\alpha )(m)\! =\! 1\wedge\forall l\! <\! m\ \ \rho_0^{\eta}(\alpha )(l)\! =\! 0
\wedge (m)_0\!\equiv\!\varepsilon ~(\mbox{mod }2)\big\}.$$
Then we set $\mathbb{S}^d_\xi\! :=\! S^d_{C_\xi}$ and 
$\mathbb{S}^\varepsilon_\xi\! :=\! S^d_{C^\varepsilon_\xi}$.\bigskip

 Theorem 3.1 is a corollary of Proposition 2.2, Lemmas 2.3 and 2.6, and of the following lemma.
 
\begin{lem} Let $1\!\leq\!\xi\! <\!\omega_1$.\smallskip

\noindent (1) The set $C_\xi$ is a $\boraxi$-complete ccs set.\smallskip 

\noindent (2) The sets $C^0_\xi$, $C^1_\xi\!\in\!\boraxi$, are disjoint, ccs, and not separable by a $\borxi$ set.\end{lem}

\noindent\bf Proof.\rm ~(1) $C_\xi$ is $\boraxi$-complete since $H_\xi$ is $\bormxi$-complete.\bigskip

\noindent $\bullet$ Assume that $\alpha_0$, $F$ satisfy the conclusion of Lemma 2.4.(b). Let us prove that 
$$\rho^\eta_0(\alpha )\! =\!\rho^\eta_0\Big( {\cal S}\big(\alpha_0\Delta F_0(\alpha )\big)\Big)\mbox{,}$$ 
for each $1\!\leq\!\eta\! <\!\omega_1$ and $\alpha\!\in\! 2^\omega$. For $\eta\! =\! 1$ we 
apply the conclusion of Lemma 2.4.(b) to $s\!\in\!\omega$. Then we have, by induction, 
$\rho^{\theta +1}_0(\alpha )\! =\!\rho_0\big(\rho^{\theta }_0(\alpha )\big)\! =\! 
\rho_0\Bigg(\rho^{\theta }_0\Big( {\cal S}\big(\alpha_0\Delta F_0(\alpha )\big)\Big)\Bigg)\! =
\!\rho^{\theta +1}_0\Big( {\cal S}\big(\alpha_0\Delta F_0(\alpha )\big)\Big)$. From this we deduce, for 
$\lambda\! >\! 0$ limit, by induction again, that
$$\rho^{(0,1)}_0(\alpha )\! =\!\rho^{\theta^\lambda_0}_0(\alpha )\! =\! 
\rho^{\theta^\lambda_0}_0\Big( {\cal S}\big(\alpha_0\Delta F_0(\alpha )\big)\Big)\! =\! 
\rho^{(0,1)}_0\Big( {\cal S}\big(\alpha_0\Delta F_0(\alpha )\big)\Big).$$
Thus $\rho^{(0,m+1)}_0(\alpha )\! =\!\rho^{(0,m+1)}_0\Big( 
{\cal S}\big(\alpha_0\Delta F_0(\alpha )\big)\Big)$, and
$$\rho^{\lambda}_0(\alpha )(m)\! =\! 
\rho^{(0,m+1)}_0(\alpha )(m)\! =\! 
\rho^{(0,m+1)}_0\Big( {\cal S}\big(\alpha_0\Delta F_0(\alpha )\big)\Big)(m)\! =\! 
\rho^{\lambda}_0\Big( {\cal S}\big(\alpha_0\Delta F_0(\alpha )\big)\Big)(m).$$
$\bullet$ If we apply the previous point, or the conclusion of Lemma 2.4.(b) to $s\! :=\!\emptyset$, then we get
$$\alpha\!\in\! C_\xi
\Leftrightarrow\exists m\!\in\!\omega\ \ \rho_0^{\eta}(\alpha )(m)\! =\! 1
\Leftrightarrow\exists m'\!\in\!\omega\ \ 
\rho_0^{\eta}\Big({\cal S}\big(\alpha_0\Delta F_0(\alpha )\big)\Big) (m')\! =\! 1
\Leftrightarrow {\cal S}\big(\alpha_0\Delta F_0(\alpha )\big)\!\in\! C_\xi .$$
Thus $C_\xi$ is ccs.

\vfill\eject

\noindent (2) Note first that $C^0_\xi$, $C^1_\xi\!\in\!\boraxi$ since $\rho_0^\eta$ is 
${\bf\Sigma}^0_{1+\eta}$-measurable, are clearly disjoint, and are ccs as in (1) since 
$(m)_0\! =\!\big( B_\alpha (m)\big)_0$ in Lemma 2.4.(b).\bigskip

\noindent $\bullet$ We set, for $\varepsilon\!\in\! 2$, $V_\varepsilon\! :=\!\big\{\alpha\!\in\! 2^\omega\mid\exists m\!\in\!\omega\ \ \rho_0^{\eta}(\alpha )(m)\! =\! 1\ \mbox{ and }\ 
(m)_0\!\equiv\!\varepsilon\ (\mbox{mod }2)\big\}$. Then $V_\varepsilon$ is a $\boraxi$ set since 
$\rho_0^\eta$ is ${\bf\Sigma}^0_{1+\eta}$-measurable. Let us prove that $V_\varepsilon$ is $\boraxi$-complete.\bigskip

\noindent - If $\eta\! =\! 0$, then $0^\infty\!\in\!\overline{V_\varepsilon}\!\setminus\! V_\varepsilon$, so that $V_\varepsilon$ is $\boraone$-complete.\bigskip

\noindent - If $\eta\! =\!\theta\! +\! 1$, then $\rho_0^{\eta}$ is an independent $\eta$-function. Let 
$(A_m)_{m\in\omega}$ be a sequence of $\bormpz (2^\omega )$ sets. Choose a continuous map 
$f_m\! :\! 2^\omega\!\rightarrow\! 2^\omega$ such that $A_m\! =\! f_m^{-1}({\cal C}_m)$. We define 
$f\! :\! 2^\omega\!\rightarrow\! 2^\omega$ by $f(\alpha )(k)\! :=\! f_m(\alpha )(k)$ if $\pi_{\eta}(k)\! =\! m$, and $f$ is continuous. Moreover,
$$\alpha\!\in\! A_m\Leftrightarrow f_m(\alpha )\!\in\! {\cal C}_n\Leftrightarrow 
f(\alpha )\!\in\! {\cal C}_m\mbox{,}$$
so that $\bigcup_{m\in\omega ,(m)_0\equiv\varepsilon\ (\mbox{mod }2)}\ A_m\! =\! f^{-1}(V_\varepsilon )$. Thus $V_\varepsilon$ is $\boraxi$-complete.\bigskip

\noindent - If $\eta$ is the limit of the $\theta_m$'s, then $\rho_0^{\eta}$ is an independent $\eta$-function. We argue as in the successor case to see that $V_\varepsilon$ is $\boraxi$-complete.\bigskip

\noindent $\bullet$ We argue by contradiction, which gives $D\!\in\!\borxi$ separating $C^0_\xi$ from $C^1_\xi$. Let $v_0$, $v_1$ be disjoint $\boraxi$ subsets of $2^\omega$. Then we can find a continuous map $f_\varepsilon\! :\! 2^\omega\!\rightarrow\! 2^\omega$ such that 
$v_\varepsilon\! =\! f_\varepsilon^{-1}(V_\varepsilon )$. As $\rho_0^\eta$ is an independent $\eta$-function, we get $\pi_{\eta}\! :\!\omega\!\rightarrow\!\omega$. We define a map 
$f\! :\! 2^\omega\!\rightarrow\! 2^\omega$ by 
$f(\alpha )(k)\! :=\! f_\varepsilon (\alpha )(k)$ if 
$\big(\pi_{\eta}(k)\big)_0\!\equiv\!\varepsilon\ (\mbox{mod }2)$, and $f$ is continuous. Note that 
$\alpha\!\in\! v_\varepsilon\Leftrightarrow f_\varepsilon (\alpha )\!\in\! V_\varepsilon\Leftrightarrow 
f(\alpha )\!\in\! V_\varepsilon$, so that $v_\varepsilon\! =\! f^{-1}(V_\varepsilon )$. Thus 
$\alpha\!\in\! v_0\Leftrightarrow f(\alpha )\!\in\! V_0\Leftrightarrow 
f(\alpha )\!\in\! V_0\!\setminus\! V_1\!\subseteq\! C^0_\xi\!\subseteq\! D$ since $v_0$ is disjoint from $v_1$. Similarly, 
$\alpha\!\in\! v_1\Leftrightarrow f(\alpha )\!\in\! V_1\!\setminus\! V_0\!\subseteq\! C^1_\xi\!\subseteq\!\neg D$. Thus $f^{-1}(D)$ separates $v_0$ from $v_1$. As $f^{-1}(D)\!\in\!\borxi$, this implies that $\boraxi$ has the separation property, which contradicts 22.C in [K].\hfill{$\square$}

\section{$\!\!\!\!\!\!$ The proof of Theorem 1.8 for the Borel classes}\indent 

 The full versions of Theorem 1.8 and Corollary 1.9 for the Borel classes are as follows:

\begin{thm} Let $T_d$ be a tree with suitable levels, $1\!\leq\!\xi\! <\!\omega_1$, $(X_i)_{i\in d}$ a sequence of Polish spaces, and $A_0$, $A_1$ disjoint analytic subsets of $\Pi_{i\in d}\ X_i$.\medskip 

\noindent (1) Let $S\!\in\!\boraxi (\lceil T_d\rceil)$. Then one of the following holds:\smallskip

\noindent (a) The set $A_0$ is separable from $A_1$ by a $\mbox{pot}(\bormxi )$ set.\smallskip

\noindent (b) The inequality $\big( (d^\omega )_{i\in d}, S,\lceil T_d\rceil\!\setminus\! S\big)\leq
\big( (X_i)_{i\in d}, A_0, A_1\big)$ holds.\smallskip

 If we moreover assume that $S$ is not separable from $\lceil T_d\rceil\!\setminus\! S$ by a 
$\mbox{pot}(\bormxi )$ set, then this is a dichotomy.\medskip 

\noindent (2) Let $S^0, S^1\!\in\!\boraxi (\lceil T_d\rceil)$ disjoint. Then one of the following holds:\smallskip

\noindent (a) The set $A_0$ is separable from $A_1$ by a $\mbox{pot}(\borxi )$ set.\smallskip

\noindent (b) The inequality $\big( (d^\omega )_{i\in d}, S^0,S^1\big)\leq
\big( (X_i)_{i\in d}, A_0, A_1\big)$ holds.\smallskip

 If we moreover assume that $S^0$ is not separable from $S^1$ by a $\mbox{pot}(\borxi )$ set, then this is a dichotomy.\end{thm}

\begin{cor} Let $\bf\Gamma$ be Borel class. Then there are Borel subsets 
$\mathbb{S}^0$, $\mathbb{S}^1$ of $(d^\omega )^d$ such that for any sequence of Polish spaces 
$(X_i)_{i\in d}$, and for any disjoint analytic subsets $A_0$, $A_1$ of 
$\Pi_{i\in d}\ X_i$, exactly one of the following holds:\smallskip

\noindent (a) The set $A_0$ is separable from $A_1$ by a $\hbox{\it pot}({\bf\Gamma})$ set.\smallskip

\noindent (b) The inequality 
$\big( (d^\omega )_{i\in d},\mathbb{S}^0,\mathbb{S}^1\big)\leq\big( (X_i)_{i\in d}, A_0, A_1\big)$ holds.
\end{cor}

\subsection{$\!\!\!\!\!\!$ Acyclicity}\indent

 In this subsection we prove a result that will be used later to show Theorem 4.1. This is the place where the essence of the notion of a finite one-sided almost acyclic set is really used.
  
\begin{lemm} Assume that ${\cal T}\!\subseteq\! {\cal X}^d$ is finite. Then the following are equivalent:\medskip

\noindent (a) The set ${\cal T}$ is one-sided and almost acyclic.\medskip

\noindent (b) For each $\overrightarrow{x^{0}}\!\in\! {\cal T}$, there is an integer 
$0\!\not=\! {\cal L}\! <\! d\! +\! 2$ and a partition $(M_j)_{j\in {\cal L}}$ of 
${\cal T}\!\setminus\!\{\overrightarrow{x^{0}}\}$ with\medskip

(1) $\forall i\!\in\! d\ \ \forall j\!\not=\! k\!\in\! {\cal L}\ \ \Pi_i[M_j]\cap\Pi_i[M_k]\! =\!\emptyset$.\medskip

(2) $\forall i\!\in\! d\ \ \forall j\!\in\! {\cal L}\ \ \forall\vec x\!\in\! M_j\ \ x_i\! =\! x^0_i\ \Rightarrow\ i\! =\! j$.
\end{lemm}

\noindent\bf Proof.\rm\ (a) $\Rightarrow$ (b) If $\vec y\!\not=\!\vec z\!\in\! {\cal T}$ and 
$\left(\overrightarrow{y^{j}}\right)_{j\leq l}$ is a walk in $G^{\cal T}$ with $\overrightarrow{y^{0}}\! =\!\vec y$ and $\overrightarrow{y^{l}}\! =\!\vec z$, then we choose such a walk of minimal length, and we call it 
$w_{\vec y,\vec z}$. We will define a partition of ${\cal T}$. We put, for $j\!\in\! d$,
$$\begin{array}{ll}
N\!\!\!\! 
& :=\{\ \vec  x\!\in\! {\cal T}\mid\vec x\!\not=\!\overrightarrow{x^{0}}\wedge 
w_{\vec x,\overrightarrow{x^{0}}}\mbox{ does not exist}\ \}\mbox{,}\cr
L_j\!\!\!\! 
& :=\{\ \vec x\!\in\! {\cal T}\mid\vec x\!\not=\!\overrightarrow{x^{0}}\wedge  
\big( w_{\vec x,\overrightarrow{x^{0}}}(|w_{\vec x,\overrightarrow{x^{0}}}|\! -\! 2)\big)_j\! =\! x^0_{j}\ \}.
\end{array}$$
So we defined a partition $\big(N,(L_j)_{j\in d}\big)$ of 
${\cal T}\!\setminus\!\{\overrightarrow{x^{0}}\}$ since $\cal T$ is one-sided. As ${\cal T}$ is finite, there is $j_0\!\in\! d$ minimal such that $L_{j}\! =\!\emptyset$ if $j\!\ >\! j_0$. We set $M_j\! :=\! L_j$ if 
$j\!\leq\! j_0$, $M_{j_0+1}\! :=\! N$ and ${\cal L}\! :=\! j_0\! +\! 2$.\bigskip

\noindent (1) Let us prove that $\Pi_{i}[L_j]\cap\Pi_{i}[N]\! =\!\emptyset$, for each $i,j\!\in\! d$. We argue by contradiction. This gives $x_{i}\!\in\!\Pi_{i}[L_j]\cap\Pi_{i}[N]$, $\vec x\!\in\! L_j$, and also $\vec y\!\in\! N$ such that $x_i\! =\! y_i$. As $\vec x,\vec y\!\in\! {\cal T}$ and $L_j\cap N\! =\!\emptyset$, $\vec x\!\not=\!\vec y$ and $\vec x,\vec y$ are $G^{\cal T}$-related. Note that $w_{\vec y,\overrightarrow{x^{0}}}$ does not exist, and that $w_{\vec x,\overrightarrow{x^{0}}}$ exists. Now the sequence 
$\left(\vec y,\vec x,...,\overrightarrow{x^{0}}\right)$ shows the existence of 
$w_{\vec y,\overrightarrow{x^{0}}}$, which is absurd.\bigskip

 It remains to see that $\Pi_{i}[L_j]\cap\Pi_{i}[L_k]\! =\!\emptyset$, for each $i,j,k\!\in\! d$ with $j\!\not=\! k$. We argue by contradiction. This gives $x_{i}\!\in\!\Pi_{i}[L_j]\cap\Pi_{i}[L_k]$, $\vec x\!\in\! L_j$, and also 
$\vec y\!\in\! L_k$ such that $x_i\! =\! y_i$. As $\vec x,\vec y\!\in\! {\cal T}$ and $j\!\not=\! k$, 
$\vec x\!\not=\!\vec y$ and $\vec x,\vec y$ are $G^{\cal T}$-related. Let us denote 
$w_{\vec x,\overrightarrow{x^{0}}}\! :=\!\left(\overrightarrow{z^n}\right)_{n\leq l+1}$ and 
$w_{\vec y,\overrightarrow{x^{0}}}\! :=\!\left(\overrightarrow{y^n}\right)_{n\leq l'+1}$. Note that 
$\overrightarrow{z^l}\!\not=\!\overrightarrow{y^{l'}}$ since 
$z^{l}_j\! =\! x^0_{j}$ and $y^{l'}_j\!\not=\! x^0_{j}$, since otherwise $\overrightarrow{y^{l'}}$, 
$\overrightarrow{x^{0}}\!\in\! {\cal T}$, $\overrightarrow{y^{l'}}\!\not=\!\overrightarrow{x^{0}}$ and 
$y^{l'}_j\! =\! x^0_{j}$, $y^{l'}_k\! =\! x^0_{k}$, which contradicts the fact that ${\cal T}$ is one-sided.\bigskip

 We denote by $W\! :=\!\left(\overrightarrow{w^n}\right)_{n\leq L}$ the following $G^{\cal T}$-walk: 
$\overrightarrow{z^l}, \overrightarrow{z^{l-1}},...,\overrightarrow{z^0},\overrightarrow{y^0},
\overrightarrow{y^1},...,\overrightarrow{y^{l'}}$. If there are $k\! <\! n\!\leq\! L$ with 
$\overrightarrow{w^k}\! =\!\overrightarrow{w^n}$, then we put 
$W'\! :=\!\overrightarrow{w^0},...,\overrightarrow{w^k},\overrightarrow{w^{n+1}},...,\overrightarrow{w^L}$. If we iterate this construction, then we get a $G^{\cal T}$-walk without repetition 
$V\! :=\!\left(\overrightarrow{v^n}\right)_{n\leq L'}$ from $\overrightarrow{w^0}$ to $\overrightarrow{w^L}$.\bigskip 

 If there are $i\!\in\! d$ and $k\! +\! 1\! <\! n\!\leq\! L'$ with $v^k_i\! =\! v^n_i$, then we put 
$V'\! :=\!\overrightarrow{v^0},...,\overrightarrow{v^k},\overrightarrow{v^{n}},...,\overrightarrow{v^{L'}}$.  
If we iterate this construction, then we get a $G^{\cal T}$-walk without repetition 
$U\! :=\!\left(\overrightarrow{u^n}\right)_{n\leq L''}$ from $\overrightarrow{w^0}$ to $\overrightarrow{w^L}$ for which it is not possible to find $i\!\in\! d$ and $k\! +\! 1\! <\! n\!\leq\! L''$ with $u^k_i\! =\! u^n_i$.\bigskip

 Now $\overrightarrow{x^{0}},\overrightarrow{u^0},...,\overrightarrow{u^{L''}},\overrightarrow{x^{0}}$ is a $G^{\cal T}$-cycle contradicting the almost acyclicity of $\cal T$.\bigskip

\noindent (2) If $\vec x\!\in\! N$, then $w_{\vec x,\overrightarrow{x^{0}}}$ does not exist. This implies that 
$x_i\!\not=\! x^0_i$ for each $i\!\in\! d$, since otherwise $\vec x$ and $\overrightarrow{x^{0}}$ would be 
$G^{\cal T}$-related, which contradicts the non-existence of $w_{\vec x,\overrightarrow{x^{0}}}$.\bigskip

 If $\vec x\!\in\! L_j$, then $i$ is the only coordinate for which $x_i\! =\! x^0_i$ since $\cal T$ is 
one-sided. Note that $w_{\vec x,\overrightarrow{x^{0}}}\! =\!\left(\vec x,\overrightarrow{x^{0}}\right)$. As 
$\vec x\!\in\! L_j$, we get 
$\big(w_{\vec x,\overrightarrow{x^{0}}}(|w_{\vec x,\overrightarrow{x^{0}}}|\! -\! 2)\big)_j\! =\! x^0_{j}$. But 
$w_{\vec x,\overrightarrow{x^{0}}}(|w_{\vec x,\overrightarrow{x^{0}}}|\! -\! 2)\! =\!\vec x$. Thus 
$x_j\! =\! x^0_j$ and $i\! =\! j$.\bigskip

\noindent (b) $\Rightarrow$ (a) Let $\overrightarrow{x^0}\!\not=\!\vec x\!\in\! {\cal T}$, $i,j\!\in\! d$ such that $x^0_i\! =\! x_i$ and $x^0_j\! =\! x_j$, and $k\!\in\! {\cal L}$ such that $\vec x\!\in\! M_k$. By (2) we get 
$i\! =\! k\! =\! j$ and $\cal T$ is one-sided. Now consider a $G^{\cal T}$-cycle 
$(\overrightarrow{x^n})_{n\leq L}$. By (1) there is $j\!\in\! {\cal L}$ such that $\overrightarrow{x^n}\!\in\! M_j$ for each $0\! <\! n\! <\! L$. Then by (2) we get $x^0_j\! =\! x^1_j\! =\! x^{L-1}_j$ and $\cal T$ is almost acyclic.
\hfill{$\square$}\bigskip

 Definition 4.1.2 and Lemma 4.1.3 below are essentially due to G. Debs (see Subsection 2.1 in [L7]):
  
\begin{defin} (Debs) Let $\Theta\! :\! {\cal X}^d\!\rightarrow\! 2^{(\omega^\omega )^d}$,  
${\cal T}\!\subseteq\! {\cal X}^d$. We say that the map 
$\theta\! =\!\Pi_{i\in d}\ \theta_{i}\!\in\!\big( (\omega^\omega )^{\cal X}\big)^d$ is a 
$\pi\mbox{-}selector~on~{\cal T}~for~\Theta$ if\smallskip

\noindent (a) $\theta (\vec x)\! =\!\big(\theta_{i} (x_{i})\big)_{i\in d}$ for each 
$\vec x\!\in\! {\cal X}^d$.\smallskip

\noindent (b) $\theta (\vec x)\!\in\!\Theta (\vec x)$ for each $\vec x\!\in\! {\cal T}$.\end{defin}

\begin{lemm} (Debs) Let $l$ be an integer, ${\cal X}\! :=\! d^{l+1}$, ${\cal T}\!\subseteq\! {\cal X}^d$ be $\Borel$, finite, one-sided, and almost acyclic, 
$\Theta\! :\! {\cal X}^d\!\rightarrow\!\Ana\big( (\omega^\omega )^d\big)$, and  
$\overline{\Theta}\! :\! {\cal X}^d\!\rightarrow\!\Ana\big( (\omega^\omega )^d\big)$ defined by 
$\overline{\Theta}(\vec x)\! :=\!\overline{\Theta (\vec x)}^{\tau_{1}}$. Then $\Theta$ 
admits a $\pi$-selector on $\cal T$ if $\overline{\Theta}$ does.\end{lemm}

\noindent\bf Proof.\rm ~(a) Let $\overrightarrow{x^{0}}\!\in\! {\cal T}$, and 
$\Psi\! :\! {\cal X}^d\!\rightarrow\!\Ana\big( (\omega^\omega )^d\big)$. We 
assume that $\Psi (\vec x)\! =\!\Theta (\vec x)$ if $\vec x\!\not=\!\overrightarrow{x^{0}}$, 
and that $\Psi \left(\overrightarrow{x^{0}}\right)\!\subseteq\!
\overline{\Theta \left(\overrightarrow{x^{0}}\right)}^{\tau_{1}}$. We first prove that $\Theta$ admits a $\pi$-selector on $\cal T$ if $\Psi$ does.\bigskip

\noindent $\bullet$ Lemma 4.1.1 gives a finite partition $(M_j)_{j\in {\cal L}}$ of 
${\cal T}\!\setminus\!\{\overrightarrow{x^{0}}\}$. Fix a $\pi$-selector $\tilde\psi$ on $\cal T$ for $\Psi$, and let $M\! :=\!\mbox{max}\ (d\cap {\cal L})$. We define $\Ana$ sets $U_{i}$, for $i\!\leq\! M$, by
$$U_{i} :=\big\{\alpha\!\in\!\omega^\omega\mid\exists\psi\!\in\!\big( (\omega^\omega )^{\cal X}\big)^d
\ \ \alpha\! =\!\psi_{i}(x^0_{i})\wedge\forall\vec x\!\in\! {\cal T}~\ \psi (\vec x)\!\in\!\Psi (\vec x)\big\}.$$

 As $\tilde\psi\left(\overrightarrow{x^{0}}\right)\! =\!\big(\tilde\psi_{i}(x^0_{i})\big)_{i\in d}\!\in\!
\Psi\left(\overrightarrow{x^{0}}\right)\cap\big( (\Pi_{i\leq M}\ U_{i})\!\times\! (\omega^\omega )^{d-M-1}\big)$ we get 
$$\emptyset\!\not=\!\Psi\left(\overrightarrow{x^{0}}\right)\cap
\big( (\Pi_{i\leq M}\ U_{i})\!\times\! (\omega^\omega )^{d-M-1}\big)\!\subseteq\! 
\overline{\Theta\left(\overrightarrow{x^{0}}\right)}^{\tau_{1}}\cap
\big( (\Pi_{i\leq M}\ U_{i})\!\times\! (\omega^\omega )^{d-M-1}\big).$$ 

 By the separation theorem this implies that 
$\Theta\left(\overrightarrow{x^{0}}\right)\cap\big( (\Pi_{i\leq M}\ U_{i})\!\times\! (\omega^\omega )^{d-M-1}\big)$ is not empty and contains some point $\vec\alpha$. Fix $i\!\leq\! M$. As $\alpha_{i}\!\in\! U_{i}$ there is $\psi^i\!\in\!\big( (\omega^\omega )^{\cal X}\big)^d$ such that $\alpha_{i}\! =\!\psi^i_{i}(x^0_{i})$ and 
$\psi^i (\vec x)\!\in\!\Psi (\vec x)$ if $\vec x\!\in\! {\cal T}$.\bigskip
 
\noindent $\bullet$ Now we can define 
$\theta_{i}\! :\! {\cal X}\!\rightarrow\!\omega^\omega$, for each $i\!\in\! d$. We put 
$$\begin{array}{ll}
\theta_{i}(x_{i})\!\!\!\! & :=\left\{\!\!\!\!\!\!\!
\begin{array}{ll}
& \alpha_{i}\mbox{ if }x_i\! =\! x^0_i\mbox{,}\cr\cr
& \psi^{j}_{i}(x_{i})\mbox{ if }x_i\!\in\!\Pi_{i}[M_j]\!\setminus\!\{ x^0_i\}\wedge j\!\leq\! M 
\mbox{,}\cr\cr
& \psi^{0}_{i}(x_{i})\mbox{ otherwise.}
\end{array}\right.
\end{array}$$
Then we set $\theta (\vec x)(i)\! :=\!\theta_{i}(x_{i})$ if $i\!\in\! d$.\bigskip

\noindent $\bullet$ It remains to see that $\theta (\vec x)\!\in\!\Theta (\vec x)$ for each 
$\vec x\!\in\! {\cal T}$.\bigskip 

 Note that 
$\theta\left(\overrightarrow{x^{0}}\right)\! =\!\vec\alpha\!\in\!\Theta\left(\overrightarrow{x^{0}}\right)$. So we may assume that $\vec x\!\not=\!\overrightarrow{x^{0}}$. So let $j\!\in\! {\cal L}$ with $\vec x\!\in\! M_j$.\bigskip

\noindent - If $x_i\!\not=\! x^0_i$ for each $i\!\in\! d$ and $j\!\leq\! M$, then 
$\theta (\vec x)\! =\!\big(\theta_{i}(x_{i})\big)_{i\in d}\! =\!\psi^{j}(\vec x)\!\in\!\Psi (\vec x)\! =\!\Theta (\vec x)$.\bigskip

\noindent - Similarly, if $x_i\!\not=\! x^0_i$ for each $i\!\in\! d$ and $j\! >\! M$, then 
$\theta (\vec x)\! =\!\big(\theta_{i}(x_{i})\big)_{i\in d}\! =\!\psi^{0}(\vec x)\!\in\!\Psi (\vec x)\! =\!\Theta (\vec x)$.\bigskip

\noindent - If $x_{i}\! =\! x^0_{i}$ for some $i\!\in\! d$, then $i\! =\! j\!\leq\! M$. This implies that 
$\theta_j(x_j)\! =\!\alpha_j\! =\!\psi^j_j(x^0_j)\! =\!\psi^j_j(x_j)$ and
$$\theta (\vec x)\! =\!\big(\theta_{i}(x_{i})\big)_{i\in d}\! =\!\psi^{j}(\vec x)\!\in\!\Psi (\vec x)\! =\!\Theta (\vec x).$$
(b) Write ${\cal T}\!:=\!\left\{\overrightarrow{x^{1}},\ldots ,\overrightarrow{x^{n}}\right\}$, and set 
$\Psi_{0}\!:=\!\overline{\Theta}$. We define 
${\Psi_{j+1}\! :\! {\cal X}^d\!\rightarrow\!\Ana\big( (\omega^\omega )^d\big)}$ as follows. We put 
$\Psi_{j+1}(\vec x)\! :=\!\Psi_{j}(\vec x)$ if $\vec x\!\not=\!\overrightarrow{x^{j+1}}$, and   
${\Psi_{j+1}\left(\overrightarrow{x^{j+1}}\right)\! :=\!\Theta\left(\overrightarrow{x^{j+1}}\right)}$, for $j\! <\! n$. The result now follows from an iterative application of (a).\hfill{$\square$} 

\subsection{$\!\!\!\!\!\!$ The topologies}\indent

 In this subsection we prove two other results that will be used to show Theorem 4.1. We use tools of effective descriptive set theory (the reader should see [M] for the basic notions). We first recall a classical result in the spirit of Theorem 3.3.1 in [H-K-Lo].\bigskip
 
\noindent\bf Notation.\rm ~Let $X$ be a recursively presented Polish space. Using the bijection between $\omega$ and $\omega^2$ defined before Definition 2.1, we can build a bijection 
$(x_n)\!\mapsto <x_n>$ between $(X^\omega )^\omega$ and $X^\omega$ by the formula 
$<x_n>(l)\! :=\! x_{(l)_0}\big( (l)_1\big)$. The inverse map $x\!\mapsto\! \big( (x)_n\big)$ is given by $(x)_n(p)\! :=\! x(<n,p>)$. These bijections are recursive.
 
\begin{lemm} Let $X$ be a recursively presented Polish space. Then there are $\Ca$ sets 
$W^{X}\!\subseteq\!\omega^\omega$, $C^{X}\!\subseteq\!\omega^\omega\!\times\! X$ with 
$\{(\alpha ,x)\!\in\!\omega^\omega\!\times\!\! X\mid\alpha\!\in\! W^{X}\mbox{ and }
x\!\notin\! C_{\alpha}^{X}\}\!\in\!\Ca$, $\Borel (X)\! =\!\{ C^{X}_{\alpha}\mid\alpha\!\in\!\Borel\cap W^{X}\}$, and $\borel (X)\! =\!\{ C^{X}_{\alpha}\mid\alpha\!\in\! W^{X}\}$.\end{lemm}
 
\noindent\bf Proof.\rm\ By 3E.2, 3F.6 and 3H.1 in [M], there is 
${\cal U}^X\!\in\!\Ca (\omega^\omega\!\times\! X)$ which is universal for $\ca (X)$ and satisfies the two following properties:\smallskip

\noindent - A subset $P$ of $X$ is $\Ca$ if and only if there is $\alpha\!\in\!\omega^\omega$ recursive with $P\! =\! {\cal U}^X_\alpha$.\smallskip

\noindent - There is $S^X\! :\!\omega^\omega\!\times\!\omega^\omega\!\rightarrow\!\omega^\omega$ recursive such that $(\alpha ,\beta ,x)\!\in\! {\cal U}^{\omega^\omega\times X}\Leftrightarrow
\big( S^X(\alpha ,\beta ),x\big)\!\in\! {\cal U}^X$.\bigskip

 We set, for $\varepsilon\!\in\! 2$, 
$U_\varepsilon\! :=\!\{ (\alpha ,x)\!\in\!\omega^\omega\!\times\! X\mid
\big( (\alpha )_\varepsilon ,x\big)\!\in\! {\cal U}^X\}$. Then $U_\varepsilon\!\in\!\Ca$. By 4B.10 in [M], $\Ca$ has the reduction property, which gives $U'_0,U'_1\!\in\!\Ca$ disjoint with 
$U'_\varepsilon\!\subseteq\! U_\varepsilon$ and $U'_0\cup U'_1\! =\! U_0\cup U_1$. We set  
$W^X\! :=\!\{\alpha\!\in\!\omega^\omega\mid (U'_0)_\alpha\cup (U'_1)_\alpha\! =\! X\}$ and 
$C^X\! :=\! U'_0$, which defines $\Ca$ sets. Moreover,
$$\alpha\!\in\! W^{X}\wedge x\!\notin\! C_{\alpha}^{X}\Leftrightarrow\alpha\!\in\! W^{X}\wedge 
(\alpha ,x)\!\in\! U'_1$$
is $\Ca$ in $(\alpha ,x)$. Assume that $A\!\in\!\Borel (X)$, which gives 
$\alpha_0,\alpha_1\!\in\!\omega^\omega$ recursive with $A\! =\! {\cal U}^X_{\alpha_0}$ (resp., 
$\neg A\! =\! {\cal U}^X_{\alpha_1}$). We define $\alpha\!\in\!\omega^\omega$ by 
$(\alpha )_\varepsilon\! :=\!\alpha_\varepsilon$, so that $\alpha$ is recursive. We get
$$\begin{array}{ll}
& x\!\in\! A\Leftrightarrow (\alpha_0,x)\!\in\! {\cal U}^X\Leftrightarrow (\alpha ,x)\!\in\! U_0\Leftrightarrow (\alpha ,x)\!\in\! U_0\!\setminus\! U_1\Leftrightarrow (\alpha ,x)\!\in\! U'_0\mbox{,}\cr
& x\!\notin\! A\Leftrightarrow (\alpha_1,x)\!\in\! {\cal U}^X\Leftrightarrow (\alpha ,x)\!\in\! U_1
\Leftrightarrow (\alpha ,x)\!\in\! U_1\!\setminus\! U_0\Leftrightarrow (\alpha ,x)\!\in\! U'_1\mbox{,}
\end{array}$$
so that $\alpha\!\in\! W^X$ and $C^X_\alpha\! =\! A$. This also proves that 
$\borel (X)\!\subseteq\!\{ C^{X}_{\alpha}\mid\alpha\!\in\! W^{X}\}$.\bigskip

 Conversely, let $\alpha\in\!\Borel\cap W^X$. Then $C^X_\alpha\!\in\!\Ca$, and 
$x\!\notin\! C^X_\alpha\Leftrightarrow\alpha\!\in\! W^{X}\mbox{ and }x\!\notin\! C_{\alpha}^{X}$, so that 
$\neg C^X_\alpha\!\in\!\Ca$ and $C^X_\alpha\!\in\!\Borel$. Note that this also proves that 
$\borel (X)\!\supseteq\!\{ C^{X}_{\alpha}\mid\alpha\!\in\! W^{X}\}$.\hfill{$\square$}\bigskip
 
 We now give some notation to state an effective version of Theorem 4.1.\bigskip
 
\noindent\bf Notation.\rm ~Let $X$ be a recursively presented Polish space.\bigskip

\noindent $\bullet$ We will use the Gandy-Harrington topology ${\it\Sigma}_{X}$ on $X$ generated by $\Ana (X)$. Recall that the set 
$\Omega_{X}\! :=\!\{x\!\in\! X\mid\omega_1^x\! =\!\omega^{\hbox{\rm CK}}_{1}\}$ is Borel and $\Ana$, that $(\Omega_{X},{\it\Sigma}_{X})$ is a $0$-dimensional Polish space (the intersection of $\Omega_{X}$ with any nonempty $\Ana$ set is a nonempty clopen subset of $(\Omega_{X},{\it\Sigma}_{X})$) (see [L8]).\bigskip 

\noindent $\bullet$ Recall the topology $\tau_1$ defined before Theorem 1.9. We will also consider some topologies between $\tau_{1}$ and ${\it\Sigma}_{(\omega^\omega )^d}$. Let 
$2\!\leq\!\xi\! <\!\omega^{\hbox{\rm CK}}_{1}$. The topology $\tau_{\xi}$ is generated by 
${\Ana\big( (\omega^\omega )^d\big)\cap\bormlxi (\tau_{1})}$. We have 
${\boraone(\tau_{\xi})\!\subseteq\!\boraxi (\tau_{1})}$, so that 
${\bormone(\tau_{\xi})\!\subseteq\!\bormxi (\tau_{1})}$. These topologies are similar to the ones considered in [Lo2] (see Definition 1.5).\bigskip

\noindent $\bullet$ We set $\hbox{\rm pot}(\bormz)\! :=\!\{\Pi_{i\in d}\ A_i\mid 
A_i\!\in\!\borel (\omega^\omega )\mbox{, and }A_i\! =\!\omega^\omega\mbox{ for almost every }i\!\in\! d\}$. We also set $W\! :=\! W^{(\omega^\omega )^d}$ and $C\! :=\! C^{(\omega^\omega )^d}$ (see Lemma 4.2.1). We will define specifically, for $\xi\! <\!\omega_1$,
$$\big\{ (\beta ,\gamma )\!\in\!\omega^\omega\!\times\! W\mid 
\beta\mbox{ codes a }\mbox{pot}(\bormxi)\mbox{ set and }C_\gamma
\mbox{ is the set coded by }\beta\big\}.$$ 
The way we will do it is not the simplest possible (we can in fact forget $\beta$, and work with $\gamma$ integer instead of real, see [L7]). We do it this way to start to give the flavor of what is going on with the Wadge classes.

\vfill\eject

\noindent $\bullet$ To do this, we set\bigskip
 
\leftline{$V_0\! :=\!
\Bigg\{\ (\beta ,\gamma )\!\in\!\omega^\omega\!\times\! W\mid
\forall i\! <\!\beta (0)\ \ (\beta^*)_i\!\in\! W^{\omega^\omega}\wedge\gamma\!\in\!\Borel (\beta )
\wedge$}\smallskip
 
\rightline{$\left[\!\!\!\!\!\!
\begin{array}{ll}
& \beta (0)\! =\! d\wedge C_\gamma\! =\!
\Pi_{i<\beta (0)}\ C^{\omega^\omega}_{(\beta^*)_i}\mbox{ if }d\! <\!\omega\cr\cr
& C_\gamma\! =\!
\left(\Pi_{i<\beta (0)}\ C^{\omega^\omega}_{(\beta^*)_i}\right)\!\times\! (\omega^\omega )^\omega
\mbox{ if }d\! =\!\omega
\end{array}
\right.
\Bigg\}.$}\bigskip

 We define an inductive operator $\Phi$ over $\omega^\omega\!\times\!\omega^\omega$ (see [C]) as follows:\bigskip

\leftline{$\Phi (A)\! :=\! A\cup V_0\ \cup\big\{  (\beta ,\gamma )\!\in\!\omega^\omega\!\times\! W\mid
\gamma\!\in\!\Borel (\beta )\wedge$}\smallskip

\rightline{$\exists\gamma'\!\in\!\Borel (\beta )\ ~\forall n\!\in\!\omega ~\ \ \big( (\beta )_n,(\gamma')_n\big)
\!\in\! A\wedge\neg C_\gamma\! =\!\bigcup_{n\in\omega}\ C_{(\gamma')_n}\big\}.$}\bigskip
\noindent Then $\Phi$ is clearly a $\Ca$ monotone inductive operator. We set, for any ordinal $\xi$, 
$V_\xi\! :=\!\Phi^\xi$ (which is coherent with the definition of $V_0$). We also set 
${V_{<\xi}\! :=\!\bigcup_{\eta <\xi}~V_{\eta}}$. The effective version of Theorem 4.1, which is the specific version of Theorem 1.9 for the Borel classes, is as follows:

\begin{them} Let $T_d$ be a tree with $\Borel$ suitable levels, 
$1\!\leq\!\xi\! <\!\omega_1^{\mbox{CK}}$, and $A_0$, $A_1$ disjoint $\Ana$ subsets of $(\omega^\omega )^d$.\medskip

\noindent (1) Assume that $S\!\in\!\boraxi (\lceil T_d\rceil )$ is not separable from 
$\lceil T_d\rceil \!\setminus\! S$ by a $\mbox{pot}(\bormxi )$ set. Then the following are equivalent:\smallskip

\noindent (a) The set $A_0$ is not separable from $A_1$ by a $\hbox{\it pot}(\bormxi )$ set.\smallskip

\noindent (b) The set $A_0$ is not separable from $A_1$ by a $\Borel\cap\mbox{pot}(\bormxi )$ set.\smallskip

\noindent (c) $\neg\big(\exists (\beta ,\gamma )\!\in\! (\Borel\!\times\!\Borel )\cap V_\xi\ \ A_0\!\subseteq\! 
C_\gamma\!\subseteq\!\neg A_1\big)$.\smallskip 

\noindent (d) The set $A_0$ is not separable from $A_1$ by a $\bormxi (\tau_1)$ set.\smallskip

\noindent (e) $\overline{A_0}^{\tau_\xi}\cap A_1\!\not=\!\emptyset$.\smallskip

\noindent (f) The inequality $\big( (d^\omega )_{i\in d}, S,\lceil T_d\rceil\!\setminus\! S\big)\leq
\big( (\omega^\omega )_{i\in d}, A_0, A_1\big)$ holds.\bigskip

\noindent (2) The sets $V_{\xi}$ and $V_{<\xi}$ are $\Ca$.\bigskip

\noindent (3) Assume that $S^0,S^1\!\in\!\boraxi (\lceil T_d\rceil )$ are disjoint and not separable by a 
$\mbox{pot}(\borxi )$ set. Then the following are equivalent:\smallskip

\noindent (a) The set $A_0$ is not separable from $A_1$ by a $\hbox{\it pot}(\borxi )$ set.\smallskip

\noindent (b) The set $A_0$ is not separable from $A_1$ by a $\Borel\cap\hbox{\it pot}(\borxi )$ set.\smallskip

\noindent (c) $\neg\big(\exists (\beta ,\gamma ), (\beta',\gamma')\!\in\! (\Borel\!\times\!\Borel )\cap V_\xi\ \ 
C_{\gamma'}\! =\!\neg C_\gamma\ \mbox{ and }\ 
A_0\!\subseteq\! C_\gamma\!\subseteq\!\neg A_1\big)$.\smallskip 

\noindent (d) The set $A_0$ is not separable from $A_1$ by a $\borxi (\tau_1)$ set.\smallskip

\noindent (e) $\overline{A_0}^{\tau_\xi}\cap\overline{A_1}^{\tau_\xi}\!\not=\!\emptyset$.\smallskip

\noindent (f) The inequality $\big( (d^\omega )_{i\in d}, S^0,S^1\big)\leq
\big( (\omega^\omega )_{i\in d}, A_0, A_1\big)$ holds.\end{them}

 The proofs of Theorems 4.1 and 4.2.2 will be by induction on $\xi$. This appears in the statement of the following lemma.

\begin{lemm} (1) The set $V_{0}$ is $\Ca$.\medskip

\noindent (2) Let $1\!\leq\!\xi\! <\!\omega^{\hbox{\it CK}}_{1}$. We assume that Theorem 4.2.2 is proved for $\eta\! <\!\xi$.\smallskip

\noindent (a) The set $V_{<\xi}$ is $\Ca$.\smallskip

\noindent (b) Fix $A\!\in\!\Ana\big( (\omega^\omega )^d\big)$. Then 
$\overline{A}^{\tau_{\xi}}\!\in\!\Ana\big( (\omega^\omega )^d)$.\smallskip

\noindent (c) Let $n\!\geq\! 1$, $1\!\leq\!\xi_{1}\! <\!\xi_{2}\! <\!\ldots\! <\!
\xi_{n}\!\leq\!\xi$, and $S_{1}$, $\ldots$, $S_{n}$ be $\Ana$ sets. 
Assume that $S_{i}\!\subseteq\!\overline{S_{i+1}}^{\tau_{\xi_{i}+1}}$ for 
$1\!\leq\! i\! <\! n$. Then 
${S_{n}\cap\bigcap_{1\leq i<n}~\overline{S_{i}}^{\tau_{\xi_{i}}}}$ is 
$\tau_{1}$-dense in $\overline{S_{1}}^{\tau_{1}}$.\end{lemm}

\noindent\bf Proof.\rm ~(1) The set $V_{0}$ is clearly $\Ca$.\bigskip

\noindent (2).(a) The proof is contained in the proof of Theorem 4.1 in [L7]. It is a consequence of  Lemma 4.8 in [C].\bigskip

\noindent (b) The proof is essentially the proof of Lemma 2.2.2.(a) in [L7].\bigskip

\noindent (c) The proof is essentially the proof of Lemma 2.2.2.(b) in [L7].\hfill{$\square$}

\begin{lemm} Let $S,T\!\in\!\Ana\big( (\omega^\omega )^d\big)$ such that $S$ is $\tau_1$-dense in $T$, 
$(X_i)_{i\in d}$ a sequence of $\Ana$ subsets of $\omega^\omega$ such that 
$X_i\! =\!\omega^\omega$ if $i\!\geq\! i_0$. Then $S\cap (\Pi_{i\in d}\ X_i)$ is $\tau_1$-dense in 
$T\cap (\Pi_{i\in d}\ X_i)$.\end{lemm}

\noindent\bf Proof.\rm\ Let $(\Delta_i)_{i\in d}$ be a sequence of $\Borel$ subsets of 
$\omega^\omega$ such that $\Delta_i\! =\!\omega^\omega$ if $i\!\geq\! j_0\!\geq\! i_0$, and also 
$T\cap (\Pi_{i\in d}\ I_i)\!\not=\!\emptyset$, where $I_i\! :=\! X_i\cap\Delta_i$. We have to see that 
$S\cap (\Pi_{i\in d}\ I_i)\!\not=\!\emptyset$. We argue by contradiction. This gives a sequence 
$(D_i)_{i\in d}$ of $\Borel$ subsets of $\omega^\omega$ such that $I_i\!\subseteq\! D_i$ if 
$i\!\in\! d$, and $S\cap (\Pi_{i\in d}\ D_i)\! =\!\emptyset$, by $j_0$ applications of the separation theorem. But $T\cap (\Pi_{i\in d}\ D_i)\!\not=\!\emptyset$, and $D_i\! =\!\omega^\omega$ if 
$i\!\geq\! j_0$. So $S\cap (\Pi_{i\in d}\ D_i)\!\not=\!\emptyset$, by $\tau_1$-density of $S$ in $T$, which is absurd.\hfill{$\square$}

\subsection{$\!\!\!\!\!\!$ Representation of Borel sets}\indent

 Now we come to the representation theorem of Borel sets by G. 
Debs and J. Saint Raymond (see [D-SR]). It specifies the classical 
result of Lusin asserting that any Borel set in a Polish space is the 
bijective continuous image of a closed subset of the Baire space. The material in this Subsection can be found in Subsection 2.3 of [L7], but we recall most of it since it will be used iteratively in the case of Wadge classes. The following definition can be found in [D-SR].

\begin{defin} (Debs-Saint Raymond) Let $c$ be a countable set. A partial order relation $R$ on 
$c^{<\omega}$ is a $tree~relation$ if, for ${t\!\in\! c^{<\omega}}$,\smallskip

\noindent (a) $\emptyset ~R~t$.\smallskip

\noindent (b) The set $P_{R}(t)\! :=\!\{s\!\in\! c^{<\omega}\mid s~R~t\}$ is finite 
and linearly ordered by $R$.\smallskip

 For instance, the non strict extension relation $\subseteq$ is a tree relation.\smallskip

\noindent $\bullet$ Let $R$ be a tree relation. An $R\mbox{-}branch$ 
is an $\subseteq$-maximal subset of $c^{<\omega}$ linearly ordered by $R$. We 
denote by $[R]$ the set of all infinite $R$-branches.\smallskip

 We equip $(c^{<\omega})^\omega$ with the product of the discrete topology on 
$c^{<\omega}$. If $R$ is a tree relation, then the space 
$[R]\!\subseteq\! (c^{<\omega})^\omega$ is equipped with the topology induced by 
that of $(c^{<\omega})^\omega$. The map 
$h\! :\! c^\omega\!\rightarrow\! [\subseteq ]$ defined by 
$h(\gamma )\! :=\! (\gamma\vert j)_{j\in\omega}$ is an homeomorphism.\smallskip

\noindent $\bullet$ Let $R$, $S$ be tree relations with $R\!\subseteq\! S$. The 
$canonical\ map$ $\Pi\! :\! [R]\!\rightarrow\! [S]$ is defined by
$$\Pi ({\cal B})\! :=\!\hbox{\it the unique $S$-branch containing ${\cal B}$.}$$
$\bullet$ Let $S$ be a tree relation. We say that $R\!\subseteq\! S$ is 
$distinguished$ in $S$ if
$$\forall s,t,u\!\in\! c^{<\omega}\ \ \left.
\begin{array}{ll}
& s~S~t~S~u\cr
& \ \ s~R~u
\end{array}\!\!\right\}\Rightarrow s~R~t.$$
For example, let $C$ be a closed subset of $c^\omega$, and define
$$s~R~t\ \Leftrightarrow\ s\!\subseteq\! t\wedge N_s\cap C\!\not=\!\emptyset .$$ 
Then $R$ is distinguished in $\subseteq$.\smallskip

\noindent $\bullet$ Let $\eta\! <\!\omega_{1}$. A family $(R^{(\rho )})_{\rho\leq\eta}$ of tree 
relations is a $resolution~family$ if\smallskip

\noindent (a) $R^{(\rho +1)}$ is a distinguished subtree of $R^{(\rho )}$, for 
all $\rho\! <\!\eta$.\smallskip

\noindent (b) $R^{(\lambda )}\! =\!\bigcap_{\rho <\lambda}~R^{(\rho )}$, for 
all limit $\lambda\!\leq\!\eta$.\end{defin}

 We will use the following extension of the property of distinction:

\begin{lemm} Let $\eta\! <\!\omega_{1}$, $(R^{(\rho )})_{\rho\leq\eta}$ a resolution family, and 
$\rho\! <\!\eta$. Assume that $s~R^{(0)}~s'~R^{(\rho )}~s''$ and $s~R^{(\rho +1)}~s''$. Then 
$s~R^{(\rho +1)}~s'$.\end{lemm}

\noindent\bf Notation.\rm ~Let $\eta\! <\!\omega_{1}$, 
$(R^{(\rho )})_{\rho\leq\eta}$ a resolution family such that $R^{(0)}$ is a subrelation of $\subseteq$, 
$\rho\!\leq\!\eta$ and $t\!\in\! c^{<\omega}\!\setminus\!\{\emptyset\}$. We set
$$t^{\rho}\! :=\! t~\vert ~\hbox{\rm max}\{ r\! <\! |t|\mid t\vert r~R^{(\rho )}~t\}.$$
We enumerate $\{ t^{\rho}\mid\rho\!\leq\!\eta\}$ by 
$\{ t^{\xi_{i}}\mid 1\!\leq\! i\!\leq\! n\}$, where $1\!\leq\! n\!\in\!\omega$ and 
${\xi_{1}\! <\!\ldots\! <\!\xi_{n}\! =\!\eta}$. We can write 
$t^{\xi_{n}}\!\subset_{\not=}\! 
t^{\xi_{n-1}}\!\subset_{\not=}\!\ldots\!\subset_{\not=}\! t^{\xi_{2}}\!\subset_{\not=}
\! t^{\xi_{1}}\!\subset_{\not=}\! t$. By Lemma 4.3.2 we have 
$t^{\xi_{i+1}}~R^{(\xi_{i}+1)}~t^{\xi_{i}}$ for each $1\!\leq\! i\! <\! n$.

\begin{lemm} Let $\eta\! <\!\omega_{1}$, $(R^{(\rho )})_{\rho\leq\eta}$ a resolution family such that 
$R^{(0)}$ is a subrelation of $\subseteq$, $t$ in $c^{<\omega}\!\setminus\!\{\emptyset\}$ and 
$1\!\leq\! i\! <\! n$.\smallskip

\noindent (a) Set $\eta_{i}\! :=\!\{\rho\!\leq\!\eta\mid t^{\xi_{i}}\!\subseteq\! t^{\rho}\}$. Then $\eta_{i}$ is a successor ordinal.\smallskip

\noindent (b) We may assume that $t^{\xi_{i}+1}\!\subset_{\not=}\! t^{\xi_{i}}$.\end{lemm}

 The following is part of Theorem I-6.6 in [D-SR].

\begin{them} (Debs-Saint Raymond) Let $\eta\! <\!\omega_{1}$, $R$ a tree relation,  $(I_n)_{n\in\omega}$ a sequence of $\bormep$ subsets of $[R]$. Then there is a resolution family $(R^{(\rho )})_{\rho\leq\eta}$ with\smallskip

\noindent (a) $R^{(0)}=R$.\smallskip

\noindent (b) The canonical map $\Pi\! :\! [R^{(\eta )}]\!\rightarrow\! [R]$ is a continuous bijection.\smallskip

\noindent (c) The set $\Pi^{-1}(I_n)$ is a closed subset of $[R^{(\eta )}]$ for each integer $n$.\end{them}

 Now we come to the actual proof of Theorem 4.1. 

\subsection{$\!\!\!\!\!\!$ Proof of Theorem 4.1}\indent

 The next result is essentially Theorem 2.4.1 in [L7]. But we give its proof since it is the basis for further generalizations.

\begin{them} Let $T_d$ be a tree with $\Borel$ suitable levels, 
$\xi\! <\!\omega^{\hbox{\it CK}}_{1}$ a successor ordinal, $S\!\in\!\boraxi (\lceil T_d\rceil )$, and $A_0$, $A_1$ disjoint $\Ana$ subsets of $(\omega^\omega )^d$. We assume that Theorem 4.2.2 is proved for $\eta\! <\!\xi$. Then one of the following holds:\smallskip

\noindent (a) $\overline{A_0}^{\tau_\xi}\cap A_1\! =\!\emptyset$.\smallskip

\noindent (b) The inequality $\big( (\Pi_i''\lceil T_d\rceil )_{i\in d}, S,\lceil T_d\rceil\!\setminus\! S\big)\leq
\big( (\omega^\omega )_{i\in d}, A_0, A_1\big)$ holds.\end{them}

\noindent\bf Proof.\rm ~Fix $\eta\! <\!\omega^{\hbox{\rm CK}}_{1}$ with $\xi\! =\!\eta\! +\! 1$.\bigskip

\noindent $\bullet$ Recall the finite sets $c_l$ defined at the end of the proof of Proposition 2.2 (we only used the fact that $T_d$ has finite levels to see that they are finite). With the notation of Definition 4.3.1, we put $c\! :=\!\bigcup_{l\in\omega}\ c_l$, so that $c$ is countable. The set 
$I\!:=\! h[\lceil T_d\rceil\!\setminus\! S]$ is a $\bormep$ subset of $[\subseteq ]$. Theorem 4.3.4 provides a resolution family. We put
$$D\! :=\!\big\{ \vec s \!\in\!T_d\mid\vec s\! =\!\vec\emptyset\vee 
\exists {\cal B}\!\in\!\Pi^{-1}(I)\ \ \vec s\!\in\! {\cal B}\big\}.$$
$\bullet$ Assume that $\overline{A_0}^{\tau_\xi}\cap A_1$ is not empty. Recall that $(\Omega_{X},{\it\Sigma}_{X})$ is a Polish space (see the notation at the beginning of Section 4.2). We fix a complete metric $d_X$ on $(\Omega_{X},{\it\Sigma}_{X})$.\bigskip

\noindent $\bullet$ We construct\smallskip

\noindent - $(\alpha^i_{s})_{i\in d,s\in\Pi_i''T_d}\!\subseteq\!\omega^\omega$,\smallskip

\noindent - $(O^i_{s})_{i\leq\vert s\vert ,i\in d,s\in\Pi_i''T_d}\!\subseteq\!\Ana (\omega^\omega )$,\smallskip

\noindent - $(U_{\vec s  })_{\vec s \in T_d}\!\subseteq\!
\Ana\big((\omega^\omega )^d\big)$.\bigskip

 We want these objects to satisfy the following conditions.
$$\begin{array}{ll}
& \!\!\!\! (1)\ \alpha^i_{s}\!\in\! O^i_{s}\!\subseteq\!\Omega_{\omega^\omega}\wedge 
(\alpha^i_{s_i })_{i\in d}\!\in\! U_{\vec s  }\!\subseteq\!
\Omega_{(\omega^\omega )^d}\mbox{,}\cr\cr 
& \!\!\!\! (2)\ O^i_{sq}\!\subseteq\! O^i_{s}\mbox{,}\cr\cr
& \!\!\!\! (3)\ \mbox{diam}_{d_{\omega^\omega}}(O^i_{s})\!\leq\! 2^{-\vert s\vert}\wedge
\mbox{diam}_{d_{(\omega^\omega )^d}}(U_{\vec s  })\!\leq\! 2^{-\vert\vec s\vert}\mbox{,}\cr\cr
& \!\!\!\! (4)\ U_{\vec s  }\!\subseteq\! \overline{A_0}^{\tau_\xi}\cap A_1\mbox{ if }\vec s \!\in\! D\mbox{,}
\cr\cr  
& \!\!\!\! (5)\ U_{\vec s  }\!\subseteq\! A_0\mbox{ if }\vec s \!\notin\! D\mbox{,}\cr\cr
& \!\!\!\! (6)\ \big( 1\!\leq\!\rho\!\leq\!\eta\wedge\vec s~R^{(\rho )}~\vec t\ \big)\Rightarrow 
U_{\vec t }\!\subseteq\!\overline{U_{\vec s  }}^{\tau_{\rho }}\mbox{,}\cr\cr
& \!\!\!\! (7)\ \big( (\vec s,\vec t\!\in\! D\vee\vec s,\vec t \!\notin\! D)\wedge
\vec s ~R^{(\eta )}~\vec t\ \big)\Rightarrow U_{\vec t }\!\subseteq\! U_{\vec s }.
\end{array}$$
$\bullet$ Let us prove that this construction is sufficient to get the theorem.\bigskip

\noindent - Fix $\vec\beta\!\in\! \lceil T_d\rceil $. Then we can define  
${(j_{k})_{k\in\omega}\! :=\! (j^{\vec\beta}_{k})_{k\in\omega}}$ by 
${\Pi^{-1}\big( (\vec\beta\vert j)_{j\in\omega}\big)\! =\! (\vec\beta\vert j_{k})_{k\in\omega}}$, with the inequalities $j_{k}\! <\! j_{k+1}$. In particular, ${\vec\beta\vert j_{k}~R^{(\eta )}~\vec\beta\vert j_{k+1}}$. We have 
$$\vec\beta\!\notin\! S\ \Leftrightarrow ~
h(\vec\beta )\! =\! (\vec\beta\vert j)_{j\in\omega}\!\in\! I\ \Leftrightarrow ~
(\vec\beta\vert j_{k})_{k\in\omega}\!\in\!\Pi^{-1}(I)\ \Leftrightarrow ~
\forall k\!\geq\! k_0\! :=\! 0~~\vec\beta\vert j_{k}\!\in\! D$$ 
since $\Pi^{-1}(I)$ is a closed subset of $[R^{(\eta )}]$. Similarly, 
$\vec\beta\!\in\! S$ is equivalent to the existence of $k_{0}\!\in\!\omega$ such that 
$\vec\beta\vert j_{k}\!\notin\! D$ for each $k\!\geq\! k_{0}$.\bigskip

 This implies that $(U_{\vec\beta\vert j_k})_{k\geq k_0}$ is a non-increasing sequence of nonempty clopen subsets of the space $(\Omega_{(\omega^\omega)^d} ,{\it\Sigma}_{(\omega^\omega)^d})$ whose 
$d_{(\omega^\omega)^d}$-diameters tend to zero, and we can define 
$${\{ {\cal F}(\vec\beta )\}\! :=\!\bigcap_{k\geq k_0}~U_{\vec\beta\vert j_k}
\!\subseteq\!\Omega_{(\omega^\omega)^d}}.$$ 
Note that ${\cal F}(\vec\beta )$ is the limit of 
$\big( (\alpha^i_{\beta_i\vert j_k})_{i\in d}\big)_{k\in\omega}$.\bigskip

\noindent - Now let $\gamma\!\in\!\Pi_i''\lceil T_d\rceil$, and $\vec\beta\!\in\! \lceil T_d\rceil $ such that 
$\beta_i\! =\!\gamma$. We set $f_i(\gamma )\! :=\! {\cal F}_i(\vec\beta )$. This defines 
$f_i\! :\!\Pi_i''\lceil T_d\rceil\!\rightarrow\!\omega^\omega$.\bigskip

 Note that $f_i(\gamma )$ is the limit of $(\alpha^i_{\gamma\vert j})_{j\in\omega}$. Indeed, 
$f_i(\gamma )$ is the limit of $(\alpha^i_{\gamma\vert j_k})_{k\in\omega}$. If $j\!\geq\! i$, then 
$\alpha^i_{\gamma\vert j}\!\in\! O^i_{\gamma\vert j}$, and the sequence 
$(O^i_{\gamma\vert j})_{j\geq i}$ is decreasing. Fix $\varepsilon\! >\! 0$, 
$k\!\geq\! i$ such that $2^{-k}\! <\!\varepsilon$. Then we get, if $j\!\geq\! k$, 
$d_{\omega^\omega}\big( f_i(\gamma ),\alpha^i_{\gamma\vert j}\big)\!\leq\! 
\mbox{diam}_{d_{\omega^\omega}}(O^i_{\gamma\vert j})\!\leq\! 2^{-j}\!\leq\! 2^{-k}\! <\!\varepsilon$. 
In particular, $f_i(\gamma )$ does not depend on the choice of $\vec\beta$. This also proves that $f_i$ is continuous on $\Pi_i''\lceil T_d\rceil$.\bigskip

\noindent - Note that ${\cal F}_i(\vec\beta )$ is the limit of some subsequence of 
$(\alpha^i_{\beta_i\vert j})_{j\in\omega}$, by continuity of the projections. Thus 
${\cal F}_i(\vec\beta )\! =\! f_i(\beta_i)$, and ${\cal F}(\vec\beta )\! =\! (\Pi_{i\in d}\ f_i)(\vec\beta )$. This implies that the inclusions $S\!\subseteq\! (\Pi_{i\in d}\ f_i)^{-1}(A_0)$ and 
$\lceil T_d\rceil \!\setminus\! S\!\subseteq\! (\Pi_{i\in d}\ f_i)^{-1}(A_1)$ hold.\bigskip

\noindent $\bullet$ So let us prove that the construction is possible.\bigskip

\noindent - Let $(\alpha^i_{\emptyset})_{i\in d}\!\in\!\overline{A_0}^{\tau_\xi}\cap A_1\cap
\Omega_{(\omega^\omega )^d}$, which is nonempty since 
$\overline{A_0}^{\tau_\xi}\cap A_1\!\not=\!\emptyset$ is $\Ana$, by Lemma 4.2.3.(2).(b). Then we choose a 
$\Ana$ subset $U_{\vec\emptyset}$ of $(\omega^\omega )^d$, with $d_{(\omega^\omega )^d}$-diameter at most $1$, such that
$$(\alpha^i_{\emptyset})_{i\in d}\!\in\! U_{\vec\emptyset}\!\subseteq\! 
\overline{A_0}^{\tau_\xi}\cap A_1\cap\Omega_{(\omega^\omega )^d}.$$
We choose a $\Ana$ subset $O^0_{\emptyset}$ of $\omega^\omega$, with $d_{\omega^\omega}$-diameter at most $1$, with 
$\alpha^0_{\emptyset}\!\in\! O^0_{\emptyset}\!\subseteq\!\Omega_{\omega^\omega}$, which is possible since $\Omega_{(\omega^\omega )^d}\!\subseteq\!\Omega_{\omega^\omega}^d$. Assume that 
$(\alpha^i_{s})_{|s|\leq l}$, $(O^i_{s})_{|s|\leq l}$ and $(U_{\vec s  })_{|\vec s|\leq l}$ satisfying conditions (1)-(7) have been constructed, which is the case for $l\! =\! 0$.

\vfill\eject

\noindent - Let $\overrightarrow{tm}\!\in\!T_d\cap (d^{l+1})^d$. Note that 
$\overrightarrow{tm}^\eta\!\in\! D$ if $\overrightarrow{tm}^\eta\!\in\! D$ is not equivalent to  
$\overrightarrow{tm}\!\in\! D$ (see the notation before Lemma 4.3.3).\bigskip

\noindent - The conclusions in the assertions (a), (b) and (c) of the following claim do not really depend on their respective assumptions, but we will use these assertions later in this form. We define 
$X_i\! :=\! O^i_{t_i}$ if $i\!\leq\! l$, and $\omega^\omega$ if $i\! >\! l$.\bigskip

\noindent\bf Claim.\rm ~Assume that $\eta\! >\! 0$.\medskip

\noindent (a) The set $A_0\cap\bigcap_{1\leq\rho\leq\eta}~\overline{U_{\overrightarrow{tm}^{\rho}}}^{\tau_{\rho}}\cap (\Pi_{i\in d}\ X_i)$ is $\tau_{1}$-dense in 
$\overline{U_{\overrightarrow{tm}^1}}^{\tau_{1}}\cap (\Pi_{i\in d}\ X_i)$ if 
$\overrightarrow{tm}^{\eta}\!\in\! D$ and $\overrightarrow{tm}\!\notin\! D$.\medskip

\noindent (b) The set $U_{\overrightarrow{tm}^{\eta}}\cap\bigcap_{1\leq\rho <\eta}~
\overline{U_{\overrightarrow{tm}^{\rho}}}^{\tau_{\rho}}\cap (\Pi_{i\in d}\ X_i)$ is $\tau_{1}$-dense in $\overline{U_{\overrightarrow{tm}^1}}^{\tau_{1}}\cap (\Pi_{i\in d}\ X_i)$ if 
$\overrightarrow{tm}^{\eta},\overrightarrow{tm}\!\in\! D$ or $\overrightarrow{tm}^{\eta},\overrightarrow{tm}\!\notin\! D$.\bigskip

 Indeed, let us forget $\Pi_{i\in d}\ X_i$ for the moment. We may assume that 
$\overrightarrow{tm}^{\xi_{i}+1}\!\subset_{\not=}\! \overrightarrow{tm}^{\xi_{i}}$ if $1\!\leq\! i\! <\! n$, by Lemma 4.3.3. We set ${S_{i}\! :=\! U_{\overrightarrow{tm}^{\xi_{i}}}}$, for ${1\!\leq\!\xi_{i}\!\leq\!\eta}$. As 
$\overrightarrow{tm}^{\xi_{i+1}}\ R^{(\xi_i+1)}\ \overrightarrow{tm}^{\xi_i}$, we can write 
$S_i\!\subseteq\!\overline{S_{i+1}}^{\tau_{\xi_i+1}}$, for $1\!\leq\!\xi_{i}\! <\!\eta$, by induction assumption. If $\overrightarrow{tm}^{\eta}\!\in\! D$ and 
$\overrightarrow{tm}\!\notin\! D$, then 
$S_n\!\subseteq\!\overline{A_0}^{\tau_{\eta +1}}$. Thus 
$A_0\cap\bigcap_{1\leq\xi_i\leq\eta}~\overline{U_{\overrightarrow{tm}^{\xi_i}}}^{\tau_{\xi_i}}$ 
and $U_{\overrightarrow{tm}^{\eta}}\cap\bigcap_{1\leq\xi_i<\eta}~
\overline{U_{\overrightarrow{tm}^{\xi_i}}}^{\tau_{\xi_i}}$ are $\tau_{1}$-dense in 
$\overline{U_{\overrightarrow{tm}^1}}^{\tau_{1}}$, by Lemma 4.2.3.(2).(c).\bigskip

 But if $1\!\leq\!\rho\!\leq\!\eta$, then there is $1\!\leq\! i\!\leq\! n$ with 
$\overrightarrow{tm}^{\rho}\! =\! \overrightarrow{tm}^{\xi_i}$. And $\rho\!\leq\!\xi_i$ since we have $\overrightarrow{tm}^{\xi_{i}+1}\!\subset_{\not=}\! \overrightarrow{tm}^{\xi_{i}}$ if ${1\!\leq\! i\! <\! n}$. 
We are done since 
$\bigcap_{1\leq\rho\leq\eta}~\overline{U_{\overrightarrow{tm}^{\rho}}}^{\tau_{\rho}}\! =\! 
\bigcap_{1\leq\xi_i\leq\eta}~\overline{U_{\overrightarrow{tm}^{\xi_i}}}^{\tau_{\xi_i}}$ and 
$$U_{\overrightarrow{tm}^{\eta}}\cap\bigcap_{1\leq\rho <\eta}~
\overline{U_{\overrightarrow{tm}^{\rho}}}^{\tau_{\rho}}\! =\! 
U_{\overrightarrow{tm}^{\eta}}\cap\bigcap_{1\leq\xi_i <\eta}~
\overline{U_{\overrightarrow{tm}^{\xi_i}}}^{\tau_{\xi_i}}$$
The claim now comes from Lemma 4.2.4.\hfill{$\diamond$}\bigskip

\noindent - Let ${\cal X}\! :=\! d^{l+1}$. The map 
$\Theta\! :\! {\cal X}^d\!\rightarrow\!\Ana\big( (\omega^\omega )^d\big)$ is defined on ${\cal T}^{l+1}$ by
$$\Theta (\overrightarrow{tm})\! :=\!\left\{\!\!\!\!\!\!\!
\begin{array}{ll} 
& A_0\cap\bigcap_{1\leq\rho\leq\eta}~\overline{U_{\overrightarrow{tm}^{\rho}}}^{\tau_{\rho}}\cap (\Pi_{i\in d}\ X_i)\cap\Omega_{(\omega^\omega )^d}\mbox{ if }\overrightarrow{tm}^{\eta}\!\in\! D
\wedge\overrightarrow{tm}\!\notin\! D\mbox{,}\cr\cr 
& U_{\overrightarrow{tm}^{\eta}}\cap\bigcap_{1\leq\rho <\eta}~
\overline{U_{\overrightarrow{tm}^{\rho}}}^{\tau_{\rho}}\cap (\Pi_{i\in d}\ X_i)\mbox{ if }
\overrightarrow{tm}^{\eta},\overrightarrow{tm}\!\in\! D\vee
\overrightarrow{tm}^{\eta},\overrightarrow{tm}\!\notin\! D.
\end{array}\right.$$
By the claim, $\Theta (\overrightarrow{tm})$ is $\tau_{1}$-dense in 
$\overline{U_{\overrightarrow{tm}^1}}^{\tau_{1}}\cap (\Pi_{i\in d}\ X_i)$ if $\eta\! >\! 0$. As 
$\overrightarrow{tm}^1\!\subseteq\! \vec t \!\subseteq\!\overrightarrow{tm}$ 
and $R^{(1)}$ is distinguished in $\subseteq$ we get 
$\overrightarrow{tm}^1\ R^{(1)}\ \vec t$ and 
$U_{\vec t }\!\subseteq\!\overline{U_{\overrightarrow{tm}^1}}^{\tau_1}$, by induction 
assumption. Therefore $\overline{U_{\vec t }}^{\tau_1}\cap (\Pi_{i\in d}\ X_i)\!\subseteq\!
\overline{U_{\overrightarrow{tm}^1}}^{\tau_1}\cap (\Pi_{i\in d}\ X_i)\!\subseteq\!
\overline{\Theta} (\overrightarrow{tm})$, and 
$(\alpha^{i}_{t_i})_{i\in d}\!\in\! U_{\vec t }\cap (\Pi_{i\in d}\ X_i)\!\subseteq\!
\overline{\Theta} (\overrightarrow{tm})$ (even if $\eta\! =\! 0$). Therefore $\overline{\Theta}$ admits a 
$\pi$-selector on ${\cal T}^{l+1}$. Indeed, we define, for each $i\!\in\! d$, 
$\overline{\theta}_i\! :\! {\cal X}\!\rightarrow\!\omega^\omega$ by 
$\overline{\theta}_i(t_im_i)\! :=\!\alpha^i_{t_i}$ if $t_i\!\in\!\Pi_i''T_d$, $0^\infty$ otherwise.\bigskip
 
\noindent - As $T_d$ is a tree with $\Borel$ suitable levels, we can apply Lemma 4.1.3. Thus 
$\Theta$ admits a $\pi$-selector $\theta$ on ${\cal T}^{l+1}$. We set, for $s\!\in\!\Pi_i[{\cal T}^{l+1}]$, 
$\alpha^{i}_{s}\! :=\!\theta_{i}(s)$.\bigskip

\noindent - We choose $\Ana$ sets $U_{\overrightarrow{tm}}$ with $d_{(\omega^\omega )^d}$-diameter at most $2^{-l-1}$ such that $\theta (\overrightarrow{tm})\!\in\! U_{\overrightarrow{tm}}\!\subseteq\!
\Theta (\overrightarrow{tm})$ if $\overrightarrow{tm}\!\in\! {\cal T}^{l+1}$.\bigskip

\noindent - Finally, we choose the $O^i_{sq}$'s. We first prove that $\alpha^i_{sq}\!\in\! O^i_{s}$ if 
$sq\!\in\!\Pi_i[{\cal T}^{l+1}]$, $i\!\in\! d$ and $i\!\leq\! l$.\bigskip

 Let $\overrightarrow{tm}\!\in\! {\cal T}^{l+1}$ such that $sq\! =\! t_im_i$. Then
$\alpha^i_{sq}\! =\!\theta_i(sq)\! =\!\theta_i(t_im_i)$. As 
$\theta (\overrightarrow{tm})\!\in\!\Theta (\overrightarrow{tm})$ and $i\!\leq\! l$, we get 
$\alpha^i_{sq}\!\in\! O^i_{t_i}\! =\! O^i_{s}$.\bigskip
  
 Now we can define the $O^i_{sq}$'s. If $sq\!\in\!\Pi_i[{\cal T}^{l+1}]$, then we choose a $\Ana$ set 
$O^i_{sq}$, with $d_{\omega^\omega}$-diameter at most $2^{-l-1}$, such that 
$$\alpha^i_{sq}\!\in\! O^i_{sq}\!\subseteq\!\left\{\!\!\!\!\!\!\! 
\begin{array}{ll}
& O^i_{s}~\mbox{ if }~i\!\leq\! l\mbox{,}\cr
& \Omega_{\omega^\omega}\mbox{ otherwise.}
\end{array}
\right.$$ 
- This finishes the proof since $\vec u\ R^{(\rho )}\ \overrightarrow{tm}\ \mbox{ and }\ 
\vec u\!\not=\!\overrightarrow{tm}\ \Rightarrow\ \vec u\ R^{(\rho )}\ \overrightarrow{tm}^\rho\ R^{(\rho )}\ 
\overrightarrow{tm}$, by Lemma 4.3.2.\hfill{$\square$}\bigskip

 Now we come to the ambiguous classes.

\begin{them} Let $T_d$ be a tree with $\Borel$ suitable levels, 
$\xi\! <\!\omega^{\hbox{\it CK}}_{1}$ a successor ordinal, $S^0, S^1$ in 
$\boraxi (\lceil T_d\rceil )$ disjoint, and $A_0$, $A_1$ disjoint $\Ana$ subsets of 
$(\omega^\omega )^d$. We assume that Theorem 4.2.2 is proved for $\eta\! <\!\xi$. Then one of the following holds:\smallskip

\noindent (a) $\overline{A_0}^{\tau_\xi}\cap\overline{A_1}^{\tau_\xi}\! =\!\emptyset$.\smallskip

\noindent (b) The inequality $\big( (\Pi_i''\lceil T_d\rceil )_{i\in d}, S^0,S^1\big)\leq
\big( (\omega^\omega )_{i\in d}, A_0, A_1\big)$ holds.\end{them}

\noindent\bf Proof.\rm ~Let us indicate the differences with the proof of Theorem 4.4.1. Assume that 
$\overline{A_0}^{\tau_\xi}\cap\overline{A_1}^{\tau_\xi}\!\not=\!\emptyset$. We set 
$I^\varepsilon\! :=\! h[\lceil T_d\rceil\!\setminus\! S^\varepsilon ]$, so that 
$I^\varepsilon$ is a $\bormxi$ subset of $[\subseteq ]$. We also set, for $\varepsilon\!\in\! 2$,
$$D^\varepsilon_1\! :=\!\big\{ \vec s \!\in\!T_d\mid\vec s\! =\!\vec\emptyset\vee
\exists {\cal B}\!\in\!\Pi^{-1}(I^\varepsilon )\ \ \vec s \!\in\! {\cal B}\big\}\mbox{,}$$
and $D^\varepsilon_0\! :=\! T_d\!\setminus\! D^\varepsilon_1$. We set, for $\theta_0,\theta_1\!\in\! 2$, $D_{\theta_0,\theta_1}\! :=\! D^0_{\theta_0}\cap D^1_{\theta_1}$. For example, 
$\vec\emptyset\!\in\! D_{1,1}$.\bigskip

\noindent $\bullet$ Conditions (4), (5), and (7) become the following:\bigskip
 
\leftline{$\begin{array}{ll}
& \!\!\!\! (4)\ U_{\vec s }\!\subseteq\!\overline{A_0}^{\tau_\xi}\cap\overline{A_1}^{\tau_\xi}\mbox{ if }
\vec s \!\in\! D_{1,1}\mbox{,}\cr\cr  
& \!\!\!\! (5)\ U_{\vec s }\!\subseteq\! 
A_\varepsilon\mbox{ if }\vec s \!\in\! D_{\varepsilon ,1-\varepsilon}\mbox{,}\cr\cr
& \!\!\!\! (7)\ (\vec s,\vec t\!\in\! D_{\varepsilon ,1-\varepsilon}\wedge\vec s~R^{(\eta )}~\vec t\ )
\Rightarrow U_{\vec t}\!\subseteq\! U_{\vec s}.
\end{array}$}\bigskip
  
\noindent $\bullet$ Fix $\vec\alpha \!\in\!\lceil T_d\rceil$. There are 
$(\theta_0,\theta_1)\!\in\! 2^2$ and $k_0\!\in\!\omega$ such that, for $k\!\geq\! k_{0}$, 
$\vec\alpha\vert j_k\!\in\! D_{\theta_0,\theta_1}$. Thus 
$S^\varepsilon\!\subseteq\! (\Pi_{i\in d}\ f_i)^{-1}(A_\varepsilon )$.\bigskip

\noindent $\bullet$ Let $(\alpha^i_{\emptyset})_{i\in d}\!\in\! 
\overline{A_0}^{\tau_\xi}\cap\overline{A_1}^{\tau_\xi}\cap\Omega_{(\omega^\omega )^d}$, which is nonempty since $\overline{A_0}^{\tau_\xi}\cap\overline{A_1}^{\tau_\xi}\!\not=\!\emptyset$ is $\Ana$. We choose $U_{\vec\emptyset}$ with $(\alpha^i_{\emptyset})_{i\in d}\!\in\! U_{\vec\emptyset}
\!\subseteq\!\overline{A_0}^{\tau_\xi}\cap\overline{A_1}^{\tau_\xi}\cap\Omega_{(\omega^\omega )^d}$.

\vfill\eject
 
\noindent $\bullet$ The statement of the claim is now as follows:\bigskip

\noindent\bf Claim.\rm ~Assume that $\eta\! >\! 0$.\medskip

\noindent (a) $A_\varepsilon\cap\bigcap_{1\leq\rho\leq\eta}~\overline{U_{\overrightarrow{tm}^{\rho}}}^{\tau_{\rho}}
\cap (\Pi_{i\in d}\ X_i)$ is $\tau_{1}$-dense in 
$\overline{U_{\overrightarrow{tm}^1}}^{\tau_{1}}\cap (\Pi_{i\in d}\ X_i)$ if 
$\overrightarrow{tm}^{\eta}\!\notin\! D_{\varepsilon ,1-\varepsilon}$ and 
$\overrightarrow{tm}\!\in\! D_{\varepsilon ,1-\varepsilon}$.\medskip

\noindent (b) $U_{\overrightarrow{tm}^{\eta}}\cap\bigcap_{1\leq\rho <\eta}~\overline{U_{\overrightarrow{tm}^{\rho}}}^{\tau_{\rho}}\cap (\Pi_{i\in d}\ X_i)$ is $\tau_{1}$-dense in $\overline{U_{\overrightarrow{tm}^1}}^{\tau_{1}}\cap 
(\Pi_{i\in d}\ X_i)$ otherwise.\bigskip

 The point is that $\overrightarrow{tm}^{\eta}\!\in\! D_{1,1}$ if 
$\overrightarrow{tm}^{\eta}\!\notin\! D_{\varepsilon ,1-\varepsilon}$ since 
$\overrightarrow{tm}^{\eta}\!\in\! D_{\theta_0,\theta_1}$ with $\varepsilon\!\leq\!\theta_0$ and 
$1\! -\!\varepsilon\!\leq\!\theta_1$.\bigskip

\noindent $\bullet$ In the same fashion, $\Theta (\overrightarrow{tm})$ is now defined as follows:
$$\Theta (\overrightarrow{tm})\! :=\!\left\{\!\!\!\!\!\!
\begin{array}{ll} 
& A_\varepsilon\cap\bigcap_{1\leq\rho\leq\eta}~\overline{U_{\overrightarrow{tm}^{\rho}}}^{\tau_{\rho}}
\cap (\Pi_{i\in d}\ X_i)\cap\Omega_{(\omega^\omega )^d}\mbox{ if }
\overrightarrow{tm}^{\eta}\!\notin\! D_{\varepsilon ,1-\varepsilon}
\wedge\overrightarrow{tm}\!\in\! D_{\varepsilon ,1-\varepsilon}\mbox{,}\cr\cr 
& U_{\overrightarrow{tm}^{\eta}}\cap\bigcap_{1\leq\rho <\eta}~
\overline{U_{\overrightarrow{tm}^{\rho}}}^{\tau_{\rho}}\cap (\Pi_{i\in d}\ X_i)\mbox{ otherwise.}
\end{array}\right.$$
We conclude as in the proof of Theorem 4.4.1.\hfill{$\square$}\bigskip

 Now we come to the limit case. We need some more definitions that can be found in [D-SR].

\begin{defin} (Debs-Saint Raymond) Let $R$ be a tree relation on $c^{<\omega}$. If 
$t\!\in\! c^{<\omega}$, then $h_{R}(t)$ is the number of strict $R$-predecessors 
of $t$. So we have ${h_{R}(t)\! =\!\hbox{\it Card}\big( P_R(t)\big)\! -\! 1}$.\smallskip

 Let $\xi\! <\!\omega_{1}$ be an infinite limit ordinal. We say that a resolution family 
$(R^{(\rho)})_{\rho\leq\xi}$ is $uniform$ if
$$\forall k\!\in\!\omega\ \exists\eta_k\! <\!\xi\ 
\forall s,t\!\in\! c^{<\omega}\ \ 
\big(\hbox{\it min}\big( h_{R^{(\xi)}}(s),h_{R^{(\xi)}}(t)\big)\!\leq\! k\ 
\wedge s\ R^{(\eta_k)}\ t\big)\Rightarrow s\ R^{(\xi)}\ t.$$
We may (and will) assume that $\eta_k\!\geq\! 2$.\end{defin}

 The following is part of Theorem I-6.6 in [D-SR].
 
\begin{them} (Debs-Saint Raymond) Let $\xi\! <\!\omega_{1}$ be an infinite limit ordinal, $R$ a tree relation, $(I_n)_{n\in\omega}$ a sequence of $\bormxi$ subsets of $[R]$. Then there is a uniform resolution family $(R^{(\rho )})_{\rho\leq\xi}$ with\smallskip

\noindent (a) $R^{(0)}=R$.\smallskip

\noindent (b) The canonical map $\Pi\! :\! [R^{(\xi )}]\!\rightarrow\! [R]$ 
is a continuous bijection.\smallskip

\noindent (c) The set $\Pi^{-1}(I_n)$ is a closed subset of $[R^{(\xi )}]$ for each integer $n$.
\end{them}

 Here again, the next result is essentially in [L7] (see Theorem 2.4.4). 
 
\begin{them} Let $T_d$ be a tree with $\Borel$ suitable levels, $\xi\! <\!\omega^{\mbox{CK}}_{1}$ an infinite limit ordinal, $S$ in $\boraxi (\lceil T_d\rceil )$, and $A_0$, $A_1$ disjoint $\Ana$ subsets of 
$(\omega^\omega )^d$. We assume that Theorem 4.2.2 is proved for $\eta\! <\!\xi$. Then one of the following holds:\smallskip

\noindent (a) $\overline{A_0}^{\tau_\xi}\cap A_1\! =\!\emptyset$.\smallskip

\noindent (b) The inequality $\big( (\Pi_i''\lceil T_d\rceil )_{i\in d}, S,\lceil T_d\rceil\!\setminus\! S\big)\leq
\big( (\omega^\omega )_{i\in d}, A_0, A_1\big)$ holds.\end{them}

\noindent\bf Proof.\rm ~Let us indicate the differences with the 
proof of Theorem 4.4.1.\bigskip

\noindent $\bullet$ The set $I\!:=\! h[\lceil T_d\rceil\!\setminus\! S]$ is 
$\bormxi ([\subseteq ])$. Theorem 4.4.4 provides a uniform resolution family.\bigskip

\noindent $\bullet$ If $\vec t\!\in\! c^{<\omega}$ then we set 
$\eta (\vec t\ )\! :=\!\hbox{\rm max}\{\eta_{h_{R^{(\xi )}}(\vec s)+1}\mid\vec s\!\subseteq\!\vec t\ \}$. 
Note that $\eta (\vec s~)\!\leq\!\eta (\vec t~)$ if $\vec s\!\subseteq\!\vec t$.\bigskip

\noindent $\bullet$ Conditions (6) and (7) become\bigskip

\leftline{$
\begin{array}{ll}
& \!\!\!\! (6)\ ( 1\!\leq\!\rho\!\leq\!\eta\big( \vec s\ \big)\wedge
\vec s ~R^{(\rho )}~\vec t\ )\Rightarrow U_{\vec t }\!\subseteq\!
\overline{U_{\vec s }}^{\tau_{\rho }}\mbox{,}\cr\cr
& \!\!\!\! (7)\ \big( (\vec s,\vec t\!\in\! D\vee\vec s,\vec t\!\notin\! D)\wedge
\vec s ~R^{(\xi )}~\vec t\ \big)\Rightarrow U_{\vec t}\!\subseteq\! U_{\vec s}.
\end{array}$}\bigskip

\noindent\bf Claim 1.\rm ~Assume that $\overrightarrow{tm}^\rho\!\not=\! \overrightarrow{tm}^\xi$. Then 
$\rho\! +\! 1\!\leq\!\eta (\overrightarrow{tm}^{\rho +1})$.\bigskip 

 We argue by contradiction. We get $\rho\! +\! 1\! >\!\rho\!\geq\!\eta (\overrightarrow{tm}^{\rho +1})\!\geq\!
\eta_{h_{R^{(\xi )}}(\overrightarrow{tm}^{\xi})+1}\! =\!\eta_{h_{R^{(\xi )}}(\overrightarrow{tm})}$. As $\overrightarrow{tm}^\rho\ R^{(\rho)}\ \overrightarrow{tm}$ we get 
$\overrightarrow{tm}^\rho\ R^{(\xi)}\ \overrightarrow{tm}$, and also 
$\overrightarrow{tm}^\rho\! =\! \overrightarrow{tm}^\xi$, which is absurd.\hfill{$\diamond$}\bigskip

 Note that $\xi_{n-1}\! <\!\xi_{n-1}\! +\! 1\!\leq\!\eta (\overrightarrow{tm}^{\xi_{n-1}+1})\!\leq\!
\eta (\overrightarrow{tm})$. This implies that 
$\overrightarrow{tm}^{\eta (\overrightarrow{tm})}\! =\! \overrightarrow{tm}^{\xi}$.\bigskip

\noindent\bf Claim 2.\rm ~(a) The set $A_0\cap\bigcap_{1\leq\rho\leq\eta (\overrightarrow{tm})}~
\overline{U_{\overrightarrow{tm}^{\rho}}}^{\tau_{\rho}}\cap (\Pi_{i\in d}\ X_i)$ is $\tau_{1}$-dense in  
$\overline{U_{\overrightarrow{tm}^1}}^{\tau_{1}}\cap (\Pi_{i\in d}\ X_i)$ if 
$\overrightarrow{tm}^{\eta}\!\in\! D$ and $\overrightarrow{tm}\!\notin\! D$.\medskip

\noindent (b) The set $U_{\overrightarrow{tm}^{\xi}}\cap
\bigcap_{1\leq\rho <\eta (\overrightarrow{tm})}~
\overline{U_{\overrightarrow{tm}^{\rho}}}^{\tau_{\rho}}\cap (\Pi_{i\in d}\ X_i)$ is $\tau_{1}$-dense in 
$\overline{U_{\overrightarrow{tm}^1}}^{\tau_{1}}\cap (\Pi_{i\in d}\ X_i)$ if 
$\overrightarrow{tm}^{\xi},\overrightarrow{tm}\!\in\! D$ or 
$\overrightarrow{tm}^{\xi},\overrightarrow{tm}\!\notin\! D$.\bigskip 

 Indeed, we set $S_{i}\! :=\! U_{\overrightarrow{tm}^{\xi_{i}}}$, for $1\!\leq\!\xi_{i}\!\leq\!\xi$. 
By Claim 1 we can apply Lemma 4.2.3.(2).(c) and we are done.\hfill{$\diamond$}\bigskip

\noindent $\bullet$ The map 
$\Theta\! :\! {\cal X}^d\!\rightarrow\!\Ana\big( (\omega^\omega )^d\big)$ is defined on ${\cal T}^{l+1}$ by
$$\Theta (\overrightarrow{tm})\! :=\!\left\{\!\!\!\!\!\!
\begin{array}{ll} 
& A_0\cap\bigcap_{1\leq\rho\leq\eta (\overrightarrow{tm})}~
\overline{U_{\overrightarrow{tm}^{\rho}}}^{\tau_{\rho}}\cap (\Pi_{i\in d}\ X_i)\cap
\Omega_{(\omega^\omega )^d}\mbox{ if }\overrightarrow{tm}^{\eta}\!\in\! D
\wedge\overrightarrow{tm}\!\notin\! D\mbox{,}\cr\cr 
& U_{\overrightarrow{tm}^{\xi}}\!\cap\!\bigcap_{1\leq\rho <\eta (\overrightarrow{tm})}
\overline{U_{\overrightarrow{tm}^{\rho}}}^{\tau_{\rho}}\!
\cap\! (\Pi_{i\in d}\ X_i)\mbox{ if }
\overrightarrow{tm}^{\xi},\overrightarrow{tm}\!\in\! D\vee
\overrightarrow{tm}^{\xi},\overrightarrow{tm}\!\notin\! D.
\end{array}\right.$$
We conclude as in the proof of Theorem 4.4.1, using the facts that 
$\eta_k\!\geq\! 1$ and $\eta (.)$ is increasing.\hfill{$\square$}\bigskip

 Now we come to the ambiguous classes.

\begin{them} Let $T$ be a tree with $\Borel$ suitable levels, 
$\xi\! <\!\omega^{\hbox{\it CK}}_{1}$ an infinite limit ordinal, $S^0, S^1$ in 
$\boraxi (\lceil T_d\rceil )$ disjoint, and $A_0$, $A_1$ disjoint $\Ana$ subsets of $(\omega^\omega )^d$. We assume that Theorem 4.2.2 is proved for $\eta\! <\!\xi$. Then one of the following holds:\smallskip

\noindent (a) $\overline{A_0}^{\tau_\xi}\cap\overline{A_1}^{\tau_\xi}\! =\!\emptyset$.\smallskip

\noindent (b) The inequality $\big( (\Pi_i''\lceil T_d\rceil )_{i\in d}, S^0,S^1\big)\leq
\big( (\omega^\omega )_{i\in d}, A_0, A_1\big)$ holds.\end{them}

\noindent\bf Proof.\rm ~Let us indicate the differences with the proofs of Theorems 4.4.1, 4.4.2 and 4.4.5.\bigskip

\noindent $\bullet$ The set $I^\varepsilon\!:=\! h[\lceil T_d\rceil\!\setminus\! S^\varepsilon ]$ is 
$\bormxi ([\subseteq ])$.\bigskip

\noindent $\bullet$ The statement of Claim 2 is now as follows.\bigskip

\noindent\bf Claim 2.\rm ~(a) $A_\varepsilon\cap
\bigcap_{1\leq\rho\leq\eta (\overrightarrow{tm})}~
\overline{U_{\overrightarrow{tm}^{\rho}}}^{\tau_{\rho}}\cap (\Pi_{i\in d}\ X_i)$ is 
$\tau_{1}$-dense in $\overline{U_{\overrightarrow{tm}^1}}^{\tau_{1}}\cap (\Pi_{i\in d}\ X_i)$  
if $\overrightarrow{tm}^{\xi}\!\notin\! D_{\varepsilon ,1-\varepsilon}$ and 
$\overrightarrow{tm}\!\in\! D_{\varepsilon ,1-\varepsilon}$.\medskip

\noindent (b) $U_{\overrightarrow{tm}^{\xi}}\cap
\bigcap_{1\leq\rho <\eta (\overrightarrow{tm})}~
\overline{U_{\overrightarrow{tm}^{\rho}}}^{\tau_{\rho}}\cap (\Pi_{i\in d}\ X_i)$ is $\tau_{1}$-dense in $\overline{U_{\overrightarrow{tm}^1}}^{\tau_{1}}\cap 
(\Pi_{i\in d}\ X_i)$ otherwise.\bigskip

\noindent $\bullet$ In the same fashion, $\Theta (\overrightarrow{tm})$ is now defined as follows:
$$\Theta (\overrightarrow{tm})\! :=\!\left\{\!\!\!\!\!\!
\begin{array}{ll} 
& A_\varepsilon\cap\bigcap_{1\leq\rho\leq\eta (\overrightarrow{tm})}~
\overline{U_{\overrightarrow{tm}^{\rho}}}^{\tau_{\rho}}\cap (\Pi_{i\in d}\ X_i)\cap
\Omega_{(\omega^\omega )^d}\mbox{ if }
\overrightarrow{tm}^{\xi}\!\notin\! D_{\varepsilon ,1-\varepsilon}\wedge
\overrightarrow{tm}\!\in\! D_{\varepsilon ,1-\varepsilon}\mbox{,}\cr\cr 
& U_{\overrightarrow{tm}^{\xi}}\cap\bigcap_{1\leq\rho <\eta (\overrightarrow{tm})}~
\overline{U_{\overrightarrow{tm}^{\rho}}}^{\tau_{\rho}}\cap 
(\Pi_{i\in d}\ X_i)\mbox{ otherwise.}
\end{array}\right.$$
We conclude as in the proof of Theorem 4.4.5.\hfill{$\square$}

\begin{lemm} Let ${\bf\Gamma}$ be a Wadge class of Borel sets. Then the class of  $\mbox{pot}({\bf\Gamma})$ sets is closed under pre-images by products of continuous maps.\end{lemm}

\noindent\bf Proof.\rm\ Assume that $A\!\in\!\mbox{pot}({\bf\Gamma})$, 
$A\!\subseteq\!\Pi_{i\in d}\ Y_i$, and $f_i\! :\! X_i\!\rightarrow\! Y_i$ is continuous. Let $\tau_i$ be a finer $0$-dimensional Polish topology on $Y_i$ such that 
$A\!\in\! {\bf\Gamma}\big( \Pi_{i\in d}\ (Y_i,\tau_i)\big)$. As $f_i\! :\! X_i\!\rightarrow\! (Y_i,\tau_i)$ is Borel, there is a finer $0$-dimensional Polish topology $\sigma_i$ on $X_i$ such that 
$f_i\! :\! (X_i,\sigma_i)\!\rightarrow\! (Y_i,\tau_i)$ is continuous. Thus 
$(\Pi_{i\in d}\ f_i)^{-1}(A)\!\in\! {\bf\Gamma}\big( \Pi_{i\in d}\ (X_i,\sigma_i)\big)$ and 
$(\Pi_{i\in d}\ f_i)^{-1}(A)\!\in\!\mbox{pot}({\bf\Gamma})$.$\hfill{\square}$\bigskip

\noindent\bf Proof of Theorem 4.1 for $\xi$, assuming that Theorem 4.2.2 is proved for 
$\eta\! <\!\xi$.\rm\bigskip 

\noindent (1) We assume that (a) does not hold. This implies that the $X_i$'s are not empty.\bigskip

\noindent - We first prove that we may assume that $X_i\! =\!\omega^\omega$ for each $i\!\in\! d$.\bigskip

 By 13.5 in [K], there is a finer zero-dimensional Polish topology $\tau_i$ on $X_i$, and, by 7.8 in 
 [K], $(X_i,\tau_i)$ is homeomorphic to a closed subset $F_i$ of $\omega^\omega$, via a map 
 $\varphi_i$. By 2.8 in [K], there is a continuous retraction $r_i\! :\!\omega^\omega\!\rightarrow\! F_i$. Let $A'_\varepsilon$ be the intersection of $\Pi_{i\in d}\ F_i$ with the pre-image of $A_\varepsilon$ by 
$\Pi_{i\in d}\ (\varphi_i^{-1}\circ r_i)$. Then $A'_0$ and $A'_1$ are disjoint analytic subsets of 
$(\omega^\omega)^d$. Moreover, $A'_0$ is not separable from $A'_1$ by a 
$\mbox{pot}(\bormxi )$ set, since otherwise (a) would hold.\bigskip

 This gives $g_i\! :\! d^\omega\!\rightarrow\!\omega^\omega$ continuous with 
$S\!\subseteq\! (\Pi_{i\in d}\ g_i)^{-1}(A'_0)$ and 
$\lceil T_d\rceil \!\setminus\! S\!\subseteq\! (\Pi_{i\in d}\ g_i)^{-1}(A'_1)$. It remains to set 
$f_i(\alpha )\! :=\! (\varphi_i^{-1}\circ r_i\circ g_i)(\alpha )$ if $\alpha\!\in\! d^\omega$.\bigskip

\noindent - To simplify the notation, we may assume that $T_d$ has $\Borel$ levels, 
$\xi\! <\!\omega^{\hbox{\rm CK}}_{1}$ and $A_0$, $A_1$ are $\Ana\big( (\omega^\omega )^d\big)$. Notice that $\overline{A_0}^{\tau_\xi}\cap A_1$ is not empty, since otherwise $A_0$ would be separable from $A_1$ by a set in $\bormone (\tau_\xi )\!\subseteq\!\bormxi (\tau_1)\!\subseteq\!\mbox{pot}(\bormxi )$ set, which is absurd. So (b) holds, by Theorems 4.4.1 and 4.4.5 (as $\Pi_i''\lceil T_d\rceil$ is compact, we just have to compose with continuous retractions to get functions defined on $d^\omega$). So (a) or (b) holds.

\vfill\eject

\noindent $\bullet$ If $P\!\in\!\hbox{\rm pot}(\bormxi )$ separates $A_0$ from $A_1$ and (b) holds, then 
$S\!\subseteq\! (\Pi_{i\in d}\ f_i)^{-1}(P)\!\subseteq\!\neg (\lceil T_d\rceil \!\setminus\! S)$. This implies that $S$ is separable from $\lceil T_d\rceil \!\setminus\! S$ by a $\mbox{pot}(\bormxi )$ set, by Lemma 
4.4.7.\bigskip

\noindent (2) We argue as in the proof of (1). Here we consider 
$\overline{A_0}^{\tau_\xi}\cap\overline{A_1}^{\tau_\xi}$, and we apply Theorems 4.4.2 and 4.4.6.
This finishes the proof.\hfill{$\square$}\bigskip

\noindent\bf Proof of Theorem 4.2.2.\rm\ We assume that Theorem 4.1 is proved for $\xi$, and that Theorem 4.2.2 is proved for $\eta\! <\!\xi$.\bigskip

\noindent (1) By Lemma 4.2.3, $V_{0}$ and $V_{<\xi}$ are $\Ca$.\bigskip

\noindent (a) $\Rightarrow$ (b) and (a) $\Rightarrow$ (d) are clear since 
${\it\Delta}_{\omega^\omega}$ is Polish.\bigskip

\noindent (b) $\Rightarrow$ (c) We argue by contradiction. As $\gamma\!\in\!\Borel$ we get 
$C_\gamma\!\in\!\Borel$. If $(\beta ,\gamma)\!\in\! V_{<\xi}$, then $C_\gamma\!\in\!\mbox{pot}(\bormlxi )$, which is absurd. If $(\beta ,\gamma)\!\in\! V_0$, then 
$C_\gamma\!\in\!\mbox{pot}(\bormz )\!\subseteq\!\mbox{pot}(\bormxi )$, which is absurd. If 
$(\beta ,\gamma)\!\notin\! V_{<\xi}\cup V_0$, then we get $\gamma'\!\in\!\Borel$ (see the definition of $\Phi$ before Theorem 4.2.2). As $\big( (\beta )_n,(\gamma')_n\big)\!\in\! V_{<\xi}$, we get 
$C_{(\gamma')_n}\!\in\!\mbox{pot}(\bormlxi )$. Now the equality  
$\neg C_\gamma\! =\!\bigcup_{n\in\omega}\ C_{(\gamma')_n}$ implies that  
$C_\gamma\!\in\!\mbox{pot}(\bormxi )$, which is absurd.\bigskip

\noindent (d) $\Rightarrow$ (e) This comes from the proof of Theorem 4.1.(1).\bigskip

\noindent (e) $\Rightarrow$ (f) This comes from Theorems 4.4.1 and 4.4.5.\bigskip

\noindent (f) $\Rightarrow$ (a) This comes from Theorem 4.1.(1).\bigskip

\noindent (c) $\Rightarrow$ (e) We argue by contradiction, so that $\overline{A_0}^{\tau_{\xi}}$ separates $A_0$ from $A_1$.\bigskip

 If $\xi\! =\! 1$, then for each $\vec\delta\!\in\! A_1$ there is 
$(\tilde\beta ,\tilde\gamma )\!\in\! (\Borel\!\times\!\Borel )\cap V_0$ such that 
$\vec\delta\!\in\! C_{\tilde\gamma}\!\subseteq\!\neg A_0$. The first reflection theorem gives 
$\beta ,\gamma' \!\in\!\Borel$ such that $\big( (\beta )_n,(\gamma')_n\big)\!\in\! V_0$ for each integer $n$ and $A_1\!\subseteq\! U\! :=\!\bigcup_{n\in\omega}\ C_{(\gamma')_n}\!\subseteq\!\neg A_0$. We choose 
$\gamma\!\in\!\Borel\cap W$ with $\neg C_\gamma\! =\! U$, and $(\beta ,\gamma )$ contradicts (c).\bigskip

 If $\xi\!\geq\! 2$, then by induction assumption and the first reflection theorem there are 
$\beta ,\gamma'\!\in\!\Borel$ with $\big( (\beta )_n,(\gamma')_n\big)\!\in\! V_{<\xi}$ and 
$C_{(\gamma')_n}\!\subseteq\!\neg A_0$, for each integer $n$, and 
$A_1\!\subseteq\! U\! :=\!\bigcup_{n}\ C_{(\gamma')_n}$. But $U$ is ${\Borel\cap\hbox{\rm pot}(\boraxi )}$ and separates $A_1$ from $A_0$. So let $\gamma\!\in\!\Borel\cap W$ with $\neg C_\gamma\! =\! U$. We have $(\beta ,\gamma )\!\in\! V_\xi$ and $C_\gamma$ separates $A_0$ from $A_1$, which is absurd.\bigskip
 
\noindent (2) It is clear that $V_{\xi}$ is $\Ca$.\bigskip

\noindent (3) We argue as in the proof of (1), except for the implication (c) $\Rightarrow$ (e) (for the implication (e) $\Rightarrow$ (f) we use Theorems 4.4.2 and 4.4.6).\bigskip

\noindent (c) $\Rightarrow$ (e) We argue by contradiction. By 4D.2 in [M], there are 
$W\!\in\!\Ca (\omega )$ and a partial function ${\bf d}\! :\!\omega\!\rightarrow\!\omega^\omega$, 
$\Ca$-recursive on $W$, such that ${\bf d}''W$ is the set of $\Borel$ points of $\omega^\omega$. We define 
$$\Pi_{A_\varepsilon}\! :=\!\big\{ n\!\in\!\omega\mid (n)_0,(n)_1\!\in\! W\wedge
\big( {\bf d}\big( (n)_0\big) ,{\bf d}\big( (n)_1\big)\big)\!\in\! V_{<\xi}
\wedge C_{{\bf d}((n)_1)}\cap A_\varepsilon\! =\!\emptyset\big\}.$$ 
Then $\Pi_{A_\varepsilon}\!\in\!\Ca$ and 
$\forall \vec\beta \!\in\! (\omega^\omega)^d\ \ \exists n\!\in\!\Pi_{A_0}\cup\Pi_{A_1}\ \ 
\vec\beta \!\in\! C_{{\bf d}((n)_1)}$ 
since $\overline{A_0}^{\tau_\xi}\cap\overline{A_1}^{\tau_\xi}\! =\!\emptyset$ (we use the induction assumption). By the first reflection theorem there is $D\!\in\!\Borel (\omega )$ such that 
$D\!\subseteq\!\Pi_{A_0}\cup\Pi_{A_1}$ and 
$\forall \vec\beta \!\in\! (\omega^\omega)^d\ \ \exists n\!\in\! D\ \ \vec\beta \!\in\! C_{{\bf d}((n)_1)}$.

\vfill\eject

 As $\Ca$ has the reduction property, we can find 
$\Pi'_{A_\varepsilon}\!\in\!\Ca$ disjoint such that $\Pi'_{A_\varepsilon}\!\subseteq\!\Pi_{A_\varepsilon}$ and $\Pi'_{A_0}\cup\Pi'_{A_1}\! =\!\Pi_{A_0}\cup\Pi_{A_1}$. We set 
$\Delta\! :=\!\bigcup_{n\in D\cap\Pi'_{A_1}}\ 
C_{{\bf d}((n)_1)}\!\setminus\! (\bigcup_{q<n}\ 
C_{{\bf d}((q)_1)})$. Then 
$$\neg\Delta\! =\!\bigcup_{n\in D\cap\Pi'_{A_0}}\ C^{(\omega^\omega)^d}_{{\bf d}((n)_1)}\!\setminus\! (\bigcup_{q<n}\ C^{(\omega^\omega)^d}_{{\bf d}((q)_1)})\mbox{,}$$ 
which proves that $\Delta\!\in\!\Borel\cap\mbox{pot}(\borxi)$, and separates $A_0$ from $A_1$. Let 
$(\beta ,\gamma ),(\beta',\gamma')\!\in\! (\Borel\!\times\!\Borel )\cap V_\xi$ with $\Delta\! =\! C_\gamma$ and 
$\neg\Delta\! =\! C_{\gamma'}$. Then we get a contradiction with (c).
\hfill{$\square$}\bigskip
 
\noindent\bf Remarks.\rm\ The assertions 4.2.3.(2).(a) and 4.2.3.(2).(b) admit uniform versions in the following sense. By 3E.2, 3F.6 and 3H.1 in [M], there is 
$S\!:\!\omega^\omega\!\times\!\omega^\omega\!\rightarrow\!\omega^\omega$ recursive such that 
for each recursively presented Polish space $X$ there is a universal set 
${\cal U}^X\!\in\!\Ca\big( (\omega^\omega )^d\big)$ satisfying the following properties:\smallskip

\noindent - $\ca (X)\! =\!\{ {\cal U}^X_\alpha\mid\alpha\!\in\!\omega^\omega\}$,\smallskip

\noindent - $\Ca (X)\! =\!\{ {\cal U}^X_\alpha\mid\alpha\!\in\!\omega^\omega
\mbox{ recursive}\}$,\smallskip

\noindent - $(\alpha ,\beta ,x)\!\in\! {\cal U}^{\omega^\omega\times X}\Leftrightarrow
\big( S(\alpha ,\beta ),x\big)\!\in\! {\cal U}^X$.\bigskip

 We set ${\cal U}\! :=\! {\cal U}^{(\omega^\omega )^d}$. The following relations are $\Ca$:
$$\begin{array}{ll}
Q(\alpha ,\beta ,\gamma )\!\!\!\! 
& \Leftrightarrow\alpha\!\in\!\mbox{WO}\wedge (\beta ,\gamma )\!\in\! V_{\vert\alpha\vert}\mbox{,}\cr
R(\alpha ,\beta ,\vec\delta )\!\!\!\! 
& \Leftrightarrow\alpha\!\in\!\Borel\cap\mbox{WO}\wedge\vert\alpha\vert\!\geq\! 1\wedge
\vec\delta\!\notin\!\overline{\neg {\cal U}_\beta}^{\tau_{\vert\alpha\vert}}.
\end{array}$$
Indeed, this comes from the proof of Lemma 4.2.3.\bigskip

\noindent $\bullet$ One can give simpler examples $\mathbb{S}^0,\mathbb{S}^1$ for which Corollary 4.2 is fullfilled when ${\bf\Gamma}\! =\!\bormone$. Indeed, recall the map $b_\omega$ defined before Lemma 2.3. As $\vert b_\omega (n)\vert\!\leq\! n$ for each integer $n$, we can define the sequence 
$s^\omega_n\! :=\! b_\omega (n)0^{n-\vert b_\omega (n)\vert}$. We set 
$\mathbb{S}^1\! :=\!\overline{\mathbb{S}^0}\!\setminus\!\mathbb{S}^0$, where
$$\mathbb{S}^0\! :=\!\Big\{\big( 0s^\omega_n0\gamma ,...,0s^\omega_nn\gamma ,
(n\! +\! 1)s^\omega_n(n\! +\! 1)\gamma ,(n\! +\! 1)s^\omega_n(n\! +\! 2)\gamma ,...\big)\mid 
(n,\gamma )\!\in\!\omega\!\times\!\omega^\omega\Big\}$$
(we do not really need $T_\omega$ when ${\bf\Gamma}\! =\!\bormone$). We have 
$\mathbb{S}^0\! =\! (\Pi_{i\in d}\ f_i)^{-1}(A_0)\cap\overline{\mathbb{S}^0}$ if (b) holds. Let us denote this by $\mathbb{S}^0\leq A_0$ (we have a quasi-order, by continuity of the $f_i$'s).\bigskip

\noindent $\bullet$ The fact that $T_d$ has finite levels was used to give a proof of Corollary 4.2 as simple as possible.  The tree $T_d$ has finite levels when $d\! <\!\omega$, and not always when $d\! =\!\omega$. This is one of the main new points in the case of the infinite dimension. Let us specify this.\bigskip

\noindent (a) We saw in the proof of Proposition 2.2 that the tree $\tilde T_d$ generated by an effective frame is a tree with one-sided almost acyclic levels. As before Lemma 2.6, we can define 
$$\tilde S^\omega_{C_1}\! :=\!\{\vec\alpha\!\in\!\lceil\tilde T_d\rceil\mid 
{\cal S}(\alpha_0\Delta\alpha_1)\!\in\! C_1\}\mbox{,}$$ 
which is not separable from $\lceil\tilde T_d\rceil\!\setminus\!\tilde S^\omega_{C_1}$ by a potentially closed set, since otherwise $S^\omega_{C_1}$ would be separable from 
$\lceil T_d\rceil\!\setminus\! S^\omega_{C_1}$ by a potentially closed set, which would contradict Lemmas 2.6 and 3.4.\bigskip

 But $\mathbb{A}_0\! :=\!\{ 0^{1+n}(1\! +\! n)^\infty\mid n\!\in\!\omega\}\!\subseteq\!\omega^\omega$ is not 
potentially closed since $0^\infty\!\in\!\overline{\mathbb{A}_0}\!\setminus\!\mathbb{A}_0$ and the topology on $\omega$ is discrete. And one can prove, in a straightforward way, that 
$\tilde S^\omega_{C_1}\not\leq\mathbb{A}_0$ and 
$\mathbb{A}_0\not\leq\tilde S^\omega_{C_1}$. This proves that the finiteness of the levels of $T_d$ is useful. But we will see that it is not necessary.

\vfill\eject

\noindent (b) We define 
$o\! :\!\{ s\!\in\! 2^{<\omega}\mid 0\!\not\subseteq\! s\}\!\rightarrow\!\omega^{<\omega}$ such that 
$\vert o(s)\vert\! =\!\vert s\vert$ by
$$o(10^{n_0}10^{n_1}...10^{n_l})\! :=\! 0^{1+n_0}(1\! +\! n_0)^{1+n_1}...
\big( (1\! +\! n_0)\! +\! ...\! +\! (1\! +\! n_{l-1})\big)^{1+n_l}.$$
In other words, we have $o(s)(i)\! =\! i$ if $s(i)\! =\! 1$, $o(s)(i)\! =\! o(s)(i\! -\! 1)$ if $s(i)\! =\! 0$. Note that $o$ is an injective homomorphism, in the sense that $o(s)\!\subseteq\! o(t)$ if $s\!\subseteq\! t$. This implies that we can extend $o$ to a continuous map from the basic clopen set $N_1$ into $\omega^\omega$ by the formula $o(\alpha)\! :=\!\mbox{sup}_{m\in\omega}\ o(\alpha\vert m)$.\bigskip

 We set $F_\omega\! :=\!\big\{ (m_i\alpha_i)_{i\in\omega}\!\in\! (\omega^\omega )^{\omega}\mid 
\vec\alpha \!\in\! \lceil\tilde T_\omega\rceil\ \mbox{ and }\ \forall i\!\in\!\omega\ \ 
m_i\! =\! o(\alpha_0\Delta\alpha_1)(i)\big\}$, and we put $\underline{S}^\omega_{C_\xi}\! :=\!
\{ (m_i\alpha_i)_{i\in\omega}\!\in\! F_\omega\mid {\cal S}(\alpha_0\Delta\alpha_1)\!\in\! C_\xi\}$. One can take $\mathbb{S}^\omega_\xi\! =\!\underline{S}^\omega_{C_\xi}$, and the proof is much more complicated than the one we gave. But the tree associated with $\overline{\underline{S}^\omega_{C_\xi}}\! =\! F_\omega$ is
$$\big\{\vec\emptyset\big\}\cup
\big\{ (m_is_i)_{i\in\omega}\!\in\! (\omega^\omega )^{<\omega}\mid 
(m_i)_{i\in\omega}\!\in\! o''[N_1]\ \mbox{ and }\ 
\vec s\!\in\!\tilde T_\omega\ \mbox{ and }\ \forall i\! <\!\vert\vec s\ \vert\ \ 
m_i\! =\! o(s_0\Delta s_1)(i)\big\}\mbox{,}$$
and has infinite levels. This proves that the finiteness of the levels of the tree associated with 
$\overline{\mathbb{S}^\omega_\xi}$ is not necessary.\bigskip

\noindent (c) In [L8], an extension to any dimension of the Kechris-Solecki-Todor\v cevi\'c dichotomy about ana-lytic graphs is proved. In [L5], it is proved that Corollary 4.2 is a consequence of the Kechris-Solecki-Todor\v cevi\'c dichotomy when ${\bf\Gamma}\! =\!\bormone$. This works as well when 
$d\! <\!\omega$, but not when $d\! =\!\omega$. More specifically, let 
$\mathbb{G}\! :=\!\{\alpha\!\in\!\omega^\omega\mid
\forall m\!\in\!\omega\ \exists n\!\geq\! m\ \ s^\omega_n0\!\subseteq\!\alpha\}$ and 
$$\mathbb{A}_\omega\! :=\!\{ (s^\omega_ni\gamma)_{i\in\omega}\mid n\!\in\!\omega\wedge
\gamma\!\in\!\omega^\omega\}.$$ 
 Then the extension to the case where $d\! =\!\omega$ of the Kechris-Solecki-Todor\v cevi\'c dichotomy works with $\mathbb{G}^\omega\cap\mathbb{A}_\omega$ (see [L8]). But one can prove the following result:

\begin{them} Let $X$ be a recursively presented Polish space, $\sigma_X$ the topology on $X^\omega$ generated by $\{\Pi_{i\in\omega}\ C_i\mid C\!\in\!\Borel (\omega\!\times\! X)\}$, and $A$ a 
$\Borel$ subset of $X^\omega$. Then exactly one of the following holds:\smallskip

\noindent (a) $\overline{A}^{\sigma_X}\!\setminus\! A\! =\!\emptyset$.\smallskip

\noindent (b) $\mathbb{G}^\omega\cap\mathbb{A}_\omega\leq A$.\end{them}

 In particular, $\mathbb{G}^\omega\cap\mathbb{A}_\omega\not\leq\mathbb{A}_0$ and we cannot take 
$\mathbb{S}^\omega_1\! =\!\mathbb{G}^\omega\cap\mathbb{A}_\omega$.

\section{$\!\!\!\!\!\!$ The proof of Theorem 1.7}

\subsection{$\!\!\!\!\!\!$ Some material in dimension one}\indent

 The material in this subsection is due to A. Louveau and J. Saint Raymond, and can be found in [Lo-SR1] or [Lo-SR2]. However, some changes are needed for our purposes, and moreover some proofs are left to the reader in these papers. So we will sometimes give some proofs. The following definition can be found in [Lo-SR2] (see Definition 1.5).
 
\begin{defin} Let $1\!\leq\!\xi\! <\!\omega_1$, ${\bf\Gamma}$ and ${\bf\Gamma}'$ two classes. Then
$$A\!\in\! S_\xi({\bf\Gamma},{\bf\Gamma}')\ \Leftrightarrow\ 
A\! =\!\bigcup_{p\geq 1}\ (A_p\cap C_p)\cup\left( B\!\setminus\!\bigcup_{p\geq 1}\ C_p\right)$$ 
for some sequence of sets $A_p$ in ${\bf\Gamma}$, $B\!\in\! {\bf\Gamma}'$, and a sequence 
$(C_p)_{p\geq 1}$ of pairwise disjoint $\boraxi$ sets.\end{defin}

 Now we come to the definition of $second\ type\ descriptions$ of non self-dual Wadge classes of Borel sets, which are elements of $\omega_1^\omega$, sometimes identified with 
$(\omega_1^\omega )^\omega$. This definition can also be found in [Lo-SR2] (see Definition 1.6).
 
\begin{defin} The relations ``$u$ is a $second\ type\ description$" and ``$u\ describes\ {\bf\Gamma}$" 
(written $u\!\in\! {\cal D}$ and ${\bf\Gamma}_u\! =\! {\bf\Gamma}$~-~ambiguously) are the least relations satisfying\smallskip

\noindent (a) If $u\! =\! 0^\infty$, then $u\!\in\! {\cal D}$ and ${\bf\Gamma}_u\! =\!\{\emptyset\}$.\smallskip

\noindent (b) If $u\! =\!\xi^\frown 1^\frown u^*$, with $u^*\!\in\! {\cal D}$ and $u^*(0)\! =\!\xi$, then 
$u\!\in\! {\cal D}$ and ${\bf\Gamma}_u\! =\!\check {\bf\Gamma}_{u^*}$.\smallskip

\noindent (c) If $u\! =\!\xi^\frown 2^\frown\!\!\! <\! u_p\! >$ satisfies $\xi\!\geq\! 1$, $u_p\!\in\! {\cal D}$, and 
$u_p(0)\!\geq\!\xi$ or $u_p(0)\! =\! 0$, then $u\!\in\! {\cal D}$ and 
${\bf\Gamma}_u\! =\! S_\xi(\bigcup_{p\geq 1}\ {\bf\Gamma}_{u_p},{\bf\Gamma}_{u_0})$.\end{defin}

\noindent\bf Remark.\rm\ If $A\!\in\! S_\xi(\bigcup_{p\geq 1}\ {\bf\Gamma}_{u_p},{\bf\Gamma}_{u_0})$, then $A$ has a decomposition as in Definition 5.1.1, and $A_p$ is in $\bigcup_{p\geq 1}\ {\bf\Gamma}_{u_p}$. But we may assume that $A_p\!\in\! {\bf\Gamma}_{u_{(p)_0+1}}$, using the fact that $C_p$ may be empty if necessary. This remark will be useful in the sequel, since it specifies the class of $A_p$.\bigskip

 The following result can be found in [Lo-SR2] (see Section 3).
 
\begin{them} Let ${\bf\Gamma}$ be a non self-dual Wadge class of Borel sets. Then there is 
$u\!\in\! {\cal D}$ such that ${\bf\Gamma}(\omega^\omega )\! =\! {\bf\Gamma}_u(\omega^\omega )$. Conversely, 
$${\bf\Gamma}_u\! :=\!\{ f^{-1}(A)\mid f\! :\! X\!\rightarrow\!\omega^\omega\mbox{ continuous }\wedge 
X\mbox{ 0-dimensional Polish space }\wedge A\!\in\! {\bf\Gamma}_u(\omega^\omega )\}$$ 
is a non self-dual Wadge class of Borel sets if $u\!\in\! {\cal D}$.\end{them}

 If $\eta\!\leq\!\xi\! <\!\omega_1$, then $\xi\! -\!\eta$ is the unique ordinal $\theta$ with 
$\eta\! +\!\theta\! =\!\xi$. The following definition can be found in [Lo-SR2] (see Definition 1.9).

\begin{defin} Let $\eta\! <\!\omega_1$ and $u\!\in\! {\cal D}$. We define 
$u^\eta\!\in\! {\cal D}$ as follows:\smallskip

\noindent (a) If $u(0)\! =\! 0$, then $u^\eta\! :=\! u$.\smallskip

\noindent (b) If $u\! =\!\xi 1u^*$, with $\xi\!\geq\! 1$, then 
$u^\eta\! :=\!\big( 1\! +\!\eta\! +\! (\xi\! -\! 1)\big) 1(u^*)^\eta$.\smallskip

\noindent (c) If $u\! =\!\xi 2 <\! u_p\! >$, with $\xi\!\geq\! 1$, then 
$u^\eta\! :=\!\big( 1\! +\!\eta\! +\! (\xi\! -\! 1)\big) 2 <(u_p)^\eta >$.\end{defin}

 The following result can be found in [Lo-SR2] (see Proposition 1.10).
 
\begin{propo} (a) If $f\! :\!\omega^\omega\!\rightarrow\!\omega^\omega$ is ${\bf\Sigma}^0_{1+\eta}$-measurable, and $A\!\in\! {\bf\Gamma}_u(\omega^\omega )$ for some $u\!\in\! {\cal D}$, then 
$f^{-1}(A)\!\in\! {\bf\Gamma}_{u^\eta}$.\smallskip

\noindent (b) The set ${\cal D}$ is the least subset $D\!\subseteq\! {\cal D}$ such that $0^\infty\!\in\! D$, 
$u(0)1u\!\in\! D$ if $u\in\! D$, $12 <u_p\! >\in\! D$ if $u_p\in\! D$ for each $p\!\in\!\omega$, and 
$u^\eta\!\in\! D$ if $u\in\! D$, for each $\eta\! <\!\omega_1$.\end{propo}

 Recall the definition of an independent $\eta$-function (see Definition 3.3).
 
\vfill\eject

\noindent\bf Example.\rm ~Let $\tau\! :\!\omega\!\rightarrow\!\omega$ be one-to-one (in [Lo-SR2] just before Lemma 2.5, increasing maps are considered, but here we relax this condition). We define 
$\tilde\tau\! :\! 2^\omega\!\rightarrow\! 2^\omega$ by $\tilde\tau (\alpha )\! :=\!\alpha\circ\tau$. Clearly 
$\tilde\tau$ is an independent $0$-function, with $\pi (k)\! =\!\tau^{-1}(k)$ if $k$ is in the range of $\tau$, $0$ otherwise. We now describe an important instance of this situation.\bigskip

\noindent\bf Example.\rm ~Let $n$ be an integer, and $\cal S$ the shift map (see the notation before Definition 2.5). Then ${\cal S}^n$ is an independent $0$-function. Indeed, if we set 
$\tau^n(m)\! :=\! m\! +\! n$, then ${\cal S}^n\! =\!\tilde {\tau^n}$, by induction on $n$. In particular, 
$\mbox{Id}_{2^\omega}\! =\! {\cal S}^0$ is an independent $0$-function.\bigskip

 The next result is essentially Lemma 2.5 in [Lo-SR2], which is given without proof, so we give the details here.
 
\begin{lemm} Let $\tau\! :\!\omega\!\rightarrow\!\omega$ be one-to-one, $\rho$ an independent $\eta$-function. Then $\tilde\tau\circ\rho$ is an independent $\eta$-function.\end{lemm}

\noindent\bf Proof.\rm ~Let $\pi$ associated with $\rho$. We define $\pi'\! :\!\omega\!\rightarrow\!\omega$ by $\pi'(k)\! :=\!\tau^{-1}\big(\pi (k)\big)$ if $\pi (k)$ is in the range of $\tau$, $0$ otherwise, so that 
$\pi'(k)\! =\! m$ if $\pi (k)\! =\!\tau (m)$. If $m$ is an integer, then 
$(\tilde\tau\circ\rho )(\alpha )(m)\! =\!\rho (\alpha )\big(\tau (m)\big)$ depends only of the values of $\alpha$ on $\pi^{-1}\big(\{\tau (m)\}\big)\!\subseteq\! (\pi')^{-1}(\{ m\})$.\bigskip

 If $\xi\! =\! 0$ (resp., $\xi\! =\!\theta\! +\! 1$, $\xi\! =\!\mbox{sup}_{m\in\omega}\ \theta_m$), then 
${\cal C}_m\! =\!\{\alpha\!\in\! 2^\omega\mid\rho (\alpha )\big(\tau (m)\big)\! =\! 1\}$ is $\borone$-complete (resp., $\bormpz$-strategically complete, $\bormpztm$-strategically complete). We are done since 
$\xi\! =\!\mbox{sup}_{p\geq 1}\ \theta_{\tau (m_p)}$ if $\xi$ is a limit ordinal ($\tau$ is one-to-one).
\hfill{$\square$}\bigskip

 After Definition 3.3, we saw that $\rho_0^\eta$ is an independent $\eta$-function. We will actually prove more, actually a result which is essentially Theorem 2.4.(b) in [Lo-SR2].

\begin{them} Let $\eta ,\xi\! <\!\omega_1$, $\rho$ an independent $\xi$-function. Then $\rho_0^\eta\circ\rho$ is an independent $(\xi\! +\!\eta )$-function.\end{them}

\noindent\bf Proof.\rm ~Note first that if $\varepsilon\!\in\! 2$, 
$\rho^\varepsilon\! :\! 2^\omega\!\rightarrow\! 2^\omega$ is equipped with $\pi^\varepsilon$ such that 
$\rho^\varepsilon (\alpha )(m)$ depends only on the values of $\alpha$ on $(\pi^\varepsilon )^{-1}(\{ m\})$, then $(\rho^0\circ\rho^1)(\alpha )(m)$ depends only on the values of $\rho^1(\alpha )$ on 
$(\pi^0)^{-1}(\{ m\})$, so it depends only on the values of $\alpha$ on 
$(\pi^1)^{-1}\big( (\pi^0)^{-1}(\{ m\})\big)$, so that if we set $\pi\! :=\!\pi^0\circ\pi^1$, then 
$(\rho^0\circ\rho^1)(\alpha )(m)$ depends only on the values of $\alpha$ on $\pi^{-1}(\{ m\})$.\bigskip

\noindent $\bullet$ We argue by induction on $\eta$. The result is clear for $\eta\! =\! 0$. So assume that 
$\eta\! =\!\theta\! +\! 1$, so that $\rho_0^\eta\circ\rho\! =\!\rho_0\circ\rho_0^\theta\circ\rho$. The induction assumption implies that $\rho^\theta\circ\rho$ is an independent $(\xi\! +\!\theta )$-function. The fact that 
$\rho_0$ is an independent $1$-function and the previous point prove the existence of $\pi_\eta$ such that $(\rho_0^\eta\circ\rho )(\alpha )(m)$ depends only on the values of $\alpha$ on 
$\pi_\eta^{-1}(\{ m\})$.\bigskip

 We set $A_n\! :=\!\{\alpha\!\in\! 2^\omega\mid (\rho_0^\theta\circ\rho)(\alpha )(<m,n>)\! =\! 1\}$. Let us prove that $\bigcap_{n\in\omega}\ \neg A_n$ is $\bormpxz$-strategically complete.\bigskip
 
 Assume first that $\xi\! +\!\theta\!\not=\! 0$. As $\rho^\theta\circ\rho$ is an independent $(\xi\! +\!\theta )$-function, $A_n$ is ${\bf\Pi}^0_{1+\theta_n}$-strategically complete, for some $\theta_n\! <\!\xi\! +\!\theta$ satisfying $\theta_n\! +\! 1\! =\!\xi\! +\!\theta$ if $\xi\! +\!\theta$ is a successor ordinal, 
$\mbox{sup}_{n\in\omega}\ \theta_n\! =\!\xi\! +\!\theta$ if $\xi\! +\!\theta$ is a limit ordinal. Note that 
$\xi\! +\!\theta\! =\!\mbox{sup}_{n\in\omega}\ (\theta_n\! +\! 1)$. As $\rho^\theta\circ\rho$ is an independent $(\xi\! +\!\theta )$-function, there is $\pi_\theta$ such that 
$(\rho_0^\theta\circ\rho )(\alpha )(q)$ depends only on the values of $\alpha$ on $\pi_\theta^{-1}(\{ q\})$. We set $\pi (\alpha )(k)\! :=\!\big(\pi_\theta (\alpha )\big)_1$, so that the fact that $\alpha\!\in\! A_n$ depends only on the values of $\alpha$ on $\pi^{-1}(\{ n\})$. By Lemma 3.7 in [Lo-SR1], 
$\bigcap_{n\in\omega}\ \neg A_n$ is $\bormpxz$-strategically complete.

\vfill\eject

 Assume now that $\xi\! +\!\theta\! =\! 0$. Then 
$A_n\! :=\!\{\alpha\!\in\! 2^\omega\mid \rho (\alpha )(<m,n>)\! =\! 1\}$ is $\borone$-complete since $\rho$ is an independent $0$-function. Let $B$ be a closed subset of $\omega^\omega$, $(B_n)_{n\in\omega}$ a sequence of clopen subsets with $B\! =\!\bigcap_{n\in\omega}\ B_n$, and 
$g_n\! :\!\omega^\omega\!\rightarrow\! 2^\omega$ continuous with $B_n\! =\! g_n^{-1}(\neg A_n)$. As 
$\rho$ is an independent $0$-function, there is $\pi_\rho$ such that $\rho (\alpha )(q)$ depends only on the values of $\alpha$ on $\pi_\rho^{-1}(\{ q\})$. We set $\pi (\alpha )(k)\! :=\!\big(\pi_\rho (\alpha )\big)_1$, so that the fact that $\alpha\!\in\! A_n$ depends only on the values of $\alpha$ on $\pi^{-1}(\{ n\})$. We define $g\! :\!\omega^\omega\!\rightarrow\! 2^\omega$ by $g(\beta )(k)\! :=\! g_{\pi (k)}(\beta )(k)$, so that $g$ is continuous. Moreover, $\beta\!\in\! B_n\Leftrightarrow g_n(\beta )\!\notin\! A_n\Leftrightarrow g(\beta )\!\notin\! A_n$ since the fact that $\alpha\!\in\! A_n$ depends only on the values of $\alpha$ on $\pi^{-1}(\{ n\})$. Thus $B\! =\! g^{-1}(\bigcap_{n\in\omega}\ \neg A_n)$ and $\bigcap_{n\in\omega}\ \neg A_n$ is 
$\bormone$-complete. Therefore $\bigcap_{n\in\omega}\ \neg A_n$ is $\bormpxz$-strategically complete.\bigskip 

 Now note that 
$$\begin{array}{ll}
\bigcap_{n\in\omega}\ \neg A_n\!\!\!\!
& \! =\!\{\alpha\!\in\! 2^\omega\mid\forall n\!\in\!\omega\ \  (\rho_0^\theta\circ\rho)(\alpha )(<m,n>)\! =\! 0\}\cr\cr 
& \! =\!\{\alpha\!\in\! 2^\omega\mid  (\rho_0^{}\circ\rho_0^\theta\circ\rho)(\alpha )(m)\! =\! 1\}
\! =\!\{\alpha\!\in\! 2^\omega\mid  (\rho_0^\eta\circ\rho)(\alpha )(m)\! =\! 1\}.
\end{array}$$
Thus $\{\alpha\!\in\! 2^\omega\mid  (\rho_0^\eta\circ\rho)(\alpha )(m)\! =\! 1\}$ is $\bormpxz$-strategically complete for each $m$, and $\xi\! +\!\eta\! =\!\xi\! +\!\theta\! +\! 1$, so that $\rho_0^\eta\circ\rho$ is an independent $(\xi\! +\!\eta )$-function.\bigskip

\noindent $\bullet$ Assume now that $\eta$ is a limit ordinal. In the definition of $\rho_0^\eta$ we fixed a sequence $(\theta^\eta_{m})_{m\in\omega}\!\subseteq\!\eta$ of successor ordinals with 
$\Sigma_{m\in\omega}~\theta^\eta_m\! =\!\eta$. As $\rho_0^{\theta^\eta_{m}}$ is an independent $\theta^\eta_{m}$-function, we get $\pi^\eta_{m}\! :\!\omega\!\rightarrow\!\omega$. We define $\pi_{m,m+1}\! :\!\omega\!\rightarrow\!\omega$ by $\pi_{m,m+1}(k)\! :=\! k$ if $k\! <\! m$, 
$\pi^\eta_{m}(k\! -\! m)\! +\! m$ if ${k\!\geq\! m}$. Let us check that $\rho^{(m,m+1)}_{0}(\alpha )(i)$ depends only on the values of $\alpha$ on $\pi_{m,m+1}^{-1}(\{ i\})$. It is clearly the case if $i\! <\! m$. So assume that $i\!\geq\! m$. Note that $\pi_{m,m+1}(k)\! =\! i$ if 
$k\!\in\! (\pi^\eta_{m})^{-1}(\{ i\! -\! m\})\! +\! m$, and we are done. Now the first point of this proof gives ${\pi_{0,m+1}\! :\!\omega\!\rightarrow\!\omega}$ such that $\rho^{(0,m+1)}_{0}(\alpha )(i)$ depends only on the values of $\alpha$ on $\pi_{0,m+1}^{-1}(\{ i\})$. We will check that 
$\rho_0^\eta (\alpha )(m)\! :=\!\rho^{(0,m+1)}_{0}(\alpha )(m)$ depends only on the values of $\alpha$ on $E_m\! :=\!\pi_{0,m+1}^{-1}(\{ m\})\cap\bigcap_{l<m}\ \pi_{0,l+1}^{-1}(\neg (l\! +\! 1))$. We actually prove something stronger: for each integer $k$, $\rho^{(0,m+1)}_{0}(\alpha )(k\! +\! m)$ depends only on the values of $\alpha$ on 
$$\pi_{0,m+1}^{-1}(\{ k\! +\! m\})\cap\bigcap_{l<m}\ \pi_{0,l+1}^{-1}(\neg (l\! +\! 1)).$$ 
We argue by induction on $m$. For $m\! =\! 0$, the result is clear. Assume that the result is true for $m$. Note that $\rho^{(0,m+2)}_{0}(\alpha )(k\! +\! m\! +\! 1)$ depends only on the values of $\alpha$ on 
$\pi_{0,m+2}^{-1}(\{ k\! +\! m\! +\! 1\})$. But
$$\rho^{(0,m+2)}_{0}(\alpha )(k\! +\! m\! +\! 1)
\! =\!\rho^{(m+1,m+2)}_{0}\big(\rho^{(0,m+1)}_{0}(\alpha )\big)(k\! +\! m\! +\! 1)
\! =\!\rho_{0}^{\theta^\eta_{m+1}}\Big( {\cal S}^{m+1}\big(\rho^{(0,m+1)}_{0}(\alpha )\big)\Big)(k)\mbox{,}$$
and we are done since $\rho^{(0,m+2)}_{0}(\alpha )(\alpha )(k\! +\! m\! +\! 1)$ depends only on the values of 
${\cal S}^{m+1}\big(\rho^{(0,m+1)}_{0}(\alpha )\big)$, which depends only on the values of $\alpha$ on 
$\pi_{0,m+1}^{-1}(\neg (m\! +\! 1))\cap\bigcap_{l<m}\ \pi_{0,l+1}^{-1}(\neg (l\! +\! 1))$.\bigskip

 As the $E_m$'s are pairwise disjoint, we can define a map $\pi^\eta\! :\!\omega\!\rightarrow\!\omega$ by 
$\pi^\eta (k)\! :=\! m$ if $k\!\in\! E_m$, and $0$ if $k\!\notin\!\bigcup_{m\in\omega}\ E_m$. Now it is clear that 
$\rho_0^\eta (\alpha )(m)$ depends only on the values of $\alpha$ on $(\pi^\eta )^{-1}(\{ m\})$. The first point of this proof gives $\pi_\eta\! :\!\omega\!\rightarrow\!\omega$ such that 
$(\rho_0^\eta\circ\rho )(\alpha )(m)$ depends only on the values of $\alpha$ on $\pi_\eta^{-1}(\{ m\})$.\bigskip

 Let $\zeta^\eta_{m}$ such that $\theta^\eta_{m}\! :=\!\zeta^\eta_{m}\! +\! 1$, and 
$\theta_m\! :=\!\xi\! +\!\Sigma_{l<m}\ \theta^\eta_{l}\! +\!\zeta^\eta_{m}$, so that $\theta_m\! <\!\xi\! +\!\eta$ and $\mbox{sup}_{p\geq 1}\ \theta_{m_p}\! =\!\xi\! +\!\eta$ for each one-to-one sequence 
$(m_p)_{p\geq 1}$ of integers. It remains to see that 
$${\cal C}_m\! :=\!\{\alpha\!\in\! 2^\omega\mid (\rho_0^\eta\circ\rho )(\alpha )(m)\! =\! 1\}$$ 
is $\bormpzm$-strategically complete for each integer $m$.\bigskip

 Let us check that ${\cal S}^m\circ\rho^{(0,m+1)}_{0}
\! =\!\rho_0^{\theta^\eta_{m}}\circ\ \circ_{l<m}\ ({\cal S}\circ\rho_0^{\theta^\eta_{m-l-1}})$ 
for each integer $m$. We argue by induction on $m$. For $m\! =\! 0$ the property is clear since 
$\rho^{(0,1)}_{0}\! =\!\rho_0^{\theta^\eta_0}$. Assume that the property is true for $m$. Then
$$\begin{array}{ll}
{\cal S}^{m+1}\circ\rho_0^{(0,m+2)}\!\!\!\!
& \! =\!\rho_0^{\theta^\eta_{m+1}}\circ {\cal S}^{m+1}\circ\rho^{(0,m+1)}_{0}
\! =\!\rho_0^{\theta^\eta_{m+1}}\circ {\cal S}\circ {\cal S}^{m}\circ\rho^{(0,m+1)}_{0}\cr\cr 
& \! =\!\rho_0^{\theta^\eta_{m+1}}\circ {\cal S}\circ\rho_0^{\theta^\eta_{m}}\circ\ \circ_{l<m}\ 
({\cal S}\circ\rho_0^{\theta^\eta_{m-l-1}})
\! =\!\rho_0^{\theta^\eta_{m+1}}\circ\ \circ_{l\leq m}\ ({\cal S}\circ\rho_0^{\theta^\eta_{m-l}})
\end{array}$$
since in the last induction we proved that ${\cal S}^{m+1}\circ\rho_0^{(0,m+2)}\! =\!
\rho_0^{\theta^\eta_{m+1}}\circ {\cal S}^{m+1}\circ\rho^{(0,m+1)}_{0}$. Thus 
$$\begin{array}{ll}
{\cal C}_m\!\!\!\!
& \! =\!\{\alpha\!\in\! 2^\omega\mid\rho^{(0,m+1)}_{0}\big(\rho (\alpha )\big)(m)\! =\! 1\}
\! =\!\{\alpha\!\in\! 2^\omega\mid ({\cal S}^m\circ\rho^{(0,m+1)}_{0}\circ\rho )(\alpha )(0)\! =\! 1\}\cr\cr 
& \! =\!\{\alpha\!\in\! 2^\omega\mid 
\big(\rho_0^{\theta^\eta_{m}}\circ\ \circ_{l<m}\ ({\cal S}\circ\rho_0^{\theta^\eta_{m-l-1}})\circ\rho\big)
(\alpha )(0)\! =\! 1\}.
\end{array}$$
So it is enough to see that 
$\rho^m\! :=\!\rho_0^{\theta^\eta_{m}}\circ\ \circ_{l<m}\ ({\cal S}\circ\rho_0^{\theta^\eta_{m-l-1}})\circ\rho$ is an independent $(\theta_m\! +\! 1)$-function.\bigskip

 We argue by induction on $m$. For $m\! =\! 0$, we are done since  $\rho_0^{\theta^\eta_{0}}\circ\rho$ is by induction assumption an independent $(\xi\! +\!\theta^\eta_0)$-function, and 
$\xi\! +\!\theta^\eta_0\! =\!\xi\! +\!\zeta^\eta_0\! +\! 1\! =\!\theta_0\! +\! 1$. Assume that the property is true for $m$. Then $\rho^{m+1}\! =\!\rho_0^{\theta^\eta_{m+1}}\circ {\cal S}\circ\rho^m$. By induction assumption, 
$\rho^m$ is an independent $(\theta_m\! +\! 1)$-function. By Lemma 5.1.6 and the example just before it, 
${\cal S}\circ\rho^m$ is also an independent $(\theta_m\! +\! 1)$-function. By induction assumption, 
$\rho^{m+1}$ is an independent $(\theta_m\! +\! 1\! +\!\theta^\eta_{m+1})$-function, and 
$$\theta_m\! +\! 1\! +\!\theta^\eta_{m+1}
\! =\!\xi\! +\!\Sigma_{l<m}\ \theta^\eta_{l}\! +\!\zeta^\eta_{m}\! +\! 1\! +\!\theta^\eta_{m+1}
\! =\!\xi\! +\!\Sigma_{l\leq m}\ \theta^\eta_{l}\! +\!\zeta^\eta_{m+1}\! +\! 1\! =\!\theta_{m+1}\! +\! 1.$$ 
This finishes the proof.\hfill{$\square$}

\subsection{$\!\!\!\!\!\!$ Some complicated sets}\indent

 Now we come to the existence of complicated sets, as in the statement of Theorem 1.7. Their construction is based on Theorem 2.7 in [Lo-SR2] that we now change. The main problem is that we want to ensure the ccs conditions of Lemma 2.6. To do this, we modify the definition of the maps $\tau_i$ of Lemma 2.11 in [Lo-SR2].\bigskip

\noindent\bf Notation.\rm\ Let $i$ be an integer. We define $\tau_i\! :\!\omega\!\rightarrow\!\omega$ by
$$\tau_i(k)\! :=\!\left\{\!\!\!\!\!\!\!
\begin{array}{ll}
& <0,k>\mbox{if }i\! =\! 0\mbox{,}\cr\cr
& <<i,(k)_0>,(k)_1>\mbox{if }i\!\geq\! 1\mbox{,}
\end{array}\right.$$
so that $\tau_i$ is one-to-one. This allows us to define, for each $\alpha\!\in\! 2^\omega$, $\alpha_i\! :=\!\tilde {\tau_i}(\alpha )$. If $s\!\in\! (\omega\!\setminus\!\{ 0\})^{<\omega}$, then we set 
$\tilde\tau_s\! :=\!\tilde\tau_{s(0)}\circ ...\circ\tilde\tau_{s(\vert s\vert -1)}$.

\begin{lemm} Let ${\bf\Gamma}$ be a non self-dual Wadge class of Borel sets, and 
$H$ a ${\bf\Gamma}$-strategically complete set. Then\smallskip

\noindent (a) The set $\tilde {\tau_i}^{-1}(H)$ is ${\bf\Gamma}$-strategically complete for each integer $i$.\smallskip

\noindent (b) Assume that $\tau\! :\!\omega\!\rightarrow\!\omega$ is one-to-one such that the fact that 
$\alpha\!\in\! H$ depends only on $\alpha\circ\tau$. Then $L\! :=\!\{\alpha\circ\tau\mid\alpha\!\in\! H\}$ is 
${\bf\Gamma}$-strategically complete.\end{lemm}

\noindent\bf Proof.\rm\ (a) As $\tilde {\tau_i}$ is continuous,  
$\tilde {\tau_i}^{-1}(H)\!\in\! {\bf\Gamma}(2^\omega )$. We define a continuous map 
$f_{\tau_i}\! :\! 2^\omega\!\rightarrow\! 2^\omega$ by 
$f_{\tau_i}(\alpha )(m)\! :=\!\alpha\big(\tau_i^{-1}(m)\big)$ if $m$ is in the range of 
$\tau_i$, $0$ otherwise. Note that $\tilde {\tau_i}\big( f_{\tau_i}(\alpha )\big)\! =\!\alpha$, so that 
$H\! =\! f_{\tau_i}^{-1}\big(\tilde {\tau_i}^{-1}(H)\big)$. This implies that $\tilde {\tau_i}^{-1}(H)$ is 
${\bf\Gamma}$-strategically complete.\bigskip

\noindent (b) As in (a), we consider the continuous map $f_\tau$, so that 
$\tilde\tau\big( f_\tau (\beta )\big)\! =\!\beta$ for each $\beta\!\in\! 2^\omega$. Here again we get that $f_\tau^{-1}(H)\!\in\! {\bf\Gamma}(2^\omega )$. Let $\beta\!\in\! L$, which gives 
$\alpha\!\in\! H$ with $\beta\! =\!\alpha\circ\tau$. As 
$f_\tau(\beta )\circ\tau\! =\!\tilde\tau\big( f_\tau(\beta )\big)\! =\!\beta$, we get 
$f_\tau(\beta )\circ\tau\! =\!\alpha\circ\tau$ and $f_\tau (\beta )\!\in\! H$ by the assumption on $H$. Conversely, if $f_\tau (\beta )\!\in\! H$, then 
$\beta\! =\!\tilde\tau\big( f_\tau (\beta )\big)\! =\!f_\tau (\beta )\circ\tau\!\in\! L$. Thus $f_\tau^{-1}(H)\! =\! L$, and $L\!\in\! {\bf\Gamma}(2^\omega )$.\bigskip

 If $\alpha\!\in\! H$, then $\tilde\tau (\alpha )\! =\!\alpha\circ\tau\!\in\! L$. Conversely, assume that 
$\tilde\tau (\alpha )\!\in\! L$. Then there is $\beta\!\in\! H$ with $\beta\circ\tau\! =\!\alpha\circ\tau$. The assumption on $H$ implies that $\alpha\!\in\! H$. Thus $H\! =\!\tilde\tau^{-1}(L)$ and $L$ is 
${\bf\Gamma}$-strategically complete.$\hfill{\square}$

\begin{lemm} Let ${\bf\Gamma}$ be a Wadge class of Borel sets, and $A\!\subseteq\! 2^\omega$. Then 
$A\!\in\! {\bf\Gamma}(2^\omega )$ if and only if there is $B\!\in\! {\bf\Gamma}(\omega^\omega )$ with 
$A\! =\! B\cap 2^\omega$.\end{lemm}

\noindent\bf Proof.\rm\ $\Rightarrow$ Let $r\! :\!\omega^\omega\!\rightarrow\! 2^\omega$ be a continuous retraction. We just have to set $B\! :=\! r^{-1}(A)$.\bigskip

\noindent $\Leftarrow$ Let $i\! :\! 2^\omega\!\rightarrow\!\omega^\omega$ be the canonical injection. Then $A\! =\! i^{-1}(B)\!\in\! {\bf\Gamma}(2^\omega )$.$\hfill{\square}$\bigskip

 This lemma shows that the notation ${\bf\Gamma}_u$ in Theorem 5.1.3 will not create any trouble, since it is equivalent to the one in Definition 5.1.2.\bigskip

\noindent\bf Notation.\rm\ The following notation can essentially be found in [Lo-SR2] (after Lemma 2.5). Let $\cal R$ be the least set of functions from $2^\omega$ into itself which contains the functions 
$\rho_0^\eta$, the functions $\tilde {\tau_i}$ for $i\!\geq\! 1$, and is closed under composition. By Lemma 5.1.6 and Theorem 5.1.7, each $\rho\!\in\! {\cal R}$ is an independent $\eta$-function for some $\eta$ called the $order$ $o(\rho )$ of $\rho$.

\begin{defin} Let $u\!\in\! {\cal D}$. A set $H\!\subseteq\! 2^\omega$ is $strongly$ 
$u$-$strategically\ complete$ if for each $\rho\!\in\! {\cal R}$ of order $\eta$, $\rho^{-1}(H)$ is ${\bf\Gamma}_{u^\eta}$-strategically complete and ccs.\end{defin}
 
\begin{them} Let $u\!\in\! {\cal D}$. Then there exists a strongly $u$-strategically complete set 
$H_u\!\subseteq\! 2^\omega$. In particular, $H_u$ is ${\bf\Gamma}_u$-complete and ccs.\end{them}

\noindent\bf Proof.\rm\ We will check that the sets $H_u$ given by Theorem 2.7 in [Lo-SR2] essentially work, even if we change them.\bigskip

 The construction is by induction on $u\!\in\! {\cal D}$. Let us say that $u$ is $nice$ if it satisfies the conclusion of the theorem. By Proposition 5.1.5, it is enough to prove that $0^\infty$ is nice, that $u(0)1u$ is nice if $u$ is nice, that $u^\eta$ is nice if $u$ is nice and $\eta\! <\!\omega_1$, and that $12<u_p>$ is nice if the $u_p$'s are nice.\bigskip

\noindent $\bullet$ We set $H_{0^\infty}\! :=\!\emptyset$, which is clearly strongly $0^\infty$-strategically complete.

\vfill\eject

\noindent $\bullet$ Assume that $u$ is nice. We set $H_{u(0)1u}\! :=\!\neg H_u$, which is strongly 
$u(0)1u$-strategically complete. Indeed, if $u(0)\! =\! 0$, then 
${\bf\Gamma}_{(u(0)1u)^\eta}\! =\! {\bf\Gamma}_{u(0)1u}\! =\!\check {\bf\Gamma}_{u}\! =\!\check {\bf\Gamma}_{u^\eta}$. If $u(0)\!\geq\! 1$, then 
$${\bf\Gamma}_{(u(0)1u)^\eta}\! =\! 
{\bf\Gamma}_{(1+\eta +(u(0)-1))1u^\eta}\! =\!\check {\bf\Gamma}_{u^\eta}$$ 
since $u^\eta (0)\! =\! 1\! +\!\eta\! +\!\big( u(0)\! -\! 1\big)$.\bigskip

\noindent $\bullet$ Assume that $u$ is nice, and let $\eta\! <\!\omega_1$. We set 
$H_{u^\eta}\! :=\! (\rho_0^\eta )^{-1}(H_u)$, which is strongly $u^\eta$-strategically complete. Indeed, 
let $\rho\!\in\! {\cal R}$ of order $\xi$. Then $\rho^{-1}(H_{u^\eta})\! =\! (\rho_0^\eta\circ\rho )^{-1}(H_u)$ is 
${\bf\Gamma}_{u^{\xi +\eta}}$-strategically complete and compatible with comeager sets since $u$ is nice and $\rho_0^\eta\circ\rho$ is in $\cal R$ of order $\xi\! +\!\eta$. It remains to notice that 
$(u^\eta )^\xi\! =\! u^{\xi +\eta}$, which is clear by induction on $u$ and by definition of the ordinal subtraction.\bigskip

\noindent $\bullet$ Assume that the $u_p$'s are nice. We set
$$\alpha\!\in\! H_{12<u_p>}\ \Leftrightarrow\ \left\{\!\!\!\!\!\!\!
\begin{array}{ll}
& \alpha_0\! =\! 0^\infty\wedge\alpha_1\!\in\! H_{u_0}\cr
& \mbox{ or }\cr
& \exists m\!\in\!\omega\ \ \alpha_0(m)\! =\! 1\wedge\forall l\! <\! m\ \ \alpha_0(l)\! =\! 0\wedge
\alpha_{(m)_0+2}\!\in\! H_{u_{((m)_0+2)_0+1}}\mbox{.}
\end{array}\right.$$
- Recall that 
${\bf\Gamma}_{12<u_p>}\! =\! S_1(\bigcup_{p\geq 1}\ {\bf\Gamma}_{u_p},{\bf\Gamma}_{u_0})$. We set $H'_0\! :=\!\{\alpha\!\in\! 2^\omega\mid\alpha_1\!\in\! H_{u_0}\}\! =\!\tilde {\tau_1}^{-1}(H_{u_0})$, and for 
$n\!\geq\! 2$,
$$\begin{array}{ll}
H'_n
& \!\!\!\! :=\!\{\alpha\!\in\! 2^\omega\mid\alpha_n\!\in\! 
H_{u_{(n)_0+1}}\}\! =\!\tilde {\tau_n}^{-1}(H_{u_{(n)_0+1}})\mbox{, }\cr\cr
C_n
& \!\!\!\! :=\!\{\alpha\!\in\! 2^\omega\mid\exists m\!\in\!\omega\ \ \alpha_0 (m)\! =\! 1\ \mbox{ and }\ 
\forall l\! <\! m\ \ \alpha_0(l)\! =\! 0\ \mbox{ and }\ (m)_0\! +\! 2\! =\! n\}.
\end{array}$$
Note that $(C_n)_{n\geq 2}$ is a sequence of pairwise disjoint open sets, and 
$H'_0\!\in\! {\bf\Gamma}_{u_0}$, $H'_n\!\in\! {\bf\Gamma}_{u_{(n)_0+1}}$ if $n\!\geq\! 2$ by Lemma 
5.2.1.(a). Moreover, $H_{12<u_p>}\! =\!\bigcup_{n\geq 2}\ (H'_n\cap C_n)\cup
(H'_0\!\setminus\!\bigcup_{n\geq 2}\ C_n)\!\in\! {\bf\Gamma}_{12<u_p>}(2^\omega )$, by Lemma 5.2.2 and the reduction property for the class of open sets (see 22.16 in [K]).\bigskip

\noindent - Let $\rho\!\in\! {\cal R}$ of order $\eta$. Then 
$\rho^{-1}(H_{12<u_p>})\!\in\!  {\bf\Gamma}_{(12<u_p>)^\eta}(2^\omega )$, by Proposition 5.1.5.(a) and a retraction argument in the style of the proof of Lemma 5.2.2. Let $\pi$ be associated with $\rho$, 
$\theta_0\! :\!\omega\!\rightarrow\!\omega$ be a one-to-one enumeration of 
$\pi^{-1}\big(\mbox{Ran}(\tau_1)\big)$, and, for $n\!\geq\! 2$, $\theta_n\! :\!\omega\!\rightarrow\!\omega$ be a one-to-one enumeration of $\pi^{-1}\big(\mbox{Ran}(\tau_n)\big)$ and 
$\theta_0^n\! :\!\omega\!\rightarrow\!\omega$ be a one-to-one enumeration of 
$$\pi^{-1}\big(\big\{ j\!\in\!\mbox{Ran}(\tau_0)\mid\big(\tau_0^{-1}(j)\big)_0\! +\! 2\! =\! n\big\}\big).$$
As $\tau_i$ is one-to-one, $\mbox{Ran}(\tau_i)$ is infinite, and $\pi^{-1}\big(\mbox{Ran}(\tau_i)\big)$ is also infinite since $\pi$ is onto. This proves the existence of the $\theta_n$'s and of the $\theta_0^n$'s. Note that  the $\mbox{Ran}(\tau_i)$'s are pairwise disjoint since ${0\! = <0,0>}$. This implies that the elements of $\{\mbox{Ran}(\theta_n)\mid n\!\not=\! 1\}\cup\{\mbox{Ran}(\theta_0^n)\mid n\!\geq\! 2\}$ 
are pairwise disjoint.\bigskip

\noindent - Note that the fact that $\alpha\!\in\! H^\eta_n\! :=\!\rho^{-1}(H'_n)$ depends only on 
$\alpha\circ\theta_n$ if $n\!\not=\! 1$. We set, for $n\!\not=\! 1$, 
$$L^\eta_n\! :=\!\{\alpha\circ\theta_n\mid\alpha\!\in\! H^\eta_n\}.$$

 Note that $\rho^{-1}(H'_0)\! =\!\rho^{-1}\big(\tilde {\tau_1}^{-1}(H_{u_0})\big)\! =\! 
(\tilde {\tau_1}\circ\rho )^{-1}(H_{u_0})$ is ${\bf\Gamma}_{u_0^\eta}$-strategically complete since $u_0$ is nice and $\tilde {\tau_1}\circ\rho$ is in $\cal R$ of order $\eta$. Similarly, $\rho^{-1}(H'_n)$ is 
${\bf\Gamma}_{u_{(n)_0+1}^\eta}$-strategically complete if $n\!\geq\! 2$. By Lemma 5.2.1.(b), we get that 
$L^\eta_0$ is ${\bf\Gamma}_{u_0^\eta}$-strategically complete, and $L^\eta_n$ is 
${\bf\Gamma}_{u_{(n)_0+1}^\eta}$-strategically complete if $n\!\geq\! 2$.

\vfill\eject

\noindent - We set, for $n\!\geq\! 2$, $C^\eta_n\! :=\!\{\alpha\circ\theta_0^n\mid\exists m\!\in\!\omega\ \ \rho(\alpha )_0(m)\! =\! 1\ \mbox{ and }\ (m)_0\! +\! 2\! =\! n\}$. Let us prove that $C^\eta_n$ is 
${\bf\Sigma}^0_{1+\eta}$-strategically complete.\bigskip

 Note first that $\{\alpha\!\in\! 2^\omega\mid f(\alpha )\!\not=\! 0^\infty\}$ is ${\bf\Sigma}^0_{1+\eta}$-strategically complete if $f$ is an independent $\eta$-function. Indeed, with the notation of Definition 3.3, we can write 
$$\{\alpha\!\in\! 2^\omega\mid f(\alpha )\! =\! 0^\infty\}\! =\!\bigcap_{m\in\omega}\ \neg\ {\cal C}_m.$$ Moreover, the fact that $\alpha\!\in\! {\cal C}_m$ depends only of the values of $\alpha$ on 
$\pi_f^{-1}(\{ m\})$. \bigskip

 Assume first that $\eta\!\geq\! 1$. As $f$ is an independent $\eta$-function, ${\cal C}_m$ is 
${\bf\Pi}^0_{1+\theta_m}$-strategically complete, for some $\theta_m\! <\!\eta$ satisfying 
$\theta_m\! +\! 1\! =\!\eta$ if $\eta$ is a successor ordinal, 
$\mbox{sup}_{m\in\omega}\ \theta_m\! =\!\eta$ if $\eta$ is a limit ordinal. Note that 
$\eta\! =\!\mbox{sup}_{m\in\omega}\ (\theta_m\! +\! 1)$. By Lemma 3.7 in [Lo-SR1], 
$\{\alpha\!\in\! 2^\omega\mid f(\alpha )\! =\! 0^\infty\}$ is ${\bf\Pi}^0_{1+\eta}$-strategically complete.\bigskip

 Assume now that $\eta\! =\! 0$. As in the proof of Theorem 5.1.7 we see that 
$\{\alpha\!\in\! 2^\omega\mid f(\alpha )\! =\! 0^\infty\}$ is ${\bf\Pi}^0_{1+\eta}$-strategically complete.\bigskip

 Now we come back to the $C^\eta_n$'s. We define $\tau\! :\!\omega\!\rightarrow\!\omega$ by 
$\tau (k)\! :=\ <\! n\! -\! 2,k\! >$, so that $\tau$ is one-to-one and 
$\mbox{Ran}(\tau )\! =\!\{ m\!\in\!\omega\mid (m)_0\! =\! n\! -\! 2\}$. As $\rho$ is an independent $\eta$-function, $\tilde {\tau_0}\circ\rho$ and $\tilde \tau\circ\tilde {\tau_0}\circ\rho$ are also independent $\eta$-functions by Lemma 5.1.6. The previous point shows that 
$$L_n\! :=\!\{\alpha\!\in\! 2^\omega\mid (\tilde\tau\circ\tilde {\tau_0}\circ\rho )(\alpha )\!\not=\! 0^\infty\}$$ 
is ${\bf\Sigma}^0_{1+\eta}$-strategically complete. But
$$\begin{array}{ll}
L_n\!\!\!\!
& =\!\{\alpha\!\in\! 2^\omega\mid\exists k\!\in\!\omega\ \ \tilde\tau\big( (\tilde {\tau_0}\circ\rho )(\alpha )\big)(k)
\! =\! 1\}\! 
=\!\{\alpha\!\in\! 2^\omega\mid\exists k\!\in\!\omega\ \ (\tilde {\tau_0}\circ\rho )(\alpha )\big(\tau (k)\big)
\! =\! 1\}\cr\cr
& =\!\{\alpha\!\in\! 2^\omega\mid\exists m\!\in\!\omega\ \ (\tilde {\tau_0}\circ\rho )(\alpha )(m)\! =\! 1\ 
\mbox{ and }\ (m)_0\! +\! 2\! =\! n\}
\end{array}$$
and the fact that $\alpha\!\in\! L_n$ depends only on $\alpha\circ\theta_0^n$. By Lemma 5.2.1.(b), we get that $C^\xi_n$ is ${\bf\Sigma}^0_{1+\eta}$-strategically complete.\bigskip

\noindent - Let $H^*\!\in\! {\bf\Gamma}_{(12<u_p>)^\eta}(\omega^\omega )$, say 
$H^*\! =\!\bigcup_{n\geq 2}\ (H^*_n\cap C^*_n)\cup (H^*_0\!\setminus\!\bigcup_{n\geq 2}\ C^*_n)$, with pairwise disjoint $C^*_n\!\in\! {\bf\Sigma}^0_{1+\eta}$, $H^*_0\!\in\! {\bf\Gamma}_{u_0^\eta}$, and without loss of generality $H^*_n\!\in\! {\bf\Gamma}_{u_{(n)_0+1}^\eta}$. Then Player 2 has for each $n\!\not=\! 1$ a winning strategy $\sigma_n$ in $G(H^*_n,L^\eta_n)$, and for each $n\!\geq\! 2$ a winning strategy 
$\sigma^*_n$ in $G(C^*_n,C^\eta_n)$. Let then Player 2 plays in 
$G\big( H^*,\rho^{-1}(H_{u_{12<u_p>}})\big)$ against $\beta$ by playing his strategies $\sigma_n$, 
$\sigma^*_n$ at the right places (the ranges of $\theta_n$ and $\theta_0^n$ respectively) against this same $\beta$, independently, and by playing $0$ out of these ranges. The result is some $\alpha$ such that $\alpha\circ\theta_n$ wins against $\beta$ in $G(H^*_n,L^\eta_n)$ and $\alpha\circ\theta_0^n$ wins against $\beta$ in $G(C^*_n,C^\eta_n)$. This wins, for $\alpha\!\in\!\rho^{-1}(H'_n)$ just in case 
$\beta\!\in\! H^*_n$, and $\rho (\alpha )_0$ takes value $1$ on some $m$ with $(m)_0\! +\! 2\! =\! n$ just in case $\beta\!\in\! C^*_n$. But as the $C^*_n$ are pairwise disjoint, there is at most one $n$ in 
$\{ (m)_0\! +\! 2\mid\rho (\alpha )_0(m)\! =\! 1\}$, and $\alpha\!\in\!\rho^{-1}(C_n)$ just in case 
$\beta\!\in\! C^*_n$. Thus $\rho^{-1}(H_{12<u_p>})$ is ${\bf\Gamma}_{(12<u_p>)^\eta}$-strategically complete.\bigskip

\noindent - It remains to see that $\rho^{-1}(H_{12<u_p>})$ is ccs. So let $\alpha_0\!\in\! d^\omega$ and 
$F\! :\! 2^\omega\!\rightarrow\! (d^\omega )^{d-1}$ satisfying the conclusion of Lemma 2.4.(b).\bigskip

\noindent $\circ$ Let $N\!\geq\! 1$ and $M\!\in\!\omega$. Then $\rho (\alpha )_N\!\in\! H_{u_M}\Leftrightarrow (\tilde {\tau_N}\circ\rho )(\alpha )\!\in\! H_{u_M}\Leftrightarrow\alpha\!\in\! (\tilde {\tau_N}\circ\rho )^{-1}(H_{u_M})$. As $N\!\geq\! 1$, $\tilde {\tau_N}\circ\rho$ is in $\cal R$, and 
$(\tilde {\tau_N}\circ\rho )^{-1}(H_{u_M})$ is ccs since $u_M$ is nice. Thus 
$\rho (\alpha )_N\!\in\! H_{u_M}$ if and only if 
$\rho\Big( {\cal S}\big(\alpha_0\Delta F_0(\alpha )\big)\Big)_N\!\in\! H_{u_M}$.\bigskip

\noindent $\circ$ Recall the notation before Lemma 2.4. We define 
$q\! :\!\omega^{<\omega}\!\setminus\!\{\emptyset\}\!\rightarrow\!\omega$ as follows:
$$q(t)\! :=\!\left\{\!\!\!\!\!\!\!
\begin{array}{ll}
& t(0)\mbox{ if }\vert t\vert\! =\! 1\mbox{,}\cr\cr
& <t(\vert t\vert\! -\! 1),q(t^-)>\mbox{if }\vert t\vert\!\geq\! 2.
\end{array}
\right.$$
$\circ$ Let us prove that $\tilde\tau_s(\alpha )(n)\! =\!\alpha (<q\big( (n)_0s\big),(n)_1>)$ for each 
$s\!\in\! (\omega\!\setminus\!\{ 0\} )^{<\omega}$.\bigskip

 We argue by induction on $\vert s\vert$. So assume that the result is proved for $\vert s\vert\!\leq\! l$, which is the case for $l\! =\! 0$. Assume that $\vert s\vert\! =\! l\! +\! 1$. We get
$$\begin{array}{ll}
\tilde\tau_s(\alpha )(n)\!\!\!\!\!
& =\! \tilde\tau_{s\vert l}\big(\tilde {\tau_{s(l)}}(\alpha )\big)(n)
\! =\!\tilde {\tau_{s(l)}}(\alpha )(<\! q\big( (n)_0(s\vert l)\big),\! (n)_1\! >)
\! =\!\alpha\big(\tau_{s(l)}(<\! q\big( (n)_0(s\vert l)\big),\! (n)_1\! >)\big)\cr\cr
& =\!\alpha\big(\big<\! <\! s(l),q\big( (n)_0(s\vert l)\big)\! >,(n)_1\big>\big)
\! =\!\alpha (<\! q\big( (n)_0s\big),(n)_1\! >).
\end{array}$$
$\circ$ Let us prove that $(\rho_0\circ\tilde\tau_s)(\alpha )\! =\!
(\rho_0\circ\tilde\tau_s)\Big( {\cal S}\big(\alpha_0\Delta F_0(\alpha)\big)\Big)$ for each 
$s\!\in\! (\omega\!\setminus\!\{ 0\})^{<\omega}$ and each $\alpha\!\in\! 2^\omega$. This comes from the following equivalences:
$$\begin{array}{ll}
(\rho_0\circ\tilde\tau_s)(\alpha )(n)\! =\! 0\!\!\!
& \Leftrightarrow\exists m\!\in\!\omega\ \ \tilde\tau_s(\alpha )(<n,m>)\! =\! 1
\Leftrightarrow\exists m\!\in\!\omega\ \ \alpha (<q(ns),m>)\! =\! 1\cr\cr
& \Leftrightarrow\exists m'\!\in\!\omega\ \ 
{\cal S}\big(\alpha_0\Delta F_0(\alpha )\big)(<q(ns),m'>)\! =\! 1\cr\cr
& \Leftrightarrow (\rho_0\circ\tilde\tau_s)\Big( {\cal S}\big(\alpha_0\Delta F_0(\alpha )\big)\Big)(n)\! =\! 0. 
 \end{array}$$
$\circ$ Let us prove that $(\rho_0^\eta\circ\tilde\tau_s)(\alpha )\! =\!
(\rho_0^\eta\circ\tilde\tau_s)\Big( {\cal S}\big(\alpha_0\Delta F_0(\alpha)\big)\Big)$ for each 
$1\!\leq\!\eta\! <\!\omega_1$, each $s\!\in\! (\omega\!\setminus\!\{ 0\})^{<\omega}$ and each 
$\alpha\!\in\! 2^\omega$.\bigskip

 We argue by induction on $\eta$. For $\eta\! =\! 1$, this comes from the previous point. If $\theta\!\geq\! 1$ and $\eta\! =\!\theta\! +\! 1$, then this comes from the fact that $\rho_0^\eta\! =\!\rho_0\circ\rho_0^\theta$. If 
$\eta$ is a limit ordinal and $m$ is an integer, then 
$$\begin{array}{ll}
& \!\!\!\! (\rho_0^\eta\circ\tilde\tau_s)(\alpha )(m)\cr\cr
=
& \!\!\!\!\rho_0^\eta\big(\tilde\tau_s(\alpha )\big)(m)
\! =\!\rho_0^{(0,m+1)}\big(\tilde\tau_s(\alpha )\big)(m)\cr\cr
=
& \!\!\!\! (\rho_0^{(m,m+1)}\circ ...\circ\rho_0^{(1,2)})
\Big(\rho_0^{(0,1)}\big(\tilde\tau_s(\alpha )\big)\Big)(m)
\! =\! (\rho_0^{(m,m+1)}\circ ...\circ\rho_0^{(1,2)})
\Big(\rho_0^{\theta^\eta_0}\big(\tilde\tau_s(\alpha )\big)\Big)(m)\cr\cr
=
& \!\!\!\! (\rho_0^{(m,m+1)}\circ ...\circ\rho_0^{(1,2)})
\Bigg(\rho_0^{\theta^\eta_0}
\bigg(\tilde\tau_s\Big( {\cal S}\big(\alpha_0\Delta F_0(\alpha)\big)\Big)\bigg)\Bigg)(m)
\! =\! (\rho_0^\eta\circ\tilde\tau_s)\Big({\cal S}\big(\alpha_0\Delta F_0(\alpha)\big)\Big)(m).
\end{array}$$
$\circ$ Note that 
$\rho (\alpha )_0\! =\! 0^\infty\Leftrightarrow\alpha\!\in\! (\tilde {\tau_0}\circ\rho )^{-1}(\{ 0^\infty\} )$. Let us prove that $(\tilde {\tau_0}\circ\rho )^{-1}(\{ 0^\infty\} )$ is ccs.\bigskip

 We can write $\rho\! =\!\circ_{j\leq l}\ \rho^j$, where $l$ is an integer and each $\rho^j$ is either of the form $\rho_0^\eta$, or one of the $\tilde {\tau_i}$'s for $i\!\geq\! 1$. By the previous point, we may assume that that each $\rho^j$ is either $\rho_0^0\! =\!\mbox{Id}_{2^\omega}$, or one of the 
$\tilde {\tau_i}$'s for $i\!\geq\! 1$. So there is $s\!\in\! (\omega\!\setminus\!\{ 0\})^{<\omega}$ such that 
$\rho\! =\!\tilde\tau_s$. We get
$$\begin{array}{ll}
\alpha\!\notin\! (\tilde {\tau_0}\circ\rho )^{-1}(\{ 0^\infty\} )\!\!\!\!
& \Leftrightarrow\exists m\!\in\!\omega\ \ (\tilde {\tau_0}\circ\rho )(\alpha )(m)\! =\! 1
\Leftrightarrow\exists m\!\in\!\omega\ \ \rho (\alpha )\big(\tau_0(m)\big)\! =\! 1\cr\cr
& \Leftrightarrow\exists m\!\in\!\omega\ \ \tilde\tau_s(\alpha )(<0,m>)\! =\! 1
\Leftrightarrow\exists m\!\in\!\omega\ \ \alpha (<q(0s),m>)\! =\! 1\cr\cr
& \Leftrightarrow\exists m\!\in\!\omega\ \ \alpha\big( p(q(0s),m)\big)\! =\! 1\cr\cr
& \Leftrightarrow\exists m'\!\in\!\omega\ \ {\cal S}\big(\alpha_0\Delta F_0(\alpha )\big)
\big( p(q(0s),m')\big)\! =\! 1\cr\cr
& \Leftrightarrow {\cal S}\big(\alpha_0\Delta F_0(\alpha )\big)\!\notin\! 
(\tilde {\tau_0}\circ\rho )^{-1}(\{ 0^\infty\}).
\end{array}$$
Thus $\rho (\alpha )_0\! =\! 0^\infty\Leftrightarrow
\rho\Big( {\cal S}\big(\alpha_0\Delta F_0(\alpha )\big)\Big)_0\! =\! 0^\infty$.\bigskip

\noindent $\circ$ It remains to see that if $\rho (\alpha )_0\!\not=\! 0^\infty$ and $m_\alpha$ is minimal with 
$\rho (\alpha )_0(m_\alpha )\! =\! 1$, then 
$$(m_\alpha )_0\! =\! (m_{{\cal S}(\alpha_0\Delta F_0(\alpha ))})_0.$$
As in the previous point we may assume that there is $s\!\in\! (\omega\!\setminus\!\{ 0\})^{<\omega}$ such that $\rho\! =\!\tilde\tau_s$. The computations of the previous point show that 
$\rho (\alpha )_0(m)\! =\!\alpha (<q(0s),m>)$ for each integer $m$. Note that 
$$n_\alpha\! :=<\! q(0s),m_\alpha\! >=\!\mbox{min}\{ n\!\in\!\omega\mid\alpha (n)\! =\! 1
\wedge (n)_0\! =\! q(0s)\}$$ 
since $<\! q( 0s),.\! >$ is increasing, and similarly 
$$<\! q( 0s),m_{{\cal S}(\alpha_0\Delta F_0(\alpha ))}\! >=\!\mbox{min}\{ m\!\in\!\omega\mid {\cal S}\big(\alpha_0\Delta F_0(\alpha )\big)(m)\! =\! 1\wedge (m)_0\! =\! q(0s)\}.$$
But 
$$B_\alpha[\{ n\!\in\!\omega\mid\alpha (n)\! =\! 1\ \mbox{ and }\ (n)_0\! =\! q(0s)\}]\! =\!\{ m\!\in\!\omega\mid {\cal S}\big(\alpha_0\Delta F_0(\alpha )\big)(m)\! =\! 1\ \mbox{ and }\ (m)_0\! =\! q(0s)\}$$ 
since $B_\alpha$ is a bijection satisfying $(n)_0\! =\!\big( B_\alpha (n)\big)_0$. As $B_\alpha$ is increasing we get 
$$B_\alpha (n_\alpha )\! =<\! q( 0s),m_{{\cal S}(\alpha_0\Delta F_0(\alpha ))}\! >\! .$$
Thus $(m_{{\cal S}(\alpha_0\Delta F_0(\alpha ))})_0\! =\!\Big(\big( B_\alpha (n_\alpha )\big)_1\Big)_0\! =\!\big( (n_\alpha )_1\big)_0\! =\! (m_\alpha )_0$ and we are done.$\hfill{\square}$

\begin{coro} Let ${\bf\Gamma}$ be a non self-dual Wadge class of Borel sets. Then there 
is $C_{\bf\Gamma}\!\subseteq\! 2^\omega$ which is ${\bf\Gamma}$-complete and ccs.\end{coro}

\noindent\bf Proof.\rm\ By Theorem 5.1.3 there is $u\!\in\! {\cal D}$ such that 
${\bf\Gamma}(\omega^\omega )\! =\! {\bf\Gamma}_u(\omega^\omega )$. By Theorem 5.2.4 there is 
$H_u\!\subseteq\! 2^\omega$ which is strongly ${\bf\Gamma}_u$-strategically complete. It is clear that 
$C_{\bf\Gamma}\! :=\! H_u$ is suitable.$\hfill{\square}$\bigskip

 Now we can prove Theorem 1.7.(1). But we need some more material to prove Theorem 1.7.(2).

\begin{defin} (a) A set $U\!\subseteq\! 2^\omega$ is $strongly\ ccs$ if for each 
$s\!\in\! (\omega\!\setminus\!\{ 0\})^{<\omega}$ the set $\tilde\tau_s^{-1}(U)$ is ccs.\smallskip

\noindent (b) Let $\bf\Gamma$ be a Wadge class of Borel sets, and 
$U_0,U_1\!\in\! {\bf\Gamma}(2^\omega )$ disjoint. We say that $(U_0,U_1)\ is\ com\mbox{-}$ 
$plete\ for\ pairs\ of\ disjoint\ {\bf\Gamma}\ sets$ if for any pair $(A_0,A_1)$ of disjoint $\bf\Gamma$ subsets of $\omega^\omega$ there is $f\! :\!\omega^\omega\!\rightarrow\! 2^\omega$ continuous such that 
$A_\varepsilon\! =\! f^{-1}(U_\varepsilon )$ for each $\varepsilon\!\in\! 2$. Similarly, we can define the notion of a sequence $(U_p)_{p\geq 1}$ complete for sequences of pairwise disjoint $\bf\Gamma$ sets.\end{defin}

\begin{lemm} (a) There is $(U_1,U_2)$ complete for pairs of disjoint $\boraone$ sets with 
$U_\varepsilon$ strongly ccs, and such that for each $s\!\in\! (\omega\!\setminus\!\{ 0\})^{<\omega}$ there is a pair $(O_1,O_2)$ of ccs $\boraone$ sets reducing the pair 
${\big(\tilde\tau_{1s1}^{-1}(U_1\cup U_2),\tilde\tau_{1s2}^{-1}(U_1\cup U_2)\big)}$.\smallskip

\noindent (b) There is $(U_p)_{p\geq 1}$ complete for sequences of pairwise disjoint $\boraone$ sets with $U_p$ strongly ccs, and such that for each $s\!\in\! (\omega\!\setminus\!\{ 0\})^{<\omega}$ there is a sequence $(O^\varepsilon_p)_{\varepsilon\in\{ 1,2\},p\geq 1}$ of ccs $\boraone$ sets reducing  
$\big(\tilde\tau_{s\varepsilon}^{-1}(U_p)\big)_{\varepsilon\in\{ 1,2\},p\geq 1}$.\end{lemm}

\noindent\bf Proof.\rm\ (a) Recall the definition of $H_1$ after Definition 3.3: $H_1\! :=\!\{ 0^\infty\}$. We saw that $H_1\!\in\!\bormone (2^\omega)$ and is $\bormone$-complete. We set $U\! :=\!\neg H_1$, so that $U$ is $\boraone$-complete. Let $(A_1,A_2)$ be a pair of disjoint $\boraone$ subsets of $\omega^\omega$. As $U$ is complete there are $f_1,f_2\! :\!\omega^\omega\!\rightarrow\! 2^\omega$ continuous such that $A_\varepsilon\! =\! f_\varepsilon^{-1}(U)$ for each $\varepsilon\!\in\!\{ 1,2\}$. We define 
$f\! :\!\omega^\omega\!\rightarrow\! 2^\omega$ by 
$$f(\alpha )\big(\big< <\varepsilon ,(k)_0>,(k)_1\big>\big)\! :=\!\left\{\!\!\!\!\!\!\!
\begin{array}{ll}
& f_\varepsilon (\alpha )(k)\mbox{ if }\varepsilon\!\in\!\{ 1,2\}\mbox{,}\cr
& 0\mbox{ otherwise,}
\end{array}
\right.$$
so that $f$ is continuous and $f_\varepsilon\! =\!\tilde\tau_\varepsilon\circ f$. Now 
$A_\varepsilon\! =\! f^{-1}\big(\tilde\tau_\varepsilon^{-1}(U)\big)$ and 
$\big(\tilde\tau_1^{-1}(U),\tilde\tau_2^{-1}(U)\big)$ is complete for pairs of $\boraone$ sets (not necessarily disjoint). Note that 
$$\begin{array}{ll}
\tilde\tau_\varepsilon^{-1}(U)\!\!\!
& \! =\!\big\{\alpha\!\in\! 2^\omega\mid\exists k\!\in\!\omega\ \ 
\alpha\big(\big< <\varepsilon ,(k)_0>,(k)_1\big>\big)\! =\! 1\big\}\cr\cr
& \! =\!\big\{\alpha\!\in\! 2^\omega\mid\exists N\!\in\!\omega\ \ \big( (N)_0\big)_0\! =\!\varepsilon\wedge
\alpha (N)\! =\! 1\big\}.
\end{array}$$
We set $V_\varepsilon\! :=\!\Big\{\alpha\!\in\! 2^\omega\mid\exists N\!\in\!\omega\ \ \big( (N)_0\big)_0\! =\!\varepsilon\wedge\alpha (N)\! =\! 1\wedge\forall l\! <\! N\ \ \Big(\big( (l)_0\big)_0\!\notin\!\{ 1,2\}
\vee\alpha (l)\! =\! 0\Big)\Big\}$. Note that $V_i\!\in\!\boraone$ and $(V_1,V_2)$ reduces 
$\big(\tilde\tau_1^{-1}(U),\tilde\tau_2^{-1}(U)\big)$. Thus
$$\alpha\!\in\! A_\varepsilon\Leftrightarrow f(\alpha )\!\in\!\tilde\tau_\varepsilon^{-1}(U)
\Leftrightarrow f(\alpha )\!\in\!\tilde\tau_\varepsilon^{-1}(U)\!\setminus\! \tilde\tau_{3-\varepsilon}^{-1}(U)
\Leftrightarrow f(\alpha )\!\in\! V_\varepsilon$$
and $(V_1,V_2)$ is complete for pairs of disjoint $\boraone$ sets. Recall the definition of $\tau_0$ before 
Lemma 5.2.1. We set $U_\varepsilon\! :=\!\tilde\tau_0^{-1}(V_\varepsilon )$, which defines a pair of disjoint $\boraone$ sets. Now $g(\alpha )\!\! :=<\!\alpha ,\alpha ,...\! >$ defines 
$g\! :\! 2^\omega\!\rightarrow\! 2^\omega$ continuous. Note that 
$\alpha\!\in\! A_\varepsilon\Leftrightarrow f(\alpha )\!\in\! V_\varepsilon
\Leftrightarrow\tilde\tau_0\Big( g\big( f(\alpha)\big)\Big)\!\in\! V_\varepsilon
\Leftrightarrow g\big( f(\alpha)\big)\!\in\! U_\varepsilon$, which shows that $(U_1,U_2)$ is complete for pairs of disjoint $\boraone$ sets.\bigskip

 Fix $s\!\in\! (\omega\!\setminus\!\{ 0\})^{<\omega}$. The proof of Theorem 5.2.4 shows that 
$\tilde\tau_s(\alpha )(n)\! =\!\alpha\big(\! <\! q\big( (n)_0s),(n)_1\! >\!\big)$. We get\bigskip
 
\leftline{$\tilde\tau_s^{-1}(U_\varepsilon)\! =\!\Big\{\alpha\!\in\! 2^\omega\mid\exists N\!\in\!\omega\ \ 
\big( (N)_0\big)_0\! =\!\varepsilon\wedge\alpha (<\! q(0s),N\! >)\! =\! 1~\wedge$}\smallskip

\rightline{$\forall l\! <\! N\ \ \Big(\big( (l)_0\big)_0\!\notin\!\{ 1,2\}\vee
\alpha (<\! q(0s),l\! >)\! =\! 0\Big)\Big\}.$}\bigskip

 Thus\bigskip
 
\leftline{$\tilde\tau_s^{-1}(U_\varepsilon)\! =\!\bigg\{\alpha\!\in\! 2^\omega\mid\exists M\!\in\!\omega\ \ 
\Big(\big( (M)_1\big)_0\Big)_0\! =\!\varepsilon\wedge (M)_0\! =\! q(0s)\wedge
\alpha (M)\! =\! 1~\wedge$}\smallskip

\rightline{$\forall l\! <\! M\ \ \bigg(\Big(\big( (M)_1\big)_0\Big)_0\!\notin\!\{ 1,2\}\vee (l)_0\!\not=\! q(0s)\vee
\alpha (l)\! =\! 0\bigg)\bigg\}.$}\bigskip

\noindent Recall the conclusion of Lemma 2.4.(b). The bijection $B_\alpha$ induces an increasing  bijection between $\Big\{ M\!\in\!\omega\mid\Big(\big( (M)_1\big)_0\Big)_0\!\in\!\{ 1,2\}\wedge 
(M)_0\! =\! q(0s)\wedge\alpha (M)\! =\! 1\Big\}$ 
and 
$$\Big\{ M'\!\in\!\omega\mid\Big(\big( (M')_1\big)_0\Big)_0\!\in\!\{ 1,2\}\wedge (M')_0\! =\! q(0s)\wedge
{\cal S}\big(\alpha_0\Delta F(\alpha )\big)(M')\! =\! 1\Big\}$$ 
since $(M)_0\! =\!\big( B_\alpha (M)\big)_0$ and 
$\big( (M)_1\big)_0\! =\!\Big(\big( B_\alpha (M)\big)_1\Big)_0$. A second application of this shows that 
$\tilde\tau_s^{-1}(U_\varepsilon )$ is ccs. Thus $U_\varepsilon$ strongly ccs. Note that 
$$\tilde\tau_{1s\varepsilon}^{-1}(U_1\cup U_2)\! =\!
\Big\{\alpha\!\in\! 2^\omega\mid\exists M\!\in\!\omega\ 
\Big(\big( (M)_1\big)_0\Big)_0\!\in\!\{ 1,2\}\wedge 
(M)_0\! =\! q(01s\varepsilon )\wedge\alpha (M)\! =\! 1\Big\}.$$
We set\bigskip

\leftline{$O_\varepsilon\! :=\!\bigg\{\alpha\!\in\! 2^\omega\mid\exists M\!\in\!\omega\ 
\Big(\big( (M)_1\big)_0\Big)_0\!\in\!\{ 1,2\}\wedge (M)_0\! =\! q(01s\varepsilon )\wedge
\alpha (M)\! =\! 1~\wedge$}\smallskip

\rightline{$\forall l\! <\! M\ \ \bigg(\Big(\big( (l)_1\big)_0\Big)_0\!\notin\!\{ 1,2\}\vee 
(l)_0\!\notin\!\{ q(01s1),q(01s2)\}\vee\alpha (l)\! =\! 0\bigg)\bigg\}$}\bigskip

\noindent This defines a pair of $\boraone$ sets reducing 
$\big(\tilde\tau_{1s1}^{-1}(U_1\cup U_2),\tilde\tau_{1s2}^{-1}(U_1\cup U_2)\big)$. We check that they are ccs as for $\tilde\tau_s^{-1}(U_\varepsilon )$.\bigskip 

\noindent (b) The proof is completely similar to that of (a).$\hfill{\square}$\bigskip

 The following result is a consequence of Theorem 1.9 and Lemmas 1.11, 1.23 in [Lo1], and of Theorem 3 in [Lo-SR3]:

\begin{them} Let $\bf\Gamma$ be a self-dual Wadge class of Borel sets. Then there is a non self-dual Wadge class of Borel sets ${\bf\Gamma}'$ such that 
${\bf\Gamma}(\omega^\omega )\! =\!\Delta ({\bf\Gamma}')(\omega^\omega )$, ${\bf\Gamma}'$ does not have the separation property, and one of the following holds:\smallskip

\noindent (1) There is $\overline{u}\!\in\! {\cal D}$ such that
$${\bf\Gamma}'(\omega^\omega )\! =\!
\Big\{ (A_0\cap C_0)\cup (A_1\cap C_1)\mid A_0,\neg A_1\!\in\! {\bf\Gamma}_{\overline{u}}(\omega^\omega )\wedge C_0,C_1\!\in\!\boraone (\omega^\omega )\wedge C_0\cap C_1\! =\!\emptyset\Big\}.$$
(2) There is $\big( (u')_p\big)_{p\geq 1}\!\in\! {\cal D}^\omega$ such that 
$\big( {\bf\Gamma}_{(u')_p}(\omega^\omega )\big)_{p\geq 1}$ is strictly increasing and
$${\bf\Gamma}'(\omega^\omega )\! =\!
\Big\{\bigcup_{p\geq 1}~(A_p\cap C_p)\mid A_p\!\in\! {\bf\Gamma}_{(u')_p}(\omega^\omega )\wedge 
C_p\!\in\!\boraone (\omega^\omega )\wedge C_p\cap C_q\! =\!\emptyset\mbox{ if }p\!\neq\! q\Big\}.$$
\end{them}

\begin{lemm} Let ${\bf\Gamma}'$ be as in the statement of Theorem 5.2.8. Then there are 
$C^0,C^1\!\in\! {\bf\Gamma}'(2^\omega )$ disjoint, ccs, and not separable by a $\Delta ({\bf\Gamma}')$ set.\end{lemm}

\vfill\eject

\noindent\bf Proof.\rm\ (1) Lemma 5.2.7.(a) gives $(U_1,U_2)$ complete for pairs of disjoint $\boraone$ sets with $U_\varepsilon$ strongly ccs, and such that for each $s\!\in\! (\omega\!\setminus\!\{ 0\})^{<\omega}$ there is a pair $(O_1,O_2)$ of ccs $\boraone$ sets reducing the pair 
${\big(\tilde\tau_{1s1}^{-1}(U_1\cup U_2),\tilde\tau_{1s2}^{-1}(U_1\cup U_2)\big)}$. Theorem 5.2.4 gives 
$H_{\overline{u}}\!\subseteq\! 2^\omega$ which is ${\bf\Gamma}_{\overline{u}}$-complete and strongly ccs. We set $H\! :=\!\big(\tilde\tau_2^{-1}(H_{\overline{u}})\cap\tilde\tau_1^{-1}(U_1)\big)\cup 
\big(\tilde\tau_3^{-1}(\neg H_{\overline{u}})\!\setminus\!\tilde\tau_1^{-1}(U_2)\big)$ and, for 
$\varepsilon\!\in\!\{ 1,2\}$, $E_\varepsilon\! :=\!\tilde\tau_\varepsilon^{-1}(H)$. Finally, we set 
$C^\varepsilon\! :=\! 
(O_\varepsilon\cap E_\varepsilon )\cup (O_{3-\varepsilon}\!\setminus\! E_{3-\varepsilon})$.\bigskip

\noindent $\bullet$ We set, for $\varepsilon ,j\!\in\!\{ 1,2\}$, 
$A^\varepsilon_1\! :=\!\tilde\tau_{2\varepsilon}^{-1}(H_{\overline{u}})$, 
$A^\varepsilon_2\! :=\!\tilde\tau_{3\varepsilon}^{-1}(\neg H_{\overline{u}})$, 
$F^\varepsilon_j\! :=\!\tilde\tau_{1\varepsilon}^{-1}(U_j)$, so that 
$$E_\varepsilon\! =\! (A^\varepsilon_1\cap F^\varepsilon_1)\cup (A^\varepsilon_2\cap F^\varepsilon_2).$$ Note that 
$$\begin{array}{ll}
C^\varepsilon\! = 
& \!\!\!\!\! (A^\varepsilon_1\cap F^\varepsilon_1\cap O_\varepsilon )\cup 
(A^\varepsilon_2\cap F^\varepsilon_2\cap O_\varepsilon )\cup 
(\neg A^{3-\varepsilon}_1\cap F^{3-\varepsilon}_1\cap O_{3-\varepsilon})\cup 
(\neg A^{3-\varepsilon}_2\cap F^{3-\varepsilon}_2\cap O_{3-\varepsilon})\cr\cr
~~~~=
& \!\!\!\!\!\Big(\big( (A^\varepsilon_1\cap F^\varepsilon_1\cap O_\varepsilon )\cup 
(\neg A^{3-\varepsilon}_2\cap F^{3-\varepsilon}_2\cap O_{3-\varepsilon})\big)\cap
\big( (F^\varepsilon_1\cap O_\varepsilon )\cup (F^{3-\varepsilon}_2\cap O_{3-\varepsilon})\big)\Big)\cup\cr\cr
& \!\!\!\!\!\Big(\big( (A^\varepsilon_2\cap F^\varepsilon_2\cap O_\varepsilon )\cup 
(\neg A^{3-\varepsilon}_1\cap F^{3-\varepsilon}_1\cap O_{3-\varepsilon})\big)\cap
\big( (F^\varepsilon_2\cap O_\varepsilon )\cup (F^{3-\varepsilon}_1\cap O_{3-\varepsilon})\big)\Big)\mbox{,}
\end{array}$$
and that $F^\varepsilon_1\cap O_\varepsilon$, $F^{3-\varepsilon}_2\cap O_{3-\varepsilon}$, 
$F^\varepsilon_2\cap O_\varepsilon$, $F^{3-\varepsilon}_1\cap O_{3-\varepsilon}$ are pairwise disjoint open subsets of $2^\omega$. By Lemma 5.2.2 and the reduction property for $\boraone$ we can write 
$C^\varepsilon$ as the intersection of $2^\omega$ with 
$$\Big(\big( ({\cal A}^\varepsilon_1\cap {\cal O}^\varepsilon_1)\cup 
(\neg {\cal A}^{3-\varepsilon}_2\cap {\cal O}^{3-\varepsilon}_2)\big)\cap 
({\cal O}^\varepsilon_1\cup {\cal O}^{3-\varepsilon}_2)\Big)\cup
\Big(\big( ({\cal A}^\varepsilon_2\cap {\cal O}^\varepsilon_2)\cup 
(\neg {\cal A}^{3-\varepsilon}_1\cap {\cal O}^{3-\varepsilon}_1)\big)\cap
({\cal O}^\varepsilon_2\cup {\cal O}^{3-\varepsilon}_1)\Big)\mbox{,}$$
where ${\cal A}^\varepsilon_1,\neg {\cal A}^\varepsilon_2\!\in\! {\bf\Gamma}_{\overline{u}}(\omega^\omega )$ and ${\cal O}^\varepsilon_j$ are four pairwise disjoint open subsets of $\omega^\omega$. By Lemma 1.4.(b)  in [Lo1], $({\cal A}^\varepsilon_1\cap {\cal O}^\varepsilon_1)\cup 
(\neg {\cal A}^{3-\varepsilon}_2\cap {\cal O}^{3-\varepsilon}_2), 
\neg\big( ({\cal A}^\varepsilon_2\cap {\cal O}^\varepsilon_2)\cup 
(\neg {\cal A}^{3-\varepsilon}_1\cap {\cal O}^{3-\varepsilon}_1)\big)\!\in\! 
{\bf\Gamma}_{\overline{u}}(\omega^\omega )$, so that 
$C^\varepsilon\!\in\! {\bf\Gamma}'(2^\omega )$, by Lemma 5.2.2 again.\bigskip

\noindent $\bullet$ It is clear that $C^1$ and $C^2$ are disjoint and ccs.\bigskip

\noindent $\bullet$ Assume, towards a contradiction, that $D\!\in\!\Delta ({\bf\Gamma}')$ separates $C^1$ from $C^2$. Let $D_1,D_2\!\in\! {\bf\Gamma}'(\omega^\omega )$ disjoint. As $H$ is complete we get 
$f_\varepsilon\! :\!\omega^\omega\!\rightarrow\! 2^\omega$ continuous such that 
$D_\varepsilon\! =\! f_\varepsilon^{-1}(H)$. We define $f\! :\!\omega^\omega\!\rightarrow\! 2^\omega$ by
$$f(\alpha )\big(\big<\! <\!\varepsilon ,(k)_0\! >,(k)_1\big>\big)\! :=\!\left\{\!\!\!\!\!\!\!\!
\begin{array}{ll}
& f_\varepsilon (\alpha )(k)\mbox{ if }\varepsilon\!\in\!\{ 1,2\}\mbox{,}\cr
& 0\mbox{ otherwise,}
\end{array}
\right.$$
so that $\big( f(\alpha )\big)_\varepsilon\! =\! f_\varepsilon (\alpha )$. Then $f$ is continuous and 
$D_\varepsilon\! =\! f^{-1}(E_\varepsilon )$. Note that 
$E_\varepsilon\!\setminus\! E_{3-\varepsilon}\!\subseteq\! C^\varepsilon$. This implies that 
$\alpha\!\in\! D_1\Leftrightarrow f(\alpha )\!\in\! E_1\Leftrightarrow 
f(\alpha )\!\in\! E_1\!\setminus\! E_2\Rightarrow f(\alpha )\!\in\! C_1\!\subseteq\! D$. Similarly, 
$D_2\!\subseteq\! f^{-1}(\neg D)$, and $f^{-1}(D)\!\in\!\Delta ({\bf\Gamma}')(\omega^\omega )$ separates $D_1$ from $D_2$. Thus ${\bf\Gamma}'$ has the separation property, which is absurd.\bigskip

\noindent (2) Lemma 5.2.7.(b) gives $(U_p)_{p\geq 1}$ complete for sequences of pairwise disjoint 
$\boraone$ sets with $U_p$ strongly ccs, and such that for each 
$s\!\in\! (\omega\!\setminus\!\{ 0\})^{<\omega}$ there is a sequence 
$(O^\varepsilon_p)_{\varepsilon\in\{ 1,2\},p\geq 1}$ of ccs $\boraone$ sets reducing  
$\big(\tilde\tau_{s\varepsilon}^{-1}(U_p)\big)_{\varepsilon\in\{ 1,2\},p\geq 1}$. Theorem 5.2.4 gives 
$H_{(u')_p}\!\subseteq\! 2^\omega$ which is ${\bf\Gamma}_{(u')_p}$-complete and strongly ccs. We set 
$H\! :=\!\bigcup_{p\geq 1}~\big(\tilde\tau_{2p}^{-1}(H_{(u')_p})\cap\tilde\tau_1^{-1}(U_p)\big)$ and, for 
$\varepsilon\!\in\!\{ 1,2\}$, $E_\varepsilon\! :=\!\tilde\tau_\varepsilon^{-1}(H)$.\bigskip

 We also set $A^\varepsilon_p\! :=\!\tilde\tau_{(2p)\varepsilon}^{-1}(H_{(u')_p})$, 
$F^\varepsilon_p\! :=\!\tilde\tau_{1\varepsilon}^{-1}(U_p)$, so that 
$E_\varepsilon\! =\!\bigcup_{p\geq 1}~(A^\varepsilon_p\cap F^\varepsilon_p)$. Finally, we set 
$C^\varepsilon\! :=\! (A^\varepsilon_1\cap O^\varepsilon_1)\cup\bigcup_{p\geq 1}~
\big( (O^{3-\varepsilon}_p\!\setminus\! A^{3-\varepsilon}_p)\cup 
(A^\varepsilon_{p+1}\cap O^\varepsilon_{p+1})\big)$.\bigskip

 Note that $C^\varepsilon\!\in\! {\bf\Gamma}'(2^\omega )$ since 
$\big( {\bf\Gamma}_{(u')_p}(\omega^\omega )\big)_{p\geq 1}$ is strictly increasing, using again Lemma 5.2.2,    the generalized reduction property for $\boraone$ (see 22.16 in [K]), and Lemma 1.4.(b)  in [Lo1]. Here again, $E_\varepsilon\!\setminus\! E_{3-\varepsilon}\!\subseteq\! C^\varepsilon$ and we conclude as in (1).
$\hfill{\square}$\bigskip

\noindent\bf Proof of Theorem 1.7.\rm\ It is clear that Proposition 2.2, Lemmas 2.3, 2.6, Corollary 5.2.5,  Lemma 5.2.9 and Theorem 3.1 imply Theorem 1.7, if we set 
$\mathbb{S}^d_{\bf\Gamma}\! :=\! S^d_{C_{\bf\Gamma}}$ and 
$\mathbb{S}^\varepsilon_{\bf\Gamma}\! :=\! S^d_{C^\varepsilon}$.$\hfill{\square}$

\section{$\!\!\!\!\!\!$ The proof of Theorem 1.8}\indent 

 We first introduce an operator in the spirit of $\Phi$ defined before Theorem 4.2.2, but in dimension one. Another important difference to notice is the following. In Theorem 4.2.2, (f) for example, $S$ is in a boldface class, while $A_0$ and $A_1$ are in a lightface class. The same phenomenon will hold in the case of Wadge classes, and in the new operator we introduce we have boldface conditions (for example, we do not ask $\gamma'$ to be $\Borel (\beta )$). We code the Borel classes, and define an operator $\Phi_1$ on $\omega^\omega\!\times\!\omega^\omega$ to do it. Recall the definition of Seq before Lemma 2.3. We set\bigskip

\leftline{$W_0\! :=\!\Big\{ (\beta ,\gamma )\!\in\!\omega^\omega\!\times\! W^{\omega^\omega}\mid
\Big(\beta (0)\!\in\!\mbox{Seq}\wedge C_\gamma^{\omega^\omega}\! =\!
\big\{\delta\!\in\!\omega^\omega\mid {\cal I}^{-1}\big(\beta (0)\big)\!\subseteq\!\delta\big\}\Big)\vee$}\smallskip

\rightline{$\Big(\beta (0)\!\notin\!\mbox{Seq}\wedge C_\gamma^{\omega^\omega}\! =\!\emptyset\Big)\Big\}\mbox{,}$}\bigskip

\leftline{$\Phi_1(A)\! :=\! A\cup W_0\cup\Big\{ (\beta ,\gamma )\!\in\!\omega^\omega\!\times\! 
W^{\omega^\omega}\mid\exists\gamma'\!\in\!\omega^\omega\ \ \forall n\!\in\!\omega\ \ 
\big( (\beta )_n,(\gamma')_n\big)\!\in\! A\ \mbox{ and}$}\smallskip

\rightline{$\neg C_\gamma^{\omega^\omega}\! =\!
\bigcup_{n\in\omega}\ C_{(\gamma')_n}^{\omega^\omega}\Big\}.$}\bigskip

 In the sequel, we will denote $\Phi^{<\xi}_1\! :=\!\bigcup_{\eta <\xi}\ \Phi^\eta_1$.

\begin{lem} Let $1\!\leq\!\xi\! <\!\omega_1$ and $B\!\subseteq\!\omega^\omega$. Then $B\!\in\!\bormxi$ if and only if there is $(\beta ,\gamma )\!\in\!\Phi_1^\xi$ such that $C_\gamma^{\omega^\omega}\! =\! B$.
\end{lem}

\noindent\bf Proof.\rm\ Note first that $B\! =\! N_s\! :=\!\{\delta\!\in\!\omega^\omega\mid s\!\subseteq\!\delta\}$ for some $s\!\in\!\omega^{<\omega}$ or $B\!=\!\emptyset$ if and only if there is 
$(\beta ,\gamma )\!\in\! W_0\! =\!\Phi_1^0$ with $C_\gamma^{\omega^\omega}\! =\! B$. Then
$$\begin{array}{ll}
B\!\in\!\bormone\!\!\!\!
& \Leftrightarrow\exists (s_n)_{n\in\omega}\!\in\! (\omega^{<\omega})^\omega\ \ 
\neg B\! =\!\bigcup_{n\in\omega}\ N_{s_n}\vee\neg B\! =\!\emptyset\cr
& \Leftrightarrow\exists\beta ,\gamma'\!\in\!\omega^\omega\ \ \forall n\!\in\!\omega\ \ 
\big( (\beta )_n,(\gamma')_n\big)\!\in\!\Phi^0_1\wedge
\neg B\! =\!\bigcup_{n\in\omega}\ C_{(\gamma')_n}^{\omega^\omega}\cr
& \Leftrightarrow\exists(\beta ,\gamma )\!\in\!\Phi_1^1\ \ C_\gamma^{\omega^\omega}\! =\! B.
\end{array}$$
Assume now that the result is proved for $1\!\leq\!\eta\! <\!\xi\!\geq\! 2$. We get
$$\begin{array}{ll}
B\!\in\!\bormxi\!\!\!\!
& \Leftrightarrow\exists (B_n)_{n\in\omega}\!\in\! (\bormlxi )^\omega\ \ 
\neg B\! =\!\bigcup_{n\in\omega}\ B_n\cr
& \Leftrightarrow\exists\beta ,\gamma'\!\in\!\omega^\omega\ \ \forall n\!\in\!\omega\ \ 
\big( (\beta )_n,(\gamma')_n\big)\!\in\!\Phi^{<\xi}_1\wedge
\neg B\! =\!\bigcup_{n\in\omega}\ C_{(\gamma')_n}^{\omega^\omega}\cr
& \Leftrightarrow\exists(\beta ,\gamma )\!\in\!\Phi_1^\xi\ \ C_\gamma^{\omega^\omega}\! =\! B.
\end{array}$$
This finishes the proof.$\hfill{\square}$

\vfill\eject

 We now define a $\Ca$ coding of $\cal D$ (recall Definition 5.1.2).\bigskip
 
\noindent\bf Notation.\rm\ We define an inductive operator $\Lambda$ over $\omega^\omega$ as follows:\bigskip

\leftline{$\Lambda (D)\! :=\! D\cup\big\{\alpha\!\in\!\omega^\omega\mid\forall n\!\in\!\omega\ \ (\alpha )_n\!\in\! \mbox{WO}\wedge\vert (\alpha )_n\vert\! =\! 0\big\}\ \cup$}\smallskip

\rightline{$\big\{\alpha\!\in\!\omega^\omega\mid\forall n\!\in\!\omega\ \ (\alpha )_n\!\in\!\mbox{WO}
\wedge (\alpha )_0\! =\! (\alpha )_2\wedge\vert (\alpha )_1\vert\! =\! 1\wedge 
<(\alpha )_{2+j}>\in\! D\big\}\cup$}\smallskip

\rightline{$\big\{\alpha\!\in\!\omega^\omega\mid\forall n\!\in\!\omega\ \ (\alpha )_n\!\in\!\mbox{WO}
\wedge\vert (\alpha )_0\vert\!\geq\! 1\wedge\vert (\alpha )_1\vert\! =\! 2~\wedge$}
\smallskip

\rightline{$\forall p\!\in\!\omega\ <(\alpha )_{2+<p,q>}>\in\! D\wedge
\big(\vert (\alpha )_{2+<p,0>}\vert\!\geq\!\vert (\alpha )_0\vert\vee
\vert (\alpha )_{2+<p,0>}\vert\! =\! 0\big)\big\}.$}\bigskip
\noindent Then $\Lambda$ is a $\Ca$ monotone inductive operator, by 4A.2 in [M].\bigskip

 By 7C.1 in [M] we get $\Lambda^\infty\! :=\!\bigcup_{\xi}\ \Lambda^\xi\! =\!\Lambda (\Lambda^\infty )\! =\!
\bigcap\{ D\!\subseteq\!\omega^\omega\mid\Lambda (D)\!\subseteq\! D\}$. An easy induction on $\xi$ shows that $\Lambda^\infty\!\subseteq\!\big\{\alpha\!\in\!\omega^\omega\mid\forall n\!\in\!\omega\ \ 
(\alpha )_n\!\in\! \mbox{WO}\big\}$, so that the coding function $c$, partially defined by 
$c(\alpha )\! :=\!\big(\vert (\alpha )_n\vert\big)_{n\in\omega}$, is defined on $\Lambda^\infty$.
 
\begin{lem} The set $\Lambda^\infty$ is a $\Ca$ coding of $\cal D$, which means that 
$\Lambda^\infty\!\in\!\Ca (\omega^\omega )$ and $c[\Lambda^\infty ]\! =\! {\cal D}$.\end{lem} 
 
\noindent\bf Proof.\rm\ We first prove that $\Lambda^\infty\!\in\!\Ca (\omega^\omega )$ (see 7C in [M] for that). We define a set relation $\varphi (\alpha ,D)$ on $\omega^\omega$ by 
$\varphi (\alpha ,D)\Leftrightarrow\alpha\!\in\!\Lambda (D)$. As $\Lambda$ is monotone, $\varphi$ is operative. If $Q\!\in\!\Ca (Z\!\times\!\omega^\omega )$, then the relation 
$\varphi (\alpha ,\{\beta\!\in\!\omega^\omega\mid (z,\beta )\!\in\! Q\})$ is in $\Ca$. Thus $\varphi$ is $\Ca$ on $\Ca$. By 7C.8 in [M], $\varphi^\infty (\alpha )$ is in $\Ca$ and $\Lambda^\infty\!\in\!\Ca (\omega^\omega )$.\bigskip
 
 Let $\beta_\varepsilon\!\in\!\mbox{WO}$ such that $\vert\beta_\varepsilon\vert\! =\!\varepsilon$, for 
$\varepsilon\!\in\! 3$. Then $<\beta_0\mid n\!\in\!\omega >\in\!\Lambda^0\!\subseteq\!\Lambda^\infty$, so that $0^\infty\!\in\! c[\Lambda^\infty ]$. Let $u^*\!\in\! c[\Lambda^\infty ]$, $\alpha^*\!\in\!\Lambda^\infty$ with 
 $u^*\! =\! c(\alpha^*)$. Then 
$<\! (\alpha^*)_0,\beta_1,(\alpha^*)_0,(\alpha^*)_1,...\! >\in\!\Lambda(\Lambda^\infty )\! =\!\Lambda^\infty$, so that $u^*(0)1u^*\! =\! 
c\big(\! <\! (\alpha^*)_0,\beta_1,(\alpha^*)_0,(\alpha^*)_1,...\! >\!\big)\!\in\! c[\Lambda^\infty ]$.\bigskip

 Now let $\xi\!\geq\! 1$, $u_p\!\in\! c[\Lambda^\infty ]$ such that $u_p(0)\!\geq\!\xi$ or $u_p(0)\! =\! 0$, for each $p\!\in\!\omega$. Choose $\alpha\!\in\!\mbox{WO}$ with $\vert\alpha\vert\! =\!\xi$, and 
$\alpha^p\!\in\!\Lambda^\infty$ with $u_p\! =\! c(\alpha^p)$. Then 
$<\alpha ,\beta_2,(\alpha^{(0)_0})_{(0)_1},(\alpha^{(1)_0})_{(1)_1},...>\in\!\Lambda(\Lambda^\infty )\! =\!\Lambda^\infty$, so that 
$\xi 2<u_p> =\! c\big( <\alpha ,\beta_2,(\alpha^{(0)_0})_{(0)_1},(\alpha^{(1)_0})_{(1)_1},...>\big)
\!\in\! c[\Lambda^\infty ]$. Thus ${\cal D}\!\subseteq\! c[\Lambda^\infty ]$.\bigskip 
 
 Assume now that $\tilde {\cal D}\!\subseteq\!\omega_1^\omega$ satisfies the following properties:\smallskip
 
\noindent (a) $0^\infty\!\in\!\tilde {\cal D}$.\smallskip
 
\noindent (b) $u^*\!\in\!\tilde {\cal D}\Rightarrow u^*(0)1u^*\!\in\!\tilde {\cal D}$.\smallskip
 
\noindent (c) $\Big(\xi\!\geq\! 1\wedge\forall p\!\in\!\omega$ $\big( u_p\!\in\!\tilde {\cal D}\wedge 
(u_p(0)\!\geq\!\xi\vee u_p(0)\! =\! 0)\big)\Big)\Rightarrow\xi 2<u_p>\in\!\tilde {\cal D}$.\bigskip

 We set $D\! :=\!\{\alpha\!\in\!\omega^\omega\mid\forall n\!\in\!\omega\ \ (\alpha )_n\!\in\!\mbox{WO}
\wedge c(\alpha )\!\in\!\tilde {\cal D}\}$. It remains to see that $\Lambda (D)\!\subseteq\! D$. Indeed, this will imply that $\Lambda^\infty\!\subseteq\! D$, $c[\Lambda^\infty ]\!\subseteq\! c[D]\!\subseteq\!\tilde {\cal D}$ and $c[\Lambda^\infty ]\!\subseteq\! {\cal D}$.\bigskip

 As $0^\infty\!\in\!\tilde {\cal D}$ we get $\big\{\alpha\!\in\!\omega^\omega\mid\forall n\!\in\!\omega\ \ 
(\alpha )_n\!\in\!\mbox{WO}\wedge\vert (\alpha )_n\vert\! =\! 0\big\}\!\subseteq\! D$. Assume that 
$(\alpha )_n\!\in\!\mbox{WO}$ for each $n\!\in\!\omega$, that $(\alpha )_0\! =\! (\alpha )_2$, 
$\vert (\alpha )_1\vert\! =\! 1$ and $<(\alpha )_{2+j}>\in\! D$. Then 
$u^*\! :=\!\big(\vert (\alpha )_{2+j}\vert\big)\!\in\!\tilde {\cal D}$, and 
$\vert (\alpha )_2\vert 1u^*\!\in\!\tilde {\cal D}$. Thus $c(\alpha )\!\in\!\tilde {\cal D}$ and 
$\alpha\!\in\! D$.\bigskip

 Assume now that $(\alpha )_n\!\in\!\mbox{WO}$ for each $n\!\in\!\omega$,  
$\vert (\alpha )_0\vert\!\geq\! 1$, $\vert (\alpha )_1\vert\! =\! 2$,  
$<\! (\alpha )_{2+<p,q>}\! >\in\! D$, and 
$\vert (\alpha )_{2+<p,0>}\vert\!\geq\!\vert (\alpha )_0\vert$ or $\vert (\alpha )_{2+<p,0>}\vert\! =\! 0$ for each 
$p\!\in\!\omega$. We set $\xi\! :=\!\vert (\alpha )_0\vert$. Then we have 
$u_p\! :=\!\big(\vert (\alpha )_{2+<p,q>}\vert\big)\!\in\!\tilde {\cal D}$, and 
$\xi 2<u_p>\in\!\tilde {\cal D}$. Thus $c(\alpha )\!\in\!\tilde {\cal D}$ and $\alpha\!\in\! D$.$\hfill{\square}$

\vfill\eject 

 Note that just like Definition 5.1.2, the definition of $\Lambda$ is cut into three cases, that we will meet again later on: $\vert (\alpha )_1\vert\! =\! 0$ (or, equivalently, $\vert (\alpha )_n\vert\! =\! 0$ for each integer $n$), $\vert (\alpha )_1\vert\! =\! 1$ or $\vert (\alpha )_1\vert\! =\! 2$.\bigskip

 Even if ``$u\!\in\! {\cal D}$" is the least relation satisfying some conditions, some simplifications are possible. For example, ${\bf\Gamma}_{01010^\infty}\! =\! {\bf\Gamma}_{0^\infty}$. Some other simplifications are possible, and some of them will simplify the notation later on. This will lead to the notion of a normalized code of a description. To define it, we need to associate a tree to a code of a description. The idea is to describe the construction of a set in ${\bf\Gamma}_u$ using simpler and simpler sets, until we get the simplest set, namely the empty set. More specifically, we define 
${\mathfrak T}\! :\!\Lambda^\infty\!\rightarrow\!\{\mbox{trees on}\ \omega\!\times\!\Lambda^\infty\}$ as follows. Let $\alpha\!\in\!\Lambda^\xi\!\setminus\!\Lambda^{<\xi}$. We set
$${\mathfrak T}(\alpha )\! :=\!\left\{\!\!\!\!\!\!\!\!
\begin{array}{ll}
& \{\emptyset\}\cup\{ <\! (0,\alpha )\! >\}\mbox{ if }\vert (\alpha )_1\vert\! =\! 0
\mbox{,}\cr\cr
& \{\emptyset\}\cup\big\{ (0,\alpha )^\frown s\mid s\!\in\! 
{\mathfrak T}(<\! (\alpha )_{2+j}\! >)\big\}\mbox{ if }\vert (\alpha )_1\vert\! =\! 1
\mbox{,}\cr\cr
& \{\emptyset\}\cup\big\{ (0,\alpha )^\frown s\mid s\!\in\! 
{\mathfrak T}(<\! (\alpha )_{2+<0,q>}\! >)\big\}\cup\cr
& \ \ \ \ \ \ \ \ \ \ \ \bigcup_{p\geq 1}\ \big\{ (p,\alpha )^\frown s\mid s\!\in\! 
{\mathfrak T}(<\! (\alpha )_{2+<(p)_0+1,q>}\! >)\big\}\mbox{ if }\vert (\alpha )_1\vert\! =\! 2.
\end{array}
\right.$$
An easy induction on $\eta$ shows that ${\mathfrak T}(\alpha )$ is always a countable well-founded tree (the first coordinate of $(p,\alpha )$ ensures the well-foundedness). A sequence 
$s\!\in\! {\mathfrak T}(\alpha )$ is said to be $maximal$ if $s\!\subseteq\! t\!\in\! {\mathfrak T}(\alpha )$ implies that $s\! =\! t$. Note that $\big\vert\big( s_1(\vert s\vert\! -\! 1)\big)_1\big\vert\! =\! 0$ if $s$ is maximal. We denote by ${\cal M}_\alpha$ the set of maximal sequences in ${\mathfrak T}(\alpha )$.
 
\begin{defi} We say that $\alpha\!\in\!\Lambda^\infty$ is $normalized$ if the following holds:
$$\big( s\!\in\! {\cal M}_\alpha\wedge i\! <\!\vert s\vert\wedge
\big\vert\big( s_1(i)\big)_1\big\vert\! =\! 1\big)\Rightarrow i\! =\!\vert s\vert\! -\! 2.$$
\end{defi}

 This means that in a maximal sequence $s$ of ${\mathfrak T}(\alpha )$, 
$\big\vert\big( s_1(i)\big)_1\big\vert$ is $2$, then possibly $1$ once, and finally $0$ once. The next lemma says that we can always assume that $\alpha$ is normalized. It is based on the fact that 
$\check {S_\xi}({\bf\Gamma},{\bf\Gamma}')\! =\! S_\xi (\check {\bf\Gamma},\check {\bf\Gamma}')$. 

\begin{lem} Let $\alpha\!\in\!\Lambda^\infty$. Then there is $\alpha'\!\in\!\Lambda^\infty$ normalized with 
$(\alpha')_0\! =\! (\alpha )_0$ and ${\bf\Gamma}_{c(\alpha')}\! =\! {\bf\Gamma}_{c(\alpha)}$.\end{lem}

\noindent\bf Proof.\rm\ Assume that $\alpha\!\in\!\Lambda^\xi\!\setminus\!\Lambda^{<\xi}$. We argue by induction on $\xi$.\medskip

\noindent\bf Case 1.\rm\ $\vert (\alpha )_1\vert\! =\! 0$.\medskip

 We just set $\alpha'\! :=\!\alpha$ since $\big\vert\big( s_1(i)\big)_1\big\vert\! =\! 0$.\medskip

\noindent\bf Case 2.\rm\ $\vert (\alpha )_1\vert\! =\! 1$.\medskip

\noindent $\bullet$ We first define $N\! :\!\Lambda^\infty\!\rightarrow\!\Lambda^\infty$ as follows. We ensure that 
$\big( N(\beta )\big)_0\! =\! (\beta )_0$ and 
${\bf\Gamma}_{c(N(\beta ))}\! =\!\check {\bf\Gamma}_{c(\beta )}$. Let $\beta_1\!\in\!\mbox{WO}$ with 
$\vert\beta_1\vert\! =\! 1$. We set
$$N(\beta )\! :=\!\left\{\!\!\!\!\!\!
\begin{array}{ll}
& <\! (\beta )_0,\beta_1,(\beta )_0,(\beta )_1,(\beta )_2,...\! >\mbox{ if }\vert (\beta )_1\vert\! =\! 0
\mbox{,}\cr
& <\! (\beta )_{2+j}\! >\mbox{ if }\vert (\beta )_1\vert\! =\! 1\mbox{,}\cr
& <\! (\beta )_0,(\beta )_1,
\bigg(\Big( N\big( <(\beta )_{2+<(i-2)_0,q>}>\big)\Big)_{(i-2)_1}\bigg)_{i\geq 2}\! >\mbox{ if }
\vert (\beta )_1\vert\! =\! 2\mbox{,}
\end{array}
\right.$$
and one easily checks that $N$ is defined and suitable.\bigskip

\noindent $\bullet$ As $<(\alpha )_{2+j}>\in\!\Lambda^{<\xi}$, the induction assumption gives 
$\alpha''\!\in\!\Lambda^\infty$ normalized satisfying the equalities $(\alpha'')_0\! =\! (\alpha )_2\! =\! (\alpha )_0$ and ${\bf\Gamma}_{c(\alpha'')}\! =\!{\bf\Gamma}_{c(<(\alpha )_{2+j}>)}$. In particular, 
$${\bf\Gamma}_{c(\alpha )}\! =\!\check {\bf\Gamma}_{c(<(\alpha )_{2+j}>)}\! =\!
\check {\bf\Gamma}_{c(\alpha'')}\! =\! {\bf\Gamma}_{c(N(\alpha''))}.$$ 
So we have to find $\alpha'\!\in\!\Lambda^\infty$ normalized with $(\alpha')_0\! =\! (\alpha'')_0$ and 
${\bf\Gamma}_{c(\alpha')}\! =\! {\bf\Gamma}_{c(N(\alpha''))}$. Assume that 
$\alpha''\!\in\!\Lambda^\eta\!\setminus\!\Lambda^{<\eta}$. We argue by induction on $\eta$.\bigskip

\noindent\bf Subcase 1.\rm\ $\vert (\alpha'')_1\vert\!\leq\! 1$.\medskip

 We just set $\alpha'\! :=\! N(\alpha'')$.\medskip

\noindent\bf Subcase 2.\rm\ $\vert (\alpha'')_1\vert\! =\! 2$.\bigskip

 Note that $<\! (\alpha'')_{2+<p,q>}\! >$ is normalized since 
$(0,\alpha'')^\frown s\!\in\! {\cal M}_{\alpha''}$ (resp., $(p,\alpha'')^\frown s\!\in\! {\cal M}_{\alpha''}$) if 
$s\!\in\! {\cal M}_{(\alpha'')_{2+<0,q>}}$ (resp., $s\!\in\! {\cal M}_{(\alpha'')_{2+<(p)_0+1,q>}}$ and 
$p\!\geq\! 1$). The induction assumption gives $<\! (\alpha')_{2+<p,q>}\! >\in\!\Lambda^\infty$ normalized with $(\alpha')_{2+<p,0>}\! =\! (\alpha'')_{2+<p,0>}$ and 
$${\bf\Gamma}_{c(<(\alpha')_{2+<p,q>}>)}\! =\! {\bf\Gamma}_{c(N(<(\alpha'')_{2+<p,q>}>))}.$$ 
We set $(\alpha')_i\! :=\! (\alpha'')_i$ if $i\!\in\! 2$ and we are done.\bigskip

\noindent\bf Case 3.\rm\ $\vert (\alpha )_1\vert\! =\! 2$.\bigskip
 
 The induction assumption gives $<\! (\alpha')_{2+<p,q>}\! >\in\!\Lambda^\infty$ normalized satisfying the equalities $(\alpha')_{2+<p,0>}\! =\! (\alpha )_{2+<p,0>}$ and 
${\bf\Gamma}_{c(<(\alpha')_{2+<p,q>}>)}\! =\!{\bf\Gamma}_{c(<(\alpha )_{2+<p,q>}>)}$. We set 
$(\alpha')_i\! :=\! (\alpha )_i$ if $i\!\in\! 2$ and we are done.$\hfill{\square}$\bigskip

 Using $\Phi_1$, we will now code the non self-dual Wadge classes of Borel sets, and define an operator 
$\Upsilon_1$ on $(\omega^\omega )^3$ to do it. We set\bigskip

\leftline{$\Upsilon_1 (A)\! :=\! A\cup
\bigg\{ (\alpha ,\beta ,\gamma )\!\in\! (\omega^\omega )^2\!\times\! W^{\omega^\omega}\mid
\forall n\!\in\!\omega\ \ (\alpha )_n\!\in\!\mbox{WO}\ \wedge$}\bigskip

\rightline{$\bigg(\forall n\!\in\!\omega\ \ \vert (\alpha )_n\vert\! =\! 0\wedge\beta (0)\! =\! 0
\wedge C^{\omega^\omega}_\gamma\! =\!\emptyset\bigg)\ \vee$}\bigskip

\rightline{$\bigg(\vert (\alpha )_1\vert\! =\! 1\wedge (\alpha )_0\! =\! (\alpha )_2\wedge
\beta (0)\! =\! 1\ \wedge$}\medskip

\rightline{$\exists\gamma'\!\in\!\omega^\omega\ \ 
( <(\alpha )_{2+j}>,\beta^*,\gamma')\!\in\! A\wedge 
C^{\omega^\omega}_\gamma\! =\!\neg C^{\omega^\omega}_{\gamma'}\bigg)\ \vee$}\bigskip

\rightline{$\bigg(\vert (\alpha )_1\vert\! =\! 2\wedge\vert (\alpha )_0\vert\!\geq\! 1\wedge
\forall p\!\in\!\omega\ \ \big(\vert (\alpha )_{2+<p,0>}\vert\!\geq\!\vert (\alpha )_0\vert\vee
\vert (\alpha )_{2+<p,0>}\vert\! =\! 0\big)\ \wedge$}\medskip

\rightline{$\beta (0)\! =\! 2\wedge\exists\gamma'\!\in\!\omega^\omega~~
( <(\alpha )_{2+<0,q>}>,(\beta^*)_0,(\gamma')_0)\!\in\! A\ \wedge$}\medskip

\rightline{$\forall p\!\geq\! 1\ \ 
\big(\! <\! (\alpha )_{2+<(p)_0+1,q>}\! >,\big( (\beta^*)_p\big)_0,\big( (\gamma' )_p\big)_0\big)
\!\in\! A\wedge\big(\big( (\beta^*)_p\big)_1,\big( (\gamma' )_p\big)_1\big)
\!\in\!\Phi_1^{\vert (\alpha )_0\vert}\ \wedge$}\medskip

\rightline{$\forall p\!\not=\! q\!\geq\! 1\ \ C^{\omega^\omega}_{((\gamma')_p)_1}\cup 
C^{\omega^\omega}_{((\gamma')_q)_1}\! =\!\omega^\omega\ \wedge$}\medskip

\rightline{$C^{\omega^\omega}_\gamma\! =\!\bigcup_{p\geq 1}\ 
(C^{\omega^\omega}_{((\gamma')_p)_0}\!\setminus\! C^{\omega^\omega}_{((\gamma')_p)_1})\cup 
(C^{\omega^\omega}_{(\gamma')_0}\cap\bigcap_{p\geq 1}\ 
C^{\omega^\omega}_{((\gamma')_p)_1})\bigg)\bigg\}.$}

\vfill\eject

\begin{lem} Let $\xi$ be an ordinal.\smallskip

\noindent (a) Assume that $(\alpha ,\beta ,\gamma )\!\in\!\Upsilon_1^\xi$. Then $\alpha\!\in\!\Lambda^\xi$.\smallskip

\noindent (b) Let $\alpha\!\in\!\Lambda^\xi$ and $B\!\subseteq\!\omega^\omega$. Then 
$B\!\in\! {\bf\Gamma}_{c(\alpha )}$ if and only if there are $\beta ,\gamma\!\in\!\omega^\omega$ such that 
$(\alpha ,\beta ,\gamma )\!\in\!\Upsilon_1^\xi$ and $C_\gamma^{\omega^\omega}\! =\! B$.
\end{lem}

\noindent\bf Proof.\rm\ (a) We argue by induction on $\xi$. So let 
$\alpha\!\in\!\Upsilon_1^\xi\!\setminus\!\Upsilon_1^{<\xi}$. We may assume that 
$\vert (\alpha )_1\vert\!\geq\! 1$. If $\vert (\alpha )_1\vert\! =\! 1$, then 
$( <\! (\alpha )_{2+j}\! >,\beta^*,\gamma')\!\in\!\Upsilon_1^{<\xi}$ for some $\gamma'$ and 
$<\! (\alpha )_{2+j}\! >\in\!\Lambda^{<\xi}$ by induction assumption, so we are done. If 
$\vert (\alpha )_1\vert\! =\! 2$, then 
$$( <(\alpha )_{2+<0,q>}>,(\beta^*)_0,(\gamma')_0),
\big( <(\alpha )_{2+<(p)_0+1,q>}>,\big( (\beta^*)_p\big)_0,\big( (\gamma' )_p\big)_0\big)\!\in\!
\Upsilon_1^{<\xi}$$
for some $\gamma'$ and $<(\alpha )_{2+<p,q>}>\in\!\Lambda^{<\xi}$ by induction assumption for each integer $p$.\bigskip

\noindent (b) $\Rightarrow$ We argue by induction on $\xi$, and we may assume that 
$\alpha\!\notin\!\Lambda^{<\xi}$.\bigskip

\noindent\bf Case 1.\rm\ $\vert (\alpha )_1\vert\! =\! 0$.\bigskip

 Note that $c(\alpha )\! =\! 0^\infty$ and $B\! =\!\emptyset$. We set $\beta\! :=\! 0^\infty$, and we choose 
$\gamma\!\in\! W^{\omega^\omega}$ with $C_\gamma\! =\!\emptyset$. Then 
$(\alpha ,\beta ,\gamma )\!\in\!\Upsilon_1^0\!\subseteq\!\Upsilon_1^\xi$.\bigskip

\noindent\bf Case 2.\rm\ $\vert (\alpha )_1\vert\! =\! 1$.\bigskip
 
 Note that $<(\alpha )_{2+j}>\in\!\Lambda^{<\xi}$, and $\neg B\!\in\! {\bf\Gamma}_{c(<(\alpha )_{2+j}>)}$. By induction assumption we get $\beta',\gamma'\!\in\!\omega^\omega$ such that 
$(<(\alpha )_{2+j}>,\beta',\gamma')\!\in\!\Upsilon_1^{<\xi}$ and 
$C_{\gamma'}^{\omega^\omega}\! =\!\neg B$. We set $\beta\! :=\! 1\beta'$ and we choose 
$\gamma\!\in\! W^{\omega^\omega}$ with 
$C^{\omega^\omega}_\gamma\! =\!\neg C^{\omega^\omega}_{\gamma'}$.\bigskip

\noindent\bf Case 3.\rm\ $\vert (\alpha )_1\vert\! =\! 2$.\bigskip

 Note that $<(\alpha )_{2+<p,q>}>\in\!\Lambda^{<\xi}$ for each integer $p$. We can write 
$$B\! =\!\bigcup_{p\geq 1}\ (A_p\cap C_p)\cup (B'\!\setminus\!\bigcup_{p\geq 1}\ C_p)\mbox{,}$$
where $(C_p)_{p\geq 1}$ is a sequence of pairwise disjoint ${\bf\Sigma}^0_{\vert (\alpha )_0\vert}$ sets, 
$B'\!\in\! {\bf\Gamma}_{c(<(\alpha )_{2+<0,q>}>)}$ and 
$$A_p\!\in\! {\bf\Gamma}_{c(<(\alpha )_{2+<(p)_0+1,q>}>)}.$$
Lemma 6.1 gives 
$\big(\big( (\beta^*)_p\big)_1,\big( (\gamma')_p\big)_1\big)\!\in\!\Phi_1^{\vert (\alpha )_0\vert}$ such that 
$C_{((\gamma')_p)_1}^{\omega^\omega}\! =\!\neg C_p$. The induction assumption gives 
$(\beta^*)_0,(\gamma')_0\!\in\!\omega^\omega$ such that 
$( <\! (\alpha )_{2+<0,q>}\! >,(\beta^*)_0,(\gamma')_0)\!\in\!\Upsilon_1^{<\xi}$ and 
$C_{(\gamma')_0}^{\omega^\omega}\! =\! B'$, and 
$\big( (\beta^*)_p\big)_0,\big( (\gamma')_p\big)_0\!\in\!\omega^\omega$ such that 
$\big( <\! (\alpha )_{2+<(p)_0+1,q>}\! >,\big( (\beta^*)_p\big)_0,\big( (\gamma')_p\big)_0\big)\!\in\!
\Upsilon_1^{<\xi}$ and 
$C_{((\gamma')_p)_0}^{\omega^\omega}\! =\! A_p$. We set $\beta (0)\! :=\! 2$ and we choose 
$\gamma\!\in\! W^{\omega^\omega}$ with 
$$C^{\omega^\omega}_\gamma\! =\!\bigcup_{p\geq 1}\ 
(C^{\omega^\omega}_{((\gamma')_p)_0}\!\setminus\! C^{\omega^\omega}_{((\gamma')_p)_1})\cup 
(C^{\omega^\omega}_{(\gamma')_0}\cap\bigcap_{p\geq 1}\ 
C^{\omega^\omega}_{((\gamma')_p)_1}).$$
$\Leftarrow$ We argue by induction on $\xi$, and we may assume that 
$(\alpha ,\beta ,\gamma )\!\notin\!\Upsilon_1^{<\xi}$.\bigskip

\noindent\bf Case 1.\rm\ $\vert (\alpha )_1\vert\! =\! 0$.\bigskip

 Note that $B\! =\! C^{\omega^\omega}_\gamma\! =\!
\emptyset\!\in\! {\bf\Gamma}_{0^\infty}\! =\! {\bf\Gamma}_{c(\alpha )}$.\bigskip

\noindent\bf Case 2.\rm\ $\vert (\alpha )_1\vert\! =\! 1$.\bigskip

 Note that there is $\gamma'$ such that $( <\! (\alpha )_{2+j}\! >,\beta^*,\gamma')\!\in\!\Upsilon_1^{<\xi}$ and $C^{\omega^\omega}_\gamma\! =\!\neg C^{\omega^\omega}_{\gamma'}$, which implies that  
$B\!\in\!\check {\bf\Gamma}_{c(<(\alpha )_{2+j}>)}\! =\! {\bf\Gamma}_{c(\alpha )}$.\bigskip

\noindent\bf Case 3.\rm\ $\vert (\alpha )_1\vert\! =\! 2$.\bigskip

 We get $\gamma'$ since $(\alpha ,\beta ,\gamma )\!\in\!\Upsilon_1^\xi$. As 
$$( <(\alpha )_{2+<0,q>}>,(\beta^*)_0,(\gamma')_0),
\big( <(\alpha )_{2+<(p)_0+1,q>}>,\big( (\beta^*)_p\big)_0,\big( (\gamma')_p\big)_0\big)\!\in\!
\Upsilon_1^{<\xi}$$
we get $C^{\omega^\omega}_{(\gamma')_0}\!\in\! {\bf\Gamma}_{c(<(\alpha )_{2+<0,q>}>)}$ and 
$C^{\omega^\omega}_{((\gamma')_p)_0}\!\in\! {\bf\Gamma}_{c(<(\alpha )_{2+<(p)_0+1,q>}>)}$, by induction assumption. As 
$\big(\big( (\beta^*)_p\big)_1,\big( (\gamma' )_p\big)_1\big)\!\in\!\Phi_1^{\vert (\alpha )_0\vert}$, we get 
$C^{\omega^\omega}_{((\gamma')_p)_1}\!\in\! {\bf\Pi}^0_{\vert (\alpha )_0\vert}$ by Lemma 6.1. This implies that 
$$B\!\in\! S_{\vert (\alpha )_0\vert}
(\bigcup_{p\geq 1}\ {\bf\Gamma}_{c(<(\alpha )_{2+<p,q>}>)}, {\bf\Gamma}_{c(<(\alpha )_{2+<0,q>}>)})
\! =\! {\bf\Gamma}_{c(\alpha )}.$$
This finishes the proof.$\hfill{\square}$\bigskip

\noindent\bf Remark.\rm\ We will also consider the operator $\Upsilon$ defined just like $\Upsilon_1$, except that\smallskip

\noindent - We replace $(W^{\omega^\omega},C^{\omega^\omega})$ with $(W,C)$ (we work in 
$(\omega^\omega )^d$ instead of $\omega^\omega$).\smallskip

\noindent - We replace the condition of the form 
$(\tilde\beta ,\tilde\gamma )\!\in\!\Phi_1^{\vert (\alpha )_0\vert}$ with 
$\big( (\alpha )_0,\tilde\beta ,\tilde\gamma\big)\!\in\! Q$ (see the remark at the end of Section 4 for the definition of $Q$).\smallskip

\noindent - We ask $\beta ,\gamma ,\gamma'$ to be $\Borel (\alpha )$, so that $\Upsilon$ is a $\Ca$ monotone inductive operator.\bigskip

 To prove Theorem 1.8, we will consider some tuples 
$\vec v\! :=\! (\alpha ,a_0,a_1,\underline{a}_0,\underline{a}_1,r)$, where 
$\alpha\!\in\!\Lambda^\infty$. We will inductively define them through an inductive operator over 
$(\omega^\omega )^6$ called $\Theta$. The definition of $\Theta$ is in the spirit of that of 
$\Upsilon_1$, and is cut into three cases, depending on the value of $\vert (\alpha )_1\vert$. As the definition of $\Theta$ is long and technical, we give first some more informal explanations about its meaning. We will have $\vec v\!\in\!\Theta^\infty$. So there is an ordinal $\xi$ such that 
$\vec v\!\in\!\Theta^\xi$.\bigskip
 
\noindent - $\alpha\!\in\!\Lambda^\xi$ is a (normalized in practice) code for a description 
$u\! =\! c(\alpha )$.\bigskip
 
\noindent - $a_0, a_1\!\in\!\Borel (\alpha )$ are codes for a pair of disjoint analytic subsets of 
$(\omega^\omega )^d$. Using the good universal set $\cal U$ for $\Ca$ defined after the proof of Theorem 4.2.2, at the end of Section 4, we will actually code the complement of these analytic sets, so that we will set $A_i\! :=\!\neg {\cal U}_{a_i}$ for $i\!\in\! 2$.\bigskip

\noindent - Similarly, $\underline{a}_0,\underline{a}_1\!\in\!\Borel (\alpha )$ are codes for a pair of disjoint analytic subsets of $(\omega^\omega )^d$. In fact, we will have 
$\underline{A}_i\! :=\!\neg {\cal U}_{\underline{a}_i}\!\subseteq\! A_i$. These codes will be used to build $r$, and $\underline{a}_0,\underline{a}_1,r$ will be completely determined by $(\alpha ,a_0,a_1)$. So one should think that $\underline{a}_i\! =\!\underline{a}_i(\alpha ,a_0,a_1)\!\simeq\!\underline{a}_i(u,a_0,a_1)$, $r\! =\! r(\alpha ,a_0,a_1)\!\simeq\! r(u,a_0,a_1)$. We need the following lemma to specify their meaning.

\begin{lem} There is a recursive map $f_a\! :\! (\omega^\omega )^2\!\rightarrow\!\omega^\omega$ such that ${\cal U}_{f_a(\alpha ,r)}\! =\! {\cal U}_{(r)_0}\cup\bigcup_{p\geq 1}\ \neg
\overline{\neg {\cal U}_{(r)_p}}^{\tau_{\vert\alpha\vert}}$ if $\alpha\!\in\!\Borel\cap\mbox{WO}$ and 
$\vert\alpha\vert\!\geq\! 1$.\end{lem}

\noindent\bf Proof.\rm\ Note first that  
$P\! :=\!\{ (\beta ,\vec\delta )\!\in\!\omega^\omega\!\times\! (\omega^\omega )^d\mid 
(\beta )_0\!\in\!\Borel\cap\mbox{WO}\wedge\vert (\beta )_0\vert\!\geq\! 1\ \wedge$\medskip

\rightline{$\vec\delta\!\in\! {\cal U}_{((\beta )_1)_0}\cup
\bigcup_{p\geq 1}\ \neg\overline{\neg {\cal U}_{((\beta )_1)_p}}^{\tau_{\vert (\beta )_0\vert}}\}$}\bigskip

\noindent is a $\Ca$ set, by the remark at the end of Section 4 defining $R$. This gives 
$\gamma\!\in\!\omega^\omega$ recursive with 
$P\! =\! {\cal U}^{\omega^\omega\times (\omega^\omega )^d}_\gamma$. Let 
$\alpha\!\in\!\Borel\cap\mbox{WO}$ with $\vert\alpha\vert\!\geq\! 1$, and $r\!\in\!\omega^\omega$. We have
$$\begin{array}{ll}
\vec\delta\!\in\! {\cal U}_{(r)_0}\cup\bigcup_{p\geq 1}\ \neg
\overline{\neg {\cal U}_{(r)_p}}^{\tau_{\vert\alpha\vert}}\!\!\!\!
& \Leftrightarrow (<\alpha ,r,r,...>,\vec\delta\ )\!\in\! P\cr
& \Leftrightarrow (\gamma ,<\alpha ,r,r,...>,\vec\delta\ )\!\in\! 
{\cal U}^{\omega^\omega\times (\omega^\omega )^d}\cr
& \Leftrightarrow\big( S(\gamma ,<\alpha ,r,r,...>),\vec\delta\ \big)\!\in\! {\cal U}
\end{array}$$
We just have to set $f_a(\alpha ,r)\! :=\! S(\gamma ,<\alpha ,r,r,...>)$.$\hfill{\square}$\bigskip

 The following will hold:\smallskip

$\circ$ If $u\! =\! 0^\infty$ or $u\! =\!\xi 1u^*$, then 
$\underline{a}_i\! =\!\underline{a}_i(\alpha ,a_0,a_1)\! =\!\underline{a}_i(u,a_0,a_1)\! =\! a_i$.\smallskip

$\circ$ If $u\! =\!\xi 2<u_p>$, then there will be $a'_0,a'_1,r'\!\in\!\Borel (\alpha )$ such that, for each 
$p\!\geq\! 1$,
$$\big( <(\alpha )_{2+<(p)_0+1,q>}>,a_0,a_1,(a'_0)_p,(a'_1)_p,(r')_p\big)\!\in\!\Theta^{<\xi}.$$
We will have $\underline{a}_i\! =\!\underline{a}_i(u,a_0,a_1)\! =\! 
f_a\big( (\alpha )_0,<a_i,(r')_1,(r')_2,...>\big)\mbox{,}$ 
and $(r')_p\! =\! r(u_{(p)_0+1},a_0,a_1)$ if $p\!\geq\! 1$. In particular, 
$\underline{A}_i\! =\! A_i\cap
\bigcap_{p\geq 1}\ \overline{\neg {\cal U}_{r(u_{(p)_0+1},a_0,a_1)}}^{\tau_\xi}$.\bigskip

\noindent - $r\!\in\!\Borel (\alpha )$ is a code for an analytic subset of $(\omega^\omega )^d$ playing the role that $\overline{A_0}^{\tau_\xi}\cap A_1$ played in Theorem 4.2.2. In other words, the emptyness of this analytic set is equivalent to the possibility of separating $A_0$ from $A_1$ by a $\mbox{pot}({\bf\Gamma}_u)$ set. Here again, using $\cal U$, we will actually code the complement of this analytic set: $\neg {\cal U}_r$ is an analytic subset of $(\omega^\omega )^d$. In particular,\smallskip

$\circ$ If $u\! =\! 0^\infty$, then $r\! =\! r(\alpha ,a_0,a_1)\! =\! r(u,a_0,a_1)\! =\! a_1$.\smallskip

$\circ$ If $u\! =\!\xi 1u^*$, then $r\! =\! r(\alpha ,a_0,a_1)\! =\! r(u,a_0,a_1)\! =\! a_0$.\smallskip

$\circ$ If $u\! =\!\xi 2<u_p>$, then we there will be $a''_0,a''_1\!\in\!\Borel (\alpha )$ such that 
$$( <(\alpha )_{2+<0,q>}>,\underline{a}_0,\underline{a}_1,a''_0,a''_1,r)\!\in\!\Theta^{<\xi}.$$
In particular, $r(u,a_0,a_1)\! =\! r(u_0,\underline{a}_0,\underline{a}_1)\! =\! 
r\big( u_0,\underline{a}_0(u,a_0,a_1),\underline{a}_1(u,a_0,a_1)\big)$. We are now ready to define $\Theta$ (recall the remark at the end of Section 4 defining $Q$).

\vfill\eject

 The operator $\Theta$ is defined as follows (recall the definition of $\Lambda$):\smallskip

\leftline{$\Theta (A)\! :=\! A\cup
\bigg\{ (\alpha ,a_0,a_1,\underline{a}_0,\underline{a}_1,r)\!\in\!
\big(\omega^\omega\cap\Borel (\alpha )\big)^6\mid
\forall n\!\in\!\omega\ \ (\alpha )_n\!\in\!\mbox{WO}\ \wedge$}\smallskip

\rightline{$\bigg(\forall n\!\in\!\omega\ \ \vert (\alpha )_n\vert\! =\! 0\wedge 
{\cal U}_{a_0}\cup {\cal U}_{a_1}\! =\! (\omega^\omega )^d\wedge 
(\underline{a}_0,\underline{a}_1)\! =\! (a_0,a_1)\wedge r\! =\! a_1\bigg)\ \vee$}\smallskip

\rightline{$\bigg(\vert (\alpha )_1\vert\! =\! 1\wedge (\alpha )_0\! =\! (\alpha )_2
\wedge ( <(\alpha )_{2+j}>,a_0,a_1,\underline{a}_0,\underline{a}_1,a_1)\!\in\! A\wedge 
r\! =\! a_0\bigg)\ \vee$}\smallskip

\rightline{$\bigg(\vert (\alpha )_1\vert\! =\! 2\wedge\vert (\alpha )_0\vert\!\geq\! 1\wedge
\forall p\!\in\!\omega\ \ \big(\vert (\alpha )_{2+<p,0>}\vert\!\geq\!\vert (\alpha )_0\vert\vee
\vert (\alpha )_{2+<p,0>}\vert\! =\! 0\big)\ \wedge$}\smallskip

\rightline{$\exists a'_0,a'_1,r'\!\in\!\Borel (\alpha )\ \ 
\big( <(\alpha )_{2+<0,q>}>,a_0,a_1,(a'_0)_0,(a'_1)_0,(r')_0\big)\!\in\! A\ \wedge$}\smallskip

\rightline{$\forall p\!\geq\! 1\ \ 
\big( <(\alpha )_{2+<(p)_0+1,q>}>,a_0,a_1,(a'_0)_p,(a'_1)_p,(r')_p\big)\!\in\! A\ \wedge$}\smallskip

\rightline{$\forall i\!\in\! 2\ \ \underline{a}_i\! =\! f_a\big( (\alpha )_0,<a_i,(r')_1,(r')_2,...>\!\big)\ 
\wedge$}\smallskip

\rightline{$\exists a''_0,a''_1\!\in\!\Borel (\alpha)\ \ 
( <(\alpha )_{2+<0,q>}>,\underline{a}_0,\underline{a}_1,a''_0,a''_1,r)\!\in\! A\bigg)\bigg\}.$}\smallskip

\noindent Then $\Theta$ is a $\Ca$ monotone inductive operator.\bigskip

\noindent\bf Remark.\rm\ Let $\xi$ be an ordinal, and 
$\vec v\! :=\! (\alpha ,a_0,a_1,\underline{a}_0,\underline{a}_1,r)\!\in\!\Theta^\xi$. Then an induction on $\xi$ shows the following properties:\smallskip

\noindent - $\neg {\cal U}_{a_0}\cap\neg {\cal U}_{a_1}\! =\!\emptyset$.\smallskip

\noindent - $\neg {\cal U}_{\underline{a}_i}\!\subseteq\!\neg {\cal U}_{a_i}$ for each $i\!\in\! 2$. In particular, $\neg {\cal U}_{\underline{a}_0}\cap\neg {\cal U}_{\underline{a}_1}\! =\!\emptyset$.\smallskip

\noindent - $\underline{a}_0,\underline{a}_1,r$ are completely determined by $(\alpha ,a_0,a_1)$.
\smallskip

\noindent - If $\neg {\cal U}_{a_i}\!\subseteq\!\neg {\cal U}_{b_i}$ for each $i\!\in\! 2$, then 
$\neg {\cal U}_{\underline{a}_i}\!\subseteq\!\neg {\cal U}_{\underline{b}_i}$ for each $i\!\in\! 2$ and 
$\neg {\cal U}_{r(\alpha ,a_0,a_1)}\!\subseteq\!\neg {\cal U}_{r(\alpha ,b_0,b_1)}$.\smallskip

\noindent - There is $i\!\in\! 2$ such that $\neg {\cal U}_{r}\!\subseteq\!\neg {\cal U}_{a_i}$.

\begin{lem} (a) Let $\xi$ be an ordinal, $\alpha\!\in\!\Borel$, and 
$(\alpha ,\beta ,\gamma )\!\in\!\Upsilon^\xi$. Then 
$\alpha\!\in\!\Lambda^\xi$ and the set $C_\gamma$ is in $\Borel\cap {\bf\Gamma}_{c(\alpha )}(\tau_1)$.\smallskip

\noindent (b) Let $\alpha\!\in\!\Borel\cap\Lambda^\infty$ normalized, $a_0,a_1\!\in\!\Borel$ with 
$A_0\cap A_1\! =\!\emptyset$. Then there are  
$\underline{a}_0,\underline{a}_1,r\!\in\!\omega^\omega$ such that 
$(\alpha ,a_0,a_1,\underline{a}_0,\underline{a}_1,r)\!\in\!\Theta^\infty$.\end{lem}

\noindent\bf Proof.\rm\ (a) We argue as in the proof of Lemmas 6.5.(a) and 6.5.(b)$\Leftarrow$. The only thing to notice is that in the case $\vert (\alpha )_1\vert\! =\! 2$, 
$\big( (\alpha )_0,\big( (\beta^*)_p\big)_1,\big( (\gamma' )_p\big)_1\big)\!\in\! Q$. Proposition 2.2, Lemma 2.3 and Theorem 3.1 give a tree $T_d$ with $\Borel$ suitable levels and 
$S\!\in\! {\bf\Sigma}^0_{\vert (\alpha )_0\vert}(\lceil T_d\rceil )$ not separable from 
$\lceil T_d\rceil\!\setminus\! S$ by a $\mbox{pot}({\bf\Pi}^0_{\vert (\alpha )_0\vert})$ set. As 
$\alpha\!\in\!\Borel$, $\vert (\alpha)_0\vert\! <\!\omega_1^{\mbox{CK}}$ and Theorem 4.2.2 implies that  $C_{((\gamma')_p)_1}$ is in ${\bf\Pi}^0_{\vert (\alpha )_0\vert}(\tau_1)$. Thus 
$C_\gamma\!\in\! {\bf\Gamma}_{c(\alpha )}(\tau_1)$.\bigskip

\noindent (b) Let $\xi$ be an ordinal with $\alpha\!\in\!\Lambda^\xi$. Here again we argue by induction on 
$\xi$. So assume that $\alpha\!\notin\!\Lambda^{<\xi}$.\bigskip

\noindent\bf Case 1.\rm\ $\vert (\alpha )_1\vert\! =\! 0$.\bigskip

 Let $\underline{a}_i\! :=\! a_i$ and $r\! :=\! a_1$. Then 
$(\alpha ,a_0,a_1,\underline{a}_0,\underline{a}_1,r)\!\in\!\Theta^0\!\subseteq\!\Theta^\infty$.

\vfill\eject

\noindent\bf Case 2.\rm\ $\vert (\alpha )_1\vert\! =\! 1$.\bigskip

 As $<(\alpha )_{2+j}>\in\!\Lambda^{<\xi}$ we get, by induction assumption, $(\underline{a}_0,\underline{a}_1,r')$ with 
$$(<(\alpha )_{2+j}>,a_0,a_1,\underline{a}_0,\underline{a}_1,r')\!\in\!\Theta^\infty .$$
As $\alpha$ is normalized we get $\vert (\alpha )_{2+j}\vert\! =\! 0$ for each $j$, and $r'\! =\! a_1$. We set $r\! :=\! a_0$. Then 
$$(\alpha ,a_0,a_1,\underline{a}_0,\underline{a}_1,r)\!\in\!\Theta (\Theta^\infty )\! =\!\Theta^\infty.$$
\bf Case 3.\rm\ $\vert (\alpha )_1\vert\! =\! 2$.\bigskip

 As $<(\alpha )_{2+<p,q>}>\in\!\Lambda^{<\xi}$ we get, by induction assumption, $(a^p_0,a^p_1,r'_p)$ with
$$\big( <(\alpha )_{2+<0,q>}>,a_0,a_1,a_0^0,a_1^0,r'_0\big)\!\in\!\Theta^\infty\mbox{,}$$ 
and $(<(\alpha )_{2+<(p)_0+1,q>}>,a_0,a_1,a^p_0,a^p_1,r'_p)\!\in\!\Theta^\infty$, for each $p\!\geq\! 1$. As in the proof of Lemma 6.2 we see that $\Theta^\infty\!\in\!\Ca$. By $\Borel$-selection, we may assume that the sequences $(a^p_0)$, $(a^p_1)$ and $(r'_p)$ are $\Borel$. In particular, there is $a'_i\!\in\!\Borel$ with $(a'_i)_p\! =\! a^p_i$. We set $(r')_p\! :=\! r'_p$, and 
$$\underline{a}_i\! :=\! f_a\big( (\alpha )_0,<a_i,(r')_1,(r')_2,...>\!\big).$$ 
The induction assumption gives $a''_0,a''_1,r$ such that 
$( <(\alpha )_{2+<0,q>}>,\underline{a}_0,\underline{a}_1,a''_0,a''_1,r)\!\in\!\Theta^\infty$. We are done since $(\alpha ,a_0,a_1,\underline{a}_0,\underline{a}_1,r)\!\in\!\Theta^\infty$.\hfill{$\square$}\bigskip

 The next lemma is the crucial separation lemma announced in the presentation of $r$.
 
\begin{lem} Let 
$\vec v\! :=\! (\alpha ,a_0,a_1,\underline{a}_0,\underline{a}_1,r)\!\in\!\Theta^\infty$ with 
$\alpha\!\in\!\Borel$ normalized and $a_0,a_1\!\in\!\Borel$, 
$\Sigma$ in $\Ana\big( (\omega^\omega )^d\big)$ with $(\neg {\cal U}_r)\cap\Sigma\! =\!\emptyset$. Then there are $\beta',\gamma'\!\in\!\omega^\omega$ such that 
$(\alpha ,\beta',\gamma')\!\in\!\Upsilon^\infty$ and $C_{\gamma'}$ separates $A_1\cap\Sigma$ from $A_0\cap\Sigma$. In particular, $A_1\cap\Sigma$ is separable from $A_0\cap\Sigma$ by a 
$\Borel\cap {\bf\Gamma}_{c(\alpha )}(\tau_1)$ set.\end{lem}

\noindent\bf Proof.\rm\  The last assertion comes from Lemma 6.7.(a). Let $\eta$ be an ordinal with 
$\vec v\!\in\!\Theta^\eta$. We argue by induction on $\eta$. So assume that 
$\vec v\!\in\!\Theta^\eta\!\setminus\!\Theta^{<\eta}$.\bigskip

\noindent\bf Case 1.\rm\ $\vert (\alpha )_1\vert\! =\! 0$.\bigskip

 We set $\beta'\! :=\! 0^\infty$, and choose $\gamma'\!\in\!\Borel\cap W$ with $C_{\gamma'}\! =\!\emptyset$. We are done since $\emptyset\! =\! A_1\cap\Sigma$.\bigskip

\noindent\bf Case 2.\rm\ $\vert (\alpha )_1\vert\! =\! 1$.\bigskip

 As $\alpha$ is normalized, we get $\vert (\alpha )_{2+j}\vert\! =\! 0$ for each $j$. We set 
$\beta'\! :=\! 10^\infty$, and choose $\gamma'\!\in\!\Borel\cap W$ with 
$C_{\gamma'}\! =\! (\omega^\omega )^d$. Then $\gamma''\!\in\!\Borel\cap W$ with 
$C_{\gamma''}\! =\!\emptyset$ is a witness for the fact that $(\alpha ,\beta',\gamma')\!\in\!\Upsilon^\infty$. We are done since $r\! =\! a_0$.\bigskip

\noindent\bf Case 3.\rm\ $\vert (\alpha )_1\vert\! =\! 2$.\bigskip

 There are $a'_0,a'_1,r'\!\in\!\Borel$ with 
$\big( <(\alpha )_{2+<(p)_0+1,q>}>,a_0,a_1,(a'_0)_p,(a'_1)_p,(r')_p\big)\!\in\!\Theta^{<\eta}$, for each 
$p\!\geq\! 1$, and, for each $i\!\in\! 2$, 
$\underline{a}_i\! =\! f_a\big( (\alpha )_0,<a_i,(r')_1,(r')_2,...>\!\big)$. Moreover, there are 
$a''_0,a''_1\!\in\!\Borel$ with 
$(<(\alpha )_{2+<0,q>}>,\underline{a}_0,\underline{a}_1,a''_0,a''_1,r)\!\in\!\Theta^{<\eta}$.\bigskip

 By Lemma 6.7.(a), one of the goals is to build $C_{\gamma'}\!\in\! {\bf\Gamma}_{c(\alpha )}(\tau_1)$. The proof of Lemma 6.7.(a) shows that ${\bf\Gamma}_{c(\alpha )}\! =\! S_{\vert (\alpha )_0\vert}
(\bigcup_{p\geq 1}\ {\bf\Gamma}_{c(<(\alpha )_{2+<p,q>}>)}, {\bf\Gamma}_{c(<(\alpha )_{2+<0,q>}>)})$. This means that we want to find sequences $(C_p)_{p\geq 1}$, $(S_p)_{p\geq 1}$ and $B$ such that 
$C_{\gamma'}\! =\!\bigcup_{p\geq 1}\ (S_p\cap C_p)\cup (B\!\setminus\!\bigcup_{p\geq 1}\ C_p)$.

\vfill\eject

\noindent - Let us construct $B$.\bigskip

 The induction assumption gives $\beta''',\gamma'''\!\in\!\omega^\omega$ such that 
$(<(\alpha )_{2+<0,q>}>,\beta''',\gamma''')\!\in\!\Upsilon^\infty$ and $C_{\gamma'''}$ separates 
$\underline{A}_1\cap\Sigma$ from $\underline{A}_0\cap\Sigma$. We set $B\! :=\! C_{\gamma'''}$.\bigskip

\noindent - Let us construct the $C_p$'s.\bigskip

 We set $\xi\! :=\!\vert (\alpha )_0\vert$. Note that $\underline{A}_i\! =\!
A_i\cap\bigcap_{p\geq 1}\ \overline{\neg {\cal U}_{(r')_p}}^{\tau_\xi}$. This implies that  
$$U\! :=\!\big( C_{\gamma'''}\cap A_0\cap\Sigma\big)\cup
\big(\neg C_{\gamma'''}\cap A_1\cap\Sigma\big)
\!\subseteq\!\bigcup_{p\geq 1}\ \neg\overline{\neg {\cal U}_{(r')_p}}^{\tau_\xi}.$$ 
As in the proof of Lemma 6.6 we see that the relation 
``$\vec\delta\!\notin\!\overline{\neg {\cal U}_{(r')_p}}^{\tau_{\vert (\alpha )_0\vert}}$" is $\Ca$ in 
$(p,\alpha ,r',\vec\delta )$. By $\Borel$-selection there is a $\Borel$-recursive map 
$f\! :\! (\omega^\omega )^d\!\rightarrow\!\omega$ such that $f(\vec\delta )\!\geq\! 1$ for each 
${\vec\delta\!\in\! (\omega^\omega )^d}$ and 
$\vec\delta\!\notin\!\overline{\neg {\cal U}_{(r')_{f(\vec\delta )}}}^{\tau_\xi}$ for each $\vec\delta\!\in\! U$.\bigskip

 In particular, for each $\vec\delta\!\in\! U$ there is $P\!\in\!\Ana\cap\bormlxi (\tau_1)$ such that 
$\vec\delta\!\in\! P\!\subseteq\! {\cal U}_{(r')_{f(\vec\delta )}}$. Now $P$ and 
$\neg {\cal U}_{(r')_{f(\vec\delta )}}$ are disjoint $\Ana$ sets, and separable by a 
$\bormlxi (\tau_1)$ set. As $\alpha\!\in\!\Borel$ we get 
$1\!\leq\!\vert (\alpha )_0\vert\! <\!\omega_1^{\mbox{CK}}$. As in the proof of Lemma 6.7.(a) we get $T_d$ and $S$. By Theorem 4.2.2 we get $(\beta ,\gamma )\!\in\! (\Borel\!\times\!\Borel )\cap V_{<\xi}$ with 
$P\!\subseteq\! C_\gamma\!\subseteq\! {\cal U}_{(r')_{f(\vec\delta )}}$.\bigskip

 By Lemma 4.2.3.(2).(a) the relation ``$(\beta ,\gamma )$ is in $(\Borel\!\times\!\Borel )\cap V_{<\xi}$" is 
$\Ca$, so there is a $\Borel$-recursive map $g\! :\! (\omega^\omega )^d\!\rightarrow\!\omega\!\times\! (\omega^\omega\!\times\!\omega^\omega )$ such that
$$\forall\vec\delta\!\in\! U\ \ g_0(\vec\delta )\! =\! f(\vec\delta )\ \mbox{ and }\ 
g_1(\vec\delta )\!\in\! (\Borel\!\times\!\Borel )\cap V_{<\xi}\ \mbox{ and }\ 
\vec\delta\!\in\! C_{(g_1(\vec\delta ))_1}\!\subseteq\! {\cal U}_{(r')_{f(\vec\delta )}}\mbox{,}$$
by $\Borel$-selection. In particular, the $\Ana$ set $g[U]$ is a subset of 
$$\big\{\big( p,(\beta ,\gamma )\big)\!\in\!\omega\!\times\!\big( (\Borel\!\times\!\Borel )\cap V_{<\xi}\big)\mid 
C_\gamma\!\subseteq\! {\cal U}_{(r')_p}\big\}\mbox{,}$$
which is $\Ca$ and countable. The separation theorem gives $D\!\in\!\Borel$ between these two sets. As $D$ is countable, there are $N,\tilde\beta ,\tilde\gamma\!\in\!\Borel$ with 
$D\! =\!\Big\{\Big( N(q),\big( (\tilde\beta )_q,(\tilde\gamma )_q\big)\Big)\mid q\!\in\!\omega\Big\}$. Now we can define $C_p\! :=\!\bigcup_{q\in\omega ,N(q)=p}\ C_{(\tilde\gamma )_q}\!\setminus\! (\bigcup_{q'<q}\ C_{(\tilde\gamma )_{q'}})$.\bigskip

\noindent - We now study the properties of the $C_p$'s. We can say that\smallskip

$\circ$ The relation ``$\vec\delta\!\in\! C_p$" is $\Borel$ in $(p,\vec\delta )$.\smallskip

$\circ$ The $C_p$'s are pairwise disjoint.\smallskip

$\circ$ $C_p\!\in\!\boraxi (\tau_1)$ since 
$C_{(\tilde\gamma )_q}\!\in\!\bormlxi (\tau_1)\!\subseteq\!\borxi (\tau_1)$, by Theorem 4.2.2.\smallskip

$\circ$ We set $\tilde C\! :=\!\{ (p,\vec\delta )\!\in\!\omega\!\times\! (\omega^\omega )^d\mid
\exists q\!\in\!\omega\ \ N(q)\! =\! p\ \mbox{ and }\ \vec\delta\!\in\! C_{(\tilde\gamma )_q}\}\mbox{,}$ so that 
$\tilde C\!\in\!\Borel$ and $\tilde C_p\!\in\!\boraone (\tau_\xi )$ for each $p\!\geq\! 1$. We have 
$C_p\!\subseteq\!\tilde C_p$.\smallskip

$\circ$ $\bigcup_{p\geq 1}\ C_p\! =\!\bigcup_{p\geq 1}\ \tilde C_p$.\smallskip

$\circ$ $\tilde C_p$ separates $U\cap f^{-1}(\{ p\})$ from $\neg {\cal U}_{(r')_p}$. In particular, $U$ is a subset of the $\Borel$ set $\bigcup_{p\geq 1}\ C_p$. Moreover, 
$\bigcap_{p\geq 1}\ \overline{\neg {\cal U}_{(r')_p}}^{\tau_\xi}\!\subseteq\!
\neg (\bigcup_{p\geq 1}\ \tilde C_p)$.

\vfill\eject

\noindent - The induction assumption gives, for each $p\!\geq\! 1$, $\beta^p,\gamma^p$ with  
$(<(\alpha )_{2+<(p)_0+1,q>}>,\beta^p,\gamma^p)\!\in\!\Upsilon^\infty$ and $C_{\gamma^p}$ separates $A_1\cap\tilde C_p$ from $A_0\cap\tilde C_p$. As in the proof of Lemma 6.7.(b) we may assume that the sequences $(\beta^p)$ and $(\gamma^p)$ are $\Borel$. By $\Borel$-selection again there is a $\Borel$-recursive map $h\! :\!\omega\!\rightarrow\!\omega^\omega\!\times\!\omega^\omega$ such that 
$h(p)\!\in\! (\Borel\!\times\!\Borel )\cap V_\xi$ and $C_{h_1(p)}\! =\!\neg C_p$ for each $p\!\geq\! 1$. We set $\big( ( {\beta'}^*)_p\big)_1\! :=\! h_0(p)$ and 
$\big( (\overline{\gamma})_p\big)_1\! :=\! h_1(p)$, so that 
$\Big( (\alpha )_0,\big( ( {\beta'}^*)_p\big)_1,\big( (\overline{\gamma})_p\big)_1\Big)\!\in\! Q$ for each 
$p\!\geq\! 1$.\bigskip

 We set $\beta'(0)\! :=\! 2$, $( {\beta'}^*)_0\! :=\!\beta'''$, and $\big( ( {\beta'}^*)_p\big)_0\! :=\!\beta^p$ if 
$p\!\geq\! 1$, so that $\beta'$ is completely defined. Similarly, we set 
$(\overline{\gamma})_0\! :=\!\gamma'''$, and $\big( (\overline{\gamma})_p\big)_0\! :=\!\gamma^p$ if 
$p\!\geq\! 1$. Finally, we choose $\gamma'\!\in\!\Borel\cap W$ with $C_{\gamma'}\! =\!
\bigcup_{p\geq 1}\ (C_{\gamma^p}\!\setminus\! C_{h_1(p)})\cup (C_{(\overline{\gamma})_0}\cap\bigcap_{p\geq 1}\ 
C_{h_1(p)})$, so that $(\alpha ,\beta',\gamma')\!\in\!\Upsilon^\infty$ and 
$C_{\gamma'}$ separates $A_1\cap\Sigma$ from $A_0\cap\Sigma$.$\hfill{\square}$\bigskip

 The next result is the actual (effective) content of Theorem 1.8.(1). It is also the version of Theorem 4.4.1 for the non self-dual Wadge classes of Borel sets. Let 
$j_d\! :\! (d^\omega )^d\!\rightarrow\!\omega^\omega$ be a continuous embedding (for example we can embed $(d^\omega )^d$ into $(\omega^\omega )^d$ in the obvious way, and then use a bijection between $(\omega^\omega )^d$ and $\omega^\omega$).

\begin{thm} Let $T_d$ be a tree with $\Borel$ suitable levels, $\alpha$ in $\Borel$ normalized, 
$\beta ,\gamma$ in $\omega^\omega$ such that $(\alpha,\beta ,\gamma )\!\in\! \Upsilon_1^\infty$, 
$S\! :=\! j_d^{-1}(C^{\omega^\omega}_{\gamma})\cap\lceil T_d\rceil$, and 
$a_0,a_1,\underline{a}_0,\underline{a}_1,r\!\in\!\omega^\omega$ with 
$\vec v\! :=\! (\alpha ,a_0,a_1,\underline{a}_0,\underline{a}_1,r)\!\in\!\Theta^\infty$. Then one of the following holds:\smallskip

\noindent (a) $\neg {\cal U}_r\! =\!\emptyset$.\smallskip

\noindent (b) The inequality $\big( (\Pi_i''\lceil T_d\rceil )_{i\in d}, S,\lceil T_d\rceil\!\setminus\! S\big)\leq
\big( (\omega^\omega )_{i\in d}, A_0, A_1\big)$ holds.\end{thm}

 Now we can state the version of Theorem 4.2.2 for the non self-dual Wadge classes of Borel sets. 
 
\begin{thm} Let $T_d$ be a tree with $\Borel$ suitable levels, $\alpha$ in $\Borel$ normalized, 
$\beta ,\gamma$ in $\omega^\omega$ such that $(\alpha,\beta ,\gamma )\!\in\!\Upsilon_1^\infty$, 
$S\! :=\! j_d^{-1}(C^{\omega^\omega}_{\gamma})\cap\lceil T_d\rceil$, and 
$a_0,a_1,\underline{a}_0,\underline{a}_1,r\!\in\!\omega^\omega$ with 
$\vec v\! :=\! (\alpha ,a_0,a_1,\underline{a}_0,\underline{a}_1,r)\!\in\!\Theta^\infty$. We assume that $S$ is not separable from $\lceil T_d\rceil \!\setminus\! S$ by a 
$\mbox{pot}(\check {\bf\Gamma}_{c(\alpha )})$ set. Then the following are equivalent:\smallskip

\noindent (a) The set $A_0$ is not separable from $A_1$ by a 
$\hbox{\it pot}(\check {\bf\Gamma}_{c(\alpha )})$ set.\smallskip

\noindent (b) The set $A_0$ is not separable from $A_1$ by a 
$\Borel\cap\hbox{\it pot}(\check {\bf\Gamma}_{c(\alpha )})$ set.\smallskip

\noindent (c) $\neg\big(\exists\beta',\gamma'\!\in\!\omega^\omega$ such that 
$(\alpha ,\beta',\gamma')\!\in\!\Upsilon^\infty$ and 
$A_1\!\subseteq\! C_{\gamma'}\!\subseteq\!\neg A_0\big)$.\smallskip

\noindent (d) The set $A_0$ is not separable from $A_1$ by a 
$\check {\bf\Gamma}_{c(\alpha )}(\tau_1)$ set.\smallskip

\noindent (e) $\neg {\cal U}_r\!\not=\!\emptyset$.\smallskip

\noindent (f) The inequality $\big( (d^\omega )_{i\in d}, S,\lceil T_d\rceil\!\setminus\! S\big)\leq
\big( (\omega^\omega )_{i\in d}, A_0, A_1\big)$ holds.\end{thm}

\noindent\bf Proof.\rm\ (a) $\Rightarrow$ (b) and (a) $\Rightarrow$ (d) are clear since 
${\it\Delta}_{\omega^\omega}$ is Polish.\bigskip

\noindent (b) $\Rightarrow$ (c) This comes from Lemma 6.7.(a).\bigskip

\noindent (b) $\Rightarrow$ (e), (c) $\Rightarrow$ (e) and (d) $\Rightarrow$ (e) This comes from Lemma 6.8.\bigskip

\noindent (e) $\Rightarrow$ (f) This comes from Theorem 6.9 (as $\Pi_i''\lceil T_d\rceil$ is compact, we just have to compose with continuous retractions to get functions defined on $d^\omega$).\bigskip

\noindent (f) $\Rightarrow$ (a) If $P\!\in\!\hbox{\rm pot}(\check {\bf\Gamma}_{c(\alpha )})$ separates $A_0$ from $A_1$ and (f) holds, then 
$S\!\subseteq\! (\Pi_{i\in d}\ f_i)^{-1}(P)\!\subseteq\!\neg (\lceil T_d\rceil \!\setminus\! S)$.  
This implies that $S$ is separable from $\lceil T_d\rceil \!\setminus\! S$ by a 
$\mbox{pot}(\check {\bf\Gamma}_{c(\alpha )})$ set, by Lemma 4.4.7. But this contradicts the assumption on $S$.\hfill{$\square$}

\vfill\eject

\noindent\bf Proof of Theorem 1.8.(1).\rm\ Note first that (a) and (b) cannot hold simultaneously, as in the proof of Theorem 6.10.\bigskip

 We assume that (a) does not hold. This implies that the $X_i$'s are not empty, since otherwise 
$A_0\! =\! A_1\! =\!\emptyset$, and $\emptyset\!\in\!\check {\bf\Gamma}$ unless 
${\bf\Gamma}\! =\!\{\emptyset\}$. As in the proof of Theorem 4.1, we may assume that 
$X_i\! =\!\omega^\omega$ for each $i\!\in\! d$, by Lemma 4.4.7. By Theorem 5.1.3 there is 
$u\!\in\! {\cal D}$ with ${\bf\Gamma}(\omega^\omega )\! =\! {\bf\Gamma}_u(\omega^\omega )$. If $E$ is a $0$-dimensional Polish space, then we also have ${\bf\Gamma}(E)\! =\! {\bf\Gamma}_u(E)$, by Theorem 4.1.3 in [Lo-SR2]. It follows that $\mbox{pot}({\bf\Gamma})\! =\!\mbox{pot}({\bf\Gamma}_u)$. By Lemmas 6.2 and 6.4 we may assume that there is $\alpha\!\in\!\Lambda^\infty$ normalized with 
$c(\alpha )\! =\! u$.\bigskip

 By Theorem 4.1.3 in [Lo-SR2] there is 
$B\!\in\! {\bf\Gamma}(\omega^\omega )$ with $S\! =\! j_d^{-1}(B)\cap\lceil T_d\rceil$. To simplify the notation, we may assume that $T_d$ has $\Borel$ levels, $\alpha\!\in\!\Borel$, and 
$A_0, A_1\!\in\!\Ana\big( (\omega^\omega )^d\big)$. By Lemma 6.5 there are  
$\beta ,\gamma\!\in\!\omega^\omega$ such that $(\alpha ,\beta ,\gamma )\!\in\!\Upsilon_1^\infty$ and 
$C^{\omega^\omega}_{\gamma }\! =\! B$. Lemma 6.7.(b) gives $\underline{a}_0,\underline{a}_1,r$ with 
$(\alpha ,a_0,a_1,\underline{a}_0,\underline{a}_1,r)\!\in\!\Theta^\infty$. Lemma 6.8  implies that 
$\neg {\cal U}_r\!\not=\!\emptyset$. So (b) holds, by Theorem 6.10.\hfill{$\square$}\bigskip

 The sequel is devoted to the proof of Theorem 6.9. We have to introduce a lot of objects before we can do it. We will create some paragraphs to describe these objects. We start with a general notion. The idea is that, given a set $S$ in ${\bf\Gamma}_{c(\alpha )}(\lceil T_d\rceil)$, and with the help of the tree $\mathfrak{T}(\alpha )$, we will keep in mind all the $\boraxi$ (or equivalently $\bormxi$, passing to complements) used to build $S$. We will represent these $\bormxi$ sets, on most sequences $s$ of 
$\mathfrak{T}(\alpha )$, by induction on $\vert s\vert$, applying the Debs-Saint Raymond theorem. At each induction step, we make closed some $\bormxi$ sets of this level, but we also partially simplify the $\bormxi$ sets to come. This is why the ordinal substraction is involved (recall the definition of ordinal substraction after Theorem 5.1.3).

\begin{defi} Let $X$ be a set, $A\!\subseteq\! X$, $\cal B$ a countable family of subsets of $X$, and $\bf\Gamma$ a Borel class. We say that $A\!\in\! {\bf\Gamma}({\cal B})$ if $A\!\in\! {\bf\Gamma}(X,\tau )$ for any topology $\tau$ on $X$ containing ${\cal B}$.\end{defi}

\begin{prop} Let $X$ be a topological space.\smallskip

\noindent (a) Let $A\!\subseteq\! X$, $\cal B$ a countable family of open subsets of $X$, and $\bf\Gamma$ a Borel class. Then $A\!\in\! {\bf\Gamma}(X)$ if $A\!\in\! {\bf\Gamma}({\cal B})$.\smallskip

\noindent (b) Let $Y$ be a set, $B\!\subseteq\! Y$, $f\! :\! X\!\rightarrow\! Y$ a bijection, 
$\cal B$ a countable family of subsets of $Y$, and $\bf\Gamma$ a Borel class. Then 
$f^{-1}(B)\!\in\! {\bf\Gamma}(\{ f^{-1}(D)\mid D\!\in\! {\cal B}\})$ if $B\!\in\! {\bf\Gamma}({\cal B})$.\smallskip

\noindent (c) Let $1\!\leq\!\eta\!\leq\!\xi$ and $A\!\in\! {\bf\Pi}^0_{\xi}(X)$. We assume that $X$ is  metrizable. Then there is ${\cal B}\!\subseteq\! {\bf\Pi}^0_{\eta}(X)$ countable such that 
$A\!\in\! {\bf\Pi}^0_{1+(\xi -\eta )}(\check {\cal B})$, where 
$\check {\cal B}\! :=\!\{\neg B\mid B\!\in\! {\cal B}\}$.\end{prop}

 In practice, $X$ will be the metrizable space $[R]$ for some tree relation $R$, and $f$ will be the canonical map given by the Debs-Saint Raymond theorem.\bigskip

\noindent\bf Proof.\rm\ (a) The topology $\tau$ is simply the topology of $X$.\bigskip

\noindent (b) Let $\tau$ be a topology on $X$ containing $\{ f^{-1}(D)\mid D\!\in\! {\cal B}\}$. Then 
$\sigma\! :=\!\{ f[A]\mid A\!\in\!\tau\}$ is a topology on $Y$ containing ${\cal B}$. Thus  
$B\!\in\!{\bf\Gamma}(Y,\sigma )$ since $B\!\in\! {\bf\Gamma}({\cal B})$. Therefore 
$f^{-1}(B)\!\in\! {\bf\Gamma}(X,\tau )$ since $f\! :\! (X,\tau )\!\rightarrow\! (Y,\sigma )$ is continuous.\bigskip

\noindent (c) We argue by induction on $\xi\! -\!\eta$. The result is clear if $\xi\! -\!\eta\! =\! 0$. So assume that $\xi\! -\!\eta\!\geq\! 1$. Write $A\! =\!\bigcap_{n\in\omega}\ \neg A_n$, where $\eta_n\! <\!\xi$ and 
$A_n\!\in\! {\bf\Pi}^0_{\eta_n}(X)$. As $X$ is metrizable, we may assume that $\eta\!\leq\!\eta_n$. The induction assumption gives ${\cal B}_n\!\subseteq\! {\bf\Pi}^0_{\eta}(X)$ countable such that 
$A_n\!\in\! {\bf\Pi}^0_{1+(\eta_n-\eta )}(\check {\cal B}_n)$. It remains to set 
${\cal B}\! :=\!\bigcup_{n\in\omega}\ {\cal B}_n$.$\hfill{\square}$\bigskip

\bf\noindent (A) The witnesses\rm\bigskip

\noindent\bf Notation.\rm\ We first define a map producing witnesses for the fact that 
$\vec v\!\in\!\Theta^\infty$. More specifically, we define a map 
${\mathfrak W}\! :\!\Theta^\infty\!\rightarrow\!\Theta^\infty\cup (\Theta^\infty)^\omega$. Let 
$\vec v\! :=\! (\alpha ,a_0,a_1,\underline{a}_0,\underline{a}_1,r)\!\in\!\Theta^\xi\!\setminus\!\Theta^{<\xi}$. 
If $\vert (\alpha )_1\vert\! =\! 0$, then we set ${\mathfrak W}(\vec v)\! :=\!\vec v$. If 
$\vert (\alpha )_1\vert\! =\! 1$, then using the definition of $\Theta$ we set 
$${\mathfrak W}(\vec v)\! :=\! (<(\alpha )_{2+j}>,a_0,a_1,\underline{a}_0,\underline{a}_1,a_1).$$ 
Note that ${\mathfrak W}(\vec v)\!\in\!\Theta^{<\xi}$. If $\vert (\alpha )_1\vert\! =\! 2$, then we set
$${\mathfrak W}(\vec v)(p)\! :=\!\left\{\!\!\!\!\!\!\!
\begin{array}{ll}
& \big(<\! (\alpha )_{2+<0,q>}\! >,a_0,a_1,(a'_0)_0,(a'_1)_0,(r')_0\big)\mbox{ if }p\! =\! 0\mbox{,}\cr\cr
& \big(<\! (\alpha )_{2+<(p)_0+1,q>}\! >,a_0,a_1,(a'_0)_p,(a'_1)_p,(r')_p\big)\mbox{ if }p\!\geq\! 1.
\end{array}
\right.$$
Here again, ${\mathfrak W}(\vec v)(p)\!\in\!\Theta^{<\xi}$.\bigskip

\noindent $\bullet$ Similarly, we define a map ${\mathfrak W}^1$ witnessing that 
$\vec w\!\in\!\Upsilon_1^\infty$. Moreover, we keep in mind $\gamma'$. More specifically, we define a map 
${\mathfrak W}^1\! :\!\Upsilon_1^\infty\!\rightarrow\!
\Upsilon_1^\infty\cup (\omega^\omega\!\times\!\Upsilon_1^\infty )\cup
\big(\omega^\omega\!\times\! (\Upsilon_1^\infty)^\omega\big)$. Let $\vec w\! :=\! (\alpha ,\beta ,\gamma )$ 
in $\Upsilon_1^\xi\!\setminus\!\Upsilon_1^{<\xi}$. If $\vert (\alpha )_1\vert\! =\! 0$, then we set 
${\mathfrak W}^1(\vec w)\! :=\!\vec w$. If $\vert (\alpha )_1\vert\! =\! 1$, then using the definition of 
$\Upsilon_1$ and some choice for $\gamma'$, we set 
${\mathfrak W}^1(\vec w)\! :=\!\big(\gamma',(<(\alpha )_{2+j}>,\beta^*,\gamma')\big)$. If 
$\vert (\alpha )_1\vert\! =\! 2$, then we set 
${\mathfrak W}^1(\vec w)\! :=\!\big(\gamma',{\mathfrak W}^1_1(\vec w)\big)$, where
$${\mathfrak W}^1_1(\vec w)(p)\! :=\!\left\{\!\!\!\!\!\!\!
\begin{array}{ll}
& \big(<\! (\alpha )_{2+<0,q>}\! >,(\beta^*)_0,(\gamma')_0\big)\mbox{ if }p\! =\! 0\mbox{,}\cr\cr
& \big(<\! (\alpha )_{2+<(p)_0+1,q>}\! >,\big( (\beta^*)_p\big)_0,\big( (\gamma')_p\big)_0\big)
\mbox{ if }p\!\geq\! 1.
\end{array}
\right.$$
\bf (B) The trees associated with the codes for the non self-dual Wadge classes of Borel  sets\rm\bigskip

\noindent $\bullet$ Recall the definition of $\mathfrak T(\alpha )$ after Lemma 6.2. Similarly, we define 
${\mathfrak T}\! :\!\Upsilon_1^\infty\!\rightarrow\!\{\mbox{trees on}\ \omega\!\times\!\Upsilon_1^\infty\}$ as follows. Let $\vec w\! :=\! (\alpha ,\beta ,\gamma )\!\in\!\Upsilon_1^\xi\!\setminus\!\Upsilon_1^{<\xi}$. We set
$${\mathfrak T}(\vec w)\! :=\!\left\{\!\!\!\!\!\!\!\!
\begin{array}{ll}
& \{\emptyset\}\cup\{ <\! (0,\vec w)\! >\}\mbox{ if }\vert (\alpha )_1\vert\! =\! 0\mbox{,}\cr
& \{\emptyset\}\cup\big\{ (0,\vec w)^\frown s\mid s\!\in\! 
{\mathfrak T}\big( {\mathfrak W}^1_1(\vec w)\big)\big\}\mbox{ if }\vert (\alpha )_1\vert\! =\! 1\mbox{,}\cr
& \{\emptyset\}\cup\bigcup_{p\in\omega}\ \big\{ (p,\vec w)^\frown s\mid s\!\in\! 
{\mathfrak T}\big( {\mathfrak W}^1_1(\vec w)(p)\big)\big\}\mbox{ if }\vert (\alpha )_1\vert\! =\! 2.
\end{array}
\right.$$
Here again ${\mathfrak T}(\vec w)$ is always a countable well founded tree containing the sequence 
$<\! (0,\vec w)\! >$. The set of maximal sequences in ${\mathfrak T}(\vec w)$ is 
${\cal M}_{\vec w}\! :=\!\{ s\!\in\! {\mathfrak T}(\vec w)\mid\forall t\!\in\! {\mathfrak T}(\vec w)\ \ 
s\!\subseteq\! t\Rightarrow s\! =\! t\}$.\bigskip

\noindent $\bullet$ Fix ${\vec w\! :=\! (\alpha ,\beta ,\gamma )\!\in\!\Upsilon_1^\infty}$ with 
$\alpha\!\in\!\Borel$ normalized. In the sequel, it will be convenient to set, for 
$s\!\in\! {\mathfrak T}(\vec w)\!\setminus\! {\cal M}_{\vec w}$,
$$s_1(\vert s\vert )\! :=\!\left\{\!\!\!\!\!\!\!
\begin{array}{ll}
& \vec w\mbox{ if }s\! =\!\emptyset\mbox{,}\cr
& {\mathfrak W}^1_1\big( s_1(\vert s\vert\! -\! 1)\big)\mbox{ if }s\!\not=\!\emptyset\wedge
\big\vert\big( s_1(\vert s\vert\! -\! 1)(0)\big)_1\big\vert\! =\! 1\mbox{,}\cr
& {\mathfrak W}^1_1\big( s_1(\vert s\vert\! -\! 1)\big)\big( s_0(\vert s\vert\! -\! 1)\big)\mbox{ if }
s\!\not=\!\emptyset\wedge\big\vert\big( s_1(\vert s\vert\! -\! 1)(0)\big)_1\big\vert\! =\! 2.
\end{array}
\right.$$
$\bullet$ Let $s\!\in\! {\mathfrak T}(\vec w)$. We set 
$B_s\! :=\!\{ i\! <\!\vert s\vert\mid\vert\big( s_1(i)(0)\big)_1\vert\! =\! 2\}$. As $\alpha$ is normalized, $B_s$ is an integer. We always have $B_s\!\leq\!\vert s\vert$. If moreover 
$s\!\in\! {\mathfrak T}(\vec w)\!\setminus\! {\cal M}_{\vec w}$, then we set 
$B'_s\! :=\!\{ i\!\leq\!\vert s\vert\mid\vert\big( s_1(i)(0)\big)_1\vert\! =\! 2\}$.\bigskip

\noindent $\bullet$ The ordinals $\vert (\alpha )_0\vert$, for $\alpha\!\in\!\Borel\cap\Lambda^\infty$, will be of particular importance in the sequel. We define a function ${\cal Z}\! :\! 
{\mathfrak T}(\vec w)\!\setminus\! {\cal M}_{\vec w}\!\rightarrow\! (\omega_1^{\mbox{CK}})^{<\omega}$ satisfying $\vert {\cal Z}(s)\vert\! =\!\vert s\vert\! +\! 1$. The sequence ${\cal Z}(s)$ gives the ordinals $\xi$ of the $\bormxi$ sets coded by $s$. We set 
${\cal Z}(s)(i)\! :=\!\big\vert\big( s_1(i)(0)\big)_0\big\vert$ if $i\!\leq\!\vert s\vert$. Note the following properties of ${\cal Z}(s)$, easy to check:\smallskip

- ${\cal Z}(s)(i)$ depends only on $s\vert i$.\smallskip

- ${\cal Z}(s)\!\subseteq\! {\cal Z}(t)$ if $s\!\subseteq\! t$.\smallskip

- ${\cal Z}(s)(i\! +\! 1)\!\geq\! {\cal Z}(s)(i)$ or ${\cal Z}(s)(i\! +\! 1)\! =\! 0$ if $i\! <\!\vert s\vert$.\smallskip

- ${\cal Z}(s)(i\! +\! 1)\! =\! 0$ if ${\cal Z}(s)(i)\! =\! 0$ and $i\! <\!\vert s\vert$.\smallskip

- $\big( {\cal Z}(s)(i)\big)_{i\in B'_s}$ is a non-decreasing sequence of non zero recursive ordinals.\bigskip

\noindent\bf (C) The resolution families\rm\bigskip

\noindent $\bullet$ Fix $\vec w\! :=\! (\alpha ,\beta ,\gamma )\!\in\!\Upsilon_1^\infty$ with 
$\alpha\!\in\!\Borel$ normalized, and $p\!\geq\! 1$. We set
$$\tilde P^{\vec w}_p\! :=\!\left\{\!\!\!\!\!\!\!\!
\begin{array}{ll}
& \omega^\omega\mbox{ if }\vert (\alpha )_1\vert\!\leq\! 1\mbox{,}\cr
& C^{\omega^\omega}_{(({\mathfrak W}^1_0(\vec w))_p)_1}\mbox{ if }\vert (\alpha )_1\vert\! =\! 2.
\end{array}
\right.$$
Note that $\tilde P^{\vec w}_p\!\in\! {\bf\Pi}^0_{\vert (\alpha )_0\vert}(\omega^\omega )$ if 
$\vert (\alpha )_1\vert\! =\! 2$, by Lemma 6.1.\bigskip

\noindent $\bullet$ Recall the finite sets $c_l\!\subseteq\! d^d$ defined at the end of the proof of Proposition 2.2 (we only used the fact that $T_d$ has finite levels to see that they are finite). We put 
$c\! :=\!\bigcup_{l\in\omega}\ c_l$, so that $c$ is countable. This will be the countable set $c$ of Definition 4.3.1.\bigskip

\noindent $\bullet$ Recall the embedding $j_d$ defined before Theorem 6.9. We set 
${\cal P}^{\vec w}_p\! :=\! h[j_d^{-1}(\tilde P^{\vec w}_p)\cap c^\omega ]$, so that the union 
${\cal P}^{\vec w}_p\cup {\cal P}^{\vec w}_q\! =\! [\subseteq]$ if $p\!\not=\! q\!\geq\! 1$. Moreover, 
${\cal P}^{s_1(i)}_p\!\in\! {\bf\Pi}^0_{{\cal Z}(s)(i)}([\subseteq ])$ if 
$s\!\in\! {\mathfrak T}(\vec w)\!\setminus\! {\cal M}_{\vec w}$ and $i\!\in\! B'_s$.\bigskip

\noindent $\bullet$ If $T$ is a tree and $s\!\in\! T$, then $T_s\! :=\!\{ t\!\in\! T\mid s\!\subseteq\! t\}$.\bigskip 

\noindent $\bullet$ Fix $\vec w\! :=\! (\alpha ,\beta ,\gamma )\!\in\!\Upsilon_1^\infty$ with 
$\alpha\!\in\!\Borel$ normalized and $\vert (\alpha )_1\vert\! =\! 2$. We say that 
$s\!\in\! {\mathfrak T}({\vec w})$ is $extensible$ if there is $t\!\in\! {\mathfrak T}({\vec w})_s$ such that 
$\vert s\vert\! <\! B_t$ (which implies that $s\!\notin\! {\cal M}_{\vec w}$). We will construct, for each $s$ extensible, a resolution family $(R^{(\rho)}_{s})_{\rho\leq\eta_{s}}$. Simultaneously, we construct some ordinals $\xi_{s}$ and $\theta_{s}$. If $\theta$ is an ordinal, then we set
$$\theta^*\! :=\!\left\{\!\!\!\!\!\!\!
\begin{array}{ll}
& \eta\mbox{ if }\theta\! =\!\eta\! +\! 1\mbox{,}\cr
& \theta\mbox{ otherwise}
\end{array}
\right.$$ 
(this is what appears in the Debs-Saint Raymond theorem). We will have $\eta_{s}\! =\!\theta^*_{s}$, $\xi_{s}\! =\! {\cal Z}(s)(\vert s\vert )$ and 
$$\theta_{s}\! :=\!\left\{\!\!\!\!\!\!\!
\begin{array}{ll}
& \xi_{s}\! =\! {\cal Z}(s)(0)\! =\!\vert (\alpha )_0\vert\mbox{ if }s\! =\!\emptyset\mbox{,}\cr
& 1\! +\! (\xi_{s}\! -\!\xi_{s^-})\mbox{ if }s\!\not=\!\emptyset .
\end{array}
\right.$$ 

\vfill\eject

 We want the resolution family  to satisfy the following conditions:\bigskip

\noindent - The family $(R^{(\rho)}_{s})_{\rho\leq\eta_{s}}$ is uniform if $\theta_{s}$ is a limit ordinal.\smallskip

\noindent - $R^{(0)}_{\emptyset}\! =\subseteq$, and 
$R^{(\eta_{s^-})}_{s^-}\! =\! R^{(0)}_{s}$ if $s\!\not=\!\emptyset$.\smallskip

\noindent - $\Pi_{s}\! :\! [R^{(\eta_{s})}_{s}]\!\rightarrow\! [R^{(0)}_{s}]$ is a continuous bijection.\smallskip

\noindent - $(\Pi_{s\vert 0}\circ\Pi_{s\vert 1}\circ ...\circ \Pi_{s})^{-1}({\cal P}^{s_1(\vert s\vert )}_p)\!\in\!
{\bf\Pi}^0_1([R^{(\eta_{s})}_{s}])$ if $p\!\geq\! 1$.\smallskip

\noindent - $(\Pi_{s\vert 0}\circ\Pi_{s\vert 1}\circ ...\circ \Pi_{s})^{-1}({\cal P}^{t_1(j+1)}_p)\!\in\!
{\bf\Pi}^0_{1+({\cal Z}(t)(j+1)-\xi_{s})}([R^{(\eta_{s})}_{s}])$ if $p\!\geq\! 1$, 
$t\!\in\! {\mathfrak T}(\vec w)_{s}\!\setminus\! {\cal M}_{\vec w}$ and $\vert s\vert\! <\! j\! +\! 1\!\in\! B'_t$.\bigskip

\noindent $\bullet$ The construction is by induction on $\vert s\vert$. Assume that $s\! =\!\emptyset$, 
$p\!\geq\! 1$, $t\!\in\! {\mathfrak T}(\vec w)\!\setminus\! {\cal M}_{\vec w}$ and $j\! +\! 1\!\in\! B'_t$. Proposition 6.12.(c) gives ${\cal B}^{t,j}_p\!\subseteq\! {\bf\Pi}^0_{\theta_\emptyset}([\subseteq ])$ countable such that 
${\cal P}^{t_1(j+1)}_p\!\in\! {\bf\Pi}^0_{1+({\cal Z}(t)(j+1)-\theta_\emptyset )}(\check {\cal B}^{t,j}_p)$. This implies that $u_\emptyset\! :=\!\{ {\cal P}^{\vec w}_p\mid p\!\geq\! 1\}\cup
\bigcup_{p\geq 1,t\in {\mathfrak T}(\vec w)\setminus {\cal M}_{\vec w},j+1\in B'_t}\ {\cal B}^{t,j}_p$ is countable and made of ${\bf\Pi}^0_{\theta_\emptyset}([\subseteq ])$ sets. Theorems 4.3.4 and 4.4.4 give a  family $(R^{(\rho)}_{\emptyset})_{\rho\leq\eta_{\emptyset}}$, uniform if $\theta_{\emptyset}$ is a limit ordinal, such that\bigskip

\noindent - $R^{(0)}_{\emptyset}\! =\subseteq$.\smallskip

\noindent - $\Pi_{\emptyset}\! :\! [R^{(\eta_{\emptyset})}_{\emptyset}]\!\rightarrow\! [R^{(0)}_{\emptyset}]$ is a continuous bijection.\smallskip

\noindent - $\Pi_{\emptyset}^{-1}(Q)\!\in\!\bormone ([R^{(\eta_{\emptyset})}_{\emptyset}])$ for each 
$Q\!\in\! u_{\emptyset}$.\bigskip

 This family is suitable, by Proposition 6.12.\bigskip

\noindent $\bullet$ Assume now that $s\!\not=\!\emptyset$ is extensible, and the construction is done for the strict predecessors of $s$. Note that 
$(\Pi_{s\vert 0}\circ\Pi_{s\vert 1}\circ ...\circ \Pi_{s^-})^{-1}({\cal P}^{s_1(\vert s\vert )}_p)
\!\in\! {\bf\Pi}^0_{\theta_{s}}([R^{(\eta_{s^-})}_{s^-}])$. Assume that $p\!\geq\! 1$, 
$t\!\in\! {\mathfrak T}(\vec w)_{s}\!\setminus\! {\cal M}_{\vec w}$ and $\vert s\vert\! <\! j\! +\! 1\!\in\! B'_t$. Then Proposition 6.12.(c) gives a countable family 
${\cal C}^{t,j}_p\!\subseteq\! {\bf\Pi}^0_{\theta_{s}}([R_{s^-}^{\eta_{s^-}}])$ such that 
$(\Pi_{s\vert 0}\circ\Pi_{s\vert 1}\circ ...\circ\Pi_{s^-})^{-1}({\cal P}^{t_1(j+1)}_p)
\!\in\! {\bf\Pi}^0_{1+({\cal Z}(t)(j+1)-\xi_{s})}(\check {\cal C}^{t,j}_p)$. This implies that
$$u_{s}\! :=\!\{ (\Pi_{s\vert 0}\circ\Pi_{s\vert 1}\circ ...\circ\Pi_{s^-})^{-1}
({\cal P}^{s_1(\vert s\vert)}_p)\mid p\!\geq\! 1\}\cup\bigcup_{p\geq 1,t\in {\mathfrak T}(\vec w)_{s}\setminus {\cal M}_{\vec w},\vert s\vert <j+1\in B'_t}\ {\cal C}^{t,j}_p$$
is countable and made of 
${\bf\Pi}^0_{\theta_{s}}([R^{(\eta_{s^-})}_{s^-}])$ sets. Theorems 4.3.4 and 4.4.4 give a resolution family 
$(R^{(\rho)}_{s})_{\rho\leq\eta_{s}}$, uniform if $\theta_{s}$ is a limit ordinal, such that\bigskip

\noindent - $R^{(0)}_{s}\! =\! R^{(\eta_{s^-})}_{s^-}$.\smallskip

\noindent - $\Pi_{s}\! :\! [R^{(\eta_{s})}_{s}]\!\rightarrow\! [R^{(0)}_{s}]$ is a continuous bijection.\smallskip

\noindent - $\Pi_{s}^{-1}(Q)\!\in\!\bormone ([R^{(\eta_{s})}_{s}])$ for each $Q\!\in\! u_{s}$.\bigskip

 This family is suitable, by Proposition 6.12. This completes the construction of the families.
 
\vfill\eject
 
\bf\noindent (D) The subsets of $T_d$\rm\bigskip

 We now build some subsets of $T_d$ that will play the role that $D$ and $T_d\!\setminus\! D$ played in the proof of Theorem 4.4.1. Fix $\vec w\! :=\! (\alpha ,\beta ,\gamma )\!\in\!\Upsilon_1^\infty$ with 
$\alpha\!\in\!\Borel$ normalized and $\vert (\alpha )_1\vert\! =\! 2$. We will define a family of subsets of $T_d$ as follows. Assume that $s\!\in\! {\mathfrak T}(\vec w)$ is extensible. We set, for $q\!\geq\! 1$, 
$$\begin{array}{ll}
& P_0(s)\! :=\!\Big\{\vec s\!\in\! T_d\mid\vec s\! =\!\vec {\emptyset}\vee
\forall p\!\geq\! 1\ \ \exists {\cal B}_p\!\in\! (\Pi_{s\vert 0}\circ\Pi_{s\vert 1}\circ ...\circ\Pi_{s})^{-1}
({\cal P}^{s_1(\vert s\vert)}_p)\ \ \vec s\!\in\! {\cal B}_p\Big\}\mbox{,}\cr
& P_q(s)\! :=\!\Big\{\vec s\!\in\! T_d\mid\vec s\!\not=\!\vec {\emptyset}\wedge
\forall {\cal B}_q\!\in\! (\Pi_{s\vert 0}\circ\Pi_{s\vert 1}\circ ...\circ\Pi_{s})^{-1}({\cal P}^{s_1(\vert s\vert)}_q)\ \ 
\vec s\!\notin\! {\cal B}_q\ \wedge\cr
& \hfill{\forall p\!\in\!\omega\!\setminus\!\{ 0,q\}\ \ 
\exists {\cal B}_p\!\in\! (\Pi_{s\vert 0}\circ\Pi_{s\vert 1}\circ ...\circ\Pi_{s})^{-1}({\cal P}^{s_1(\vert s\vert)}_p)\ \ \vec s\!\in\! {\cal B}_p\Big\}.}
\end{array}$$
Note that the $P_q(s)$'s are pairwise disjoint. The next lemma associates to each $\vec t\!\in\! T_d$ a sequence $s(\vec t~)$ in ${\mathfrak T}(\vec w)$ saying in which $P_q(s)$'s the sequence $\vec t$ is.

\begin{prop} Let $\vec w\! :=\! (\alpha ,\beta ,\gamma )\!\in\!\Upsilon_1^\infty$ with $\alpha\!\in\!\Borel$ normalized and $\vert (\alpha )_1\vert\! =\! 2$, and $\vec t\!\in\! T_d$. Then there are $l\!\in\!\omega$ and 
$s(\vec t~)\!\in\! {\mathfrak T}(\vec w)$ of length $l$ such that\smallskip

\noindent (a) $\vec t\!\in\!\bigcap_{i<l}\ P_{s(\vec t~)(i)(0)}\big( s(\vec t~)\vert i\big)$.\smallskip

\noindent (b) If $s(\vec t~)$ is extensible by $t$, then $\vec t\!\notin\! P_{t(l)(0)}(t\vert l)$.\end{prop}

\noindent\bf Proof.\rm\ We actually construct, for $j\!\in\!\omega$, a sequence 
$s_j\!\in\! {\mathfrak T}(\vec w)$. We will have $s_j\!\subseteq\! s_{j+1}$, $\vert s_j\vert\! =\! j$ if $j\!\leq\! l$, $s_j\! =\! s_l$ if $j\! >\! l$, and $\vec t\!\in\!\bigcap_{i<\vert s_j\vert}\ P_{s_j(i)(0)}\big( s_j\vert i\big)$. At the end, $s(\vec t~)$ will be $s_l$. The definition of $s_j$ is by induction on $j$. Assume that $(s_k)_{k\leq j}$ are constructed satisfying these properties, which is the case for $j\! =\! 0$. We may assume that 
$\vert s_j\vert\! =\! j$.\bigskip

 If $s_j$ is not extensible or $\vec t\!\notin\! {\cal B}$ for each ${\cal B}\!\in\! [R^{(\eta_{s_j})}_{s_j}]$, then we set $s_{j+1}\! :=\! s_j$. If $\vec t\!\in\! {\cal B}$ for some ${\cal B}\!\in\! [R^{(\eta_{s_j})}_{s_j}]$, then there is a unique integer $q$ such that $\vec t\!\in\! P_q(s_j)$ since 
$$(\Pi_{s_j\vert 0}\circ\Pi_{s_j\vert 1}\circ ...\circ \Pi_{s_j})^{-1}({\cal P}^{(s_j)_1(j)}_p)\cup 
(\Pi_{s_j\vert 0}\circ\Pi_{s_j\vert 1}\circ ...\circ \Pi_{s_j})^{-1}({\cal P}^{(s_j)_1(j)}_q)\! =\! 
[R^{(\eta_{s_j})}_{s_j}]$$ if $p\!\not=\! q\!\geq\! 1$. We will have $\vert s_{j+1}\vert\! =\! j\! +\! 1$, and 
$s_{j+1}(j)(0)\! :=\! q$. Moreover,
$$s_{j+1}(j)(1)\! :=\!\left\{\!\!\!\!\!\!\!
\begin{array}{ll}
& \vec w\mbox{ if }j\! =\! 0\mbox{,}\cr
& {\mathfrak W}^1_1\big( s_j(j\! -\! 1)(1)\big)\big( s_j(j\! -\! 1)(0)\big)\mbox{ if }j\!\geq\! 1.
\end{array}
\right.$$
This completes the construction of the $s_j$'s, and they are in ${\mathfrak T}(\vec w)$. The well-foundedness of ${\mathfrak T}(\vec w)$ proves the existence of $l$, and $s(\vec t~)$ is suitable.
$\hfill{\square}$\bigskip

\noindent\bf Notation.\rm\ Proposition 6.13 associates $s(\vec t~)\!\in\! {\mathfrak T}(\vec w)$ to 
$\vec t\!\in\! T_d$. Under the same conditions, we can associate $S(\vec t~)\!\in\! {\cal M}_{\vec w}$ to 
$\vec t$. To do this, we need the following lemma:

\begin{lem} Let $\vec w\! :=\! (\alpha ,\beta ,\gamma )\!\in\!\Upsilon_1^\infty$ with $\alpha\!\in\!\Borel$ normalized and $\vert (\alpha )_1\vert\! =\! 2$, and $s\!\in\! {\mathfrak T}({\vec w})$. Then there is 
$S\!\in\! {\cal M}_{\vec w}$ extending $s$ such that $S_0(i)\! =\! 0$ for $\vert s\vert\!\leq\! i\! <\!\vert S\vert$.
\end{lem}

\noindent\bf Proof.\rm\ If $s\! =\!\emptyset$, then we set $S(0)\! :=\! (0,\vec w)$ and, if 
${\mathfrak W}^1\big( S_1(i)\big)\!\not=\! S_1(i)$, then we set
$$S(i\! +\! 1)\! :=\!\left\{\!\!\!\!\!\!\!
\begin{array}{ll}
& \Big( 0,{\mathfrak W}^1_1\big( S(i)\big)\Big)\mbox{ if }
{\mathfrak W}^1_1\big( S(i)\big)\!\in\!\Upsilon_1^\infty\mbox{,}\cr\cr
& \Big( 0,{\mathfrak W}^1_1\big( S(i)\big)(0)\Big)\mbox{ if }
{\mathfrak W}^1_1\big( S(i)\big)\!\in\! (\Upsilon_1^\infty )^\omega .
\end{array}
\right.$$
By induction, we see that $S\vert (i\! +\! 1)\!\in\! {\mathfrak T}(\vec w)$ for each $i\! <\!\vert S\vert$, which proves that the length of $S$ is finite since ${\mathfrak T}(\vec w)$ is well-founded. Thus 
$S\!\in\! {\cal M}_{\vec w}$.\bigskip

 If $s\!\not=\!\emptyset$, then $S(\vert s\vert\! -\! 1)$ is defined. We argue similarly. The only thing to change is that 
$$S(\vert s\vert )\! :=\!\Big( 0,{\mathfrak W}^1_1\big( s(\vert s\vert\! -\! 1)\big)
\big(s_0(\vert s\vert\! -\! 1)\big)\Big)$$ 
if ${\mathfrak W}^1\big( s_1(\vert s\vert\! -\! 1)\big)\!\not=\! s_1(\vert s\vert\! -\! 1)$ and 
${\mathfrak W}^1_1\big( s(\vert s\vert\! -\! 1)\big)\!\in\! (\Upsilon_1^\infty )^\omega$.
$\hfill{\square}$\bigskip

 We now associate a maximal extension $S(\vec t~)$ of $s(\vec t~)$ to any $\vec t$ in $T_d$.\bigskip

\noindent\bf Remark.\rm\ In particular, there is $S(\vec\emptyset )\!\in\! {\cal M}_{\vec w}$ with 
$\big( S(\vec\emptyset )\big)_0(i)\! =\! 0$ for $i\! <\!\vert S(\vec\emptyset )\vert$. Note that 
$s(\vec\emptyset )\!\subseteq\! S(\vec\emptyset )$. If $\vec\emptyset\!\not=\!\vec t\!\in\! T_d$, then we define 
$S(\vec t~)$ by induction on $\vert\vec t\vert$:\bigskip

\noindent - If $s(\vec t~)\! =\!\emptyset$, then $\vec t\!\not=\!\emptyset$ since 
$\vec\emptyset\!\in\! P_0(\emptyset )$, and $S(\vec t~)\! :=\! S({\vec t~}^{\eta_\emptyset}_\emptyset )$.\smallskip

\noindent - If $s(\vec t~)\!\not=\!\emptyset$ and ${\vec t~}^{\eta_{s(\vec t~)^-}}_{s(\vec t~)^-}\!\in\!
\bigcap_{i<\vert s(\vec t~)\vert}\ P_{s(\vec t~)(i)(0)}\big( s(\vec t~)\vert i\big)$, then 
$S(\vec t~)\! :=\! S({\vec t~}^{\eta_{s(\vec t~)^-}}_{s(\vec t~)^-})$.\smallskip

\noindent - If $s(\vec t~)\!\not=\!\emptyset$ and ${\vec t~}^{\eta_{s(\vec t~)^-}}_{s(\vec t~)^-}\!\notin\!
\bigcap_{i<\vert s(\vec t~)\vert}\ P_{s(\vec t~)(i)(0)}\big( s(\vec t~)\vert i\big)$, then $S(\vec t~)$ is the extension of $s(\vec t~)$ given by Lemma 6.14 applied to $s\! :=\! s(\vec t~)$.\bigskip

 Note that $S(\vec t~)\!\in\! {\cal M}_{\vec w}$ and is always an extension of $s(\vec t~)$, by induction on 
$\vert\vec t\vert$. This comes from the fact that 
$s(\vec t~)\!\subseteq\! s({\vec t~}^{\eta_{s(\vec t~)^-}}_{s(\vec t~)^-})$ in the second case.\bigskip

\bf\noindent (E) The tuples\rm\bigskip

 We now keep in mind the tuples $(\alpha ,a_0,a_1,\underline{a}_0,\underline{a}_1,r)$ along any sequence of ${\mathfrak T}({\vec w})$, using the witness map $\mathfrak W$. Fix 
$\vec w\! :=\! (\alpha ,\beta ,\gamma )\!\in\!\Upsilon_1^\infty$,  
$\vec v\! :=\! (\alpha ,a_0,a_1,\underline{a}_0,\underline{a}_1,r)\!\in\!\Theta^\infty$ with 
$\alpha\!\in\!\Borel$ normalized and $\vert (\alpha )_1\vert\! =\! 2$. We will define a map 
$V\! :\! {\mathfrak T}({\vec w})\!\rightarrow\! (\Theta^\infty )^{<\omega}$ such that 
$\vert V(s)\vert\! =\!\vert s\vert$, $V(s)(i)$ depends only on $s\vert i$ as follows. We set, for 
$i\! <\!\vert s\vert$, 
$$V(s)(i)\! :=\!\left\{\!\!\!\!\!\!\!
\begin{array}{ll}
& \vec v\mbox{ if }i\! =\! 0\mbox{,}\cr
& {\mathfrak W}\big( V(s)(i\! -\! 1)\big)\mbox{ if }i\!\geq\! 1\wedge
\big\vert\big( V(s)(i\! -\! 1)(0)\big)_1\big\vert\!\leq\! 1\mbox{,}\cr
& {\mathfrak W}\big( V(s)(i\! -\! 1)\big)\big( s_0(i\! -\! 1)\big)\mbox{ if }i\!\geq\! 1\wedge
\big\vert\big( V(s)(i\! -\! 1)(0)\big)_1\big\vert\! =\! 2.
\end{array}
\right.$$

\begin{lem} Let $\vec w\! :=\! (\alpha ,\beta ,\gamma )\!\in\!\Upsilon_1^\infty$,  
$\vec v\! :=\! (\alpha ,a_0,a_1,\underline{a}_0,\underline{a}_1,r)\!\in\!\Theta^\infty$ with 
$\alpha\!\in\!\Borel$ normalized and $\vert (\alpha )_1\vert\! =\! 2$, $s\!\in\! {\mathfrak T}({\vec w})$, and 
$i\! <\!\vert s\vert$. Then $V(s)(i)(0)\! =\! s_1(i)(0)$. In particular, $s\!\notin\! {\cal M}_{\vec w}$ and 
$i\!\leq\!\vert s\vert$ imply that ${\cal Z}(s)(i)\! =\!\vert\big( V(s)(i)(0)\big)_0\vert$.\end{lem}

\noindent\bf Proof.\rm\ The last assertion clearly comes from the first one. The proof is by induction on $i$. The assertion is clear for $i\! =\! 0$ since $V(s)(0)(0)\! =\! s_1(0)(0)\! =\!\alpha$. Assume that it holds for 
$i\! <\!\vert s\vert\! -\! 1$.\bigskip

\noindent $\bullet$ If $i\!\notin\! B_s$, then 
$\vert\big( V(s)(i)(0)\big)_1\vert\! =\!\vert\big( s_1(i)(0)\big)_1\vert\! =\! 1$. Thus\smallskip

\centerline{$V(s)(i\! +\! 1)(0)\! =\! {\mathfrak W}\big( V(s)(i)\big)(0)\! =<\! \big( V(s)(i)(0)\big)_{2+j}\! >=
<\!\big( s_1(i)(0)\big)_{2+j}\! >=\! s_1(i\! +\! 1)(0).$}\smallskip

\noindent $\bullet$ If $i\!\in\! B_s$, then 
$\big\vert\big( V(s)(i)(0)\big)_1\big\vert\! =\!\vert\big( s_1(i)(0)\big)_1\vert\! =\! 2$. If moreover $s_0(i)\! =\! 0$, then
$$V(s)(i\! +\! 1)(0)\! =<\!\big( V(s)(i)(0)\big)_{2+<0,q>}\! >=<\!\big( s_1(i)(0)\big)_{2+<0,q>}\! >
=\! s_1(i\! +\! 1)(0).$$
The argument is similar if $s_0(i)\!\geq\! 1$.$\hfill{\square}$\bigskip

 The next lemma is a preparation for Lemma 6.21, which is the crucial step to prove a version of the claim in the proof of Theorem 4.4.1 for the non self-dual Wadge classes of Borel sets.

\begin{lem} Let $\vec w\! :=\! (\alpha ,\beta ,\gamma )\!\in\!\Upsilon_1^\infty$,  
$\vec v\! :=\! (\alpha ,a_0,a_1,\underline{a}_0,\underline{a}_1,r)\!\in\!\Theta^\infty$ with 
${\alpha\!\in\!\Borel}$ normalized and $\vert (\alpha )_1\vert\! =\! 2$, $s\!\in\! {\mathfrak T}({\vec w})$, and 
$i\!\in\! B_s$.\smallskip

\noindent (a) If $s_0(i)\! =\! 0$, then $\neg {\cal U}_{V(s)(i)(5)}\!\subseteq\!\neg {\cal U}_{V(s)(i+1)(5)}$.\smallskip

\noindent (b) We have $\neg {\cal U}_{V(s)(i)(5)}\!\subseteq\!
\overline{\neg {\cal U}_{V(s)(i+1)(5)}}^{\tau_{\xi_{s\vert i}}}$.\end{lem}

\noindent\bf Proof.\rm\ (a) We have $V(s)(i\! +\! 1)\! =\! {\mathfrak W}\big( V(s)(i)\big)(0)$, by Lemma 6.15. Thus 
$$V(s)(i+1)(5)\! =\! {\mathfrak W}\big( V(s)(i)\big)(0)(5)\! =\! (r')_0$$ 
for some $r'$ for which $\neg {\cal U}_{V(s)(i)(5)}\!\subseteq\!\neg {\cal U}_{(r')_0}$, by the 2nd and the 4th remarks after the definition of $\Theta$.\bigskip

\noindent (b) We may assume that $s_0(i)\!\geq\! 1$, so that $V(s)(i+1)(5)\! =\! (r')_{s_0(i)}$, and 
$$\neg {\cal U}_{V(s)(i)(5)}\!\subseteq\!\overline{\neg {\cal U}_{V(s)(i+1)(5)}}^
{\tau_{\vert (V(s)(i)(0))_0\vert}}$$ 
by the 5th remark after the definition of $\Theta$ and the definition of $f_a$. We are done by Lemma 6.15.
$\hfill{\square}$\bigskip

\bf\noindent (F) The sequences of integers\rm\bigskip

 We have to keep in mind the integers $s_0(i)$ for $s\!\in\! {\mathfrak T}({\vec w})$. We will consider an ordering of these finite sequences of integers that will help us to prove the claim just mentioned.\bigskip

\noindent\bf Notation.\rm\ Fix $\vec w\! :=\! (\alpha ,\beta ,\gamma )\!\in\!\Upsilon_1^\infty$,  
$\vec v\! :=\! (\alpha ,a_0,a_1,\underline{a}_0,\underline{a}_1,r)\!\in\!\Theta^\infty$ with 
$\alpha\!\in\!\Borel$ normalized and $\vert (\alpha )_1\vert\! =\! 2$, and 
$s,s'\!\in\! {\mathfrak T}({\vec w})$.\bigskip

\noindent $\bullet$ If $s$ and $s'$ are not compatible, then we denote $s\wedge s'\! :=\! s\vert i\! =\! s'\vert i$, where $i$ is minimal with $s(i)\!\not=\! s'(i)$. Note that $\vert s\wedge s'\vert\!\in\! B_s$.\bigskip

\noindent $\bullet$ We define $O(s)\!\in\!\omega^{\vert s\vert}$: we set $O(s)(i)\! :=\! s_0(i)$.\bigskip

\noindent $\bullet$ We also define a partial order on $\omega^{<\omega}$ as follows:
$$O\sqsubseteq O'\ \Leftrightarrow\ O\! =\! O'\vee\exists i\! <\!\mbox{min}(\vert O\vert ,\vert O'\vert )\ \ 
\big( O\vert i\! =\! O'\vert i\wedge O(i)\! =\! 0\! <\! O'(i)\big).$$

\begin{lem} Let $\vec w\! :=\! (\alpha ,\beta ,\gamma )\!\in\!\Upsilon_1^\infty$,  
$\vec v\! :=\! (\alpha ,a_0,a_1,\underline{a}_0,\underline{a}_1,r)\!\in\!\Theta^\infty$ with 
$\alpha\!\in\!\Borel$ normalized and $\vert (\alpha )_1\vert\! =\! 2$, and $s,s'\!\in\! {\mathfrak T}({\vec w})$ uncompatible.  Assume that $\vec s\!\in\!\bigcap_{i\leq\vert s\wedge s'\vert}\ P_{s_0(i)}(s\vert i)$, 
$\vec t$ is in $\bigcap_{i\leq\vert s\wedge s'\vert}\ P_{s'_0(i)}(s'\vert i)$ and 
$\vec s\ R^{(\eta_{s\vert\vert s\wedge s'\vert})}_{s\vert\vert s\wedge s'\vert}\ \vec t$. Then 
$O(s)\sqsubseteq O(s')$.\end{lem}

\noindent\bf Proof.\rm\ As $s(\vert s\wedge s'\vert )\!\not=\! s'(\vert s\wedge s'\vert )$ and 
$s_1(\vert s\wedge s'\vert )\! =\! s'_1(\vert s\wedge s'\vert )$, we get 
$s_0(\vert s\wedge s'\vert )\!\not=\! s'_0(\vert s\wedge s'\vert )$. Recall the definition of the $P_q(s)$'s. Note the following facts. Assume that $i\!\in\! B_s$ and $\vec s\ R^{(\eta_{s\vert i})}_{s\vert i}\ \vec t$.\smallskip

\noindent - If $s_0(i)\! =\! 0$ and $\vec t\!\in\! P_0(s\vert i)$, then $\vec s\!\in\! P_0(s\vert i)$ too.\smallskip

\noindent - If $s_0(i)\!\geq\! 1$ and $\vec t\!\in\! P_{s_0(i)}(s\vert i)$, then 
$\vec s\!\in\! P_0(s\vert i)\cup P_{s_0(i)}(s\vert i)$.\bigskip

 These facts imply that $s_0(\vert s\wedge s'\vert )\! =\! 0\! <\! s'_0(\vert s\wedge s'\vert )$. Therefore
$O(s)\sqsubseteq O(s')$.$\hfill{\square}$\bigskip

\bf\noindent (G) The ranges\rm\bigskip

 The goal of this paragraph is to defiine the analytic sets $r\big( S(\vec t~)\big)$ that will contain 
$U_{\vec t}$ in the inductive construction of the proof of Theorem 6.9. They will play the role that $\overline{A_0}^{\tau_\xi}\cap A_1$ and $A_0$ played in the proof of Theorem 4.4.1, Conditions 
(4)-(5).\bigskip

\noindent\bf Notation.\rm\ Fix $\vec w\! :=\! (\alpha ,\beta ,\gamma )\!\in\!\Upsilon_1^\infty$,  
$\vec v\! :=\! (\alpha ,a_0,a_1,\underline{a}_0,\underline{a}_1,r)\!\in\!\Theta^\infty$ with 
$\alpha\!\in\!\Borel$ normalized and $\vert (\alpha )_1\vert\! =\! 2$, and 
$s\!\in\! {\mathfrak T}({\vec w})\!\setminus\!\{\emptyset\}$. We set 
$$\begin{array}{ll}
i^s
& \!\!\!\!\! :=\!\left\{\!\!\!\!\!\!\!
\begin{array}{ll}
& \vert s\vert\! -\! 1\mbox{ if }\forall j\! <\!\vert s\vert\ \ s_0(j)\!\geq\! 1\mbox{,}\cr
& \mbox{min}\{ i\! <\!\vert s\vert\mid s_0(i)\! =\! 0\}\mbox{ otherwise,}
\end{array}
\right.\cr\cr
I^s
& \!\!\!\!\! :=\!\left\{\!\!\!\!\!\!\!
\begin{array}{ll}
& \vert s\vert\! -\! 1\mbox{ if }s_0(\vert s\vert\! -\! 1)\!\geq\! 1\mbox{,}\cr
& \mbox{min}\{ i\! <\!\vert s\vert\mid\forall j\!\geq\! i\ \ s_0(j)\! =\! 0\}\mbox{ otherwise.}
\end{array}
\right.
\end{array}$$ 
Note that $i^s\!\leq\! I^s\!\leq\! B_s$. We associate, with each $i^s\!\leq\! i\! <\!\vert s\vert$, 
$\underline{a}^{s,i}_0,\underline{a}^{s,i}_1,r^{s,i}\!\in\!\omega^\omega$. The definition is by induction on 
$i$. We set $\underline{a}^{s,i^s}_\varepsilon\! :=\!\underline{a}_\varepsilon\big( V(s)(i^s)(0),a_0,a_1\big)$, 
$r^{s,i^s}\! :=\! r\big( V(s)(i^s)(0),a_0,a_1\big)\! =\! V(s)(i^s)(5)$. Then
$$\begin{array}{ll}
\underline{a}^{s,i+1}_\varepsilon
& \!\!\!\!\! :=\!\left\{\!\!\!\!\!\!\!
\begin{array}{ll}
& \underline{a}^{s,i}_\varepsilon\mbox{ if }s_0(i\! +\! 1)\!\geq\! 1\mbox{,}\cr
& \underline{a}_\varepsilon\big( V(s)(i\! +\! 1)(0),\underline{a}^{s,i}_0,\underline{a}^{s,i}_1\big)\mbox{ if }
s_0(i\! +\! 1)\!=\! 0\mbox{,}
\end{array}
\right.\cr\cr
r^{s,i+1}
& \!\!\!\!\! :=\!\left\{\!\!\!\!\!\!\!
\begin{array}{ll}
& r^{s,i}\mbox{ if }s_0(i\! +\! 1)\!\geq\! 1\mbox{,}\cr
& r\big( V(s)(i\! +\! 1)(0),\underline{a}^{s,i}_0,\underline{a}^{s,i}_1\big)\mbox{ if }
s_0(i\! +\! 1)\!=\! 0.
\end{array}
\right.
\end{array}$$ 
The $range$ of $s$ is $r(s)\! :=\!\neg {\cal U}_{r^{s,I^s}}$.

\begin{lem} Let $\vec w\! :=\! (\alpha ,\beta ,\gamma )\!\in\!\Upsilon_1^\infty$,  
$\vec v\! :=\! (\alpha ,a_0,a_1,\underline{a}_0,\underline{a}_1,r)\!\in\!\Theta^\infty$ with $\alpha\!\in\!\Borel$ normalized and $\vert (\alpha )_1\vert\! =\! 2$, $s\!\in\! {\mathfrak T}({\vec w})\!\setminus\!\{\emptyset\}$, and $i^s\!\leq\! i\! <\! B_s\! -\! 1$ with $s_0(i)\! =\! 0$. Then $r^{s,i}\! =\! r^{s,i+1}$.\end{lem}

\noindent\bf Proof.\rm\ We may assume that $s_0(i\! +\! 1)\! =\! 0$. Assume first that $i\! =\! i^s$. Then 
$$\begin{array}{ll}
r^{s,i^s}\!\!\!\!
& \! =\! r\big( V(s)(i^s)(0),a_0,a_1\big)\cr
& \! =\! r\Big( {\mathfrak W}\big( V(s)(i^s)\big)(0)(0),\underline{a}_0\big( V(s)(i^s)(0),a_0,a_1\big),
\underline{a}_1\big( V(s)(i^s)(0),a_0,a_1\big)\Big)\cr
& \! =\! r\Big( {\mathfrak W}\big( V(s)(i^s)\big)\big( s_0(i^s)\big)(0),
\underline{a}_0\big( V(s)(i^s)(0),a_0,a_1\big),\underline{a}_1\big( V(s)(i^s)(0),a_0,a_1\big)\Big)\cr
& \! =\! r\Big( V(s)(i^s\! +\! 1)(0),\underline{a}_0\big( V(s)(i^s)(0),a_0,a_1\big),
\underline{a}_1\big( V(s)(i^s)(0),a_0,a_1\big)\Big)\cr
& \! =\! r\big( V(s)(i^s\! +\! 1)(0),\underline{a}_0^{s,i^s},\underline{a}_1^{s,i^s}\big)\cr
& \! =\! r^{s,i^s+1}.
\end{array}$$
The argument is similar if $i\! >\! i^s$.$\hfill{\square}$

\begin{lem} Let $\vec w\! :=\! (\alpha ,\beta ,\gamma )\!\in\!\Upsilon_1^\infty$,  
$\vec v\! :=\! (\alpha ,a_0,a_1,\underline{a}_0,\underline{a}_1,r)\!\in\!\Theta^\infty$ with $\alpha\!\in\!\Borel$ normalized and $\vert (\alpha )_1\vert\! =\! 2$. Then there is $S(\vec\emptyset )\!\in\! {\cal M}_{\vec w}$ with 
$\vec\emptyset\!\in\!\bigcap_{i<B_{S(\vec\emptyset )}}\ P_{(S(\vec\emptyset ))_0(i)}
\big( S(\vec\emptyset )\vert i\big)$ and 
$\neg {\cal U}_r\!\subseteq\! r\big( S(\vec\emptyset )\big)$.\end{lem}
 
\noindent\bf Proof.\rm\ We set $s\! :=\! S(\vec\emptyset )$ for short. We already saw that 
$s\!\in\! {\cal M}_{\vec w}$, $\vec\emptyset\!\in\!\bigcap_{i<B_s}\ P_{s_0(i)}( s\vert i)$, and $s_0(i)\! =\! 0$ for each $i\! <\!\vert s\vert$ after Lemma 6.14. Note that $i^s\! =\! I^s\! =\! 0$. We get 
$$\neg {\cal U}_r\! =\!\neg {\cal U}_{V(s)(0)(5)}\! =\!\neg {\cal U}_{V(s)(i^s)(5)}\! =\!\neg {\cal U}_{r^{s,i^s}}
\! =\!\neg {\cal U}_{r^{s,I^s}}\! =\! r(s).$$ 
This finishes the proof.$\hfill{\square}$\bigskip

 The role of the next objects is to determine if we go to the $A_0$ side or the $A_1$ side in the inductive construction of the proof of Theorem 6.9.\bigskip

\noindent\bf Notation.\rm\ Let $\vec w\! :=\! (\alpha ,\beta ,\gamma )\!\in\!\Upsilon_1^\infty$ with 
$\alpha\!\in\!\Borel$ normalized and $\vert (\alpha )_1\vert\! =\! 2$, and $s\!\in\! {\cal M}_{\vec w}$. We set 
$\varepsilon_s\! :=\! 0$ if $B_s\! <\!\vert s\vert\! -\! 1$, $\varepsilon_s\! :=\! 1$ otherwise, i.e., if 
$B_s\! =\!\vert s\vert\! -\! 1$.

\begin{lem} Let $\vec w\! :=\! (\alpha ,\beta ,\gamma )\!\in\!\Upsilon_1^\infty$,  
$\vec v\! :=\! (\alpha ,a_0,a_1,\underline{a}_0,\underline{a}_1,r)\!\in\!\Theta^\infty$ with 
$\alpha\!\in\!\Borel$ normalized and $\vert (\alpha )_1\vert\! =\! 2$, and $s\!\in\! {\cal M}_{\vec w}$. Then 
$r(s)\!\subseteq\!\neg {\cal U}_{a_{\varepsilon_s}}$.\end{lem}
 
\noindent\bf Proof.\rm\ Note first that 
$\neg {\cal U}_{\underline{a}^{s,i}_\varepsilon}\!\subseteq\!\neg {\cal U}_{a_\varepsilon}$, by induction on $i$ and the 2nd remark after the definition of $\Theta$. This implies that 
$\neg {\cal U}_{r^{s,I^s}}\!\subseteq\!\neg {\cal U}_{r(V(s)(I^s)(0),a_0,a_1)}\! =\!\neg {\cal U}_{V(s)(I^s)(5)}$, by the 4th remark after the definition of $\Theta$. Thus 
$r(s)\! =\!\neg {\cal U}_{r^{s,I^s}}\!\subseteq\!\neg {\cal U}_{V(s)(I^s)(5)}$. Lemma 6.16 implies that 
$\neg {\cal U}_{V(s)(I^s)(5)}\!\subseteq\!\neg {\cal U}_{V(s)(B_s)(5)}$. But 
$V(s)(B_s)(5)\! =\! a_{\varepsilon_s}$, by Lemma 6.15.$\hfill{\square}$\bigskip

 Now we come to the crucial lemma for the claim mentioned earlier.

\begin{lem} Let $\vec w\! :=\! (\alpha ,\beta ,\gamma )\!\in\!\Upsilon_1^\infty$,  
$\vec v\! :=\! (\alpha ,a_0,a_1,\underline{a}_0,\underline{a}_1,r)\!\in\!\Theta^\infty$ with 
$\alpha\!\in\!\Borel$ normalized and $\vert (\alpha )_1\vert\! =\! 2$,  
$s,s'\!\in\! {\mathfrak T}({\vec w})$ with $O(s)\!\not=\! O(s')$ and $O(s)\sqsubseteq O(s')$. Then 
$r(s)\!\subseteq\!\overline{r(s')}^{\tau_{\xi_{s\vert\vert s\wedge s'\vert}}}$.\end{lem}

\noindent\bf Proof.\rm\ We can write $O(s)\! :=\! 0^{k_0}n_0...0^{k_{l-1}}n_{l-1}0^{k_l}$, with 
$l,k_i\!\in\!\omega$, and $n_i\!\geq\! 1$. Similarly, we write 
$O(s')\! :=\! 0^{k'_0}n'_0...0^{k'_{l'-1}}n'_{l'-1}0^{k'_{l'}}$. The assumption implies that $l'\!\geq\! 1$, and the existence of $j\! <\! l'$ with $(k_i,n_i)\! =\! (k'_i,n'_i)$ if $i\! <\! j$ and $k'_j\! <\! k^{}_j$. Lemma 6.14 shows the existence of $k''_{j+1}\!\geq\! 1$ and $s''\!\in\! {\cal M}_{\vec w}$ with 
$O(s'')\! =\! 0^{k'_0}n'_0...0^{k'_{j-1}}n'_{j-1}0^{k'_j}n'_j0^{k''_{j+1}}$ if $j\! <\! l'\! -\! 1$. If $j\! =\! l'\! -\! 1$, then we set $s''\! :=\! s'$.

 Note that $O(s)\!\not=\! O(s'')$, $O(s)\sqsubseteq O(s'')$, and $O(s'')\sqsubseteq O(s')$. Moreover, $O(s'')\!\not=\! O(s')$  and $\vert s\wedge s'\vert\! =\!\vert s\wedge s''\vert\! <\!\vert s'\wedge s''\vert$ if 
$j\! <\! l'\! -\! 1$. It is enough to prove that 
$r(s)\!\subseteq\!\overline{r(s'')}^{\tau_{\xi_{s\vert\vert s\wedge s''\vert}}}$. This means that we may assume that $(k_i,n_i)\! =\! (k'_i,n'_i)$ if $i\! <\! l'\! -\! 1$ and $k'_{l'-1}\! <\! k^{}_{l'-1}$. This implies that 
$I^{s'}\!\geq\! 1$, $\vert s\wedge s'\vert\! =\! I^{s'}\! -\! 1$, $s\vert (I^{s'}\! -\! 1)\! =\! s'\vert (I^{s'}\! -\! 1)$, 
$s_0(I^{s'}\! -\! 1)\! =\! 0\! <\! s'_0(I^{s'}\! -\! 1)$ and $i^s\!\leq\! I^{s'}\! -\! 1$.\bigskip

\noindent\bf Case 1.\rm\ $i^s\! =\! I^s$ and $i^{s'}\! =\! I^{s'}$.\bigskip

 Note that $r(s)\! =\!\neg {\cal U}_{r^{s,I^s}}\! =\!\neg {\cal U}_{r^{s,i^s}}\! =\!
\neg {\cal U}_{V(s)(i^s)(5)}\! =\!\neg {\cal U}_{V(s')(I^s)(5)}$. Lemma 6.16 implies that 
$$r(s)\! =\!\neg {\cal U}_{V(s')(I^s)(5)}\!\subseteq\!\neg {\cal U}_{V(s')(I^{s'}-1)(5)}
\!\subseteq\!\overline{\neg {\cal U}_{V(s')(I^{s'})(5)}}^{\tau_{\xi_{s'\vert (I^{s'}-1)}}}\! =\!
\overline{r(s')}^{\tau_{\xi_{s\vert\vert s\wedge s'\vert}}}.$$
\bf Case 2.\rm\ $i^s\! =\! I^s$ and $i^{s'}\! <\! I^{s'}$.\bigskip

 Note that $i^s\! =\! i^{s'}\! <\! I^{s'}\! -\! 1$. Lemma 6.18 implies that 
$r(s)\! =\!\neg {\cal U}_{r^{s,I^s}}\! =\!\neg {\cal U}_{r^{s,I^{s'}-1}}$. Thus
$$\begin{array}{ll}
r(s)\!\!\!\!
& \! =\!\neg {\cal U}_{r(V(s)(I^{s'}-1)(0),\underline{a}_0^{s,I^{s'}-2},\underline{a}_1^{s,I^{s'}-2})}\cr
& \! =\!\neg {\cal U}_{r(V(s')(I^{s'}-1)(0),\underline{a}_0^{s',I^{s'}-2},\underline{a}_1^{s',I^{s'}-2})}\cr
& \! =\!\neg {\cal U}_{r(V(s')(I^{s'}-1)(0),\underline{a}_0^{s',I^{s'}-1},\underline{a}_1^{s',I^{s'}-1})}\cr
& \!\subseteq\!\overline{\neg {\cal U}_{r(V(s')(I^{s'})(0),\underline{a}_0^{s',I^{s'}-1},
\underline{a}_1^{s',I^{s'}-1})}}^{\tau_{\xi_{s'\vert (I^{s'}-1)}}}\cr
& \! =\!\overline{r(s')}^{\tau_{\xi_{s\vert\vert s\wedge s'\vert}}}\mbox{,}
\end{array}$$
by Lemma 6.16.\bigskip

\noindent\bf Case 3.\rm\ $i^s\! <\! I^s\! <\! I^{s'}$.\bigskip

 We argue as in Case 2.\bigskip

\noindent\bf Case 4.\rm\ $i^s\! <\! I^s$ and $I^{s'}\!\leq\! I^s$, which implies that $I^{s'}\! <\! I^s$.\bigskip

 The 5th remark after the definition of $\Upsilon$ gives $\varepsilon\!\in\! 2$ with 
$r(s)\! =\!\neg {\cal U}_{r^{s,I^s}}\!\subseteq\!\neg {\cal U}_{\underline{a}_\varepsilon^{s,I^{s}-1}}$. Thus 
$r(s)\!\subseteq\!\neg {\cal U}_{\underline{a}_\varepsilon^{s,I^{s}-1}}\!\subseteq\! ...\!\subseteq\!
\neg {\cal U}_{\underline{a}_\varepsilon^{s,I^{s'}-1}}$. If $I^{s'}\!\geq\! 2$, then we get
$$\begin{array}{ll}
\neg {\cal U}_{\underline{a}_\varepsilon^{s,I^{s'}-1}}\!\!\!\!
& \! =\!\neg {\cal U}_{a_\varepsilon (V(s')(I^{s'}-1)(0),\underline{a}_0^{s',I^{s'}-2},\underline{a}_1^{s',I^{s'}-2})}\cr
& \!\subseteq\!\overline{
\neg {\cal U}_{r(V(s')(I^{s'})(0),\underline{a}_0^{s',I^{s'}-2},\underline{a}_1^{s',I^{s'}-2})}}^
{\tau_{\xi_{s\vert\vert s\wedge s'\vert}}}\cr
& \! =\!\overline{
\neg {\cal U}_{r(V(s')(I^{s'})(0),\underline{a}_0^{s',I^{s'}-1},\underline{a}_1^{s',I^{s'}-1})}}^
{\tau_{\xi_{s\vert\vert s\wedge s'\vert}}}\cr
& \! =\!\overline{r(s')}^{\tau_{\xi_{s\vert\vert s\wedge s'\vert}}}.
\end{array}$$

 Otherwise, we get $I^{s'}\!=\! 1$, $i^s\! =\! 0$, $i^{s'}\! =\! I^{s'}$ and
$$\neg {\cal U}_{\underline{a}_\varepsilon^{s,0}}\! =\!\neg {\cal U}_{a_\varepsilon (V(s')(0)(0),a_0,a_1)}
\!\subseteq\!\overline{\neg {\cal U}_{r(V(s')(1)(0),a_0,a_1)}}^{\tau_{\xi_{s\vert\vert s\wedge s'\vert}}}\! =\!\overline{r(s')}^{\tau_{\xi_{s\vert\vert s\wedge s'\vert}}}.$$
This finishes the proof.$\hfill{\square}$\bigskip

\bf\noindent (H) The maximal sequences\rm\bigskip

 We now associate a maximal sequence to a couple $(\vec\beta ,\vec w)$ with 
$\vec\beta\!\in\!\lceil T_d\rceil$. It is build in a way similar to that of the $s(\vec t~)$'s, but for infinite sequences instead of finite ones.\bigskip

\noindent $\bullet$ Let $\vec w\! :=\! (\alpha ,\beta ,\gamma )\!\in\!\Upsilon_1^\infty$ with 
$\alpha\!\in\!\Borel$ normalized and $\vert (\alpha )_1\vert\! =\! 2$, and $\vec\beta\!\in\!\lceil T_d\rceil$. We will define $s(\vec\beta ,\vec w)\!\in\! {\cal M}_{\vec w}$. Recall the definition of $\tilde {\cal P}^{\vec w}_p$. We set, for $s\!\in\! {\cal M}_{\vec w}$ and $i\!\in\! B_s$, 
$$E^s_i\! :=\!\left\{\!\!\!\!\!\!\!\!
\begin{array}{ll}
& \bigcap_{p\geq 1}\ \tilde {\cal P}^{s(i)(1)}_p\mbox{ if }s(i)(0)\! =\! 0
\mbox{,}\cr
& \neg\tilde {\cal P}^{s(i)(1)}_{s(i)(0)}\mbox{ if }s(i)(0)\!\geq\! 1.
\end{array}
\right.$$

 We define $s(\vec\beta ,\vec w)$ in such a way that 
$j_d(\vec\beta )\!\in\!\bigcap_{i\in B_{s(\vec\beta ,\vec w)}}\ E^{s(\vec\beta ,\vec w)}_i$. Let $\xi$ be an ordinal such that $\vec w\!\in\!\Upsilon_1^\xi\!\setminus\!\Upsilon_1^{<\xi}$. The definition of 
$s(\vec\beta ,\vec w)$ is by induction on $\xi$.\bigskip

\noindent\bf Case 1.\rm\ $\vert (\alpha )_1\vert\! =\! 0$.\bigskip

 We set $s(\vec\beta ,\vec w)\! :=<\! (0,\vec w)\! >$.\bigskip

\noindent\bf Case 2.\rm\ $\vert (\alpha )_1\vert\! =\! 1$.\bigskip

 We set $s(\vec\beta ,\vec w)\! :=\! (0,\vec w)^\frown s\big(\vec\beta ,{\mathfrak W}^1_1(\vec w)\big)$.\bigskip

\noindent\bf Case 3.\rm\ $\vert (\alpha )_1\vert\! =\! 2$.\bigskip

 We set $s(\vec\beta ,\vec w)\! :=\!\left\{\!\!\!\!\!\!\!\!
\begin{array}{ll}
& (0,\vec w)^\frown s\big(\vec\beta ,{\mathfrak W}^1_1(\vec w)(0)\big)\mbox{ if }
j_d(\vec\beta )\!\in\!\bigcap_{p\geq 1}\ \tilde {\cal P}^{\vec w}_p\mbox{,}\cr
& (p,\vec w)^\frown s\big(\vec\beta ,{\mathfrak W}^1_1(\vec w)(p)\big)\mbox{ if }
j_d(\vec\beta )\!\notin\!\tilde {\cal P}^{\vec w}_p\wedge p\!\geq\! 1.
\end{array}
\right.$\bigskip

\noindent $\bullet$ We set $(\vec\beta\vert j_k)_{k\in\omega}\! :=\! (\Pi_{s(\vec\beta ,\vec w)\vert 0}\circ ...\circ
\Pi_{s(\vec\beta ,\vec w)\vert (B_{s(\vec\beta ,\vec w)}-1)})^{-1}\big( h(\vec\beta )\big)$.\bigskip

 Recall the definition of $\varepsilon_s$ before Lemma 6.20.

\begin{lem} Let $\vec w\! :=\! (\alpha ,\beta ,\gamma )\!\in\!\Upsilon_1^\infty$ with 
$\alpha\!\in\!\Borel$ normalized and $\vert (\alpha )_1\vert\! =\! 2$, and 
$\vec\beta\!\in\!\lceil T_d\rceil$.\smallskip

\noindent (a) There is $k_0\!\in\!\omega$ such that $\vec\beta\vert j_k\!\in\! 
\bigcap_{i<B_{s(\vec\beta ,\vec w)}}\ P_{s(\vec\beta ,\vec w)(i)(0)}\big( s(\vec\beta ,\vec w)\vert i\big)$ if  $k\!\geq\! k_0$. In this case, the sequence $s(\vec\beta\vert j_k)$ given by Proposition 6.13 is 
$s(\vec\beta ,\vec w)\vert B_{s(\vec\beta ,\vec w)}$, and is not extensible.\smallskip

\noindent (b) We have $j_d(\vec\beta )\!\in\! C^{\omega^\omega}_{\gamma}$ if and only if 
$\varepsilon_{s(\vec\beta ,\vec w)}\! =\! 0$.\end{lem}

\noindent\bf Proof.\rm\ We set $s\! :=\! s(\vec\beta ,\vec w)$ for simplicity.\medskip

\noindent (a) To define $k_0$, we will define, for $i\! <\! B_s$, $k^i_0\!\in\!\omega$ and we will set 
$k_0\! :=\!\mbox{max}\{ k^i_0\mid i\! <\! B_s\}$. To do this, we set 
$(\vec\beta\vert j^i_k)_{k\in\omega}\! :=\! (\Pi_{s\vert 0}\circ ...\circ\Pi_{s\vert i})^{-1}\big( h(\vec\beta )\big)$, so that $(\vec\beta\vert j^{i+1}_k)_{k\in\omega}$ is a subsequence of $(\vec\beta\vert j^i_k)_{k\in\omega}$ if $i\! <\! B_s\! -\! 1$.

\vfill\eject

 By the choice of the $E^s_i$'s we get, for $i\! <\! B_s$,
$$\begin{array}{ll}
\ \ \ \ \ \ \ \ \ \ h(\vec\beta)
& \!\!\!\!\in\!\left\{\!\!\!\!\!\!\!\!
\begin{array}{ll}
& \bigcap_{p\geq 1}\ {\cal P}^{s_1(i)}_p\mbox{ if }s_0(i)\! =\! 0\mbox{,}\cr
& \neg {\cal P}^{s_1(i)}_{s_0(i)}\mbox{ if }s_0(i)\!\geq\! 1\mbox{,}
\end{array}
\right.\cr
(\vec\beta\vert j^i_k)_{k\in\omega}
& \!\!\!\!\in\!\left\{\!\!\!\!\!\!\!\!
\begin{array}{ll}
& \bigcap_{p\geq 1}\ (\Pi_{s\vert 0}\circ ...\circ\Pi_{s\vert i})^{-1}({\cal P}^{s_1(i)}_p)\mbox{ if }
s_0(i)\! =\! 0\mbox{,}\cr
& \neg (\Pi_{s\vert 0}\circ ...\circ\Pi_{s\vert i})^{-1}({\cal P}^{s_1(i)}_{s_0(i)})\mbox{ if }s_0(i)\!\geq\! 1.
\end{array}
\right.
\end{array}$$
Note the existence of ${\cal B}^i_p$ in $(\Pi_{s\vert 0}\circ ...\circ\Pi_{s\vert i})^{-1}({\cal P}^{s_1(i)}_p)$ such that $\vec\beta\vert j^i_k\!\in\! {\cal B}^i_p$ if $s_0(i)\! =\! 0$, $k\!\in\!\omega$ and $p\!\geq\! 1$. If 
$s_0(i)\!\geq\! 1$ and $p\!\in\!\omega\!\setminus\!\{ 0,s_0(i)\}$, then 
$(\vec\beta\vert j^i_k)_{k\in\omega}\!\in\! (\Pi_{s\vert 0}\circ ...\circ\Pi_{s\vert i})^{-1}({\cal P}^{s_1(i)}_p)$ since ${\cal P}^{s_1(i)}_p\cup {\cal P}^{s_1(i)}_{s_0(i)}\! =\! [\subseteq ]$. This implies the existence of 
${\cal B}^i_p\!\in\! (\Pi_{s\vert 0}\circ ...\circ\Pi_{s\vert i})^{-1}({\cal P}^{s_1(i)}_p)$ such that 
$\vec\beta\vert j^i_k\!\in\! {\cal B}^i_p$ if $k\!\in\!\omega$. As 
$(\Pi_{s\vert 0}\circ ...\circ\Pi_{s\vert i})^{-1}({\cal P}^{s_1(i)}_{s_0(i)})\!\in\! 
{\bf\Pi}^0_1([R^{(\eta_{s\vert i})}_{s\vert i}])$, there is $k^i_0\!\geq\! 1$ such that 
$\vec\beta\vert j^i_k\!\notin\! {\cal B}^i_{s_0(i)}$ if $s_0(i)\!\geq\! 1$, 
${\cal B}^i_{s_0(i)}\!\in\! (\Pi_{s\vert 0}\circ ...\circ\Pi_{s\vert i})^{-1}({\cal P}^{s_1(i)}_{s_0(i)})$ and 
$k\!\geq\! k^i_0$. This defines $k^i_0$ and $k_0$. It remains to check that 
$\vec\beta\vert j_k\!\in\! P_{s(i)(0)}(s\vert i)$ if $i\! <\! B_s$ and $k\!\geq\! k_0$. This comes from the fact that $j_k\! =\! j^{B_s-1}_k\! =\! j^i_{K(k)}$ for some $K(k)\!\geq\! k\!\geq\! k_0\!\geq\! k^i_0$. The last assertion comes from the construction of $s(\vec t~)$.\bigskip

\noindent (b) We define, for $i\! <\!\vert s\vert$, $\varepsilon^i_{s}\!\in\! 2$. The definition is by induction on 
$i$. We first set $\varepsilon^{0}_{s}\! :=\! 1$. Then $\varepsilon^{i+1}_{s}\! :=\! 0$ if 
$\vert s\vert\! -\! i\! -\! 2\!\notin\! B_{s}$, $\varepsilon^{i+1}_{s}\! :=\!\varepsilon^{i}_{s}$ otherwise. Note that 
$\varepsilon_s\! =\! \varepsilon^{\vert s\vert -1}_s$ ($\varepsilon_s$ is defined before Lemma 6.20). We have to see that $j_d(\vec\beta )$ is in $C^{\omega^\omega}_{s_1(0)(2)}$ is equivalent to 
$\varepsilon^{\vert s\vert -1}_s\! =\! 0$. We prove the following stronger fact: 
$j_d(\vec\beta )\!\in\! C^{\omega^\omega}_{s_1(\vert s\vert -i-1)(2)}$ is equivalent to 
$\varepsilon^i_s\! =\! 0$ if $i\! <\!\vert s\vert$. Here again we argue by induction on $i$. The result is clear for $i\! =\! 0$ since $C^{\omega^\omega}_{s_1(\vert s\vert -1)(2)}\! =\!\emptyset$. So assume that the result is true for $i\! <\!\vert s\vert\! -\! 1$.\bigskip

 If $\vert s\vert\! -\! i\! -\! 2\!\notin\! B_s$, then we are done since 
$\varepsilon^{i+1}_s\! =\! 1\! -\!\varepsilon^i_s$ and 
$C^{\omega^\omega}_{s_1(\vert s\vert -i-2)(2)}\! =\!\neg C^{\omega^\omega}_{s_1(\vert s\vert -i-1)(2)}$. If 
$\vert s\vert\! -\! i\! -\! 2\!\in\! B_s$, then $\varepsilon^{i+1}_s\! =\!\varepsilon^i_s$ and\bigskip

\leftline{$C^{\omega^\omega}_{s_1(\vert s\vert -i-2)(2)}\! =\!\bigcup_{p\geq 1}\ 
(C^{\omega^\omega}_{(({\mathfrak W}^1_0(s_1(\vert s\vert -i-2)))_p)_0}\!\setminus\! 
C^{\omega^\omega}_{(({\mathfrak W}^1_0(s_1(\vert s\vert -i-2)))_p)_1})\cup$}\smallskip

\rightline{$(C^{\omega^\omega}_{({\mathfrak W}^1_0(s_1(\vert s\vert -i-2)))_0}\cap\bigcap_{p\geq 1}\ 
C^{\omega^\omega}_{(({\mathfrak W}^1_0(s_1(\vert s\vert -i-2)))_p)_1}).$}\bigskip

 If $s_0(\vert s\vert\! -\! i\! -\! 2)\! =\! 0$, then 
$j_d(\vec\beta )\!\in\!\bigcap_{p\geq 1}\ \tilde {\cal P}^{s_1(\vert s\vert -i-2)}_p\! =\!\bigcap_{p\geq 1}\ 
C^{\omega^\omega}_{(({\mathfrak W}^1_0(s_1(\vert s\vert -i-2)))_p)_1}$. We can say that 
$j_d(\vec\beta )\!\in\! C^{\omega^\omega}_{s_1(\vert s\vert -i-2)(2)}$ is equivalent to 
$j_d(\vec\beta )\!\in\! C^{\omega^\omega}_{({\mathfrak W}^1_0(s_1(\vert s\vert -i-2)))_0}\! =\! 
C^{\omega^\omega}_{s_1(\vert s\vert -i-1)(2)}$, and we are done by induction assumption. We argue similarly when $s_0(\vert s\vert\! -\! i\! -\! 2)\!\geq\! 1$.$\hfill{\square}$\bigskip

\noindent\bf Remark.\rm\ Recall the definition of an extensible sequence at the beginning of the construction of the resolution families. If $s$ is not extensible, then $s$ admits a unique extension $M(s)$ in ${\cal M}_{\vec w}$. In particular, in Lemma 6.22.(a), 
$M\big( s(\vec\beta\vert j_k)\big)\! =\! s(\vec\beta ,\vec w)\! =\! S(\vec\beta\vert j_k)$. In Lemma 6.19, 
$s(\vec\emptyset )\! =\! s\vert B_s$ is not extensible and 
$M\big( s(\vec\emptyset )\big)\! =\! S(\vec\emptyset )$.\bigskip

\noindent\bf Notation.\rm\ Recall the construction of the resolution families, and also the proof of Theorem 4.4.5, especially the definition of $\eta (\vec t~)$. If $\theta_{s}$ is a limit ordinal, then we consider some ordinals $\eta_{s}(\vec t~)$'s, as in the proof of Theorem 4.4.5. We set 
$\rho (s,\vec s~)\! :=\!\left\{\!\!\!\!\!\!
\begin{array}{ll}
& \eta_{s}\mbox{ if }\theta_{s}\ \mbox{is a successor ordinal,}\cr
& \eta_{s}(\vec s~)\mbox{ if }\theta_{s}\ \mbox{is a limit ordinal.}
\end{array}
\right.$

\vfill\eject

 The next lemma is the final preparation for the claim mentioned earlier.

\begin{lem} Let $\vec w\! :=\! (\alpha ,\beta ,\gamma )\!\in\!\Upsilon_1^\infty$ with 
$\alpha\!\in\!\Borel$ normalized and $\vert (\alpha )_1\vert\! =\! 2$, $s\!\in\! {\mathfrak T}({\vec w})$, and 
$i\! <\! B_s$. Then $\big(\Sigma_{i'\leq i}\ \rho (s\vert i',\overrightarrow{tm})\big)\! +\! 1\!\leq\!\xi_{s\vert i}$.\end{lem}

\noindent\bf Proof.\rm\ We argue by induction on $i$. Note first that 
$\rho (s\vert 0,\overrightarrow{tm})\! +\! 1\!\leq\!\theta_{s\vert 0}\! =\!\xi_{s\vert 0}$. Then, inductively, 
$$\begin{array}{ll}
\big(\Sigma_{i'\leq i+1}\ \rho (s\vert i',\overrightarrow{tm})\big)\! +\! 1\!\!\!
& \!\leq\!\big(\Sigma_{i'\leq i}\ \rho (s\vert i',\overrightarrow{tm})\big)\! +\!\theta_{s\vert (i+1)}\cr
& \!\leq\!\big(\Sigma_{i'\leq i}\ \rho (s\vert i',\overrightarrow{tm})\big)\! +\! 1\! +\! 
(\xi_{s\vert (i+1)}\! -\!\xi_{s\vert i})\cr
& \!\leq\!\xi_{s\vert i}\! +\! (\xi_{s\vert (i+1)}\! -\!\xi_{s\vert i})\cr
& \!\leq\!\xi_{s\vert (i+1)}
\end{array}$$
This finishes the proof.$\hfill{\square}$\bigskip

\noindent\bf Proof of Theorem 6.9.\rm\ Let $\xi$ be an ordinal with 
$\vec w\! :=\! (\alpha ,\beta ,\gamma )\!\in\!\Upsilon_1^\xi$. We argue by induction on $\xi$. So assume that 
$\vec w\!\in\!\Upsilon_1^\xi\!\setminus\!\Upsilon_1^{<\xi}$.\bigskip

\noindent\bf Case 1.\rm\ $\vert (\alpha )_1\vert\! =\! 0$.\bigskip

 Lemma 6.5 implies that $C^{\omega^\omega}_\gamma\!\in\! {\bf\Gamma}_{c(\alpha )}\! =\! 
{\bf\Gamma}_{0^\infty}\! =\!\{\emptyset\}$, so that $S\! =\!\emptyset$. We also have 
$r\! =\! a_1$. Assume that (a) does not hold. Then $A_1\!\not=\!\emptyset$, so it contains some 
$\vec\alpha$. We just have to set $f_i(\beta_i)\! :=\!\alpha_i$.\bigskip

\noindent\bf Case 2.\rm\ $\vert (\alpha )_1\vert\! =\! 1$.\bigskip

 As $\vec w\!\in\!\Upsilon_1^\xi$ we get $\gamma'\!\in\!\omega^\omega$ with 
$(<(\alpha )_{2+j}>,\beta^*,\gamma')\!\in\!\Upsilon_1^{<\xi}$ and 
$C^{\omega^\omega}_\gamma\! =\!\neg C^{\omega^\omega}_{\gamma'}$ (see the definition of 
$\Upsilon_1$). As $\alpha$ is normalized, $C^{\omega^\omega}_{\gamma'}\! =\!\emptyset$, so that 
$S\! =\!\lceil T_d\rceil$. We also have $r\! =\! a_0$. Assume that (a) does not hold. Then 
$A_0\!\not=\!\emptyset$, and we argue as in Case 1.\bigskip

\noindent\bf Case 3.\rm\ $\vert (\alpha )_1\vert\! =\! 2$.\bigskip

 Assume that (a) does not hold. As for Theorems 4.4.1 and 4.4.5 we construct 
$(\alpha^i_{s})_{i\in d,s\in\Pi_i''T_d}$, $(O^i_{s})_{i\leq\vert s\vert ,i\in d,s\in\Pi_i''T_d}$, 
$(U_{\vec s  })_{\vec s \in T_d}$. We want these objects to satisfy the following conditions.
$$\begin{array}{ll}
& \!\!\!\! (1)\ \alpha^i_{s}\!\in\! O^i_{s}\!\subseteq\!\Omega_{\omega^\omega}\wedge 
(\alpha^i_{s_i })_{i\in d}\!\in\! U_{\vec s  }\!\subseteq\!
\Omega_{(\omega^\omega )^d}\mbox{,}\cr\cr 
& \!\!\!\! (2)\ O^i_{sq}\!\subseteq\! O^i_{s}\mbox{,}\cr\cr
& \!\!\!\! (3)\ \hbox{\rm diam}_{d_{\omega^\omega}}(O^i_{s})\!\leq\! 2^{-\vert s\vert}\wedge
\mbox{diam}_{d_{(\omega^\omega )^d}}(U_{\vec s  })\!\leq\! 2^{-\vert\vec s\vert}\mbox{,}\cr\cr
& \!\!\!\! (4)\ \vec t\!\in\! T_d~\Rightarrow ~U_{\vec t}\!\subseteq\! r\big( S(\vec t~)\big)\mbox{,}
\cr\cr 
& \!\!\!\! (5)\ \left(\!\!\!\!\!\!
\begin{array}{lcl}
& \vec s,\vec t \!\in\!\bigcap_{i'<i,\eta_{s\vert i'}\geq 1}\ P_{s_0(i')}(s\vert i')\cr 
& 1\!\leq\!\rho\!\leq\!\rho (s\vert i,\vec s)\cr
& \vec s ~R_{s\vert i}^{(\rho )}~\vec t
\end{array}\right)
\ \Rightarrow\ U_{\vec t}\!\subseteq\!\overline{U_{\vec s}}
^{\tau_{(\Sigma_{i'<i}\ \rho (s\vert i',\vec s))+\rho}}\mbox{,}\cr\cr
& \!\!\!\! (6)\ \left(\vec s\!\in\!\bigcap_{i<\vert s(\vec t~)\vert}\ P_{s(\vec t~)(i)(0)}\big( s(\vec t~)\vert i\big)
\wedge\vec s ~R^{(\eta_{s(\vec t~)^-})}_{s(\vec t~)^-}~\vec t~\right)\Rightarrow ~
U_{\vec t}\!\subseteq\! U_{\vec s}.
\end{array}$$
$\bullet$ Let us prove that this construction is sufficient to get the theorem.\bigskip

\noindent - Fix $\vec\beta\!\in\! \lceil T_d\rceil $. Lemma 6.22 gives $k_0\!\in\!\omega$ such that 
$\vec\beta\vert j_k\!\in\! 
\bigcap_{i<B_{s(\vec\beta ,\vec w)}}\ P_{s(\vec\beta ,\vec w)(i)(0)}\big( s(\vec\beta ,\vec w)\vert i\big)$ for each $k\!\geq\! k_0$. Proposition 6.13 gives $l$ and $s(\vec\beta\vert j_k)\!\in\! {\mathfrak T}({\vec w})$ with 
$\vec\beta\vert j_k\!\in\!\bigcap_{i<l}\ P_{s(\vec\beta\vert j_k)(i)(0)}\big( s(\vec\beta\vert j_k)\vert i\big)$, and Lemma 6.22.(a) implies that 
$s(\vec\beta\vert j_k)\! =\! s(\vec\beta ,\vec w)\vert B_{s(\vec\beta ,\vec w)}$. This implies that $(U_{\vec\beta\vert j_k})_{k\geq k_0}$ is non-increasing since 
$\vec\beta\vert j_k~R^{(\eta_{s(\beta ,\vec w)\vert (B_{s(\beta ,\vec w)}-1)})}_
{s(\beta ,\vec w)\vert (B_{s(\beta ,\vec w)}-1)}~\vec\beta\vert j_{k+1}$ for each integer $k$, by Condition (6). As in the proof of Theorem 4.4.1 we define ${\cal F}(\vec\beta )$ and $f_i$ continuous with 
${\cal F}(\vec\beta )\! =\! (\Pi_{i\in d}\ f_i)(\vec\beta )$. The inclusions 
$$S\!\subseteq\! (\Pi_{i\in d}\ f_i)^{-1}(A_0)$$ 
and $\lceil T_d\rceil \!\setminus\! S\!\subseteq\! (\Pi_{i\in d}\ f_i)^{-1}(A_1)$ hold, by Lemmas 6.20 and 6.22, since $r\big( s(\beta ,\vec w)\big)\!\subseteq\! A_{\varepsilon_{s(\beta ,\vec w)}}$.\bigskip

\noindent $\bullet$ So let us prove that the construction is possible.\bigskip

\noindent - As $\neg {\cal U}_r$ is nonempty  and $\Ana$, we can choose 
$(\alpha^i_{\emptyset})_{i\in d}\!\in\!\neg {\cal U}_r\cap\Omega_{(\omega^\omega )^d}$. Then we choose a $\Ana$ subset $U_{\vec\emptyset}$ of $(\omega^\omega )^d$, with $d_{(\omega^\omega )^d}$-diameter at most $1$, such that $(\alpha^i_{\emptyset})_{i\in d}\!\in\! U_{\vec\emptyset}\!\subseteq\! 
\neg {\cal U}_r\cap\Omega_{(\omega^\omega )^d}$. We choose a $\Ana$ subset $O^0_{\emptyset}$ of 
$\Omega_{\omega^\omega}$, with $d_{\omega^\omega}$-diameter at most $1$, with 
$\alpha^0_{\emptyset}\!\in\! O^0_{\emptyset}\!\subseteq\!\Omega_{\omega^\omega}$, which is possible since $\Omega_{(\omega^\omega )^d}\!\subseteq\!\Omega_{\omega^\omega}^d$. Assume that 
$(\alpha^i_{s})_{|s|\leq l}$, $(O^i_{s})_{|s|\leq l}$ and $(U_{\vec s  })_{|s_0|\leq l}$ satisfying conditions (1)-(6) have been constructed, which is the case for $l\! =\! 0$ by Lemma 6.19.\bigskip

\noindent - Let $\overrightarrow{tm}\!\in\!T_d\cap (d^{l+1})^d$. We define $X_i\! :=\! O^i_{t_i}$ if $i\!\leq\! l$, and $\omega^\omega$ if $i\! >\! l$.\bigskip

\noindent\bf Claim.\it\ Assume that $s\!\in\! {\mathfrak T}({\vec w})$, $i\! <\! B_s$, 
$\overrightarrow{tm}^{\eta_{s\vert i}}_{s\vert i},\overrightarrow{tm}\!\in\!\bigcap_{i'<i}\ P_{s_0(i')}(s\vert i')$, and $i_0\!\leq\! i$ is minimal with $\eta_{s\vert i_0}\!\geq\! 1$.\smallskip

\noindent (a) The set\smallskip

\leftline{$U_{\overrightarrow{tm}^{\rho (s\vert i,\overrightarrow{tm})}_{s\vert i}}
\cap\bigcap_{1\leq\rho <\rho (s\vert i,\overrightarrow{tm})}\ 
\overline{U_{\overrightarrow{tm}^{\rho}_{s\vert i}}}^
{\tau_{(\Sigma_{i'<i}\ \rho (s\vert i',\overrightarrow{tm}))+\rho}}$}\smallskip

\rightline{$\cap\bigcap_{i'<i}\ \bigcap_{1\leq\rho\leq\rho (s\vert i',\overrightarrow{tm})}\ 
\overline{U_{\overrightarrow{tm}^{\rho}_{s\vert i'}}}
^{\tau_{(\Sigma_{i''<i'}\ \rho (s\vert i'',\overrightarrow{tm}))+\rho}}\cap (\Pi_{i\in d}\ X_i)$}\smallskip

\noindent is $\tau_{1}$-dense in 
$\overline{U_{\overrightarrow{tm}^1_{s\vert i_0}}}^{\tau_{1}}\cap (\Pi_{i\in d}\ X_i)$.\smallskip

\noindent (b) Assume moreover that $s'\!\in\! {\mathfrak T}({\vec w})$, $s$ and $s'$ are uncompatible, 
$i\! :=\!\vert s\wedge s'\vert$, $\overrightarrow{tm}\!\in\! P_{s'_0(i)}(s'\vert i)$, and 
$\overrightarrow{tm}^{\eta_{s\vert i}}_{s\vert i}\!\in\! P_{s_0(i)}(s\vert i)$. Then  
$$r\big( S(\overrightarrow{tm})\big)\cap\bigcap_{i'\leq i}\ 
\bigcap_{1\leq\rho\leq\rho (s\vert i',\overrightarrow{tm})}\ 
\overline{U_{\overrightarrow{tm}^{\rho}_{s\vert i'}}}
^{\tau_{(\Sigma_{i''<i'}\ \rho (s\vert i'',\overrightarrow{tm}))+\rho}}\cap (\Pi_{i\in d}\ X_i)$$ 
is $\tau_{1}$-dense in $\overline{U_{\overrightarrow{tm}^1_{s\vert i_0}}}^{\tau_{1}}\cap (\Pi_{i\in d}\ X_i)$.\rm\bigskip

\noindent (a) Assume first that $i_0\! =\! 0$. Note that 
$\overrightarrow{tm}^{\rho +1}_{\emptyset}\ R^{(\rho +1)}_{\emptyset}\ 
\overrightarrow{tm}^{\rho}_{\emptyset}\ R^{(\rho )}_{\emptyset}\ \overrightarrow{tm}$ if 
$1\!\leq\!\rho\! <\!\rho (\emptyset,\overrightarrow{tm})$, by Lemma 4.3.2. As in the proof of Claim 2 in Theorem 4.4.5, this implies that $U_{\overrightarrow{tm}^{\rho}_{\emptyset}}\!\subseteq\!
\overline{U_{\overrightarrow{tm}^{\rho +1}_{\emptyset}}}^{\tau_{\rho +1}}$. By assumption, 
$\overrightarrow{tm}^{\eta_{s\vert i}}_{s\vert i},\overrightarrow{tm}\!\in\!\bigcap_{i'<i}\ P_{s_0(i')}(s\vert i')$. 

\vfill\eject

 Note that $\overrightarrow{tm}^{\rho}_{s\vert (i'+1)}\!\in\! P_{s_0(i'')}(s\vert i'')$ if $i''\!\leq\! i'\! <\! i$ and 
$\rho\!\leq\!\eta_{s\vert (i'+1)}$. Indeed, this comes from the fact that 
$\overrightarrow{tm}^{\eta_{s\vert i}}_{s\vert i}\ R^{(\eta_{s\vert i''})}_{s\vert i''}\ 
\overrightarrow{tm}^{\rho}_{s\vert (i'+1)}\ R^{(\eta_{s\vert i''})}_{s\vert i''}\ \overrightarrow{tm}$. As in the proof of Claim 2 in Theorem 4.4.5 again, this implies that 
$U_{\overrightarrow{tm}^{\rho}_{s\vert (i'+1)}}
\!\subseteq\!\overline{U_{\overrightarrow{tm}^{\rho +1}_{s\vert (i'+1)}}}^
{\tau_{(\Sigma_{i''<i'+1}\ \rho (s\vert i'',\overrightarrow{tm}))+\rho +1}}$ 
if $\rho\! <\!\rho (s\vert (i'+1),\overrightarrow{tm})$. Note that 
$\overrightarrow{tm}^{0}_{s\vert (i'+1)}\! =\!\overrightarrow{tm}^{\eta_{s\vert i'}}_{s\vert i'}
\! =\!\overrightarrow{tm}^{\rho (s\vert i',\overrightarrow{tm})}_{s\vert i'}$. This implies the result. We argue similarly if $i_0\! >\! 0$.\bigskip

\noindent (b) By (a) and Lemma 6.22, it is enough to see that 
$U\! :=\! U_{\overrightarrow{tm}^{\rho (s\vert i,\overrightarrow{tm})}_{s\vert i}}\!\subseteq\!
\overline{r\big( S(\overrightarrow{tm})\big)}^{\tau_{\xi_{s\vert i}}}$. The induction assumption implies that 
$U\!\subseteq\! r\big( S(\overrightarrow{tm}^{\eta_{s\vert i}}_{s\vert i})\big)$. So let us prove that 
$r\big( S(\overrightarrow{tm}^{\eta_{s\vert i}}_{s\vert i})\big)\!\subseteq\!
\overline{r\big( S(\overrightarrow{tm})\big)}^{\tau_{\xi_{s\vert i}}}$. Note that 
$s\vert (i\! +\! 1)\!\subseteq\! s(\overrightarrow{tm}^{\eta_{s\vert i}}_{s\vert i})\!\subseteq\! 
S(\overrightarrow{tm}^{\eta_{s\vert i}}_{s\vert i})$ and, similarly, 
$s'\vert (i\! +\! 1)\!\subseteq\! S(\overrightarrow{tm})$. Lemma 6.17 implies that 
$O\big( S(\overrightarrow{tm}^{\eta_{s\vert i}}_{s\vert i})\big)\sqsubseteq 
O\big( S(\overrightarrow{tm})\big)$,  and the beginning of its proof that 
$O\big( S(\overrightarrow{tm}^{\eta_{s\vert i}}_{s\vert i})\big)\!\not=\! O\big( S(\overrightarrow{tm})\big)$. Now Lemma 6.21 implies the result.\hfill{$\diamond$}\bigskip

\noindent - Let ${\cal X}\! :=\! d^{l+1}$. The map 
$\Psi\! :\! {\cal X}^d\!\rightarrow\!\Ana\big( (\omega^\omega )^d\big)$ is defined on ${\cal T}^{l+1}$ by
$$\Psi (\overrightarrow{tm})\! :=\!\left\{\!\!\!\!\!\!
\begin{array}{ll} 
& r\big( S(\overrightarrow{tm})\big)\cap\bigcap_{1\leq\rho\leq\rho (\emptyset ,\overrightarrow{tm})}\ 
\overline{U_{\overrightarrow{tm}^{\rho}_{\emptyset}}}^{\tau_\rho}\cap (\Pi_{i\in d}\ X_i)\cap
\Omega_{(\omega^\omega )^d}\mbox{ if }s(\overrightarrow{tm})\! =\!\emptyset\mbox{,}\cr\cr
& U_{\overrightarrow{tm}^{\rho (s(\overrightarrow{tm})^-,\overrightarrow{tm})}_{s(\overrightarrow{tm})^-}}
\cap\bigcap_{1\leq\rho <\rho (s(\overrightarrow{tm})^-,\overrightarrow{tm})}\ 
\overline{U_{\overrightarrow{tm}^{\rho}_{s(\overrightarrow{tm})^-}}}^
{\tau_{(\Sigma_{i'<\vert s(\overrightarrow{tm})\vert -1}\ \rho (s\vert i',\overrightarrow{tm}))+\rho}}\cr
& \hfill{\cap\bigcap_{i'<\vert s(\overrightarrow{tm})\vert -1}\ 
\bigcap_{1\leq\rho\leq\rho (s\vert i',\overrightarrow{tm})}\ \overline{U_{\overrightarrow{tm}^{\rho}_{s\vert i'}}}
^{\tau_{(\Sigma_{i''<i'}\ \rho (s\vert i'',\overrightarrow{tm}))+\rho}}\cap (\Pi_{i\in d}\ X_i)}\cr
& \hfill{\mbox{if }s(\overrightarrow{tm})\!\not=\!\emptyset\wedge
\overrightarrow{tm}^{\eta_{s(\overrightarrow{tm})^-}}_{s(\overrightarrow{tm})^-}
\!\in\!\bigcap_{i'<\vert s(\overrightarrow{tm})\vert}\ 
P_{s(\overrightarrow{tm})(i')(0)}\big( s(\overrightarrow{tm})\vert i'\big)}\cr 
& \hfill{\wedge\exists i_0\! <\! \vert s(\overrightarrow{tm})\vert ~~
\eta_{s(\overrightarrow{tm})\vert i_0}\!\geq\! 1\mbox{,}}\cr\cr 
& r\big( S(\overrightarrow{tm})\big)\cr
& \hfill{\cap\bigcap_{i'\leq i}\ \bigcap_{1\leq\rho\leq\rho (s(\overrightarrow{tm})\vert i',\overrightarrow{tm})}\ 
\overline{U_{\overrightarrow{tm}^{\rho}_{s(\overrightarrow{tm})\vert i'}}}
^{\tau_{(\Sigma_{i''<i'}\ \rho (s(\overrightarrow{tm})\vert i'',\overrightarrow{tm}))+\rho}}\cap (\Pi_{i\in d}\ X_i)\cap\Omega_{(\omega^\omega )^d}}\cr
& \hfill{\mbox{if }s(\overrightarrow{tm})\!\not=\!\emptyset}\wedge
\overrightarrow{tm}^{\eta_{s(\overrightarrow{tm})^-}}_{s(\overrightarrow{tm})^-}
\!\notin\!\bigcap_{i'<\vert s(\overrightarrow{tm})\vert}\ 
P_{s(\overrightarrow{tm})(i')(0)}\big( s(\overrightarrow{tm})\vert i'\big)\cr
& \hfill{\wedge ~i\! <\!\vert s(\overrightarrow{tm})\vert\mbox{ is\ maximal\ with }
\overrightarrow{tm}^{\eta_{s(\overrightarrow{tm})\vert i}}_{s(\overrightarrow{tm})\vert i}\!\in\!
\bigcap_{i'<i}\ P_{s(\overrightarrow{tm})(i')(0)}\big( s(\overrightarrow{tm})\vert i'\big)}\cr 
& \hfill{\wedge ~\exists i_0\!\leq\! i~~\eta_{s(\overrightarrow{tm})\vert i_0}\!\geq\! 1\mbox{,}}\cr\cr 
& U_{\vec t}\cap (\Pi_{i\in d}\ X_i)\mbox{ if }s(\overrightarrow{tm})\!\not=\!\emptyset\wedge
\overrightarrow{tm}^{\eta_{s(\overrightarrow{tm})^-}}_{s(\overrightarrow{tm})^-}\!\in\!
\bigcap_{i'<\vert s(\overrightarrow{tm})\vert}\ 
P_{s(\overrightarrow{tm})(i')(0)}\big( s(\overrightarrow{tm})\vert i'\big)\cr 
& \hfill{\wedge ~\forall i_0\! <\! \vert s(\overrightarrow{tm})\vert ~~
\eta_{s(\overrightarrow{tm})\vert i_0}\! =\! 0\mbox{,}}\cr\cr 
& r\big( S(\overrightarrow{tm})\big)\cap (\Pi_{i\in d}\ X_i)\cap\Omega_{(\omega^\omega )^d}\ \mbox{ if }~
s(\overrightarrow{tm})\!\not=\!\emptyset\cr
& \hfill{\wedge ~\overrightarrow{tm}^{\eta_{s(\overrightarrow{tm})^-}}_{s(\overrightarrow{tm})^-}
\!\notin\!\bigcap_{i'<\vert s(\overrightarrow{tm})\vert}\ 
P_{s(\overrightarrow{tm})(i')(0)}\big( s(\overrightarrow{tm})\vert i'\big)}\cr
& \hfill{\wedge ~i\! <\!\vert s(\overrightarrow{tm})\vert\mbox{ is\ maximal\ with }
\overrightarrow{tm}^{\eta_{s(\overrightarrow{tm})\vert i}}_{s(\overrightarrow{tm})\vert i}\!\in\!
\bigcap_{i'<i}\ P_{s(\overrightarrow{tm})(i')(0)}\big( s(\overrightarrow{tm})\vert i'\big)}\cr 
& \hfill{\wedge ~\forall i_0\!\leq\! i~~\eta_{s(\overrightarrow{tm})\vert i_0}\! =\! 0.}
\end{array}\right.$$
By the claim, $\Psi (\overrightarrow{tm})$ is $\tau_{1}$-dense in 
$\overline{U_{\overrightarrow{tm}^1_{s(\overrightarrow{tm})\vert i_0}}}^{\tau_{1}}\cap (\Pi_{i\in d}\ X_i)$ in the 2nd and in the 3rd cases.

\vfill\eject

 In these cases, as 
$\overrightarrow{tm}^1_{s(\overrightarrow{tm})\vert i_0}\!\subseteq\! \vec t \!\subseteq\!\overrightarrow{tm}$ 
and $R^{(1)}_{s(\overrightarrow{tm})\vert i_0}$ is distinguished in 
$R^{(0)}_{s(\overrightarrow{tm})\vert i_0}\! =\subseteq$, 
$$\overrightarrow{tm}^1_{s(\overrightarrow{tm})\vert i_0}\ R^{(1)}_{s(\overrightarrow{tm})\vert i_0}\ 
\vec t$$ 
and 
$U_{\vec t }\!\subseteq\!\overline{U_{\overrightarrow{tm}^1_{s(\overrightarrow{tm})\vert i_0}}}^{\tau_1}$, by induction assumption. Therefore 
$$U_{\vec t }\cap (\Pi_{i\in d}\ X_i)\!\subseteq\!
\overline{U_{\overrightarrow{tm}^1_{s(\overrightarrow{tm})\vert i_0}}}^{\tau_1}\cap (\Pi_{i\in d}\ X_i)
\!\subseteq\!\overline{\Psi} (\overrightarrow{tm}).$$ 
Using similar arguments, one can prove that this also holds in the last two cases.\bigskip

 Let us look at the first case. If $\eta_\emptyset\!\geq\! 1$, then using arguing as in the claim one can prove that $U_{\overrightarrow{tm}^{\rho (\emptyset ,\overrightarrow{tm})}_{\emptyset}}\cap
\bigcap_{1\leq\rho <\rho (\emptyset ,\overrightarrow{tm})}\ 
\overline{U_{\overrightarrow{tm}^{\rho}_\emptyset}}^{\tau_\rho}\cap (\Pi_{i\in d}\ X_i)$ is $\tau_{1}$-dense in $\overline{U_{\overrightarrow{tm}^1_\emptyset}}^{\tau_{1}}\cap (\Pi_{i\in d}\ X_i)$. Now we can write 
$U_{\overrightarrow{tm}^{\rho (\emptyset ,\overrightarrow{tm})}_{\emptyset}}\!\subseteq\! 
r\big( S(\overrightarrow{tm}^{\eta_\emptyset}_\emptyset )\big)\! =\! r\big( S(\overrightarrow{tm})\big)$ and we can repeat the previous argument since $i_0\! =\! 0$. If $\eta_\emptyset\! =\! 0$, then we get 
$\overrightarrow{tm}^{\eta_\emptyset}_\emptyset\! =\!\vec t$, and 
$U_{\vec t }\cap (\Pi_{i\in d}\ X_i)\!\subseteq\! r\big( S(\vec t~)\big)\cap (\Pi_{i\in d}\ X_i)\! =\! 
r\big( S(\overrightarrow{tm})\big)\cap (\Pi_{i\in d}\ X_i)$ and we are done.\bigskip

 Now we can write $(\alpha^{i}_{t_i})_{i\in d}\!\in\! U_{\vec t }\cap (\Pi_{i\in d}\ X_i)\!\subseteq\!
\overline{\Psi} (\overrightarrow{tm})$, and we conclude as in the proof of Theorem 4.4.1.\hfill{$\square$}\bigskip

 The rest of this section is devoted to the proof of Theorem 1.8.(2) when $\Delta ({\bf\Gamma})$ is a Wadge class. Recall Theorem 5.2.8. We will say that $\alpha\!\in\!\Borel\cap\Lambda^\infty$ is $suitable$ if 
$\Delta ({\bf\Gamma}_{c(\alpha )})$ is a Wadge class and one of the following holds:\bigskip
 
 \noindent (1) There is $\overline{\alpha}\!\in\!\Borel\cap\Lambda^\infty$ normalized with 
$${\bf\Gamma}_{c(\alpha )}\! =\!
\Big\{ (A_0\cap C_0)\cup (A_1\cap C_1)\mid A_0,\neg A_1\!\in\! {\bf\Gamma}_{c(\overline{\alpha})}\wedge C_0,C_1\!\in\!\boraone\wedge C_0\cap C_1\! =\!\emptyset\Big\}.$$
(2) There is $\alpha'\!\in\!\Borel$ such that $(\alpha')_p\!\in\!\Lambda^\infty$ is normalized for each 
$p\!\geq\! 1$, $\big( {\bf\Gamma}_{c((\alpha')_p)}\big)_{p\geq 1}$ is strictly increasing, and 
${\bf\Gamma}_{c(\alpha )}\! =\!
\Big\{\bigcup_{p\geq 1}~(A_p\cap C_p)\mid A_p\!\in\! {\bf\Gamma}_{c((\alpha')_p)}\wedge 
C_p\!\in\!\boraone\wedge C_p\cap C_q\! =\!\emptyset\mbox{ if }p\!\neq\! q\Big\}$.\bigskip

 Assume that $\alpha$ is suitable and $a_0,a_1\!\in\!\Borel$ satisfy $A_0\cap A_1\! =\!\emptyset$. Then Lemma 6.7.(b) gives $r(\overline{\alpha}, a_0,a_1)$ and $r(\overline{\alpha}, a_1,a_0)$, or 
$r\big( (\alpha')_p, a_0,a_1\big)$. We set 
$R(\overline{\alpha},a_0,a_1)\! :=\!\neg {\cal U}_{r(\overline{\alpha},a_0,a_1)}$ in the same fashion as before, and 
$$R'(\alpha ,a_0,a_1)\! :=\!\left\{\!\!\!\!\!\!
\begin{array}{ll}
& \overline{R(\overline{\alpha},a_0,a_1)}^{\tau_1}\cap
\overline{R(\overline{\alpha},a_1,a_0)}^{\tau_1}\mbox{ if we are in Case }(1)\mbox{,}\cr\cr
& \bigcap_{p\geq 1}~\overline{R\big( (\alpha')_p, a_0,a_1\big)}^{\tau_1}\mbox{ if we are in Case }(2).
\end{array}
\right.$$
We now give the self-dual version of Lemma 6.8.

\begin{lem} Let $\alpha$ suitable, and $a_0,a_1\!\in\!\Borel$ such that $A_0\cap A_1\! =\!\emptyset$. We assume that $R'(\alpha ,a_0,a_1)\! =\!\emptyset$. Then $A_0$ is separable from $A_1$ by a 
$\Borel\cap\Delta ({\bf\Gamma}_{c(\alpha )})(\tau_1)$ set.\end{lem}

\noindent\bf Proof.\rm\ (1) As $\overline{R(\overline{\alpha},a_0,a_1)}^{\tau_1}\cap
\overline{R(\overline{\alpha},a_1,a_0)}^{\tau_1}\! =\!\emptyset$, there is $C\!\in\!\borone (\tau_1)$ separating $R(\overline{\alpha},a_0,a_1)$ from $R(\overline{\alpha},a_1,a_0)$. As 
$R(\overline{\alpha},a_0,a_1)$ and $R(\overline{\alpha},a_1,a_0)$ are $\Ana$, we may assume that 
$C\!\in\!\Borel$, by Theorem 4.2.2. A double application of Lemmas 6.7.(b) and 6.8 gives some sets 
$B_0,B_1\!\in\!\Borel\cap {\bf\Gamma}_{c(\overline{\alpha})}(\tau_1)$ such that $B_0$ (resp., $B_1$) separates $A_0\cap C$ (resp., $A_1\!\setminus\! C$) from $A_1\cap C$ (resp., $A_0\!\setminus\! C$). Now  the set $(B_0\cap C)\cup (\neg B_1\cap\neg C)$ is suitable.\bigskip

\noindent (2) The proof is similar, but we have to make some $\Borel$-selection. As $\Theta^\infty$ is 
$\Ca$ and $r\big( (\alpha')_p, a_0,a_1\big)$ is $\Borel$ and completely determined by $(\alpha')_p$, 
$a_0$ and $a_1$, the sequence $\Big( r\big( (\alpha')_p, a_0,a_1\big)\Big)_{p\geq 1}$ is $\Borel$. As 
$\bigcap_{p\geq 1}~\overline{R\big( (\alpha')_p, a_0,a_1\big)}^{\tau_1}\! =\!\emptyset$, there is a $\Borel$-recursive map $f\! :\! (\omega^\omega )^d\!\rightarrow\!\omega$ such that 
$f(\vec\alpha )\!\geq\! 1$ and 
$\vec\alpha\!\notin\!\overline{R\big( (\alpha')_{f(\vec\alpha )}, a_0,a_1\big)}^{\tau_1}$ for each 
$\vec\alpha\!\in\! (\omega^\omega )^d$.\bigskip

 We set $U_p\! :=\! f^{-1}(\{ p\})$, so that $U_p$ and $R\big( (\alpha')_p, a_0,a_1\big)$ are disjoint $\Ana$ sets and separable by a $\tau_1$-open set. By Theorem 4.2.2, there is 
$V_p\!\in\!\Borel\cap\boraone (\tau_1)$ separating them. Moreover, we may assume that the 
sequence $(V_p)$ is $\Borel$. We reduce the sequence $(V_p)$ into a $\Borel$-sequence $(C_p)$ 
of $\Borel\cap\boraone (\tau_1)$ sets. Note that $(C_p)$ is a partition of $(\omega^\omega )^d$ into 
$\borone (\tau_1)$ sets. As $R\big( (\alpha')_p, a_0,a_1\big)\cap C_p\! =\!\emptyset$, Lemma 6.8 gives 
$\beta',\gamma'\!\in\!\omega^\omega$ such that 
$\big( (\alpha')_p,(\beta')_p,(\gamma')_p\big)\!\in\!\Upsilon^\infty$ and $C_{(\gamma')_p}$ separates $A_1\cap C_p$ from $A_0\cap C_p$ for each $p\!\geq\! 1$. Moreover, we may assume that 
$\beta',\gamma'\!\in\!\Borel$. Now the set $\bigcup_{p\geq 1}~(\neg C_{(\gamma')_p}\cap C_p)$ is suitable.$\hfill{\square}$\bigskip

 We now give the self-dual version of Theorem 6.9.

\begin{thm} Let $T_d$ be a tree with $\Borel$ suitable levels, $\alpha$ suitable, 
$\beta^\varepsilon ,\gamma^\varepsilon\!\in\!\omega^\omega$ with 
$(\alpha,\beta^\varepsilon ,\gamma^\varepsilon )\!\in\! \Upsilon_1^\infty$, 
$S^\varepsilon\! :=\! j_d^{-1}(C^{\omega^\omega}_{\gamma^\varepsilon})\cap\lceil T_d\rceil$, and 
$a_0,a_1,\underline{a}_0,\underline{a}_1,r\!\in\!\omega^\omega$ such that 
$\vec v\! :=\! (\alpha ,a_0,a_1,\underline{a}_0,\underline{a}_1,r)\!\in\!\Theta^\infty$. We assume that 
$S^0$ and $S^1$ are disjoint. Then one of the following holds:\smallskip

\noindent (a) $R'(\alpha ,a_0,a_1)\! =\!\emptyset$.\smallskip

\noindent (b) The inequality $\big( (\Pi_i''\lceil T_d\rceil )_{i\in d}, S^0,S^1\big)\leq
\big( (\omega^\omega )_{i\in d}, A_0, A_1\big)$ holds.\end{thm}

 Now we can state the version of Theorem 4.2.2 for the self-dual Wadge classes of Borel sets.
 
\begin{thm} Let $T_d$ be a tree with $\Borel$ suitable levels, $\alpha$ suitable, 
$\beta^\varepsilon ,\gamma^\varepsilon\!\in\!\omega^\omega$ with 
$(\alpha,\beta^\varepsilon ,\gamma^\varepsilon )\!\in\! \Upsilon_1^\infty$, 
$S^\varepsilon\! :=\! j_d^{-1}(C^{\omega^\omega}_{\gamma^\varepsilon})\cap\lceil T_d\rceil$, and 
$a_0,a_1,\underline{a}_0,\underline{a}_1,r\!\in\!\omega^\omega$ such that 
$\vec v\! :=\! (\alpha ,a_0,a_1,\underline{a}_0,\underline{a}_1,r)\!\in\!\Theta^\infty$. We assume that $S^0$, $S^1$ are disjoint and not separable by a $\mbox{pot}\big(\Delta ({\bf\Gamma}_{c(\alpha )})\big)$ set. Then the following are equivalent:\smallskip

\noindent (a) The set $A_0$ is not separable from $A_1$ by a 
$\hbox{\it pot}\big(\Delta ({\bf\Gamma}_{c(\alpha )})\big)$ set.\smallskip

\noindent (b) The set $A_0$ is not separable from $A_1$ by a 
$\Borel\cap\hbox{\it pot}\big(\Delta ({\bf\Gamma}_{c(\alpha )})\big)$ set.\smallskip

\noindent (c) The set $A_0$ is not separable from $A_1$ by a 
$\Delta ({\bf\Gamma}_{c(\alpha )})(\tau_1)$ set.\smallskip

\noindent (d) $R'(\alpha ,a_0,a_1)\!\not=\!\emptyset$.\smallskip

\noindent (e) The inequality $\big( (d^\omega )_{i\in d}, S^0,S^1\big)\leq
\big( (\omega^\omega )_{i\in d}, A_0, A_1\big)$ holds.\end{thm}

\noindent\bf Proof.\rm\ We argue as in the proof of Theorem 6.10, using Lemma 6.24 (resp., Theorem 6.25) instead of Lemma 6.8 (resp., Theorem 6.9).\hfill{$\square$}\bigskip

\noindent\bf Proof of Theorem 1.8.(2).\rm\ We argue as in the proof of Theorem 1.8.(1). Theorem 5.2.8 gives $\overline{u}$ or $\big( (u')_p\big)_{p\geq 1}$. The equalities in Theorem 5.2.8 hold in $\omega^\omega$, but also in any $0$-dimensional Polish space (we argue like in Lemma 5.2.2 to see it). Using Definition 5.1.2, we can build $u\!\in\! {\cal D}$ with ${\bf\Gamma}\! =\! {\bf\Gamma}_u$.

\vfill\eject
 
 Using Lemmas 6.2 and 6.4, we get $\alpha\!\in\!\Lambda^\infty$ normalized with 
${\bf\Gamma}_{c(\alpha )}\! =\! {\bf\Gamma}_u$, and $\overline{\alpha}\!\in\!\Lambda^\infty$ (or 
$\alpha'\!\in\!\Lambda^\infty$ such that $(\alpha')_p$ is) normalized with 
${\bf\Gamma}_{\overline{u}}\! =\! {\bf\Gamma}_{c(\overline{\alpha})}$ (or 
${\bf\Gamma}_{(u')_p}\! =\! {\bf\Gamma}_{c((\alpha')_p)}$).\bigskip

 By Theorem 4.1.3 in [Lo-SR2] there is $B^\varepsilon\!\in\! {\bf\Gamma}(\omega^\omega )$ with 
$S^\varepsilon\! =\! j_d^{-1}(B^\varepsilon)\cap\lceil T_d\rceil$. To simplify the notation, we may assume that $T_d$ has $\Borel$ levels, $\alpha$, as well as $\overline{\alpha}$ (or $\alpha'$), are $\Borel$, and 
$A_0, A_1$ are $\Ana$. By Lemma 6.5 there are  
$\beta^\varepsilon ,\gamma^\varepsilon\!\in\!\omega^\omega$ such that 
$(\alpha ,\beta^\varepsilon ,\gamma^\varepsilon )\!\in\!\Upsilon_1^\infty$ and 
$C^{\omega^\omega}_{\gamma^\varepsilon}\! =\! B^\varepsilon$. Lemma 6.7.(b) gives 
$\underline{a}_0,\underline{a}_1,r$ with 
$(\alpha ,a_0,a_1,\underline{a}_0,\underline{a}_1,r)\!\in\!\Theta^\infty$. Lemma 6.24 implies that 
$R'(\alpha ,a_0,a_1)\!\not=\!\emptyset$. So (b) holds, by Theorem 6.26.\hfill{$\square$}\bigskip

\noindent\bf Proof of Theorem 6.25.\rm\ (1) Let 
$C^\varepsilon_{\varepsilon'}\!\in\!\boraone (\lceil T_d\rceil )$, 
$A^\varepsilon_0\!\in\! {\bf\Gamma}_{c(\overline{\alpha})}(\lceil T_d\rceil )$, 
$A^\varepsilon_1\!\in\!\check {\bf\Gamma}_{c(\overline{\alpha})}(\lceil T_d\rceil )$ such that 
$S^\varepsilon\! =\! (A^\varepsilon_0\cap C^\varepsilon_0)\cup (A^\varepsilon_1\cap C^\varepsilon_1)$. We reduce the family $(C^0_0,C^0_1,C^1_0,C^1_1)$ into a family $(O^0_0,O^0_1,O^1_0,O^1_1)$ of open subsets of $\lceil T_d\rceil$. Note that $S^\varepsilon\!\subseteq\! T^\varepsilon\! :=\! 
(A^\varepsilon_0\cap O^\varepsilon_0)\cup (A^\varepsilon_1\cap O^\varepsilon_1)\cup 
(\neg A^{1-\varepsilon}_0\cap O^{1-\varepsilon}_0)\cup 
(\neg A^{1-\varepsilon}_1\cap O^{1-\varepsilon}_1)$. We will in fact ensure that 
$\big( (\Pi_i''\lceil T_d\rceil )_{i\in d}, T^0,T^1\big)\leq
\big( (\omega^\omega )_{i\in d}, A_0, A_1\big)$ if (a) does not hold, which will be enough.\bigskip

\noindent\bf Subcase 1.\rm\ $\vert (\alpha )_0\vert\! =\! 0$.\bigskip

 We set $o^\varepsilon_{\varepsilon'}\! :=\! h[\lceil T_d\rceil\!\setminus\! O^\varepsilon_{\varepsilon'}]$, so that $o^\varepsilon_{\varepsilon'}\!\in\!\bormone ([\subseteq ])$. We also set 
$$D\! :=\!\{\vec s\!\in\! T_d\mid\vec s\! =\!\vec\emptyset\vee\forall (\varepsilon,\varepsilon')\!\in\! 2^2~~
\exists {\cal B}\!\in\! o^\varepsilon_{\varepsilon'}\ \ \vec s\!\in\! {\cal B}\}\mbox{,}$$ 
$$D^\varepsilon_{\varepsilon'}\! :=\!\{\vec s\!\in\! T_d\mid\vec s\!\not=\!\vec\emptyset\wedge
\forall {\cal B}\!\in\! o^\varepsilon_{\varepsilon'}\ \ \vec s\!\notin\! {\cal B}\ \wedge
\forall (\varepsilon'',\varepsilon''')\!\in\! 2^2\!\setminus\!\{ (\varepsilon ,\varepsilon')\}\ \ 
\exists {\cal B}\!\in\! o^{\varepsilon''}_{\varepsilon'''}\ \ \vec s\!\in\! {\cal B}\}\mbox{,}$$
so that $(D,D^0_0,D^0_1,D^1_0,D^1_1)$ is a partition of $T_d$. The proof is very similar to the proof of Theorem 4.4.2 when $\xi\! =\! 1$. The changes to make in the conditions (1)-(7) are as follows:
$$\begin{array}{ll}
& \!\!\!\! (4)\ U_{\vec s  }\!\subseteq\! R'(\alpha ,a_0,a_1)\! =\!\overline{A_0}^{\tau_1}\cap
\overline{A_1}^{\tau_1}\mbox{ if }\vec s \!\in\! D\mbox{,}\cr\cr  
& \!\!\!\! (5)\ U_{\vec s  }\!\subseteq\! A_0\mbox{ if }\vec s \!\in\! D^0_1\cup D^1_0\mbox{,}\cr\cr
& \!\!\!\! (6)\ U_{\vec s  }\!\subseteq\! A_1\mbox{ if }\vec s \!\in\! D^0_0\cup D^1_1\mbox{,}\cr\cr
& \!\!\!\! (7)\ (\vec s,\vec t \!\in\! D\vee\vec s,\vec t \!\in\! D^\varepsilon_{\varepsilon'})\Rightarrow 
U_{\vec t }\!\subseteq\! U_{\vec s }.
\end{array}$$
We conclude as in the proof of Theorem 4.4.2.\bigskip

\noindent\bf Subcase 2.\rm\ $\vert (\alpha )_0\vert\!\geq\! 1$.\bigskip

 We will have the same scheme of construction as in the proof of Theorem 6.9. As long as 
$\vec t\!\in\! D$, we will have $U_{\vec t}\!\subseteq\! R'(\alpha ,a_0,a_1)$. If 
$\vec t\!\in\! D^\varepsilon_{\varepsilon'}$, then all the extensions of $\vec t$ will stay in 
$D^\varepsilon_{\varepsilon'}$, and we will copy the construction of the proof of Theorem 6.9, since inside the clopen set defined by $\vec t$ we want to reduce a pair $(\tilde S^0,\tilde S^1)$ to $(A_0,A_1)$.\bigskip
 
 As $A^\varepsilon_0\!\in\! {\bf\Gamma}_{c(\overline{\alpha})}(\lceil T_d\rceil )$, there is 
$B^\varepsilon_0\!\in\! {\bf\Gamma}_{c(\overline{\alpha})}(\omega^\omega )$ with 
$A^\varepsilon_0\! =\! j_d^{-1}(B^\varepsilon_0)\cap\lceil T_d\rceil$. As 
$\overline{\alpha}\!\in\!\Borel\cap\Lambda^\infty$, Lemma 6.5.(b) gives 
$\beta^\varepsilon_0,\gamma^\varepsilon_0\!\in\!\omega^\omega$ such that 
$(\overline{\alpha},\beta^\varepsilon_0,\gamma^\varepsilon_0)\!\in\!\Upsilon_1^\infty$ and 
$C^{\omega^\omega}_{\gamma^\varepsilon_0}\! =\! B^\varepsilon_0$. Similarly, there are 
$\beta^\varepsilon_1,\gamma^\varepsilon_1\!\in\!\omega^\omega$ such that 
$(\overline{\alpha},\beta^\varepsilon_1,\gamma^\varepsilon_1)\!\in\!\Upsilon_1^\infty$ and 
$A^\varepsilon_1\! =\! j_d^{-1}(\neg C^{\omega^\omega}_{\gamma^\varepsilon_1})\cap\lceil T_d\rceil$. 

\vfill\eject

 We can associate with any $(\varepsilon ,\varepsilon')\!\in\! 2^2$ the objects met before, among which the function ${\cal Z}^{\varepsilon ,\varepsilon'}$, the ordinals $\eta^{\varepsilon ,\varepsilon'}_s$, the resolution families $(R^{(\rho )}_{\varepsilon ,\varepsilon',s})_{\rho\leq\eta^{\varepsilon ,\varepsilon'}_s}$, the ordinals $\rho (\varepsilon ,\varepsilon',s,\vec s~)$. Instead of considering the set $P_q(s)$, we will consider $P^{\varepsilon ,\varepsilon'}_q(s)\cap D^\varepsilon_{\varepsilon'}$. If 
$\vec t\!\in\! D^\varepsilon_{\varepsilon'}$, then we set 
$\vec w(\vec t~)\! :=\!\vec w^\varepsilon_{\varepsilon'}$. This allows us to define 
$s(\vec t~)\!\in\! {\mathfrak T}\big(\vec w(\vec t~)\big)$ and $S(\vec t~)\!\in\! {\cal M}_{\vec w(\vec t)}$. We also set
$$\vec v(\vec t~)\! :=\!\left\{\!\!\!\!\!\!
\begin{array}{ll}
& (\overline{\alpha},a_0,a_1,\underline{a}_0,\underline{a}_1,r)\mbox{ if }\vec t\!\in\! D^0_0\cup D^1_1
\mbox{,}\cr\cr
& (\overline{\alpha},a_1,a_0,\underline{a}_0,\underline{a}_1,r)\mbox{ if }\vec t\!\in\! D^0_1\cup D^1_0.
\end{array}
\right.$$
The other modifications to make in the conditions (1)-(6) are as follows. In condition (4), we ask for the inclusion $U_{\vec t}\!\subseteq\! R\big( S(\vec t~)\big)$ only if $\vec t\!\notin\! D$. If $\vec t\!\in\! D$, then we want that $U_{\vec t}\!\subseteq\! R'(\alpha ,a_0,a_1)$. Condition (6) was described when 
$\vec s,\vec t\!\in\! D^\varepsilon_{\varepsilon'}$. If $\vec s,\vec t\!\in\! D$, then we also want that 
$U_{\vec t}\!\subseteq\! U_{\vec s}$.\bigskip

 The sequence ${\cal F}(\vec\beta )$ is defined if 
$\beta\!\in\! C^0_0\cup C^0_1\cup C^1_0\cup C^1_1$. If 
$\beta\!\notin\! C^0_0\cup C^0_1\cup C^1_0\cup C^1_1$, then $\vec\beta\vert k\!\in\! D$ for each integer $k$, and ${\cal F}(\vec\beta )$ is also defined. The definition of $\vec v(\vec t~)$ ensures that 
$T^\varepsilon\!\subseteq\! (\Pi_{i\in d}~f_i)^{-1}(A_\varepsilon )$.\bigskip

 The defintion of $\Psi (\overrightarrow{tm})$ is done if $\overrightarrow{tm}\!\notin\! D$. If 
$\overrightarrow{tm}\!\in\! D$, then we simply set 
$$\Psi (\overrightarrow{tm})\! :=\! U_{\vec t}~\cap (\Pi_{i\in d}~X_i).$$
Then we conclude as in the proof of Theorem 6.9.\bigskip

\noindent (2) Let $C^\varepsilon_p\!\in\!\boraone (\lceil T_d\rceil )$ and 
$A^\varepsilon_p\!\in\! {\bf\Gamma}_{c((\alpha')_p)}(\lceil T_d\rceil )$ such that 
$S^\varepsilon\! =\!\bigcup_{p\geq 1}~(A^\varepsilon_p\cap C^\varepsilon_p)$. We reduce the family 
$(C^0_1,C^0_2,...,C^1_1,C^1_2,...)$ into a family $(O^0_1,O^0_2,...,O^1_1,O^1_2,...)$ of open subsets of $\lceil T_d\rceil$. Note that $S^\varepsilon\!\subseteq\! T^\varepsilon\! :=\! 
(A^\varepsilon_1\cap O^\varepsilon_1)\cup\bigcup_{p\geq 1}~
\big( (\neg A^{1-\varepsilon}_p\cap O^{1-\varepsilon}_p)\cup 
(A^\varepsilon_{p+1}\cap O^\varepsilon_{p+1})\big)$. We will in fact ensure that 
$\big( (\Pi_i''\lceil T_d\rceil )_{i\in d}, T^0,T^1\big)\leq
\big( (\omega^\omega )_{i\in d}, A_0, A_1\big)$ if (a) does not hold, which will be enough.\bigskip

 The proof is similar. We can assume that $\big\vert\big( (\alpha')_p\big)_0\big\vert\!\geq\! 1$ for each 
$p\!\geq\! 1$, since $({\bf\Gamma}_{c((\alpha')_p)})_{p\geq 1}$ is strictly increasing.  So there is no Subcase 1. We set
$$\vec v(\vec t~)\! :=\!\left\{\!\!\!\!\!\!\!
\begin{array}{ll}
& (\overline{\alpha},a_0,a_1,\underline{a}_0,\underline{a}_1,r)\mbox{ if }
\vec t\!\in\!\bigcup_{p\geq 1}~D^0_p\mbox{,}\cr\cr
& (\overline{\alpha},a_1,a_0,\underline{a}_0,\underline{a}_1,r)\mbox{ if }
\vec t\!\in\!\bigcup_{p\geq 1}~D^1_p.
\end{array}
\right.$$
We conclude as in Case 1.\hfill{$\square$}

\section{$\!\!\!\!\!\!$ Injectivity complements}\indent

 In the introduction, we saw that G. Debs proved that we can have the $f_i$'s one-to-one in Theorem 1.3 when $d\! =\! 2$, ${\bf\Gamma}\!\in\!\{\bormxi ,\boraxi\}$ and $\xi\!\geq\! 3$.\bigskip
 
\noindent $\bullet$ This cannot be extended to higher dimensions, even if we replace $(d^\omega )^d$ with $\Pi_{i\in d}~Z_i$, where $Z_i$ is a sequence of Polish spaces.

\vfill\eject

 Indeed, we argue by contradiction. Recall the proof of Theorem 3.1. We saw that there is $C_\xi$ in 
$\boraxi (2^\omega )\!\setminus\!\bormxi$ such that 
${\mathbb{S}}^3_\xi\! :=\!\{\vec\alpha\!\in\!\lceil T_3\rceil\mid {\cal S}(\alpha_0\Delta\alpha_1)\!\in\! C_\xi\}$ is not separable from $\lceil T_3\rceil\!\setminus\! {\mathbb{S}}^3_\xi$ by a $\mbox{pot}(\bormxi )$ set. We set
$$\begin{array}{ll}
& B^0\! :=\!\{\vec\alpha\!\in\! 3^\omega\!\times\! 3^\omega\!\times\! 1\mid 
{\cal S}(\alpha_0\Delta\alpha_1)\!\in\! C_\xi\}\mbox{,}\cr
& B^1\! :=\!\{\vec\alpha\!\in\! 3^\omega\!\times\! 1\!\times\! 3^\omega\mid 
{\cal S}(\alpha_0\Delta\alpha_2)\!\in\! C_\xi\}\mbox{,}\cr
& B^2\! :=\!\{\vec\alpha\!\in\! 1\!\times\! 3^\omega\!\times\! 3^\omega\mid 
{\cal S}(\alpha_1\Delta\alpha_2)\!\in\! C_\xi\}.
\end{array}$$
Let $O\! :\! 3^\omega\!\rightarrow\! 1$. As ${\mathbb{S}}^3_\xi\! :=\! 
(\mbox{Id}_{3^\omega}\!\times\!\mbox{Id}_{3^\omega}\!\times\! O)^{-1}(B^0)\cap \lceil T_3\rceil$, 
$B^0\!\notin\!\mbox{pot}(\bormxi )$. Similarly, $B^1,B^2\!\notin\!\mbox{pot}(\bormxi )$. This implies that the $Z_i$'s have cardinality at most one, and ${\mathbb{S}}_0\!\in\!\borone$. Thus ${\mathbb{S}}_0$ is separable from ${\mathbb{S}}_1$ by a $\mbox{pot}(\bormxi )$ set, which is absurd.\bigskip

\noindent $\bullet$ If $d\! =\!\omega$, ${\bf\Gamma}\! =\!\bormxi$ and $\xi\!\geq\! 3$, then we cannot ensure that at least two of the $f_i$'s are one-to-one. Indeed, we again argue by contradiction. Consider 
$X_i\! :=\!\omega$, and $B_\xi\!\in\!\boraxi (\omega^\omega )\!\setminus\!\bormxi$. Then $B_\xi$ is not 
$\mbox{pot}(\bormxi )$ since the topology on $\omega$ is discrete. This implies that two of the $Z_i$'s at least are countable, say $Z_0,Z_1$ for example. Consider now $A_0\! :=\! {\mathbb{S}}^\omega_\xi$ and 
$A_1\! :=\!\lceil T_\omega\rceil\!\setminus\!  {\mathbb{S}}^\omega_\xi$. Then 
$(f_i\circ\Pi_i)[{\mathbb{S}}_0]$ is countable for each $i\!\in\! 2$. Thus $P\! :=\! (\Pi_{i\in d}~f_i)
[{\mathbb{S}}_0]\!\subseteq\! {\mathbb{S}}^\omega_\xi\!\subseteq\!\lceil T_\omega\rceil$ is countable since an element of $\lceil T_\omega\rceil$ is completely determined by two of its coordinates. Thus 
$P\!\in\!\mbox{pot}(\boratwo )\!\subseteq\!\mbox{pot}(\bormxi )$. Therefore $(\Pi_{i\in d}~f_i)^{-1}(P)$ is a 
$\mbox{pot}(\bormxi )$ set separating ${\mathbb{S}}_0$ from ${\mathbb{S}}_1$, which is absurd.\bigskip

\noindent $\bullet$ However, if ${\bf\Gamma}\!\in\!\{\bormxi ,\boraxi ,\borxi\}$ and $\xi\!\geq\! 3$, then we can ensure that $(\Pi_{i\in d}~f_i)_{\vert {\mathbb{S}}_0\cup {\mathbb{S}}_1}$ is one-to-one, using G. Debs's proof and some additional arguments.  This is also true if ${\bf\Gamma}\! =\! {\bf\Gamma}_u$ is a non self-dual Wadge class of Borel sets with $u(0)\!\geq\! 3$. This leads to the following notation. Let 
$(Z_i)_{i\in d}$, $(X_i)_{i\in d}$ be sequences of Polish spaces, and $S_0$, $S_1$ (resp., $A_0$, $A_1$) disjoint analytic subsets of $\Pi_{i\in d}\ Z_i$ (resp., $\Pi_{i\in d}\ X_i$). Then\bigskip

\leftline{$\big( (Z_i)_{i\in d}, S_0, S_1\big)\sqsubseteq\big( (X_i)_{i\in d}, A_0, A_1\big)\ \Leftrightarrow\ 
\forall i\!\in\! d\ \ \exists f_i\! :\! Z_i\!\rightarrow\! X_i\ \mbox{ continuous such that }$}\smallskip

\rightline{$(\Pi_{i\in d}\ f_i)_{\vert S_0\cup S_1}~\mbox{ is one-to-one and }\ 
\forall\varepsilon\!\in\! 2\ \ S_\varepsilon\!\subseteq\! (\Pi_{i\in d}\ f_i)^{-1}(A_\varepsilon ).$}

\begin{thm} There is no tuple $\big( (Z_i)_{i\in 2},S_0,S_1)$, where the $Z_i$'s are Polish spaces and $S_0$, $S_1$ disjoint analytic subsets of $\Pi_{i\in 2}\ Z_i$, such that for any tuple 
$\big( (X_i)_{i\in 2}, B_0, B_1\big)$ of the same type exactly one of the following holds:\smallskip

\noindent (a) The set $B_0$ is separable from $B_1$ by a $\mbox{pot}(\bormone )$ set.\smallskip

\noindent (b) The inequality $\big( (Z_i)_{i\in 2},S_0,S_1\big)\sqsubseteq
\big( (X_i)_{i\in 2}, B_0, B_1\big)$ holds.\end{thm}

 One can prove this result using the Borel digraph 
$B_0\! :=\!\bigcup_{n\in\omega}\ \mbox{Gr}({g_n}_{\vert 2^\omega\setminus M})$ considered in [L5] (see Section 3), which has countable vertical sections but is not locally countable. We give here another proof which moreover shows that we cannot hope for a positive result, even if $B_0$ is locally countable. This has to be noticed, since the locally countable sets have been considered a lot in the last decades.

\begin{lem} Let $\bf\Gamma$ be a Borel class, and $\big( (Z_i)_{i\in 2},S_0,S_1)$ as in the statement of Theorem 5.1 such that $S_0$ is not separable from $S_1$ by a $\mbox{pot}({\bf\Gamma})$ set. Then 
$S_0\cap (\Pi_0''S_1\!\times\!\Pi_1''S_1)$ is not separable from $S_1$ by a $\mbox{pot}({\bf\Gamma})$ set. Moreover, $S_0$ is not separable from $S_1\cap (\Pi_0''S_0\!\times\!\Pi_1''S_0)$ by a $\mbox{pot}({\bf\Gamma})$ set.\end{lem}

\noindent\bf Proof.\rm ~We prove the first assertion by contradiction, which gives 
$P\!\in\!\mbox{pot}({\bf\Gamma})$. The first reflection theorem gives Borel sets $C_0, C_1$ such that $\Pi_i''S_1\!\subseteq\! C_i$ and $S_0\cap (C_0\!\times\! C_1)\!\subseteq\! P$. Now 
$$S_0\!\subseteq\! P\cup (\neg C_0\!\times\! Z_1)\cup (Z_0\!\times\!\neg C_1)\!\subseteq\!\neg S_1\mbox{,}$$
which contradicts the fact that $S_0$ is not separable from $S_1$ by a $\mbox{pot}({\bf\Gamma})$ set.\bigskip

 We prove the second assertion using the first one, passing to complements.\hfill{$\square$}

\begin{lem} Let $\big( (Z_i)_{i\in 2},S_0,S_1)$ and $\big( (X_i)_{i\in 2}, B_0, B_1\big)$ be as in the statement of Theorem 5.1 such that 
$\big( (Z_i)_{i\in 2},S_0,S_1\big)\sqsubseteq\big( (X_i)_{i\in 2}, B_0, B_1\big)$, $(f_i)_{i\in 2}$ witnesses for this inequality, and $\varepsilon_0\!\in\! 2$ such that $B_{\varepsilon_0}$ is Borel locally countable. Then 
${f_i}_{\vert\Pi_i''S_{\varepsilon_0}}$ is countable-to-one for each $i\!\in\! 2$ and $S_{\varepsilon_0}$ is locally countable.\end{lem}

\noindent\bf Proof.\rm ~The inequality 
$\big( (Z_i)_{i\in 2},S_0,S_1\big)\sqsubseteq\big( (X_i)_{i\in 2}, B_0, B_1\big)$ gives 
$f_i\! :\! Z_i\!\rightarrow\! X_i$ continuous such that $(\Pi_{i\in 2}\ f_i)_{\vert S_0\cup S_1}$ is one-to-one, and also $S_\varepsilon\!\subseteq\! (\Pi_{i\in 2}\ f_i)^{-1}(B_\varepsilon )$ for each $\varepsilon\!\in\! 2$.\bigskip

\noindent $\bullet$ By the Lusin-Novikov theorem and Lemma 2.4.(a) in [L2] we can find Borel one-to-one partial functions $b_n$ with Borel domain such that 
$B_{\varepsilon_0}\! =\!\bigcup_{n\in\omega}\ \mbox{Gr}(b_n)$. Let us prove that 
$${f_i}_{\vert\Pi_i[S_{\varepsilon_0}\cap (\Pi_{i\in 2}\ f_i)^{-1}(\mbox{Gr}(b_n))]}$$ 
is one-to-one for each $i\!\in\! 2$.\bigskip

 Assume for example that $i\! =\! 0$. Let 
$z\!\not=\! z'\!\in\!\Pi_0\big[ S_{\varepsilon_0}\cap (\Pi_{i\in 2}\ f_i)^{-1}\big(\mbox{Gr}(b_n)\big)\big]$, and $y,y'\!\in\! Z_1$ such that 
$(z,y),(z',y')\!\in\! S_{\varepsilon_0}\cap (\Pi_{i\in 2}\ f_i)^{-1}\big(\mbox{Gr}(b_n)\big)$. As 
$(z,y)\!\not=\! (z',y')$, we get 
$$\big( f_0(z),f_1(y)\big)\!\not=\!\big( f_0(z'),f_1(y')\big).$$ 
But $b_n\big( f_0(z)\big)\! =\! f_1(y)$, $b_n\big( f_0(z')\big)\! =\! f_1(y')$, so that $f_0(z)\!\not=\! f_0(z')$ since $b_n$ is a partial function. If $i\! =\! 1$, then we use the fact that $b_n$ is one-to-one to see that 
${f_i}_{\vert\Pi_i[S_{\varepsilon_0}\cap (\Pi_{i\in 2}\ f_i)^{-1}(\mbox{Gr}(b_n))]}$ is also one-to-one.
\bigskip

\noindent $\bullet$ This proves that ${f_i}_{\vert\Pi_i''S_{\varepsilon_0}}$ is countable-to-one since 
$S_{\varepsilon_0}\! =\!
\bigcup_{n\in\omega}\ S_{\varepsilon_0}\cap (\Pi_{i\in 2}\ f_i)^{-1}\big(\mbox{Gr}(b_n)\big)$.\bigskip

\noindent $\bullet$ Now $S_{\varepsilon_0}$ is locally countable since $S_{\varepsilon_0}\!\subseteq\! 
(\Pi_{i\in 2}\ {f_i}_{\vert\Pi_i''S_{\varepsilon_0}})^{-1}(B_{\varepsilon_0})$, $B_{\varepsilon_0}$ is locally countable and ${f_i}_{\vert\Pi_i''S_{\varepsilon_0}}$ is countable-to-one for each $i\!\in\! 2$.
\hfill{$\square$}

\begin{lem} Let $Y$ be a Polish space, $C$ a Borel subset of $Y$ and $(m_n)_{n\in\omega}$ a sequence of Borel partial functions from a Borel subset of $C$ into $C$. We assume that 
$M\! :=\!\bigcup_{n\in\omega}\ \mbox{Gr}(m_n)$ is disjoint from $\Delta (C)$, but not separable from 
$\Delta (C)$ by a $\mbox{pot}(\bormone )$ set. Then there are integers $n\! <\! p$ and $y\!\in\! C$ such that 
$m_n(y)$ and $m_n\big( m_p(y)\big)$ are defined.\end{lem}

\noindent\bf Proof.\rm ~We may assume that $Y$ is recursively presented and $C,M$ and the $m_n$'s are $\Borel$. We put 
$$V\! :=\!\bigcup\{ D\!\in\!\Borel (Y)\mid D^2\cap M\mbox{ has finite vertical sections}\}.$$
Then $V\!\in\!\Ca (Y)$.

\vfill\eject

\noindent\bf Case 1.\rm\ $V\! =\! Y$.\bigskip

 We can find a sequence $(D_n)_{n\in\omega}$ of $\Borel$ subsets of $Y$ such that 
$Y\! =\!\bigcup_{n\in\omega}\ D_n$ and $D_n^2\cap M$ has finite vertical sections. By Theorem 3.6 in [Lo2], $D_n^2\cap M$ is $\mbox{pot}(\bormone )$, so that $D_n^2\!\setminus\! M$ is 
$\mbox{pot}(\boraone )$. Thus 
$\Delta (C)\!\subseteq\!\bigcup_{n\in\omega}\ D_n^2\!\setminus\! M\!\subseteq\!\neg M$ and $\Delta (C)$ is separable from $M$ by a $\mbox{pot}(\boraone )$ set, which is absurd.\bigskip

\noindent\bf Case 2.\rm\ $V\!\not=\! Y$.\bigskip

 The first reflection theorem proves that for each nonempty $\Ana$ subset $S$ of $Y$ contained in 
$Y\!\setminus\! V$ there is $y\!\in\! S$ such that $(S^2\cap M)_y$ is infinite. So there is an integer $n$ such that $(Y\!\setminus\! V)^2\cap\mbox{Gr}(m_n)\!\not=\!\emptyset$. In particular, 
$S\! :=\! (Y\!\setminus\! V)\cap m_n^{-1}(Y\!\setminus\! V)$ is a nonempty $\Ana$ subset of $Y$, which gives $y\!\in\! S$ such that $(S^2\cap M)_y$ is infinite. This proves the existence of $p\! >\! n$ such that 
$\big( y,m_p(y)\big)\!\in\! S^2$. Note that $y\!\in\! C$ since $Y\!\setminus\! C\!\subseteq\! V$. Now it is clear that $n$, $p$ and $y$ are suitable.\hfill{$\square$}

\begin{lem} Let $i\!\in\! 2$, $Y_i$ a Polish space, $\delta_i$ a Borel subset of $Y_i$, $c\! :\!\delta_0\!\rightarrow\!\delta_1$ a Borel isomorphism, $n\!\in\!\omega$, $c_n$ a Borel one-to-one partial function from $Y_0$ into $Y_1$ with Borel domain, and $C_0\! :=\!\bigcup_{n\in\omega}\ \mbox{Gr}(c_n)$. We assume that $C_0\cap (\delta_0\!\times\!\delta_1)$ is disjoint from $\mbox{Gr}(c)$, but not separable from 
$\mbox{Gr}(c)$ by a $\mbox{pot}(\bormone )$ set. Then there are integers $n\! <\! p$ and $y_0\!\in\! Y_0$ such that $(cc_n^{-1}c_p)(y_0)$ and $(cc_n^{-1}c)(y_0)$ are defined and different.\end{lem}

\noindent\bf Proof.\rm ~We set $c'_n\! :=\! {c_n}_{\vert\delta_0\cap c_n^{-1}(\delta_1)}$, so that 
$C_0\cap (\delta_0\!\times\!\delta_1)\! =\!\bigcup_{n\in\omega}\ \mbox{Gr}(c'_n)$. Now we consider the 
pre-images 
$$\Delta (\delta_1)\! =\!
(c^{-1}\!\times\!\mbox{Id}_{\delta_1})^{-1}\big(\mbox{Gr}(c)\big)$$ 
and $\mbox{Gr}(c''_n)\! =\!
(c^{-1}\!\times\!\mbox{Id}_{\delta_1})^{-1}\big(\mbox{Gr}(c'_n)\big)$, where 
$c''_n\! :=\! c'_n\circ c^{-1}_{\vert c[\delta_0\cap c_n^{-1}(\delta_1)]}$. Note that $c''_n$ is a Borel one-to-one partial function with Borel domain and that 
$C''_0\! :=\!\bigcup_{n\in\omega}\ \mbox{Gr}(c''_n)$ is not separable from $\Delta (\delta_1)$ by a 
$\mbox{pot}(\bormone )$ set. This implies that $\bigcup_{n\in\omega}\ \mbox{Gr}\big( (c''_n)^{-1}\big)$ is not separable from $\Delta (\delta_1)$ by a $\mbox{pot}(\bormone )$ set.\bigskip

 By Lemma 7.4 there are integers $n\! <\! p$ and $y_1\!\in\!\delta_1$ such that $(c''_n)^{-1}(y_1)$ and 
$(c''_n)^{-1}\big( (c''_p)^{-1}(y_1)\big)$ are defined. We set $y_0\! :=\! (c'_p)^{-1}(y_1)$, so that 
$\big( c(c'_n)^{-1}c'_p\big)(y_0)$ and $\big( c(c'_n)^{-1}c\big)(y_0)$ are defined and equal respectively to 
$\big( cc_n^{-1}c_p\big)(y_0)$ and $\big( cc_n^{-1}c\big)(y_0)$. Now note that 
$y_1\!\not=\! (c''_p)^{-1}(y_1)$ for each $y_1$ in the range of $c''_p$. This implies that 
$(c''_n)^{-1}(y_1)\!\not=\! (c''_n)^{-1}\big( (c''_p)^{-1}(y_1)\big)$, 
$$\big( c(c'_n)^{-1}\big)(y_1)\!\not=\!\big( c(c'_n)^{-1}c(c'_p)^{-1}\big)(y_1)\mbox{,}$$ 
$\big( c(c'_n)^{-1}c'_p\big)(y_0)\!\not=\!\big( c(c'_n)^{-1}c\big)(y_0)$ and 
$\big( cc_n^{-1}c_p\big)(y_0)\!\not=\!\big( cc_n^{-1}c\big)(y_0)$.\hfill{$\square$}

\begin{lem} Let $Y$ be a Polish space, $n\!\in\!\omega$, $c$ and $c_n$ continuous open partial functions from $Y$ into $Y$ with open domain, $\varepsilon\!\in\! 2$, 
$C^\varepsilon\! :=\!\bigcup_{n\in\omega}\ \mbox{Gr}(c_{2n+\varepsilon})$. We assume that 
$C^0$ is disjoint from $C^1\cup\mbox{Gr}(c)$, but 
$\emptyset\!\not=\!\mbox{Gr}(c)\!\subseteq\!\overline{C^0}\cap\overline{C^1}$. Then 
$C^0$ is not separable from $C^1$ by a $\mbox{pot}(\borone )$ set, and $C^0$ is not separable from 
$\mbox{Gr}(c)$ by a $\mbox{pot}(\bormone )$ set. If moreover the domains $\mbox{Dom}(c_n)$ are dense, then $C^0\cap (\bigcap_{n\in\omega}\ \mbox{Dom}(c_n)\!\times\! 2^\omega )$ is not separable from 
$C^1\cap (\bigcap_{n\in\omega}\ \mbox{Dom}(c_n)\!\times\! 2^\omega )$ by a $\mbox{pot}(\borone )$ set.\end{lem}

\noindent\bf Proof.\rm ~We argue by contradiction, which gives $P\!\in\!\mbox{pot}(\borone )$. Let $G_i$ be a dense $G_\delta$ subset of $Y_i$ such that $P\cap (G_0\!\times\! G_1)\!\in\!\borone (G_0\!\times\! G_1)$. The proof of Lemma 3.5 in [L1] shows the inclusion $\mbox{Gr}(c)\!\subseteq\!\overline{\mbox{Gr}(c)\cap (G_0\!\times\! G_1)}$, and similarly with $c_n$. Thus
$$\begin{array}{ll}
\mbox{Gr}(c)\!\!\!\! 
& \!\subseteq\!\overline{\overline{C^0}\cap\overline{C^1}\cap (G_0\!\times\! G_1)}\!\subseteq\!
\overline{\overline{C^0\cap (G_0\!\times\! G_1)}\cap\overline{C^1\cap (G_0\!\times\! G_1)}\cap 
(G_0\!\times\! G_1)}\cr\cr
& \!\subseteq\!\overline{\big( P\cap (G_0\!\times\! G_1)\big)\!\setminus\!\big( P\cap (G_0\!\times\! G_1)\big)}
\! =\!\emptyset\mbox{,}
\end{array}$$
which is absurd. The last assertion follows since we may assume that 
$G_0\!\subseteq\!\bigcap_{n\in\omega}\ \mbox{Dom}(c_n)$. The proof of the second assertion is similar and simpler.\hfill{$\square$}

\begin{lem} There is a tuple $\big( (Y_i)_{i\in 2}, C_0,C_1\big)$ such that\smallskip

\noindent (a) $Y_0$ and $Y_1$ are Polish spaces.\smallskip

\noindent (b) $C_0\! =\!\bigcup_{n\in\omega}\ \mbox{Gr}(c_n)\!\subseteq\!\Pi_{i\in 2}\ Y_i$, for some 
Borel one-to-one partial functions $c_n$ with Borel domain.\smallskip

\noindent (c) $C_1\! =\!\mbox{Gr}(c)$, for some Borel function $c\! :\! Y_0\!\rightarrow\! Y_1$.\smallskip

\noindent (d) $C_0$ is disjoint from $C_1$, but not separable from $C_1$ by a $\mbox{pot}(\bormone )$ set.\smallskip

\noindent (e) We set $C^\varepsilon_0\! :=\!\big(\bigcup_{n\in\omega}\ \mbox{Gr}(c_{2n+\varepsilon})\big)\cap (\bigcap_{n\in\omega}\ \mbox{Dom}(c_n)\!\times\! 2^\omega )$, for $\varepsilon\!\in\! 2$. Then 
$C^0_0$ is disjoint from $C^1_0$, but not separable from $C^1_0$ by a $\mbox{pot}(\borone )$ set, and 
$\overline{C^0_0}\cap\overline{C^1_0}\cap (\bigcap_{n\in\omega}\ \mbox{Dom}(c_n)\!\times\! 2^\omega )
\!\subseteq\!\mbox{Gr}(c)$.\smallskip

\noindent (f) The equality $\big( cc_n^{-1}c_p\big)(y_0)\! =\!\big( cc_n^{-1}c\big)(y_0)$ holds as soon as the two members of the equality  are defined and $n\! <\! p$.\end{lem}

\noindent\bf Proof.\rm ~We set $Y_i\! :=\! 2^\omega$ and $c(\alpha )(k)\! :=\!\alpha (2k)$.\bigskip

\noindent $\bullet$ We first build an increasing sequence $(S_n)_{n\in\omega}$ of co-infinite subsets of 
$\omega$, a sequence $(\psi_n)_{n\in\omega}$ of bijections, and a sequence $(h_n)_{n\in\omega}$ of homeomorphisms of $2^\omega$ onto itself. We do it by induction on $n$. We set $S_0\! :=\!\emptyset$, $\psi_0\! :=\!\mbox{Id}_\omega$ and $h_0\! :=\!\mbox{Id}_{2^\omega}$. Assume that $(S_q)_{q\leq n}$, 
$(\psi_q)_{q\leq n}$ and $(h_q)_{q\leq n}$ are constructed, which is the case for $n\! =\! 0$. We define a map $\varphi_n\! :\!\omega\!\rightarrow\!\omega$ by 
$$\varphi_n(k)\! :=\!\left\{\!\!\!\!\!\!
\begin{array}{ll}
& \psi_n^{-1}(k)\mbox{ if }k\!\notin\! 2S_n\mbox{,}\cr\cr
& \frac{k}{2}\mbox{ if }k\!\in\! 2S_n.
\end{array}
\right.$$
Note that $\varphi_n$ is a bijection. We set $S_{n+1}\! :=\!\varphi_n[2\omega ]\cup (n\! +\! 1)$, which is co-infinite. The sequence $(S_n)_{n\in\omega}$ is increasing since 
$S_n\! =\!\varphi_n[2S_n]\!\subseteq\! S_{n+1}$. As $S_{n+1}$ is co-infinite we can build the bijection  
$\psi_{n+1}\! :\!\omega\!\setminus\! S_{n+1}\!\rightarrow\!\omega\!\setminus\! 2S_{n+1}$ in such a way that 
$\psi_{n+1}(k)\!\not=\!\psi_q(k)$ for infinitely many $k\!\notin\! S_{n+1}$, for each $q\!\leq\! n$. We set 
$$h_{n+1}(\alpha )(k)\! :=\!\left\{\!\!\!\!\!\!
\begin{array}{ll}
& c(\alpha )(k)\mbox{ if }k\!\in\! S_{n+1}\mbox{,}\cr\cr
& \alpha\big(\psi_{n+1}(k)\big)\mbox{ if }k\!\notin\! S_{n+1}.
\end{array}
\right.$$
As $h_{n+1}$ permutes the coordinates, it is an homeomorphism.\bigskip

\noindent $\bullet$ We set $D_n\! :=\!\{\alpha\!\in\! 2^\omega\mid c(\alpha )\!\not=\! h_n(\alpha)\wedge\forall q\! <\! n\ \ h_n(\alpha )\!\not=\! h_q(\alpha )\}$, so that $D_n$ is an open subset of $2^\omega$. We set $c_n\! :=\! {h_n}_{\vert D_n}$, so that 
$c_n$ is an homeomorphism from $D_n$ onto its open range, $C_0$ is disjoint from $C_1$, and $C^0_0$ is disjoint from $C^1_0$.\bigskip

 Let us prove that $D_n$ is dense for each integer $n$. Note that 
$D_0\! =\!\{\alpha\!\in\! 2^\omega\mid\exists k\!\in\!\omega\ \ \alpha (2k)\!\not=\!\alpha (k)\}$, which is clearly dense. Now $D_{n+1}$ contains 
$$\{\alpha\!\in\! 2^\omega\mid\exists k\!\notin\! S_{n+1}\ \ \alpha (2k)\!\not=\!\alpha\big(\psi_{n+1}(k)\big)\}\cap\bigcap_{q<n}\ \{\alpha\!\in\! 2^\omega\mid\exists k\!\notin\! S_{n+1}\ \ \alpha\big(\psi_{n+1}(k)\big)\!\not=\!\alpha\big(\psi_{q}(k)\big)\}.$$
The set 
$\{\alpha\!\in\! 2^\omega\mid\exists k\!\notin\! S_{n+1}\ \ \alpha (2k)\!\not=\!\alpha\big(\psi_{n+1}(k)\big)\}$ is open dense since the odd integers are in $\psi_{n+1}[\omega\!\setminus\! S_{n+1}]$. The set 
$\{\alpha\!\in\! 2^\omega\mid\exists k\!\notin\! S_{n+1}\ \ \alpha\big(\psi_{n+1}(k)\big)\!\not=\!
\alpha\big(\psi_{q}(k)\big)\}$ is open dense by construction of $\psi_{n+1}$. This proves that $D_{n+1}$ is dense.\bigskip

\noindent $\bullet$ Note that $\mbox{Gr}(c)\!\subseteq\!\overline{C^0_0}\cap\overline{C^1_0}$ since 
$c(\alpha )\vert n\! =\! h_n(\alpha )\vert n$, $D_n$ is dense and $c$ is continuous. Lemma 7.6 proves the non-separation assertions. We also have 
$\overline{C^0_0}\cap\overline{C^1_0}\cap (\bigcap_{n\in\omega}\ \mbox{Dom}(c_n)\!\times\! 2^\omega )
\!\subseteq\!\mbox{Gr}(c)$ since $c(\alpha )\vert n\! =\! h_n(\alpha )\vert n$ and $c_n$ is continuous.\bigskip

\noindent $\bullet$ Now it is enough to prove that $ch_n^{-1}h_p\! =\! ch_n^{-1}c$ if $n\! <\! p$. We have
$$h_{n}^{-1}(\beta )(j)\! :=\!\left\{\!\!\!\!\!\!
\begin{array}{ll}
& \beta (k)\mbox{ if }j\! =\! 2k\!\in\! 2S_{n}\mbox{,}\cr\cr
& \beta\big(\psi_{n}^{-1}(j)\big)\mbox{ if }j\!\notin\! 2S_{n}.
\end{array}
\right.$$
Thus 
$$(ch_n^{-1}c)(\alpha )(k)\! =\! c\big( (h_n^{-1}c)(\alpha )\big)(k)\! =\! (h_n^{-1}c)(\alpha )(2k)\! =\!\left\{\!\!\!\!\!\!
\begin{array}{ll} 
& c(\alpha )(k)\mbox{ if }k\!\in\! S_n\mbox{,}\cr\cr
& c(\alpha )\big(\psi_n^{-1}(2k)\big)\mbox{ if }k\!\notin\! S_n.
\end{array}
\right.$$
Similarly,
$$(ch_n^{-1}h_p)(\alpha )(k)\! =\!\left\{\!\!\!\!\!\!
\begin{array}{ll} 
& h_p(\alpha )(k)\mbox{ if }k\!\in\! S_n\mbox{,}\cr\cr
& h_p(\alpha )\big(\psi_n^{-1}(2k)\big)\mbox{ if }k\!\notin\! S_n.
\end{array}
\right.$$
Note that $S_n\! =\!\varphi_n[2S_n]\!\subseteq\! S_{n+1}$, so that $S_n\!\subseteq\! S_p$. Thus 
$(ch_n^{-1}h_p)(\alpha )(k)\! =\! (ch_n^{-1}c)(\alpha )(k)$ if $k\!\in\! S_n$. If $k\!\notin\! S_n$, then 
$2k\!\notin\! 2S_n$ and $\varphi_n(2k)\! =\!\psi_n^{-1}(2k)\!\in\! S_{n+1}\!\subseteq\! S_p$. Thus 
$$(ch_n^{-1}h_p)(\alpha )(k)\! =\! h_p(\alpha )\big(\psi_n^{-1}(2k)\big)\! =\! c(\alpha )\big(\psi_n^{-1}(2k)\big)
\! =\! (ch_n^{-1}c)(\alpha )(k).$$
This finishes the proof.\hfill{$\square$}\bigskip

\noindent\bf Proof of Theorem 7.1.\rm ~We argue by contradiction. Note that $S_0$ is not separable from $S_1$ by a $\mbox{pot}(\bormone )$ set since (b) holds. By Lemma 7.2 we may assume that 
the inequality $S_1\!\subseteq\!\Pi_0''S_0\!\times\!\Pi_1''S_0$ holds.\bigskip

\noindent $\bullet$ Recall the digraph $A_1$ in [L5]. If we take $X_i\! :=\! 2^\omega$, $B_0\! :=\! A_1$ and $B_1\! :=\!\Delta (2^\omega )$, then by Corollary 12 in [L5], $B_0$ is Borel locally countable, not 
$\mbox{pot}(\bormone )$, and $B_1\! =\!\overline{B_0}\!\setminus\! B_0$. It follows that $B_0$ is not separable from $B_1$ by a $\mbox{pot}(\bormone )$ set $Q$, since otherwise we would have 
$B_0\! =\! Q\cap\overline{B_0}\!\in\!\mbox{pot}(\bormone )$. This implies that 
$\big( (X_i)_{i\in 2}, B_0,B_1\big)$ satisfies condition (b) in Theorem 7.1. By Lemma 7.3, 
${f_i}_{\vert\Pi_i''S_0}$ is countable-to-one for each $i\!\in\! 2$ and $S_0$ is locally countable.\bigskip

\noindent $\bullet$ Lemma 7.7 gives a tuple $\big( (Y_i)_{i\in 2}, C_0,C_1\big)$. Note that 
$\big( (Y_i)_{i\in 2}, C_0,C_1\big)$ satisfies condition (b) in Theorem 7.1, which gives 
$g_i\! :\! Z_i\!\rightarrow\! Y_i$. Lemma 7.3 implies that ${g_i}_{\vert\Pi_i''S_0}$ is countable-to-one for each $i\!\in\! 2$. The first reflection theorem gives a Borel set $O_i\!\supseteq\!\Pi_i''S_0$ such that 
${f_i}_{\vert O_i}$ and ${g_i}_{\vert O_i}$ are countable-to-one, for each $i\!\in\! 2$. By Lemma 2.4.(a) in [L2] we can find a partition $(O^i_n)_{n\in\omega}$ of $O_i$ into Borel sets such that 
${f_i}_{\vert O^i_n}$ and ${g_i}_{\vert O^i_n}$ are one-to-one, for each $i\!\in\! 2$.\bigskip 

\noindent $\bullet$ We set $S''_\varepsilon\! :=\! (\Pi_{i\in 2}\ {f_i}_{\vert O_i})^{-1}(B_\varepsilon )\cap 
(\Pi_{i\in 2}\ {g_i})^{-1}(C_\varepsilon )$, for each $\varepsilon\!\in\! 2$, so that $S''_\varepsilon$ is a Borel subset of $\Pi_{i\in 2}\ Z_i$ containing $S_\varepsilon$. In particular, $S''_0$ is not separable from $S''_1$ by a $\mbox{pot}(\bormone )$ set. We choose integers  $n_0$ and $n_1$ such that 
$S''_0\cap (\Pi_{i\in 2}\ O^i_{n_i})$ is not separable from $S''_1\cap (\Pi_{i\in 2}\ O^i_{n_i})$ by a 
$\mbox{pot}(\bormone )$ set. We set 
$D_\varepsilon\! :=\! (\Pi_{i\in 2}\ {f_i}_{\vert O^i_{n_i}})[S''_\varepsilon\cap (\Pi_{i\in 2}\ O^i_{n_i})]$, so that $D_0$ is a Borel subset of $B_0$ which is not separable from $D_1$ by a $\mbox{pot}(\bormone )$ set. Note that $D_1$  is a Borel subset of $B_1\! =\!\Delta (2^\omega )$. In particular, there is a Borel subset $D$ of $2^\omega$ such that $D_1\! =\!\Delta (D)$. By Lemma 7.2, $D_0\cap D^2$ is not separable from $D_1$ by a $\mbox{pot}(\bormone )$ set. Let $h_i\! :\! D\!\rightarrow\! Y_i$ be defined by 
$h_i(\alpha )\! :=\! (g_i\circ {f_i}_{\vert O^i_{n_i}}^{-1})(\alpha)$. Then $h_i$ is Borel, one-to-one, and 
$D_\varepsilon\cap D^2\!\subseteq\! B_\varepsilon\cap (\Pi_{i\in 2}\ h_i)^{-1}(C_\varepsilon )$.\bigskip 

\noindent $\bullet$ Note that $(\Pi_{i\in 2}\ h_i)[\Delta (D)]$ is a Borel subset of $C_1$, which proves the existence of a Borel subset $\delta$ of $Y_0$ such that 
$(\Pi_{i\in 2}\ h_i)[\Delta (D)]\! =\!\mbox{Gr}(c_{\vert\delta})$. If 
$y\!\not=\! y'\!\in\!\delta$, then $\big( y,c(y)\big)\! =\!\big( h_0(d),h_1(d)\big)$ and 
$$\big( y',c(y')\big)\! =\!\big( h_0(d'),h_1(d')\big)$$ 
for some $d\!\not=\! d'\!\in\! D$. As $h_1$ is one-to-one we get $c(y)\!\not=\! c(y')$, $c_{\vert\delta}$ is 
one-to-one and $c''\delta$ is Borel.\bigskip

 As $D_0\cap D^2\!\subseteq\! (\Pi_{i\in 2}\ h_i)^{-1}(C_0)$ and 
$D_1\!\subseteq\! (\Pi_{i\in 2}\ h_i)^{-1}\big(\mbox{Gr}(c_{\vert\delta})\big)$, $C_0$ is not separable from 
$\mbox{Gr}(c_{\vert\delta})$ by a $\mbox{pot}(\bormone )$ set. By Lemma 7.2, 
$C'_0\! :=\! C_0\cap (\delta\!\times\! c''\delta )$ is not separable from $\mbox{Gr}(c_{\vert\delta})$ by a 
$\mbox{pot}(\bormone )$ set.\bigskip

\noindent $\bullet$ By Lemma 7.5 applied to $\delta_0\! :=\!\delta$ and $\delta_1\! :=\! c''\delta$ there are 
$n\! <\! p$ and $y_0\!\in\! Y_0$ such that $(cc_n^{-1}c_p)(y_0)$ and $(cc_n^{-1}c)(y_0)$ are defined and different, which contradicts Lemma 7.7.(f).\hfill{$\square$}\bigskip

\noindent\bf Remark.\rm ~We recover the algebraic relation ``$g_n\! =\! g_n\circ g_p$ if $n\! <\! p$" that was already present in Section 3 of [L5] mentioned just after the statement of Theorem 7.1.

\begin{thm} There is no tuple $\big( (Z_i)_{i\in 2},S_0,S_1)$, where the $Z_i$'s are Polish spaces and $S_0$, $S_1$ disjoint analytic subsets of $\Pi_{i\in 2}\ Z_i$, such that for any tuple 
$\big( (X_i)_{i\in 2}, B_0, B_1\big)$ of the same type exactly one of the following holds:\smallskip

\noindent (a) The set $B_0$ is separable from $B_1$ by a $\mbox{pot}(\borone )$ set.\smallskip

\noindent (b) The inequality $\big( (Z_i)_{i\in 2},S_0,S_1\big)\sqsubseteq
\big( (X_i)_{i\in 2}, B_0, B_1\big)$ holds.\end{thm}

\noindent\bf Proof.\rm ~Let us indicate the differences with the proof of Theorem 7.1. This time, $S_0$ is not separable from $S_1$ by a $\mbox{pot}(\borone )$ set.\bigskip

\noindent $\bullet$ Note that $A_1\! =\!\bigcup_{n\in\omega}\ \mbox{Gr}(H_n)$, where 
$H_n\! :\! N_{s_n0}\!\rightarrow\! N_{s_n1}$ is a partial homeomorphism with clopen domain and range. The crucial properties of $(s_n)_{n\in\omega}\!\subseteq\! 2^{<\omega}$ is that it is dense and 
$\vert s_n\vert\! =\! n$. We can easily ensure this in such a way that $(s_{2n})_{n\in\omega}$ and 
$(s_{2n+1})_{n\in\omega}$ are dense. We set 
$B_\varepsilon\! :=\!\bigcup_{n\in\omega}\ \mbox{Gr}(H_{2n+\varepsilon})$. The previous remark implies that $\Delta (2^\omega )\! =\!\overline{B_\varepsilon}\!\setminus\! B_\varepsilon$. By Lemma 7.6, $B_0$ is not separable from $B_1$ by a $\mbox{pot}(\borone )$ set. So here again ${f_i}_{\vert\Pi_i''S_0}$ is countable-to-one for each $i\!\in\! 2$, and $S_0$, $S_1$ are locally countable by Lemma 7.3.

\vfill\eject

\noindent $\bullet$ Lemma 7.7 gives a tuple 
$\Big(\big(\bigcap_{n\in\omega}\mbox{ Gr}(c_n), 2^\omega\big), C^0_0,C^1_0\Big)$. Note that 
$\Big(\big(\bigcap_{n\in\omega}\mbox{ Gr}(c_n), 2^\omega\big), C^0_0,C^1_0\Big)$ satisfies condition (b) in Theorem 7.8.\bigskip

\noindent $\bullet$ We change the topology on $2^\omega$ into a finer Polish topology $\tau$ so that the sets $f_i''O^i_{n_i}$ become clopen and the maps $({f_i}_{\vert O^i_{n_i}})^{-1}$ become continuous. Now $$\overline{D_0}^{\tau^2}\cap\overline{D_1}^{\tau^2}\!\subseteq\!\overline{B_0}\cap\overline{B_1}\! =\!
\big( B_0\cup\Delta (2^\omega )\big)\cap\big( B_1\cup\Delta (2^\omega )\big)\! =\!\Delta (2^\omega ).$$
So there is a Borel subset $D$ of $2^\omega$ such that 
$\overline{D_0}^{\tau^2}\cap\overline{D_1}^{\tau^2}\! =\!\Delta (D)$, and 
$D\!\subseteq\!\bigcap_{i\in 2}\ f_i''O^i_{n_i}$.\bigskip

\noindent $\bullet$ Let us prove that $D_0\cap D^2$ is not separable from $D_1\cap D^2$ by a 
$\mbox{pot}(\borone )$ set.\bigskip

 We argue by contradiction, which gives $P\!\in\!\mbox{pot}(\borone )$ such that 
$D_0\cap D^2\!\subseteq\! P\!\subseteq\! D^2\!\setminus\! D_1$. The sets 
$\overline{D_0}^{\tau^2}\cap (\neg D\!\times\! 2^\omega )$ and 
$\overline{D_1}^{\tau^2}\cap (\neg D\!\times\! 2^\omega )$ are disjoint, $\mbox{pot}(\bormone )$, so that they are separable by $\Delta_l$ in $\mbox{pot}(\borone )$. Similarly, there is 
$\Delta_r\!\in\!\mbox{pot}(\borone )$ which separates 
$\overline{D_0}^{\tau^2}\cap (2^\omega\!\times\!\neg D)$ from 
$\overline{D_1}^{\tau^2}\cap (2^\omega\!\times\!\neg D)$. Now 
$$D_0\!\subseteq\! P\cup\big( D_0\cap (\neg D\!\times\! 2^\omega )\big)\cup
\big( D_0\cap (2^\omega\!\times\!\neg D)\big)\!\subseteq\! P\cup
\big(\Delta_l\cap (\neg D\!\times\! 2^\omega )\big)\cup
\big(\Delta_r\cap (2^\omega\!\times\!\neg D)\big)\!\subseteq\!\neg D_1$$
which is absurd since $P\cup\big(\Delta_l\cap (\neg D\!\times\! 2^\omega )\big)\cup
\big(\Delta_r\cap (2^\omega\!\times\!\neg D)\big)\!\in\!\mbox{pot}(\borone )$.\bigskip

\noindent $\bullet$ Let us prove that $D_0\cap D^2$ is not separable from $\Delta (D)$ by a 
$\mbox{pot}(\bormone )$ set.\bigskip

 We argue by contradiction, which gives $Q\!\in\!\mbox{pot}(\bormone )$ such that 
$D_0\cap D^2\!\subseteq\! Q\!\subseteq\! D^2\!\setminus\!\Delta (D)$. The sets 
$Q$ and $\Delta (D)$ are disjoint, $\mbox{pot}(\bormone )$, so that there is $R$ in $\mbox{pot}(\borone )$ such that $Q\!\subseteq\! R\!\subseteq\! D^2\!\setminus\!\Delta (D)$. The sets 
$\overline{D_0}^{\tau^2}\cap R$ and $\overline{D_1}^{\tau^2}\cap R$ are disjoint, 
$\mbox{pot}(\bormone )$, so that there is $S$ in $\mbox{pot}(\borone )$ such that 
$\overline{D_0}^{\tau^2}\cap R\!\subseteq\! S\!\subseteq\! R\!\setminus\!\overline{D_1}^{\tau^2}$. But $S$ separates $D_0\cap D^2$ from $D_1\cap D^2$, which contradicts the previous point.\bigskip

\noindent $\bullet$ Note that 
$(\Pi_{i\in 2}\ h_i)[\Delta (D)]\!\subseteq\!\overline{C^0_0}\cap\overline{C^1_0}\cap 
(\bigcap_{n\in\omega}\ \mbox{Dom}(c_n)\!\times\! 2^\omega )\!\subseteq\!\mbox{Gr}(c)$. We conclude as in the proof of Theorem 7.1.\hfill{$\square$}

\vfill\eject

\section{$\!\!\!\!\!\!$ References}

\noindent [B]\ \ B. Bollob\'as,~\it Modern graph theory,~\rm 
Springer-Verlag, New York, 1998

\noindent [C]\ \ D. Cenzer, Monotone inductive definitions over the continuum,~
\it J. Symbolic Logic\rm ~41 (1976), 188-198

\noindent [D-SR]\ \ G. Debs and J. Saint Raymond, Borel liftings of Borel sets: 
some decidable and undecidable statements,~\it Mem. Amer. Math. Soc.\rm ~187, 876 (2007)

\noindent [H-K-Lo]\ \ L. A. Harrington, A. S. Kechris and A. Louveau, A Glimm-Effros 
dichotomy for Borel equivalence relations,~\it J. Amer. Math. Soc.\rm ~3 (1990), 903-928

\noindent [K]\ \ A. S. Kechris,~\it Classical Descriptive Set Theory,~\rm 
Springer-Verlag, 1995

\noindent [K-S-T]\ \ A. S. Kechris, S. Solecki and S. Todor\v cevi\'c, Borel chromatic numbers,\ \it 
Adv. Math.\rm\ 141 (1999), 1-44

\noindent [L1]\ \ D. Lecomte, Classes de Wadge potentielles et 
th\'eor\`emes d'uniformisation partielle,\it ~Fund. Math.~\rm 143 (1993), 231-258

\noindent [L2]\ \ D. Lecomte, Uniformisations partielles et crit\`eres \`a la 
Hurewicz dans le plan,~\it Trans. Amer. Math. Soc.\rm ~347, 11 (1995), 4433-4460

\noindent [L3]\ \ D. Lecomte, Tests \`a la Hurewicz dans le plan,\it ~Fund. Math.~
\rm 156 (1998), 131-165

\noindent [L4]\ \ D. Lecomte, Complexit\'e des bor\'eliens~\`a coupes 
d\'enombrables,\ \it Fund. Math.~\rm 165 (2000), 139-174

\noindent [L5]\ \ D. Lecomte, On minimal non potentially closed subsets of the plane,
\ \it Topology Appl.~\rm 154, 1 (2007) 241-262

\noindent [L6]\ \ D. Lecomte, Hurewicz-like tests for Borel subsets of the plane,
~\it Electron. Res. Announc. Amer. Math. Soc.\rm\ 11 (2005)

\noindent [L7]\ \ D. Lecomte, How can we recognize potentially $\bormxi$ subsets of the plane?,
~\it to appear in J. Math. Log. (see arXiv)\rm\ 

\noindent [L8]\ \ D. Lecomte, A dichotomy characterizing analytic digraphs of uncountable Borel chromatic number in any dimension,~\it preprint (see arXiv)\rm

\noindent [Lo1]\ \ A. Louveau, Some results in the Wadge hierarchy of Borel sets,\ \it Cabal Sem. 
79-81, Lect. Notes in Math.\ \rm 1019 (1983), 28-55

\noindent [Lo2]\ \ A. Louveau, A separation theorem for $\Ana$ sets,\ \it Trans. 
Amer. Math. Soc.\ \rm 260 (1980), 363-378

\noindent [Lo3]\ \ A. Louveau, Ensembles analytiques et bor\'eliens dans les 
espaces produit,~\it Ast\'erisque (S. M. F.)\ \rm 78 (1980)

\noindent [Lo-SR1]\ \ A. Louveau and J. Saint Raymond, Borel classes and closed games: 
Wadge-type and Hurewicz-type results,\ \it Trans. Amer. Math. Soc.\ \rm 304 (1987), 431-467

\noindent [Lo-SR2]\ \ A. Louveau and J. Saint Raymond, The strength of Borel Wadge determinacy,\ \it 
Cabal Seminar 81-85, Lecture Notes in Math.\ \rm 1333 (1988), 1-30

\noindent [Lo-SR3]\ \ A. Louveau et J. Saint Raymond, Les propri\'et\'es de 
r\'eduction et de norme pour les classes de bor\'eliens,\ \it Fund. Math.\ \rm 
131 (1988), no. 3, 223-243

\noindent [M]\ \ Y. N. Moschovakis,~\it Descriptive set theory,~\rm North-Holland, 1980

\noindent [S]\ \ J. R. Steel, Determinateness and the separation property,
\ \it J. Symbolic Logic~\rm 46 (1981) 41-44

\noindent [vW]\ \ R. van Wesep, Separation principles and the axiom of determinateness,
\ \it J. Symbolic Logic~\rm 43 (1978) 77-81
 
\end{document}